# NEW CLASSES OF CODES FOR CRYPTOLOGISTS AND COMPUTER SCIENTISTS

**W. B. Vasantha Kandasamy**
e-mail: vasanthakandasamy@gmail.com
web: http://mat.iitm.ac.in/~wbv
www.vasantha.net

**Florentin Smarandache**
e-mail: smarand@unm.edu

**2008**

# NEW CLASSES OF CODES FOR CRYPTOLOGISTS AND COMPUTER SCIENTISTS

**W. B. Vasantha Kandasamy**
**Florentin Smarandache**

**2008**

# CONTENTS

# PREFACE

Historically a code refers to a cryptosystem that deals with linguistic units: words, phrases etc. We do not discuss such codes in this book. Here codes are message carriers or information storages or information transmitters which in time of need should not be decoded or read by an enemy or an intruder. When we use very abstract mathematics in using a specific code, it is difficult for non-mathematicians to make use of it. At the same time, one cannot compromise with the capacity of the codes. So the authors in this book have introduced several classes of codes which are explained very non-technically so that a strong foundation in higher mathematics is not needed. The authors also give an easy method to detect and correct errors that occur during transmission. Further some of the codes are so constructed to mislead the intruder. False n-codes, whole n-codes can serve this purpose.

These codes can be used by computer scientists in networking and safe transmission of identity thus giving least opportunity to the hackers. These codes will be a boon to cryptologists as very less mathematical background is needed.

To honour Periyar on his 125$^{th}$ birth anniversary and to recognize his services to humanity the authors have named a few new classes of codes in his name. This book has three chapters. Chapter one is introductory in nature. The notion of bicodes and their generalization, n-codes are introduced in chapter two. Periyar linear codes are introduced in chapter three.

Many examples are given for the reader to understand these new notions. We mainly use the two methods, viz. pseudo best n-approximations and n-coset leader properties to detect and correct errors.

The authors deeply acknowledge the unflinching support of Dr.K.Kandasamy, Meena and Kama.

W.B.VASANTHA KANDASAMY
FLORENTIN SMARANDACHE

**Chapter One**

# BASIC CONCEPTS

In this chapter we introduce the basic concepts about linear codes which forms the first section. Section two recalls the notion of bimatrices and n-matrices (n a positive integer) and some of their properties.

## 1.1 Introduction of Linear Code and its basic Properties

In this section we just recall the definition of linear code and enumerate a few important properties about them. We begin by describing a simple model of a communication transmission system given by the figure 1.1.

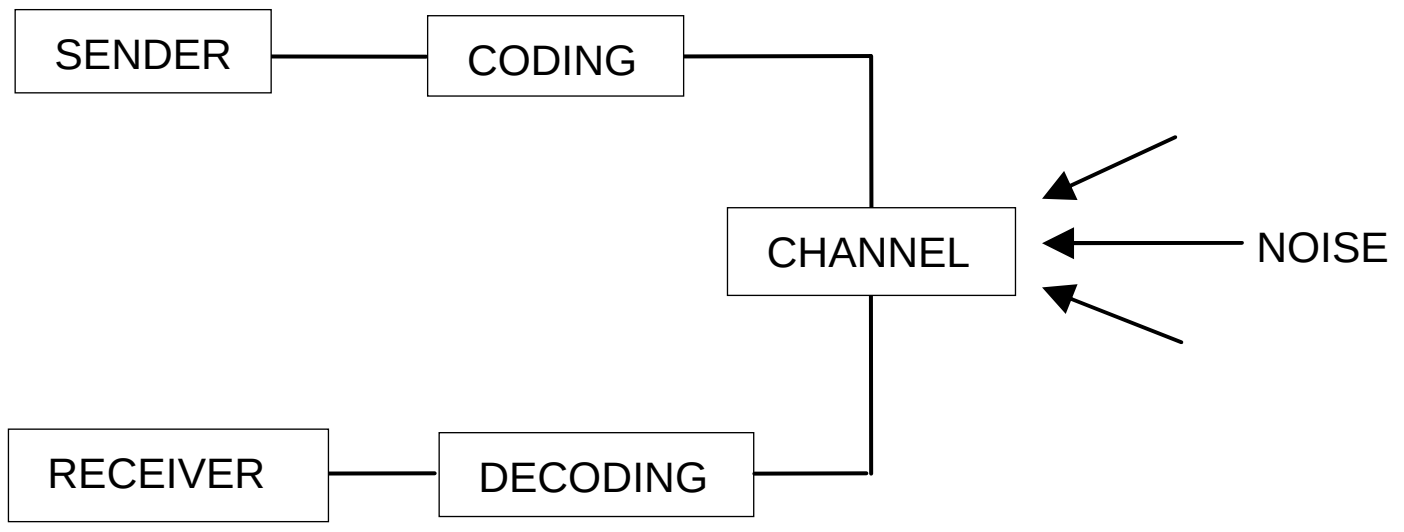

**Figure 1.1**

Messages go through the system starting from the source (sender). We shall only consider senders with a finite number of discrete signals (eg. Telegraph) in contrast to continuous sources (eg. Radio). In most systems the signals emanating from the source cannot be transmitted directly by the channel. For instance, a binary channel cannot transmit words in the usual Latin alphabet. Therefore an encoder performs the important task of data reduction and suitably transforms the message into usable form. Accordingly one distinguishes between source encoding the channel encoding. The former reduces the message to its essential(recognizable) parts, the latter adds redundant information to enable detection and correction of possible errors in the transmission. Similarly on the receiving end one distinguishes between channel decoding and source decoding, which invert the corresponding channel and source encoding besides detecting and correcting errors.

One of the main aims of coding theory is to design methods for transmitting messages error free cheap and as fast as possible. There is of course the possibility of repeating the message. However this is time consuming, inefficient and crude. We also note that the possibility of errors increases with an increase in the length of messages. We want to find efficient algebraic methods (codes) to improve the reliability of the transmission of messages. There are many types of algebraic codes; here we give a few of them.

Throughout this book we assume that only finite fields represent the underlying alphabet for coding. Coding consists of transforming a block of k message symbols $a_1, a_2, \ldots, a_k$; $a_i \in F_q$ into a code word $x = x_1 x_2 \ldots x_n$; $x_i \in F_q$, where $n \geq k$. Here the first $k_i$ symbols are the message symbols i.e., $x_i = a_i$; $1 \leq i \leq k$; the remaining $n - k$ elements $x_{k+1}, x_{k+2}, \ldots, x_n$ are check symbols or control symbols. Code words will be written in one of the forms x; $x_1, x_2, \ldots, x_n$ or $(x_1 x_2 \ldots x_n)$ or $x_1 x_2 \ldots x_n$. The check symbols can be obtained from the message symbols in such a way that the code words x satisfy a system of linear equations; $Hx^T = (0)$ where H is the given $(n - k) \times n$ matrix with elements in $F_q = Z_{p^n}$ $(q = p^n)$. A standard form for H is $(A, I_{n-k})$ with $n - k \times k$ matrix and $I_{n-k}$, the $n - k \times n - k$ identity matrix.

We illustrate this by the following example.

***Example 1.1.1:*** Let us consider $Z_2 = \{0, 1\}$. Take n = 7, k = 3. The message $a_1\ a_2\ a_3$ is encoded as the code word $x = a_1\ a_2\ a_3\ x_4\ x_5\ x_6\ x_7$. Here the check symbols $x_4\ x_5\ x_6\ x_7$ are such that for this given matrix

$$H = \begin{bmatrix} 0 & 1 & 0 & 1 & 0 & 0 & 0 \\ 1 & 0 & 1 & 0 & 1 & 0 & 0 \\ 0 & 0 & 1 & 0 & 0 & 1 & 0 \\ 0 & 0 & 1 & 0 & 0 & 0 & 1 \end{bmatrix} = (A; I_4);$$

we have $Hx^T = (0)$ where

$x = a_1\ a_2\ a_3\ x_4\ x_5\ x_6\ x_7$.

$$a_2 + x_4 = 0$$
$$a_1 + a_3 + x_5 = 0$$
$$a_3 + x_6 = 0$$
$$a_3 + x_7 = 0.$$

Thus the check symbols $x_4\ x_5\ x_6\ x_7$ are determined by $a_1\ a_2\ a_3$. The equation $Hx^T = (0)$ are also called check equations. If the message a = 1 0 0 then, $x_4 = 0$, $x_5 = 1$, $x_6 = 0$ and $x_7 = 0$. The code word x is 1 0 0 0 1 0 0. If the message a = 1 1 0 then $x_4 =1$, $x_5 = 1$, $x_6 = 1 = x_7$. Thus the code word x = 1 1 0 1 1 0 0.

We will have altogether $2^3$ code words given by

| | |
|---|---|
| 0 0 0 0 0 0 0 | 1 1 0 1 1 0 0 |
| 1 0 0 0 1 0 0 | 1 0 1 0 0 1 1 |
| 0 1 0 1 1 0 0 | 0 1 1 1 1 1 1 |
| 0 0 1 0 1 1 1 | 1 1 1 1 0 1 1 |

**DEFINITION 1.1.1:** *Let H be an n – k × n matrix with elements in $Z_q$. The set of all n-dimensional vectors satisfying $Hx^T = (0)$ over $Z_q$ is called a linear code(block code) C over $Z_q$ of block* *length n. The matrix H is called the parity check matrix of the code C. C is also called a linear(n, k) code.*

*If H is of the form(A, $I_{n-k}$) then the k-symbols of the code word $x$ is called massage(or information) symbols and the last n – k symbols in $x$ are the check symbols. C is then also called a systematic linear(n, k) code. If q = 2, then C is a binary code. k/n is called transmission (or information) rate.*

*The set C of solutions of $x$ of $Hx^T = (0)$. i.e., the solution space of this system of equations, forms a subspace of this system of equations, forms a subspace of $Z_q^n$ of dimension k. Since the code words form an additive group, C is also called a group code. C can also be regarded as the null space of the matrix H.*

***Example 1.1.2:*** (Repetition Code) If each codeword of a code consists of only one message symbol $a_1 \in Z_2$ and (n – 1) check symbols $x_2 = x_3 = \ldots = x_n$ are all equal to $a_1$ ($a_1$ is repeated n – 1 times) then we obtain a binary (n, 1) code with parity check matrix

$$H = \begin{bmatrix} 1 & 1 & 0 & 0 & \ldots & 1 \\ 0 & 0 & 1 & 0 & \ldots & 0 \\ 0 & 0 & 0 & 1 & \ldots & 0 \\ \vdots & \vdots & \vdots & \vdots & & \vdots \\ 1 & 0 & 0 & 0 & \ldots & 1 \end{bmatrix}.$$

There are only two code words in this code namely 0 0 … 0 and 1 1 …1.

If is often impracticable, impossible or too expensive to send the original message more than once. Especially in the transmission of information from satellite or other spacecraft, it is impossible to repeat such messages owing to severe time limitations. One such cases is the photograph from spacecraft as it is moving it may not be in a position to retrace its path. In such cases it is impossible to send the original message more than once. In repetition codes we can of course also consider code words with more than one message symbol.

***Example 1.1.3: (Parity-Check Code):*** This is a binary (n, n – 1) code with parity-check matrix to be H = (1 1 … 1). Each code word has one check symbol and all code words are given by all binary vectors of length n with an even number of ones. Thus if sum of the ones of a code word which is received is odd then atleast one error must have occurred in the transmission.

Such codes find its use in banking. The last digit of the account number usually is a control digit.

**DEFINITION 1.1.2:** *The matrix $G = (I_k, -A^T)$ is called a canonical generator matrix (or canonical basic matrix or encoding matrix) of a linear (n, k) code with parity check matrix $H = (A, I_{n-k})$. In this case we have $GH^T = (0)$.*

***Example 1.1.4:*** Let

$$G = \begin{bmatrix} 1 & 0 & 0 & 0 & 1 & 0 & 0 \\ 0 & 1 & 0 & 1 & 0 & 0 & 0 \\ 0 & 0 & 1 & 0 & 1 & 1 & 1 \end{bmatrix}$$

be the canonical generator matrix of the code given in example 1.1.1. The $2^3$ code words x of the binary code can be obtained from x = aG with $a = a_1 a_2 a_3$, $a_i \in Z_2$, $1 \le i \le 3$. We have the set of $a = a_1 a_2 a_3$ which correspond to the message symbols which is as follows:

[0 0 0], [1 0 0], [0 1 0], [0 0 1],
[1 1 0], [1 0 1], [0 1 1] and [1 1 1].

$$\begin{aligned} x &= \begin{bmatrix} 0 & 0 & 0 \end{bmatrix} \begin{bmatrix} 1 & 0 & 0 & 0 & 1 & 0 & 0 \\ 0 & 1 & 0 & 1 & 0 & 0 & 0 \\ 0 & 0 & 1 & 0 & 1 & 1 & 1 \end{bmatrix} \\ &= \begin{bmatrix} 0 & 0 & 0 & 0 & 0 & 0 & 0 \end{bmatrix} \end{aligned}$$

$$x = \begin{bmatrix}1 & 0 & 0\end{bmatrix}\begin{bmatrix}1 & 0 & 0 & 0 & 1 & 0 & 0\\0 & 1 & 0 & 1 & 0 & 0 & 0\\0 & 0 & 1 & 0 & 1 & 1 & 1\end{bmatrix}$$
$$= \begin{bmatrix}1 & 0 & 0 & 0 & 1 & 0 & 0\end{bmatrix}$$

$$x = \begin{bmatrix}0 & 1 & 0\end{bmatrix}\begin{bmatrix}1 & 0 & 0 & 0 & 1 & 0 & 0\\0 & 1 & 0 & 1 & 0 & 0 & 0\\0 & 0 & 1 & 0 & 1 & 1 & 1\end{bmatrix}$$
$$= \begin{bmatrix}0 & 1 & 0 & 1 & 0 & 0 & 0\end{bmatrix}$$

$$x = \begin{bmatrix}0 & 0 & 1\end{bmatrix}\begin{bmatrix}1 & 0 & 0 & 0 & 1 & 0 & 0\\0 & 1 & 0 & 1 & 0 & 0 & 0\\0 & 0 & 1 & 0 & 1 & 1 & 1\end{bmatrix}$$
$$= \begin{bmatrix}0 & 0 & 1 & 0 & 1 & 1 & 1\end{bmatrix}$$

$$x = \begin{bmatrix}1 & 1 & 0\end{bmatrix}\begin{bmatrix}1 & 0 & 0 & 0 & 1 & 0 & 0\\0 & 1 & 0 & 1 & 0 & 0 & 0\\0 & 0 & 1 & 0 & 1 & 1 & 1\end{bmatrix}$$
$$= \begin{bmatrix}1 & 1 & 0 & 1 & 1 & 0 & 0\end{bmatrix}$$

$$x = \begin{bmatrix}1 & 0 & 1\end{bmatrix}\begin{bmatrix}1 & 0 & 0 & 0 & 1 & 0 & 0\\0 & 1 & 0 & 1 & 0 & 0 & 0\\0 & 0 & 1 & 0 & 1 & 1 & 1\end{bmatrix}$$
$$= \begin{bmatrix}1 & 0 & 1 & 0 & 0 & 1 & 1\end{bmatrix}$$

$$x = \begin{bmatrix}0 & 1 & 1\end{bmatrix}\begin{bmatrix}1 & 0 & 0 & 0 & 1 & 0 & 0\\0 & 1 & 0 & 1 & 0 & 0 & 0\\0 & 0 & 1 & 0 & 1 & 1 & 1\end{bmatrix}$$
$$= \begin{bmatrix}0 & 1 & 1 & 1 & 1 & 1 & 1\end{bmatrix}$$

$$x = \begin{bmatrix}1 & 1 & 1\end{bmatrix}\begin{bmatrix}1 & 0 & 0 & 0 & 1 & 0 & 0 \\ 0 & 1 & 0 & 1 & 0 & 0 & 0 \\ 0 & 0 & 1 & 0 & 1 & 1 & 1\end{bmatrix}$$

$$= \begin{bmatrix}1 & 1 & 1 & 1 & 0 & 1 & 1\end{bmatrix}.$$

The set of codes words generated by this G are

(0 0 0 0 0 0 0), (1 0 0 0 1 0 0), (0 1 0 1 0 0 0), (0 0 1 0 1 1 1), (1 1 0 1 1 0 0), (1 0 1 0 0 1 1), (0 1 1 1 1 1 1) and (1 1 1 1 0 1 1).

The corresponding parity check matrix H obtained from this G is given by

$$H = \begin{bmatrix}0 & 1 & 0 & 1 & 0 & 0 & 0 \\ 1 & 0 & 1 & 0 & 1 & 0 & 0 \\ 0 & 0 & 1 & 0 & 0 & 1 & 0 \\ 0 & 0 & 1 & 0 & 0 & 0 & 1\end{bmatrix}.$$

Now

$$GH^T = \begin{bmatrix}1 & 0 & 0 & 0 & 1 & 0 & 0 \\ 0 & 1 & 0 & 1 & 0 & 0 & 0 \\ 0 & 0 & 1 & 0 & 1 & 1 & 1\end{bmatrix}\begin{bmatrix}0 & 1 & 0 & 0 \\ 1 & 0 & 0 & 0 \\ 0 & 1 & 1 & 1 \\ 1 & 0 & 0 & 0 \\ 0 & 1 & 0 & 0 \\ 0 & 0 & 1 & 0 \\ 0 & 0 & 0 & 1\end{bmatrix}$$

$$= \begin{bmatrix}0 & 0 & 0 & 0 \\ 0 & 0 & 0 & 0 \\ 0 & 0 & 0 & 0\end{bmatrix}.$$

We recall just the definition of Hamming distance and Hamming weight between two vectors. This notion is applied to codes to find errors between the sent message and the received message. As finding error in the received message happens to be one of the difficult problems more so is the correction of errors and retrieving the correct message from the received message.

**DEFINITION 1.1.3:** *The Hamming distance d(x, y) between two vectors* $x = x_1 x_2 \ldots x_n$ *and* $y = y_1 y_2 \ldots y_n$ *in* $F_q^n$ *is the number of coordinates in which x and y differ. The Hamming weight* $\omega(x)$ *of a vector* $x = x_1 x_2 \ldots x_n$ *in* $F_q^n$ *is the number of non zero co ordinates in* $x_i$. *In short* $\omega(x) = d(x, 0)$.

We just illustrate this by a simple example.
Suppose x = [1 0 1 1 1 1 0] and y ∈ [0 1 1 1 1 0 1 ] belong to $F_2^7$ then D(x, y) = (x ~ y) = (1 0 1 1 1 1 0) ~ (0 1 1 1 1 0 1) = (1~0, 0~1, 1~1, 1~1, 1~1, 1~0, 0~1) = (1 1 0 0 0 1 1) = 4. Now Hamming weight ω of x is ω(x) = d(x, 0) = 5 and ω(y) = d(y, 0) = 5.

**DEFINITION 1.1.4***:* *Let C be any linear code then the minimum distance* $d_{min}$ *of a linear code C is given as*

$$d_{\min} = \min_{\substack{u,v \in C \\ u \neq v}} d(u,v) .$$

*For linear codes we have*

$$d(u, v) = d(u - v, 0) = \omega(u - v).$$

Thus it is easily seen minimum distance of C is equal to the least weight of all non zero code words. A general code C of length n with k message symbols is denoted by C(n, k) or by a binary (n, k) code. Thus a parity check code is a binary (n, n – 1) code and a repetition code is a binary (n, 1) code.

If $H = (A, I_{n-k})$ be a parity check matrix in the standard form then $G = (I_k, -A^T)$ is the canonical generator matrix of the linear (n, k) code.

The check equations $(A, I_{n-k})\, x^T = (0)$ yield

$$\begin{bmatrix} x_{k+1} \\ x_{k+2} \\ \vdots \\ x_n \end{bmatrix} = -A \begin{bmatrix} x_1 \\ x_2 \\ \vdots \\ x_k \end{bmatrix} = -A \begin{bmatrix} a_1 \\ a_2 \\ \vdots \\ a_k \end{bmatrix}.$$

Thus we obtain

$$\begin{bmatrix} x_1 \\ x_2 \\ \vdots \\ x_n \end{bmatrix} = \begin{bmatrix} I_k \\ -A \end{bmatrix} \begin{bmatrix} a_1 \\ a_2 \\ \vdots \\ a_k \end{bmatrix}.$$

We transpose and denote this equation as

$$(x_1\ x_2\ \ldots\ x_n) = (a_1\ a_2\ \ldots\ a_k)\ (I_k, -A^7)$$
$$= (a_1\ a_2\ \ldots\ a_k)\ G.$$

We have just seen that minimum distance

$$d_{\min} = \min_{\substack{u,v \in C \\ u \neq v}} d(u,v).$$

If d is the minimum distance of a linear code C then the linear code of length n, dimension k and minimum distance d is called an (n, k, d) code.

Now having sent a message or vector x and if y is the received message or vector a simple decoding rule is to find the code word closest to x with respect to Hamming distance, i.e., one chooses an error vector e with the least weight. The decoding method is called "nearest neighbour decoding" and amounts to comparing y with all $q^k$ code words and choosing the closest among them. The nearest neighbour decoding is the maximum likelihood decoding if the probability p for correct transmission is > ½.

Obviously before, this procedure is impossible for large k but with the advent of computers one can easily run a program in few seconds and arrive at the result.

We recall the definition of sphere of radius r. The set $S_r(x) = \{y \in F_q^n \,/\, d(x, y) \leq r\}$ is called the sphere of radius r about $x \in F_q^n$.

In decoding we distinguish between the detection and the correction of error. We can say a code can correct t errors and can detect t + s, s ≥ 0 errors, if the structure of the code makes it possible to correct up to t errors and to detect t + j, 0 < j ≤ s errors which occurred during transmission over a channel.

A mathematical criteria for this, given in the linear code is ; A linear code C with minimum distance $d_{min}$ can correct upto t errors and can detect t + j, 0 < j ≤ s, errors if and only if zt + s ≤ $d_{min}$ or equivalently we can say "A linear code C with minimum distance d can correct t errors if and only if $t = \left[\frac{(d-1)}{2}\right]$. The real problem of coding theory is not merely to minimize errors but to do so without reducing the transmission rate unnecessarily. Errors can be corrected by lengthening the code blocks, but this reduces the number of message symbols that can be sent per second. To maximize the transmission rate we want code blocks which are numerous enough to encode a given message alphabet, but at the same time no longer than is necessary to achieve a given Hamming distance. One of the main problems of coding theory is "Given block length n and Hamming distance d, find the maximum number, A(n, d) of binary blocks of length n which are at distances ≥ d from each other".

Let $u = (u_1, u_2, \ldots, u_n)$ and $v = (v_1, v_2, \ldots, v_n)$ be vectors in $F_q^n$ and let $u.v = u_1v_1 + u_2v_2 + \ldots + u_nv_n$ denote the dot product of u and v over $F_q^n$. If u.v = 0 then u and v are called orthogonal.

**DEFINITION 1.1.5:** *Let C be a linear (n, k) code over $F_q$. The dual(or orthogonal)code $C^{\perp}$ = {u | u.v = 0 for all v ∈ C}, u ∈ $F_q^n$. If C is a k-dimensional subspace of the n-dimensional vector space $F_q^n$ the orthogonal complement is of dimension n – k and an (n, n – k) code. It can be shown that if the code C has a generator matrix G and parity check matrix H then $C^{\perp}$ has generator matrix H and parity check matrix G.*

Orthogonality of two codes can be expressed by $GH^T = (0)$.

***Example 1.1 .5:*** Let us consider the parity check matrix H of a (7, 3) code where

$$H = \begin{bmatrix} 1 & 0 & 0 & 1 & 0 & 0 & 0 \\ 0 & 0 & 1 & 0 & 1 & 0 & 0 \\ 1 & 1 & 0 & 0 & 0 & 1 & 0 \\ 1 & 0 & 1 & 0 & 0 & 0 & 1 \end{bmatrix}.$$

The code words got using H are as follows

$$\begin{matrix} 0\ 0\ 0\ 0\ 0\ 0\ 0 & \\ 1\ 0\ 0\ 1\ 0\ 1\ 1 & 0\ 1\ 1\ 0\ 1\ 1\ 1 \\ 0\ 1\ 0\ 0\ 0\ 1\ 0 & 1\ 0\ 1\ 1\ 1\ 1\ 0 \\ 0\ 0\ 1\ 0\ 1\ 0\ 1 & 1\ 1\ 1\ 1\ 1\ 0\ 0 \\ 1\ 1\ 0\ 1\ 0\ 0\ 1 & \end{matrix}.$$

Now for the orthogonal code the parity check matrix H of the code happens to be generator matrix,

$$G = \begin{bmatrix} 1 & 0 & 0 & 1 & 0 & 0 & 0 \\ 0 & 0 & 1 & 0 & 1 & 0 & 0 \\ 1 & 1 & 0 & 0 & 0 & 1 & 0 \\ 1 & 0 & 1 & 0 & 0 & 0 & 1 \end{bmatrix}.$$

$$x = \begin{bmatrix} 0 \\ 0 \\ 0 \\ 0 \end{bmatrix}^T \begin{bmatrix} 1 & 0 & 0 & 1 & 0 & 0 & 0 \\ 0 & 0 & 1 & 0 & 1 & 0 & 0 \\ 1 & 1 & 0 & 0 & 0 & 1 & 0 \\ 1 & 0 & 1 & 0 & 0 & 0 & 1 \end{bmatrix} = [0\ 0\ 0\ 0\ 0\ 0\ 0].$$

$$x = \begin{bmatrix} 1 \\ 0 \\ 0 \\ 0 \end{bmatrix}^{T} \begin{bmatrix} 1 & 0 & 0 & 1 & 0 & 0 & 0 \\ 0 & 0 & 1 & 0 & 1 & 0 & 0 \\ 1 & 1 & 0 & 0 & 0 & 1 & 0 \\ 1 & 0 & 1 & 0 & 0 & 0 & 1 \end{bmatrix} = [1\ 0\ 0\ 1\ 0\ 0\ 0].$$

$$x = \begin{bmatrix} 0 & 1 & 0 & 0 \end{bmatrix} \begin{bmatrix} 1 & 0 & 0 & 1 & 0 & 0 & 0 \\ 0 & 0 & 1 & 0 & 1 & 0 & 0 \\ 1 & 1 & 0 & 0 & 0 & 1 & 0 \\ 1 & 0 & 1 & 0 & 0 & 0 & 1 \end{bmatrix} = [0\ 0\ 1\ 0\ 1\ 0\ 0]$$

$$x = \begin{bmatrix} 0 & 0 & 1 & 0 \end{bmatrix} \begin{bmatrix} 1 & 0 & 0 & 1 & 0 & 0 & 0 \\ 0 & 0 & 1 & 0 & 1 & 0 & 0 \\ 1 & 1 & 0 & 0 & 0 & 1 & 0 \\ 1 & 0 & 1 & 0 & 0 & 0 & 1 \end{bmatrix} = [1\ 1\ 0\ 0\ 0\ 1\ 0]$$

$$x = \begin{bmatrix} 0 & 0 & 0 & 1 \end{bmatrix} \begin{bmatrix} 1 & 0 & 0 & 1 & 0 & 0 & 0 \\ 0 & 0 & 1 & 0 & 1 & 0 & 0 \\ 1 & 1 & 0 & 0 & 0 & 1 & 0 \\ 1 & 0 & 1 & 0 & 0 & 0 & 1 \end{bmatrix} = [1\ 0\ 1\ 0\ 0\ 0\ 1]$$

$$x = \begin{bmatrix} 1 & 1 & 0 & 0 \end{bmatrix} \begin{bmatrix} 1 & 0 & 0 & 1 & 0 & 0 & 0 \\ 0 & 0 & 1 & 0 & 1 & 0 & 0 \\ 1 & 1 & 0 & 0 & 0 & 1 & 0 \\ 1 & 0 & 1 & 0 & 0 & 0 & 1 \end{bmatrix} = [1\ 0\ 1\ 1\ 1\ 0\ 0]$$

$$x = \begin{bmatrix} 1 & 0 & 1 & 0 \end{bmatrix} \begin{bmatrix} 1 & 0 & 0 & 1 & 0 & 0 & 0 \\ 0 & 0 & 1 & 0 & 1 & 0 & 0 \\ 1 & 1 & 0 & 0 & 0 & 1 & 0 \\ 1 & 0 & 1 & 0 & 0 & 0 & 1 \end{bmatrix} = [0\ 1\ 0\ 1\ 0\ 1\ 0]$$

$$x=\begin{bmatrix}1 & 0 & 0 & 1\end{bmatrix}\begin{bmatrix}1 & 0 & 0 & 1 & 0 & 0 & 0\\0 & 0 & 1 & 0 & 1 & 0 & 0\\1 & 1 & 0 & 0 & 0 & 1 & 0\\1 & 0 & 1 & 0 & 0 & 0 & 1\end{bmatrix}=[0\ 0\ 1\ 1\ 0\ 0\ 1]$$

$$x=\begin{bmatrix}0 & 1 & 1 & 0\end{bmatrix}\begin{bmatrix}1 & 0 & 0 & 1 & 0 & 0 & 0\\0 & 0 & 1 & 0 & 1 & 0 & 0\\1 & 1 & 0 & 0 & 0 & 1 & 0\\1 & 0 & 1 & 0 & 0 & 0 & 1\end{bmatrix}=[1\ 1\ 1\ 0\ 1\ 1\ 0]$$

$$x=\begin{bmatrix}0 & 1 & 0 & 1\end{bmatrix}\begin{bmatrix}1 & 0 & 0 & 1 & 0 & 0 & 0\\0 & 0 & 1 & 0 & 1 & 0 & 0\\1 & 1 & 0 & 0 & 0 & 1 & 0\\1 & 0 & 1 & 0 & 0 & 0 & 1\end{bmatrix}=[1\ 0\ 0\ 0\ 1\ 0\ 1]$$

$$x=\begin{bmatrix}0 & 0 & 1 & 1\end{bmatrix}\begin{bmatrix}1 & 0 & 0 & 1 & 0 & 0 & 0\\0 & 0 & 1 & 0 & 1 & 0 & 0\\1 & 1 & 0 & 0 & 0 & 1 & 0\\1 & 0 & 1 & 0 & 0 & 0 & 1\end{bmatrix}=[0\ 1\ 1\ 0\ 0\ 1\ 1]$$

$$x=\begin{bmatrix}1 & 1 & 1 & 0\end{bmatrix}\begin{bmatrix}1 & 0 & 0 & 1 & 0 & 0 & 0\\0 & 0 & 1 & 0 & 1 & 0 & 0\\1 & 1 & 0 & 0 & 0 & 1 & 0\\1 & 0 & 1 & 0 & 0 & 0 & 1\end{bmatrix}=[0\ 1\ 1\ 1\ 1\ 1\ 0]$$

$$x=\begin{bmatrix}1 & 1 & 0 & 1\end{bmatrix}\begin{bmatrix}1 & 0 & 0 & 1 & 0 & 0 & 0\\0 & 0 & 1 & 0 & 1 & 0 & 0\\1 & 1 & 0 & 0 & 0 & 1 & 0\\1 & 0 & 1 & 0 & 0 & 0 & 1\end{bmatrix}=[0\ 0\ 0\ 1\ 1\ 0\ 1]$$

$$x = \begin{bmatrix}1 & 0 & 1 & 1\end{bmatrix}\begin{bmatrix}1 & 0 & 0 & 1 & 0 & 0 & 0\\ 0 & 0 & 1 & 0 & 1 & 0 & 0\\ 1 & 1 & 0 & 0 & 0 & 1 & 0\\ 1 & 0 & 1 & 0 & 0 & 0 & 1\end{bmatrix} = [1\ 1\ 1\ 1\ 0\ 1\ 1]$$

$$x = \begin{bmatrix}0 & 1 & 1 & 1\end{bmatrix}\begin{bmatrix}1 & 0 & 0 & 1 & 0 & 0 & 0\\ 0 & 0 & 1 & 0 & 1 & 0 & 0\\ 1 & 1 & 0 & 0 & 0 & 1 & 0\\ 1 & 0 & 1 & 0 & 0 & 0 & 1\end{bmatrix} = [0\ 1\ 0\ 0\ 1\ 1\ 1]$$

$$x = \begin{bmatrix}1 & 1 & 1 & 1\end{bmatrix}\begin{bmatrix}1 & 0 & 0 & 1 & 0 & 0 & 0\\ 0 & 0 & 1 & 0 & 1 & 0 & 0\\ 1 & 1 & 0 & 0 & 0 & 1 & 0\\ 1 & 0 & 1 & 0 & 0 & 0 & 1\end{bmatrix} = [1\ 1\ 0\ 1\ 1\ 1\ 1]\ .$$

The code words of C(7, 4) i.e., the orthogonal code of C(7, 3) are

{(0 0 0 0 0 0 0), (1 0 0 1 0 0 0), (0 0 1 0 1 0 0), (1 1 0 0 0 1 0), (1 0 1 0 0 0 1), (1 0 1 1 1 0 0), (0 1 0 1 0 1 0), (0 0 1 1 0 0 1), (1 1 1 0 1 1 0), (1 0 0 0 1 0 1), (0 1 1 0 0 1 1), (0 1 1 1 1 1 0), (0 0 0 1 1 0 1), (1 1 1 1 0 1 1), (0 1 0 0 1 1 1), (1 1 0 1 1 1 1)}

Thus we have found the orthogonal code for the given code. Now we just recall the definition of the cosets of a code C.

**DEFINITION 1.1.6:** *For $a \in F_q^n$ we have $a + C = \{a + x \,/ x \in C\}$. Clearly each coset contains $q^k$ vectors. There is a partition of $F_q^n$ of the form $F_q^n = C \cup \{a^{(1)} + C\} \cup \{a^{(2)} + C\} \cup \ldots \cup \{a^t + C\}$ for $t = q^{n-k} - 1$. If y is a received vector then y must be an element of one of these cosets say $a^i + C$. If the code word $x^{(1)}$ has been transmitted then the error vector*

$$e = y - x^{(1)} \in a^{(i)} + C - x^{(1)} = a^{(i)} + C.$$

*Now we give the decoding rule which is as follows.*

*If a vector y is received then the possible error vectors e are the vectors in the coset containing y. The most likely error is the vector* $\overline{e}$ *with minimum weight in the coset of y. Thus y is decoded as* $\overline{x} = y - \overline{e}$ *. [23, 4]*

*Now we show how to find the coset of y and describe the above method. The vector of minimum weight in a coset is called the coset leader.*

*If there are several such vectors then we arbitrarily choose one of them as coset leader. Let* $a^{(1)}, a^{(2)}, \ldots, a^{(t)}$ *be the coset leaders. We first establish the following table*

| $x^{(1)} = 0$ | $x^{(2)} = 0$ | … | $x^{(q^k)}$ | *code words in C* |
|---|---|---|---|---|
| $a^{(1)} + x^{(1)}$ | $a^{(1)} + x^{(2)}$ | … | $a^{(1)} + x^{(q^k)}$ | *other cosets* |
| ⋮ | ⋮ | | ⋮ | |
| $a^{(t)} + x^{(1)}$ | $a^{(t)} + x^{(2)}$ | … | $a^{(t)} + x^{(q^k)}$ | |
| *coset leaders* | | | | |

*If a vector y is received then we have to find y in the table. Let* $y = a^{(i)} + x^{(j)}$*; then the decoder decides that the error* $\overline{e}$ *is the coset leader* $a^{(i)}$*. Thus y is decoded as the code word* $\overline{x} = y - \overline{e} = x^{(j)}$ *. The code word* $\overline{x}$ *occurs as the first element in the column of y. The coset of y can be found by evaluating the so called syndrome.*

*Let H be parity check matrix of a linear (n, k) code. Then the vector* $S(y) = Hy^T$ *of length n–k is called syndrome of y. Clearly S(y) = (0) if and only if* $y \in C$*.*
$S(y^{(1)}) = S(y^{(2)})$ *if and only if* $y^{(1)} + C = y^{(2)} + C$ *.*

*We have the decoding algorithm as follows:*

*If* $y \in F_q^n$ *is a received vector find S(y), and the coset leader* $\overline{e}$ *with syndrome S(y). Then the most likely transmitted code word is* $\overline{x} = y - \overline{e}$ *we have* $d(\overline{x}, y) .= \min\{d(x, y)/x \in C\}$*.*

We illustrate this by the following example.

***Example 1.1.6:*** Let C be a (5, 3) code where the parity check matrix H is given by

$$H = \begin{bmatrix} 1 & 0 & 1 & 1 & 0 \\ 1 & 1 & 0 & 0 & 1 \end{bmatrix}$$

and

$$G = \begin{bmatrix} 1 & 0 & 0 & 1 & 1 \\ 0 & 1 & 0 & 0 & 1 \\ 0 & 0 & 1 & 1 & 0 \end{bmatrix}.$$

The code words of C are

{(0 0 0  0 0), (1 0 0 1 1), (0 1 0 0 1), (0 0 1 1 0), (1 1 0 1 0), (1 0 1 0 1), (0 1 1 1 1), (1 1 1 0 0)}.

The corresponding coset table is

| Message | 000 | 100 | 010 | 001 | 110 | 101 | 011 | 111 |
|---|---|---|---|---|---|---|---|---|
| code words | 00000 | 10011 | 01001 | 00110 | 11010 | 10101 | 01111 | 11100 |
| other cosets | 10000 | 00011 | 11001 | 10110 | 01010 | 00101 | 11111 | 01100 |
| | 01000 | 11011 | 00001 | 01110 | 10010 | 11101 | 00111 | 10100 |
| | 00100 | 10111 | 01101 | 00010 | 11110 | 10001 | 01011 | 11000 |

coset leaders

If y = (1 1 1 1 0) is received, then y is found in the coset with the coset leader (0 0 1 0 0)
y + (0 0 1 0 0) = (1 1 1 1 0) + (0 0 1 0 0 ) = (1 1 0 1 0) is the corresponding message.

Now with the advent of computers it is easy to find the real message or the sent word by using this decoding algorithm.

A binary code $C_m$ of length $n = 2^m - 1$, $m \geq 2$ with $m \times 2^m - 1$ parity check matrix H whose columns consists of all non zero binary vectors of length m is called a binary Hamming code.

We give example of them.

***Example 1.1.7:*** Let

$$H = \begin{bmatrix} 1 & 0 & 1 & 1 & 1 & 1 & 0 & 0 & 1 & 0 & 1 & 1 & 0 & 0 & 0 \\ 1 & 1 & 0 & 1 & 1 & 1 & 1 & 0 & 0 & 1 & 0 & 0 & 1 & 0 & 0 \\ 1 & 1 & 1 & 0 & 1 & 0 & 1 & 1 & 0 & 0 & 1 & 0 & 0 & 1 & 0 \\ 1 & 1 & 1 & 1 & 0 & 0 & 0 & 1 & 1 & 1 & 0 & 0 & 0 & 0 & 1 \end{bmatrix}$$

which gives a $C_4(15, 11, 4)$ Hamming code.

Cyclic codes are codes which have been studied extensively.

Let us consider the vector space $F_q^n$ over $F_q$. The mapping

$$Z: F_q^n \to F_q^n$$

where Z is a linear mapping called a "cyclic shift" if $Z(a_0, a_1, \ldots, a_{n-1}) = (a_{n-1}, a_0, \ldots, a_{n-2})$

$A = (F_q[x], +, ., .)$ is a linear algebra in a vector space over $F_q$. We define a subspace $V_n$ of this vector space by

$$\begin{aligned} V_n &= \{v \in F_q[x] \,/\, \text{degree } v < n\} \\ &= \{v_0 + v_1x + v_2x^2 + \ldots + v_{n-1}x^{n-1} \,/\, v_i \in F_q;\ 0 \le i \le n-1\}. \end{aligned}$$

We see that $V_n \cong F_q^n$ as both are vector spaces defined over the same field $F_q$. Let Γ be an isomorphism

$\Gamma(v_0, v_1, \ldots, v_{n-1}) \to \{v_0 + v_1x + v_2x^2 + \ldots + v_{n-1}x^{n-1}\}$.

w: $F_q^n \hookrightarrow F_q[x] / x^n - 1$

i.e., $w(v_0, v_1, \ldots, v_{n-1}) = v_0 + v_1x + \ldots + v_{n-1}x^{n-1}$.

Now we proceed onto define the notion of a cyclic code.

**DEFINITION 1.1.7:** *A k-dimensional subspace C of $F_q^n$ is called a cyclic code if Z(v) $\in$ C for all v $\in$ C that is v = $v_0, v_1, \ldots, v_{n-1}$ $\in$ C implies $(v_{n-1}, v_0, \ldots, v_{n-2}) \in C$ for v $\in F_q^n$.*

We just give an example of a cyclic code.

***Example 1.1.8:*** Let C $\subseteq F_2^7$ be defined by the generator matrix

$$G = \begin{bmatrix} 1 & 1 & 1 & 0 & 1 & 0 & 0 \\ 0 & 1 & 1 & 1 & 0 & 1 & 0 \\ 0 & 0 & 1 & 1 & 1 & 0 & 1 \end{bmatrix} = \begin{bmatrix} g^{(1)} \\ g^{(2)} \\ g^{(3)} \end{bmatrix}.$$

The code words generated by G are {(0 0 0 0 0 0 0), (1 1 1 0 1 0 0), (0 1 1 1 0 1 0), (0 0 1 1 1 0 1), (1 0 0 1 1 1 0), (1 1 0 1 0 0 1), (0 1 0 0 1 1 1), (1 0 1 0 0 1 1)}.

Clearly one can check the collection of all code words in C satisfies the rule if $(a_0 \ldots a_5) \in C$ then $(a_5\ a_0 \ldots a_4) \in C$ i.e., the codes are cyclic. Thus we get a cyclic code.

Now we see how the code words of the Hamming codes looks like.

***Example 1.1.9:*** Let

$$H = \begin{bmatrix} 1 & 0 & 0 & 1 & 1 & 0 & 1 \\ 0 & 1 & 0 & 1 & 0 & 1 & 1 \\ 0 & 0 & 1 & 0 & 1 & 1 & 1 \end{bmatrix}$$

be the parity check matrix of the Hamming (7, 4) code.

Now we can obtain the elements of a Hamming(7,4) code.

We proceed on to define parity check matrix of a cyclic code given by a polynomial matrix equation given by defining the generator polynomial and the parity check polynomial.

**DEFINITION 1.1.8:** *A linear code C in $V_n = \{v_0 + v_1x + \ldots + v_{n-1}x^{n-1} \mid v_i \in F_q, 0 \leq i \leq n-1\}$ is cyclic if and only if C is a principal ideal generated by $g \in C$.*

*The polynomial g in C can be assumed to be monic. Suppose in addition that $g / x^n - 1$ then g is uniquely determined and is called the generator polynomial of C. The elements of C are called code words, code polynomials or code vectors.*

*Let $g = g_0 + g_1x + \ldots + g_mx^m \in V_n$, $g / x^n - 1$ and deg $g = m < n$. Let C be a linear (n, k) code, with $k = n - m$ defined by the generator matrix,*

$$G = \begin{bmatrix} g_0 & g_1 & \cdots & g_m & 0 & \cdots & 0 \\ 0 & g_0 & \cdots & g_{m-1} & g_m & \cdots & 0 \\ \vdots & \vdots & & & & & \\ 0 & 0 & & g_0 & g_1 & & g_m \end{bmatrix} = \begin{bmatrix} g \\ xg \\ \\ x^{k-1}g \end{bmatrix}.$$

*Then C is cyclic. The rows of G are linearly independent and rank G = k, the dimension of C.*

***Example 1.1.10:*** Let $g = x^3 + x^2 + 1$ be the generator polynomial having a generator matrix of the cyclic(7,4) code with generator matrix

$$G = \begin{bmatrix} 1 & 0 & 1 & 1 & 0 & 0 & 0 \\ 0 & 1 & 0 & 1 & 1 & 0 & 0 \\ 0 & 0 & 1 & 0 & 1 & 1 & 0 \\ 0 & 0 & 0 & 1 & 0 & 1 & 1 \end{bmatrix}.$$

The codes words associated with the generator matrix is

0000000, 1011000, 0101100, 0010110, 0001011, 1110100, 1001110, 1010011, 0111010, 0100111, 0011101, 1100010, 1111111, 1000101, 0110001, 1101001.

The parity check polynomial is defined to be

$$h = \frac{x^7 - 1}{g}$$

$$h = \frac{x^7 - 1}{x^3 + x^2 + 1} = x^4 + x^3 + x^2 + 1.$$

$$\text{If } \frac{x^n - 1}{g} = h_0 + h_1x + \ldots + h_kx^k.$$

the parity check matrix H related with the generator polynomial g is given by

$$H = \begin{bmatrix} 0 & \dots & 0 & h_k & \dots & h_1 & h_0 \\ 0 & \dots & h_k & h_{k-1} & & h_0 & 0 \\ \vdots & \vdots & \vdots & \vdots & & & \vdots \\ h_k & \dots & h_1 & h_0 & & \dots & 0 \end{bmatrix}.$$

For the generator polynomial $g = x^3 + x^2 + 1$ the parity check matrix

$$H = \begin{bmatrix} 0 & 0 & 1 & 1 & 1 & 0 & 1 \\ 0 & 1 & 1 & 1 & 0 & 1 & 0 \\ 1 & 1 & 1 & 0 & 1 & 0 & 0 \end{bmatrix}$$

where the parity check polynomial is given by $x^4 + x^3 + x^2 + 1 = \frac{x^7 - 1}{x^3 + x^2 + 1}$. It is left for the reader to verify that the parity check matrix gives the same set of cyclic codes.

We now proceed on to give yet another new method of decoding procedure using the method of best approximations.

We just recall this definition given by [4, 23, 39]. We just give the basic concepts needed to define this notion. We know that $F_q^n$ is a finite dimensional vector space over $F_q$. If we take $Z_2 = (0, 1)$ the finite field of characteristic two. $Z_2^5 = Z_2 \times Z_2 \times Z_2 \times Z_2 \times Z_2$ is a 5 dimensional vector space over $Z_2$. Infact {(1 0 0 0 0), (0 1 0 0 0), (0 0 1 0 0), (0 0 0 1 0), (0 0 0 0 1)} is a basis of $Z_2^5$. $Z_2^5$ has only $2^5 = 32$ elements in it. Let F be a field of real numbers and V a vector space over F. An inner product on V is a function which assigns to each ordered pair of vectors α, β in V a scalar $\langle \alpha / \beta \rangle$ in F in such a way that for all α, β, γ in V and for all scalars c in F.

(a) $\langle \alpha + \beta / \gamma \rangle = \langle \alpha/\gamma \rangle + \langle \beta/\gamma \rangle$

(b) $\langle c\alpha / \beta\rangle = c\langle \alpha/\beta\rangle$
(c) $\langle \beta/\alpha\rangle = \langle \alpha/\beta\rangle$
(d) $\langle \alpha/\alpha\rangle > 0$ if $\alpha \neq 0$.

On V there is an inner product which we call the standard inner product. Let $\alpha = (x_1, x_2, \ldots, x_n)$ and $\beta = (y_1, y_2, \ldots, y_n)$

$$\langle \alpha / \beta\rangle = \sum_i x_i y_i \; .$$

This is called as the standard inner product. $\langle \alpha/\alpha\rangle$ is defined as norm and it is denoted by $||\alpha||$. We have the Gram-Schmidt orthogonalization process which states that if V is a vector space endowed with an inner product and if $\beta_1, \beta_2, \ldots, \beta_n$ be any set of linearly independent vectors in V; then one may construct a set of orthogonal vectors $\alpha_1, \alpha_2, \ldots, \alpha_n$ in V such that for each $k = 1, 2, \ldots, n$ the set $\{\alpha_1, \ldots, \alpha_k\}$ is a basis for the subspace spanned by $\beta_1, \beta_2, \ldots, \beta_k$ where $\alpha_1 = \beta_1$.

$$\alpha_2 = \beta_2 - \frac{\langle \beta_1 / \alpha_1 \rangle}{|| \alpha_1 ||^2} \; \alpha_1$$

$$\alpha_3 = \beta_3 - \frac{\langle \beta_3 / \alpha_1 \rangle}{|| \alpha_1 ||^2} \; \alpha_1 - \frac{\langle \beta_3 / \alpha_2 \rangle}{|| \alpha_2 ||^2} \; \alpha_2$$

and so on.

Further it is left as an exercise for the reader to verify that if a vector $\beta$ is a linear combination of an orthogonal sequence of non-zero vectors $\alpha_1, \ldots, \alpha_m$, then $\beta$ is the particular linear combination, i.e.,

$$\beta = \sum_{k=1}^{m} \frac{\langle \beta / \alpha_k \rangle}{|| \alpha_k ||^2} \; \alpha_k \; .$$

In fact this property that will be made use of in the best approximations.

We just proceed on to give an example.

***Example 1.1.11:*** Let us consider the set of vectors $\beta_1$ = (2, 0, 3), $\beta_2$ = (–1, 0, 5) and $\beta_3$ = (1, 9, 2) in the space $R^3$ equipped with the standard inner product.
Define $\alpha_1$ = (2, 0, 3)

$$\alpha_2 = (-1, 0, 5) - \frac{\langle(-1, 0, 5)/(2, 0, 3)\rangle}{13}(2, 0, 3)$$
$$= (-1, 0, 5) - \frac{13}{13}(2, 0, 3) = (-3, 0, 2)$$
$$\alpha_3 = (1,9,2) - \frac{\langle(-1, 9, 2)/(2, 0, 3)\rangle}{13}(2, 0, 3)$$
$$- \frac{\langle(1, 9, 2)/(-3, 0, 2)\rangle}{13}(-3,0,2)$$
$$= (1,9,2) - \frac{8}{13}(2,0,3) - \frac{1}{13}(-3,0,2)$$
$$= (1,9,2) - \left(\frac{16}{13}, 0, \frac{24}{13}\right) - \left(\frac{3}{13}, 0, \frac{2}{13}\right)$$
$$= (1,9,2) - \left\{\left(\frac{16-3}{13}, 0, \frac{24+2}{13}\right)\right\}$$
$$= (1, 9, 2) - (1, 0, 2)$$
$$= (0, 9, 0).$$

Clearly the set {(2, 0, 3), (–3, 0, 2), (0, 9, 0)} is an orthogonal set of vectors.

Now we proceed on to define the notion of a best approximation to a vector β in V by vectors of a subspace W where $\beta \notin W$. Suppose W is a subspace of an inner product space V and let β be an arbitrary vector in V. The problem is to find a best possible approximation to β by vectors in W. This means we want to find a vector α for which $||\beta - \alpha||$ is as small as possible subject to the restriction that α should belong to W. To be precisely in mathematical terms: A best approximation to β by vectors in W is a vector α in W such that $||\beta - \alpha|| \leq ||\beta - \gamma||$ for every vector γ in W ; W a subspace of V.

By looking at this problem in $R^2$ or in $R^3$ one sees intuitively that a best approximation to β by vectors in W ought to be a vector $\alpha$ in W such that $\beta - \alpha$ is perpendicular (orthogonal) to W and that there ought to be exactly one such $\alpha$. These intuitive ideas are correct for some finite dimensional subspaces, but not for all infinite dimensional subspaces.

We just enumerate some of the properties related with best approximation.

Let W be a subspace of an inner product space V and let $\beta$ be a vector in V.

(i) The vector $\alpha$ in W is a best approximation to $\beta$ by vectors in W if and only if $\beta - \alpha$ is orthogonal to every vector in W.

(ii) If a best approximation to $\beta$ by vectors in W exists, it is unique.

(iii) If W is finite-dimensional and $\{\alpha_1, \alpha_2, \ldots, \alpha_n\}$ is any orthonormal basis for W, then the vector $\alpha = \sum_k \frac{\langle \beta / \alpha_k \rangle}{\|\alpha_k\|^2} \alpha_k$, where $\alpha$ is the (unique)best approximation to $\beta$ by vectors in W.

Now this notion of best approximation for the first time is used in coding theory to find the best approximated sent code after receiving a message which is not in the set of codes used. Further we use for coding theory only finite fields $F_q$. i.e., $|F_q| < \infty$ . If C is a code of length n; C is a vector space over $F_q$ and C $\cong F_q^k \subseteq F_q^n$ , k the number of message symbols in the code, i.e., C is a C(n, k) code. While defining the notion of inner product on vector spaces over finite fields we see all axiom of inner product defined over fields as reals or complex in general is not true. The main property which is not true is if $0 \neq x \in V$; the inner product of x with itself i.e., $\langle x / x \rangle = \langle x, x \rangle \neq 0$ if $x \neq 0$ is not true i.e., $\langle x / x \rangle = 0$ does not imply $x = 0$ .

To overcome this problem we define for the first time the new notion of pseudo inner product in case of vector spaces defined over finite characteristic fields [4, 23].

**DEFINITION 1.1.9:** *Let V be a vector space over a finite field $F_p$ of characteristic p, p a prime. Then the pseudo inner product on*

*V is a map $\langle , \rangle_p : V \times V \to F_p$ satisfying the following conditions.*

1. *$\langle x, x\rangle_p \geq 0$ for all $x \in V$ and $\langle x, x\rangle_p = 0$ does not in general imply $x = 0$.*
2. *$\langle x, y\rangle_p = \langle y, x\rangle_p$ for all $x, y \in V$.*
3. *$\langle x + y, z\rangle_p = \langle x, z\rangle_p + \langle y, z\rangle_p$ for all $x, y, z \in V$.*
4. *$\langle x, y + z\rangle_p = \langle x, y\rangle_p + \langle x, z\rangle_p$ for all $x, y, z \in V$.*
5. *$\langle \alpha.x, y\rangle_p = \alpha \langle x, y\rangle_p$ and*
6. *$\langle x, \beta.y\rangle_p = \beta\langle x, y\rangle_p$ for all $x, y, \in V$ and $\alpha, \beta \in F_p$.*

*Let V be a vector space over a field $F_p$ of characteristic p, p is a prime; then V is said to be a pseudo inner product space if there is a pseudo inner product $\langle , \rangle_p$ defined on V. We denote the pseudo inner product space by $(V, \langle , \rangle_p)$.*

Now using this pseudo inner product space $(V, \langle , \rangle_p)$ we proceed on to define pseudo-best approximation.

**DEFINITION 1.1.10:** *Let V be a vector space defined over the finite field $F_p$ (or $Z_p$). Let W be a subspace of V. For $\beta \in V$ and for a set of basis $\{\alpha_1, \ldots, \alpha_k\}$ of the subspace W the pseudo best approximation to $\beta$, if it exists is given by $\sum_{i=1}^{k} \langle \beta, \alpha_i \rangle_p \alpha_i$. If $\sum_{i=1}^{k} \langle \beta, \alpha_i \rangle_p \alpha_i = 0$, then we say the pseudo best approximation does not exist for this set of basis $\{\alpha_1, \alpha_2, \ldots, \alpha_k\}$. In this case we choose another set of basis for W say $\{\gamma_1, \gamma_2, \ldots, \gamma_k\}$ and calculate $\sum_{i=1}^{k} \langle \beta, \gamma_i \rangle_p \gamma_i$ and $\sum_{i=1}^{k} \langle \beta, \gamma_i \rangle_p \gamma_i$ is called a pseudo best approximation to $\beta$.*

*Note:* We need to see the difference even in defining our pseudo best approximation with the definition of the best approximation. Secondly as we aim to use it in coding theory and most of our linear codes take only their values from the field of characteristic two we do not need $\langle x, x\rangle$ or the norm to

be divided by the pseudo inner product in the summation of finding the pseudo best approximation.

Now first we illustrate the pseudo inner product by an example.

***Example 1.1.12:*** Let V = $Z_2 \times Z_2 \times Z_2 \times Z_2$ be a vector space over $Z_2$. Define $\langle , \rangle_p$ to be the standard pseudo inner product on V; so if x = (1 0 1 1) and y = (1 1 1 1) are in V then the pseudo inner product of

$$\langle x, y\rangle_p = \langle(1\ 0\ 1\ 1), (1\ 1\ 1\ 1)\rangle_p = 1 + 0 + 1 + 1 = 1.$$

Now consider

$$\langle x, x\rangle_p = \langle(1\ 0\ 1\ 1), (1\ 0\ 1\ 1)\rangle_p = 1 + 0 + 1 + 1 \neq 0$$

but

$$\langle y, y\rangle_p = \langle(1\ 1\ 1\ 1), (1\ 1\ 1\ 1)\rangle_p = 1 + 1 + 1 + 1 = 0.$$

We see clearly y ≠ 0, yet the pseudo inner product is zero.

Now having seen an example of the pseudo inner product we proceed on to illustrate by an example the notion of pseudo best approximation.

***Example 1.1.13:*** Let

$$V = Z_2^8 = \underbrace{Z_2 \times Z_2 \times \ldots \times Z_2}_{8 \text{ times}}$$

be a vector space over $Z_2$. Now
W = {0 0 0 0 0 0 0 0), (1 0 0 0 1 0 11), (0 1 0 0 1 1 0 0), (0 0 1 0 0 1 1 1), (0 0 0 1 1 1 0 1), (1 1 0 0 0 0 1 0), (0 1 1 0 1 1 1 0), (0 0 1 1 1 0 1 0), (0 1 0 1 0 1 0 0), (1 0 1 0 1 1 0 0), (1 0 0 1 0 1 1 0), (1 1 1 0 0 1 0 1), (0 1 1 1 0 0 1 1), (1 1 0 1 1 1 1 1), (1 0 1 1 0 0 0 1), (1 1 1 1 1 0 0 0)}

be a subspace of V. Choose a basis of W as B = {$\alpha_1, \alpha_2, \alpha_3, \alpha_4$} where

$$\alpha_1 = (0\ 1\ 0\ 0\ 1\ 0\ 0\ 1),$$
$$\alpha_2 = (1\ 1\ 0\ 0\ 0\ 0\ 1\ 0),$$
$$\alpha_3 = (1\ 1\ 1\ 0\ 0\ 1\ 0\ 1)$$

and

$$\alpha_4 = (1\ 1\ 1\ 1\ 1\ 0\ 0\ 0).$$

Suppose $\beta$ = (1 1 1 1 1 1 1 1) is a vector in V using pseudo best approximations find a vector in W close to $\beta$. This is given by $\alpha$ relative to the basis B of W where

$$\alpha = \sum_{k=-1}^{4} \langle \beta, \alpha_k \rangle_p \alpha_k$$

$= \langle(1\ 1\ 1\ 1\ 1\ 1\ 1\ 1), (0\ 1\ 0\ 0\ 1\ 0\ 0\ 1)\rangle_p \alpha_1 +$
$\langle(1\ 1\ 1\ 1\ 1\ 1\ 1\ 1), (1\ 1\ 0\ 0\ 0\ 0\ 1\ 0)\rangle_p \alpha_2 +$
$\langle(1\ 1\ 1\ 1\ 1\ 1\ 1\ 1), (1\ 1\ 1\ 0\ 0\ 1\ 0\ 1)\rangle_p \alpha_3 +$
$\langle(1\ 1\ 1\ 1\ 1\ 1\ 1\ 1), (1\ 1\ 1\ 1\ 1\ 0\ 0\ 0)\rangle_p \alpha_4$.

$= 1.\alpha_1 + 1.\alpha_2 + 1.\alpha_3 + 1.\alpha_4$.

$=$ (0 1 0 0 1 0 0 1) + (1 1 0 0 0 0 1 0) + (1 1 1 0 0 1 0 1) + (1 1 1 1 1 0 0 0)

$=$ (1 0 0 1 0 1 1 0) $\in$ W.

Now having illustrated how the pseudo best approximation of a vector $\beta$ in V relative to a subspace W of V is determined, now we illustrate how the approximately the nearest code word is obtained.

***Example 1.1.14:*** Let C = C(4, 2) be a code obtained from the parity check matrix

$$H = \begin{bmatrix} 1 & 0 & 1 & 0 \\ 1 & 1 & 0 & 1 \end{bmatrix}.$$

The message symbols associated with the code C are {(0, 0), (1, 0), (1, 0), (1, 1)}. The code words associated with H are C = {(0 0 0 0), (1 0 1 1), (0 1 0 1), (1 1 1 0)}. The chosen basis for C is B = {$\alpha_1$, $\alpha_2$} where $\alpha_1$ = (0 1 0 1) and $\alpha_2$ = (1 0 1 1). Suppose the received message is $\beta$ = (1 1 1 1), consider $H\beta^T$ = (0 1) $\neq$ (0) so $\beta \notin$ C. Let $\alpha$ be the pseudo best approximation to $\beta$ relative to the basis B given as

$$\alpha = \sum_{k=1}^{2} \langle \beta, \alpha_k \rangle_p \alpha_k = \langle (1\ 1\ 1\ 1), (0\ 1\ 0\ 1)\rangle_p \alpha_1 + \langle (1\ 1\ 1\ 1), (1\ 0\ 1\ 1)\rangle_p \alpha_2.$$

$$= (1\ 0\ 1\ 1)\ .$$

Thus the approximated code word is (1 0 1 1).
This method could be implemented in case of algebraic linear bicodes and in general to algebraic linear n-codes; $n \geq 3$.

Now having seen some simple properties of codes we now proceed on to recall some very basic properties about bimatrices and their generalization, n-matrices ($n \geq 3$).

## 1.2 Bimatrices and their Generalizations

In this section we recall some of the basic properties of bimatrices and their generalizations which will be useful for us in the definition of linear bicodes and linear n-codes respectively.

In this section we recall the notion of bimatrix and illustrate them with examples and define some of basic operations on them.

**DEFINITION 1.2.1:** *A bimatrix $A_B$ is defined as the union of two rectangular array of numbers $A_1$ and $A_2$ arranged into rows and columns. It is written as follows $A_B = A_1 \cup A_2$ where $A_1 \neq A_2$ with*

$$A_1 = \begin{bmatrix} a_{11}^1 & a_{12}^1 & \cdots & a_{1n}^1 \\ a_{21}^1 & a_{22}^1 & \cdots & a_{2n}^1 \\ \vdots & & & \\ a_{m1}^1 & a_{m2}^1 & \cdots & a_{mn}^1 \end{bmatrix}$$

*and*

$$A_2 = \begin{bmatrix} a^2_{11} & a^2_{12} & \cdots & a^2_{1n} \\ a^2_{21} & a^2_{22} & \cdots & a^2_{2n} \\ \vdots & & & \\ a^2_{m1} & a^2_{m2} & \cdots & a^2_{mn} \end{bmatrix}$$

*'$\cup$' is just the notational convenience (symbol) only.*

The above array is called a m by n bimatrix (written as B(m × n) since each of $A_i$ (i = 1, 2) has m rows and n columns). It is to be noted a bimatrix has no numerical value associated with it. It is only a convenient way of representing a pair of array of numbers.

*Note:* If $A_1 = A_2$ then $A_B = A_1 \cup A_2$ is not a bimatrix. A bimatrix $A_B$ is denoted by $\left(a^1_{ij}\right) \cup \left(a^2_{ij}\right)$. If both $A_1$ and $A_2$ are m × n matrices then the bimatrix $A_B$ is called the m × n rectangular bimatrix.

But we make an assumption the zero bimatrix is a union of two zero matrices even if $A_1$ and $A_2$ are one and the same; i.e., $A_1 = A_2 = (0)$.

***Example 1.2.1:*** The following are bimatrices

i. $$A_B = \begin{bmatrix} 3 & 0 & 1 \\ -1 & 2 & 1 \end{bmatrix} \cup \begin{bmatrix} 0 & 2 & -1 \\ 1 & 1 & 0 \end{bmatrix}$$

is a rectangular 2 × 3 bimatrix.

ii. $$A'_B = \begin{bmatrix} 3 \\ 1 \\ 2 \end{bmatrix} \cup \begin{bmatrix} 0 \\ -1 \\ 0 \end{bmatrix}$$

is a column bimatrix.

iii. $A''_B = (3, -2, 0, 1, 1) \cup (1, 1, -1, 1, 2)$

is a row bimatrix.

In a bimatrix $A_B = A_1 \cup A_2$ if both $A_1$ and $A_2$ are $m \times n$ rectangular matrices then the bimatrix $A_B$ is called the rectangular $m \times n$ bimatrix.

**DEFINITION 1.2.2:** *Let $A_B = A_1 \cup A_2$ be a bimatrix. If both $A_1$ and $A_2$ are square matrices then $A_B$ is called the square bimatrix.*

*If one of the matrices in the bimatrix $A_B = A_1 \cup A_2$ is a square matrix and other is a rectangular matrix or if both $A_1$ and $A_2$ are rectangular matrices say $m_1 \times n_1$ and $m_2 \times n_2$ with $m_1 \neq m_2$ or $n_1 \neq n_2$ then we say $A_B$ is a mixed bimatrix.*

The following are examples of a square bimatrix and the mixed bimatrix.

***Example 1.2.2:*** Given

$$A_B = \begin{bmatrix} 3 & 0 & 1 \\ 2 & 1 & 1 \\ -1 & 1 & 0 \end{bmatrix} \cup \begin{bmatrix} 4 & 1 & 1 \\ 2 & 1 & 0 \\ 0 & 0 & 1 \end{bmatrix}$$

is a $3 \times 3$ square bimatrix.

$$A'_B = \begin{bmatrix} 1 & 1 & 0 & 0 \\ 2 & 0 & 0 & 1 \\ 0 & 0 & 0 & 3 \\ 1 & 0 & 1 & 2 \end{bmatrix} \cup \begin{bmatrix} 2 & 0 & 0 & -1 \\ -1 & 0 & 1 & 0 \\ 0 & -1 & 0 & 3 \\ -3 & -2 & 0 & 0 \end{bmatrix}$$

is a $4 \times 4$ square bimatrix.

***Example 1.2.3:*** Let

$$A_B = \begin{bmatrix} 3 & 0 & 1 & 2 \\ 0 & 0 & 1 & 1 \\ 2 & 1 & 0 & 0 \\ 1 & 0 & 1 & 0 \end{bmatrix} \cup \begin{bmatrix} 1 & 1 & 2 \\ 0 & 2 & 1 \\ 0 & 0 & 4 \end{bmatrix}$$

then $A_B$ is a mixed square bimatrix.
Let

$$A'_B = \begin{bmatrix} 2 & 0 & 1 & 1 \\ 0 & 1 & 0 & 1 \\ -1 & 0 & 2 & 1 \end{bmatrix} \cup \begin{bmatrix} 2 & 0 \\ 4 & -3 \end{bmatrix},$$

$A'_B$ is a mixed bimatrix.

Now we proceed on to give the operations on bimatrices.

Let $A_B = A_1 \cup A_2$ and $C_B = C_1 \cup C_2$ be two bimatrices we say $A_B$ and $C_B$ are equal written as $A_B = C_B$ if and only if $A_1$ and $C_1$ are identical and $A_2$ and $C_2$ are identical i.e., $A_1 = C_1$ and $A_2 = C_2$.

If $A_B = A_1 \cup A_2$ and $C_B = C_1 \cup C_2$, we say $A_B$ is not equal to $C_B$, we write $A_B \neq C_B$ if and only if $A_1 \neq C_1$ or $A_2 \neq C_2$.

***Example 1.2.4:*** Let

$$A_B = \begin{bmatrix} 3 & 2 & 0 \\ 2 & 1 & 1 \end{bmatrix} \cup \begin{bmatrix} 0 & -1 & 2 \\ 0 & 0 & 1 \end{bmatrix}$$

and

$$C_B = \begin{bmatrix} 1 & 1 & 1 \\ 0 & 0 & 0 \end{bmatrix} \cup \begin{bmatrix} 2 & 0 & 1 \\ 1 & 0 & 2 \end{bmatrix}$$

clearly $A_B \neq C_B$. Let

$$A_B = \begin{bmatrix} 0 & 0 & 1 \\ 1 & 1 & 2 \end{bmatrix} \cup \begin{bmatrix} 0 & 4 & -2 \\ -3 & 0 & 0 \end{bmatrix}$$

$$C_B = \begin{bmatrix} 0 & 0 & 1 \\ 1 & 1 & 2 \end{bmatrix} \cup \begin{bmatrix} 0 & 0 & 0 \\ 1 & 0 & 1 \end{bmatrix}$$

clearly $A_B \neq C_B$.

If $A_B = C_B$ then we have $C_B = A_B$.

We now proceed on to define multiplication by a scalar. Given a bimatrix $A_B = A_1 \cup B_1$ and a scalar $\lambda$, the product of $\lambda$ and $A_B$ written $\lambda A_B$ is defined to be

$$\lambda A_B = \begin{bmatrix} \lambda a_{11} & \cdots & \lambda a_{1n} \\ \vdots & & \vdots \\ \lambda a_{m1} & \cdots & \lambda a_{mn} \end{bmatrix} \cup \begin{bmatrix} \lambda b_{11} & \cdots & \lambda b_{1n} \\ \vdots & & \vdots \\ \lambda b_{m1} & \cdots & \lambda b_{mn} \end{bmatrix}$$

each element of $A_1$ and $B_1$ are multiplied by $\lambda$. The product $\lambda$ $A_B$ is then another bimatrix having m rows and n columns if $A_B$ has m rows and n columns.

We write

$$\begin{aligned} \lambda A_B &= \left[\lambda a_{ij}\right] \cup \left[\lambda b_{ij}\right] \\ &= \left[a_{ij}\lambda\right] \cup \left[b_{ij}\lambda\right] \\ &= A_B \lambda. \end{aligned}$$

***Example 1.2.5:*** Let

$$A_B = \begin{bmatrix} 2 & 0 & 1 \\ 3 & 3 & -1 \end{bmatrix} \cup \begin{bmatrix} 0 & 1 & -1 \\ 2 & 1 & 0 \end{bmatrix}$$

and $\lambda = 3$ then

$$3A_B = \begin{bmatrix} 6 & 0 & 3 \\ 9 & 9 & -3 \end{bmatrix} \cup \begin{bmatrix} 0 & 3 & -3 \\ 6 & 3 & 0 \end{bmatrix}.$$

If $\lambda = -2$ for

$$\begin{aligned} A_B &= [3\ 1\ 2\ -4] \cup [0\ 1\ -1\ 0], \\ \lambda A_B &= [-6\ -2\ -4\ 8] \cup [0\ -2\ 2\ 0]. \end{aligned}$$

Let $A_B = A_1 \cup B_1$ and $C_B = A_2 \cup B_2$ be any two $m \times n$ bimatrices. The sum $D_B$ of the bimatrices $A_B$ and $C_B$ is defined as $D_B = A_B + C_B = [A_1 \cup B_1] + [A_2 \cup B_2] = (A_1 + A_2) \cup [B_2 + B_2]$; where $A_1 + A_2$ and $B_1 + B_2$ are the usual addition of matrices i.e., if

$$A_B = \left(a^1_{ij}\right) \cup \left(b^1_{ij}\right)$$

and

$$C_B = \left(a^2_{ij}\right) \cup \left(b^2_{ij}\right)$$

then

$$A_B + C_B = D_B = \left(a^1_{ij} + a^2_{ij}\right) \cup \left(b^1_{ij} + b^2_{ij}\right) \; \left(\forall ij\right).$$

If we write in detail

$$A_B = \begin{bmatrix} a^1_{11} & \cdots & a^1_{1n} \\ \vdots & & \vdots \\ a^1_{m1} & \cdots & a^1_{mn} \end{bmatrix} \cup \begin{bmatrix} b^1_{11} & \cdots & b^1_{1n} \\ \vdots & & \vdots \\ b^1_{m1} & \cdots & b^1_{mn} \end{bmatrix}$$

$$C_B = \begin{bmatrix} a^2_{11} & \cdots & a^2_{1n} \\ \vdots & & \vdots \\ a^2_{m1} & \cdots & a^2_{mn} \end{bmatrix} \cup \begin{bmatrix} b^2_{11} & \cdots & b^2_{1n} \\ \vdots & & \vdots \\ b^2_{m1} & \cdots & b^2_{mn} \end{bmatrix}$$

$$A_B + C_B =$$

$$\begin{bmatrix} a^1_{11} + a^2_{11} & \dots & a^1_{1n} + a^2_{1n} \\ \vdots & & \vdots \\ a^1_{m1} + a^2_{m1} & \dots & a^1_{mn} + a^2_{mn} \end{bmatrix} \cup \begin{bmatrix} b^1_{11} + b^2_{11} & \dots & b^1_{1n} + b^2_{1n} \\ \vdots & & \vdots \\ b^1_{m1} + b^2_{m1} & \dots & b^1_{mn} + b^2_{mn} \end{bmatrix}.$$

The expression is abbreviated to

$$\begin{aligned} D_B &= A_B + C_B \\ &= (A_1 \cup B_1) + (A_2 \cup B_2) \\ &= (A_1 + A_2) \cup (B_1 + B_2). \end{aligned}$$

Thus two bimatrices are added by adding the corresponding elements only when compatibility of usual matrix addition exists.

*Note*: If $A_B = A^1 \cup A^2$ be a bimatrix we call $A^1$ and $A^2$ as the components of $A_B$ or component matrices of the bimatrix $A_B$.

***Example 1.2.6:***

(i) Let

$$A_B = \begin{bmatrix} 3 & 1 & 1 \\ -1 & 0 & 2 \end{bmatrix} \cup \begin{bmatrix} 4 & 0 & -1 \\ 0 & 1 & -2 \end{bmatrix}$$

and

$$C_B = \begin{bmatrix} -1 & 0 & 1 \\ 2 & 2 & -1 \end{bmatrix} \cup \begin{bmatrix} 3 & 3 & 1 \\ 0 & 2 & -1 \end{bmatrix},$$

then,

$$D_B = A_B + C_B$$

$$= \begin{bmatrix} 3 & 1 & 1 \\ -1 & 0 & 2 \end{bmatrix} + \begin{bmatrix} -1 & 0 & 1 \\ 2 & 2 & -1 \end{bmatrix} \cup$$

$$\begin{bmatrix} 4 & 0 & -1 \\ 0 & 1 & -2 \end{bmatrix} + \begin{bmatrix} 3 & 3 & 1 \\ 0 & 2 & -1 \end{bmatrix}$$

$$= \begin{bmatrix} 2 & 1 & 2 \\ 1 & 2 & 1 \end{bmatrix} \cup \begin{bmatrix} 7 & 3 & 0 \\ 0 & 3 & -3 \end{bmatrix}.$$

(ii) Let

$$A_B = (3\ 2\ -1\ 0\ 1) \cup (0\ 1\ 1\ 0\ -1)$$

and

$$C_B = (1\ 1\ 1\ 1\ 1) \cup (5\ -1\ 2\ 0\ 3),$$

$$A_B + C_B = (4\ 3\ 0\ 1\ 2) \cup (5\ 0\ 3\ 0\ 2).$$

***Example 1.2.7:*** Let

$$A_B = \begin{bmatrix} 6 & -1 \\ 2 & 2 \\ 1 & -1 \end{bmatrix} \cup \begin{bmatrix} 3 & 1 \\ 0 & 2 \\ -1 & 3 \end{bmatrix}$$

and

$$C_B = \begin{bmatrix} 2 & -4 \\ 4 & -1 \\ 3 & 0 \end{bmatrix} \cup \begin{bmatrix} 1 & 4 \\ 2 & 1 \\ 3 & 1 \end{bmatrix}.$$

$$A_B + A_B = \begin{bmatrix} 12 & -2 \\ 4 & 4 \\ 2 & -2 \end{bmatrix} \cup \begin{bmatrix} 6 & 2 \\ 0 & 4 \\ -2 & 6 \end{bmatrix} = 2A_B$$

$$C_B + C_B = \begin{bmatrix} 4 & -8 \\ 8 & -2 \\ 6 & 0 \end{bmatrix} \cup \begin{bmatrix} 2 & 8 \\ 4 & 2 \\ 6 & 2 \end{bmatrix} = 2C_B.$$

Similarly we can add

$$A_B + A_B + A_B = 3A_B = \begin{bmatrix} 18 & -3 \\ 6 & 6 \\ 3 & -3 \end{bmatrix} \cup \begin{bmatrix} 9 & 3 \\ 0 & 6 \\ -3 & 9 \end{bmatrix}.$$

*Note:* Addition of bimatrices are defined if and only if both the bimatrices are $m \times n$ bimatrices.

Let

$$A_B = \begin{bmatrix} 3 & 0 & 1 \\ 1 & 2 & 0 \end{bmatrix} \cup \begin{bmatrix} 1 & 1 & 1 \\ 0 & 2 & -1 \end{bmatrix}$$

and

$$C_B = \begin{bmatrix} 3 & 1 \\ 2 & 1 \\ 0 & 0 \end{bmatrix} \cup \begin{bmatrix} 1 & 1 \\ 2 & -1 \\ 3 & 0 \end{bmatrix}.$$

The addition of $A_B$ with $C_B$ is not defined for $A_B$ is a $2 \times 3$ bimatrix where as $C_B$ is a $3 \times 2$ bimatrix.

Clearly $A_B + C_B = C_B + A_B$ when both $A_B$ and $C_B$ are $m \times n$ matrices.

Also if $A_B$, $C_B$, $D_B$ be any three $m \times n$ bimatrices then $A_B + (C_B + D_B) = (A_B + C_B) + D_B$.

Subtraction is defined in terms of operations already considered for if

$$A_B = A_1 \cup A_2$$

and

$$B_B = B_1 \cup B_2$$

then

$$\begin{aligned} A_B - B_B &= A_B + (-B_B) \\ &= (A_1 \cup A_2) + (-B_1 \cup -B_2) \\ &= (A_1 - B_1) \cup (A_2 - B_2) \\ &= [A_1 + (-B_1)] \cup [A_2 + (-B_2)]. \end{aligned}$$

***Example 1.2.8:***

i. Let

$$A_B = \begin{bmatrix} 3 & 1 \\ -1 & 2 \\ 0 & 3 \end{bmatrix} \cup \begin{bmatrix} 5 & -2 \\ 1 & 1 \\ 3 & -2 \end{bmatrix}$$

and

$$B_B = \begin{bmatrix} 8 & -1 \\ 4 & 2 \\ -1 & 3 \end{bmatrix} \cup \begin{bmatrix} 9 & 2 \\ 2 & 9 \\ -1 & 1 \end{bmatrix}$$

$$A_B - B_B = A_B + (-B_B).$$

$$= \left\{ \begin{bmatrix} 3 & 1 \\ -1 & 2 \\ 0 & 3 \end{bmatrix} \cup \begin{bmatrix} 5 & -2 \\ 1 & 1 \\ 3 & -2 \end{bmatrix} \right\} + - \left\{ \begin{bmatrix} 8 & -1 \\ 4 & 2 \\ -1 & 3 \end{bmatrix} \cup \begin{bmatrix} 9 & 2 \\ 2 & 9 \\ -1 & 1 \end{bmatrix} \right\}$$

$$= \left\{ \begin{bmatrix} 3 & 1 \\ -1 & 2 \\ 0 & 3 \end{bmatrix} - \begin{bmatrix} 8 & -1 \\ 4 & 2 \\ -1 & 3 \end{bmatrix} \right\} \cup \left\{ \begin{bmatrix} 5 & -2 \\ 1 & 1 \\ 3 & -2 \end{bmatrix} - \begin{bmatrix} 9 & 2 \\ 2 & 9 \\ -1 & 1 \end{bmatrix} \right\}$$

$$= \begin{bmatrix} -5 & 2 \\ -5 & 0 \\ 1 & 0 \end{bmatrix} \cup \begin{bmatrix} 4 & -4 \\ -1 & -8 \\ 4 & -3 \end{bmatrix}.$$

ii. Let

$$A_B = (1, 2, 3, -1, 2, 1) \cup (3, -1, 2, 0, 3, 1)$$

and

$$B_B = (-1, 1, 1, 1, 1, 0) \cup (2, 0, -2, 0, 3, 0)$$

then

$$A_B + (-B_B) = (2, 1, 2, -2, 1, 1) \cup (1, -1, 4, 0, 0, 1).$$

Now we have defined addition and subtraction of bimatrices. Unlike in matrices we cannot say if we add two bimatrices the sum will be a bimatrix.

Now we proceed onto define the notion of n-matrices.

**DEFINITION 1.2.3:** *A n matrix A is defined to be the union of n rectangular array of numbers $A_1, \ldots, A_n$ arranged into rows and columns. It is written as $A = A_1 \cup \ldots \cup A_n$ where $A_i \neq A_j$ with*

$$A_i = \begin{bmatrix} a^i_{11} & a^i_{12} & \cdots & a^i_{1p} \\ a^i_{21} & a^i_{22} & \cdots & a^i_{2p} \\ \vdots & \vdots & & \vdots \\ a^i_{m1} & a^i_{m2} & \cdots & a^i_{mp} \end{bmatrix}$$

*i = 1, 2, …, n.*

'∪' is just the notional convenience (symbol) only ($n \geq 3$).

**Note**: If n = 2 we get the bimatrix.

***Example 1.2.9:*** Let

$$A = \begin{bmatrix} 3 & 1 & 0 & 1 \\ 0 & 0 & 1 & 1 \end{bmatrix} \cup \begin{bmatrix} 2 & 1 & 1 & 0 \\ 0 & 1 & 1 & 0 \end{bmatrix} \cup$$

$$\begin{bmatrix} 1 & 0 & 0 & 1 \\ 0 & 1 & 0 & 1 \end{bmatrix} \cup \begin{bmatrix} 5 & 1 & 0 & 2 \\ 7 & -1 & 0 & 3 \end{bmatrix}$$

A is a 4-matrix.

***Example 1.2.10:*** Let

$$A = A_1 \cup A_2 \cup A_3 \cup A_4 \cup A_5$$

$$= \quad [1\ 0\ 0] \cup \begin{bmatrix} 1 \\ 2 \\ -1 \\ 0 \\ 0 \end{bmatrix} \cup \begin{bmatrix} 3 & 1 & 2 \\ 0 & 1 & 1 \\ 9 & 7 & -8 \end{bmatrix}$$

$$\cup \begin{bmatrix} 2 & 1 & 3 & 5 \\ 0 & 1 & 0 & 2 \end{bmatrix} \cup \begin{bmatrix} 7 & 9 & 8 & 11 & 0 \\ 1 & -2 & 0 & 9 & 7 \\ 0 & 5 & 7 & -1 & 8 \\ -4 & -6 & 6 & 0 & 1 \end{bmatrix};$$

A is a 5-matrix. Infact A is a mixed 5-matrix.

***Example 1.2.11:*** Consider the 7-matrix

$$A = \begin{bmatrix} 2 & 0 \\ 1 & 1 \end{bmatrix} \cup \begin{bmatrix} 1 & 1 \\ 0 & 1 \end{bmatrix} \cup \begin{bmatrix} 1 & -1 \\ 0 & 4 \end{bmatrix} \cup$$

$$\begin{bmatrix} 3 & 1 \\ 1 & 0 \\ 0 & 1 \\ 2 & 1 \end{bmatrix} \cup \begin{bmatrix} 5 \\ 6 \\ 7 \\ 8 \\ -1 \\ 3 \\ 2 \end{bmatrix} \cup [3 \;\; 7 \;\; 8 \;\; 1 \;\; 0] \cup \begin{bmatrix} 2 & 1 & 0 & 0 \\ 1 & 0 & 2 & -1 \\ -1 & 1 & 0 & 0 \\ 0 & 2 & 0 & 1 \\ 2 & 0 & 0 & 1 \end{bmatrix}.$$

$= A_1 \cup A_2 \cup \ldots \cup A_7$. A is a mixed 7-matrix.

**Chapter Two**

# BICODES AND THEIR GENERALIZATIONS

This chapter has four sections. We introduce in this chapter four classes of new codes. Section 1 contains bicodes and the new class of n-codes ($n \geq 3$) and their properties. Section two deals with error correction using n-coset leaders and pseudo best n approximations. The notion of false n-matrix and pseudo false n-matrix are introduced in section three. Section four introduces 4 new classes of codes viz false n-codes, m-pseudo false n-codes, (t, t) pseudo false n-codes and (t, m) pseudo false n-codes which has lot of applications.

## 2.1 Bicodes and n-codes

In this section we introduce the notion of bicodes and n-codes and describe some of their basic properties. Now for the first time we give the use of bimatrices in the field of coding theory. Very recently we have defined the notion of bicodes [43, 44].

A bicode is a pair of sets of markers or labels built up from a finite "alphabet" which can be attached to some or all of the entities of a system of study. A mathematical bicode is a bicode with some mathematical structure. In many cases as in case of codes the mathematical structure of a bicodes refers to the fact that it forms a bivector space over finite field and in this case the bicode is linear. Let $Z_q(m) = Z_{q_1}(m_2) \cup Z_{q_2}(m_2)$. Denote the m-dimensional linear bicode in the bivector space V = $V_1 \cup V_2$ where $V_1$ is a $m_1$ dimensional vector space over $Z_p$ where $q_1$ is a power of the prime p and $V_2$ is a $m_2$ – dimensional vector space over $Z_p$ where $q_2$ is a power of the prime p, then $Z_q$ (m) = $Z_{q_1}(m_1) \cup Z_{q_2}(m_2)$ is a bicode over $Z_p$.

Thus a bicode can also be defined as a 'union' of two codes $C_1$ and $C_2$ where union is just the symbol.

For example we can have a bicode over $Z_2$.

***Example 2.1.1:*** Let C = $C_1 \cup C_2$ be a bicode over $Z_2$ given by

C = $C_1 \cup C_2$
= {( 0 0 0 0 0 0), (0 1 1 0 1 1 ), (1 1 0 1 1 0), (0 0 1 1 1 0), (1 0 0 0 1 1), (1 1 1 0 0 0), (0 1 0 1 0 1), (1 0 1 1 0 1)} $\cup$ {(0 0 0 0), (1 1 1 0), (1 0 0 1), (0 1 1 1), (0 1 0 1), (0 0 1 0), (1 1 0 0), (1 0 1 1)} over $Z_2$.

These codes are subbispaces of the bivector space over $Z_2$ of dimension (6, 4). Now the bicodes are generated by bimatrices and we have parity check bimatrices to play a role in finding the check symbols. Thus we see we have the applications of linear bialgebra / bivector spaces in the study of bicodes.

$$C (n_1 \cup n_2, k_1, k_2) = C_1 (n_1, k_1) \cup C_2 (n_2, k_2)$$

is a linear bicode if and only if both $C_1$ $(n_1, k_1)$ and $C_2$ $(n_2, k_2)$ are linear codes of length $n_1$ and $n_2$ with $k_1$ and $k_2$ message symbols respectively with entries from the same field $Z_q$. The check symbols can be obtained from the $k_1$ and $k_2$ messages in such a way that the bicode words $x = x^1 \cup x^2$ satisfy the system of linear biequations.

i.e. $$Hx^T = (0)$$
i.e. $$(H_1 \cup H_2) (x^1 \cup x^2)^T = (0) \cup (0)$$
i.e. $$H_1 (x^1)^T \cup H_2 (x^2)^T = (0) \cup (0)$$
where $H = H_1 \cup H_2$ is a given mixed bimatrix of order $(n_1 - k_1 \times n_1, n_2 - k_2 \times n_2)$ with entries from the same field $Z_q$, $H = H_1 \cup H_2$ is the parity check bimatrix of the bicode C.

The standard form of H is $(A_1, I_{n_1-k_1}) \cup (A_2, I_{n_2-k_2})$ with $A_1$ a $n_1 - k_1 \times k_1$ matrix and $A_2$ a $n_2 - k_2 \times k_2$ matrix. $I_{n_1-k_1}$ and $I_{n_2-k_2}$ are $n_1 - k_1 \times n_1 - k_1$ and $n_2 - k_2 \times n_2 - k_2$ identity matrices respectively. The bimatrix H is called the parity check bimatrix of the bicode $C = C_1 \cup C_2$. C is also called the linear $(n, k) = (n_1 \cup n_2, k_1 \cup k_2) = (n_1, k_1) \cup (n_2, k_2)$ bicode.

***Example 2.1.2:*** Let $C (n, k) = C (6 \cup 7, 3 \cup 4) = C_1 (6, 3) \cup C_2 (7, 4)$ be a bicode got by the parity check bimatrix, $H = H_1 \cup H_2$ where

$$H = \begin{bmatrix} 0 & 1 & 1 & 1 & 0 & 0 \\ 1 & 0 & 1 & 0 & 1 & 0 \\ 1 & 1 & 0 & 0 & 0 & 1 \end{bmatrix} \cup \begin{bmatrix} 1 & 1 & 1 & 0 & 1 & 0 & 0 \\ 0 & 1 & 1 & 1 & 0 & 1 & 0 \\ 1 & 1 & 0 & 1 & 0 & 0 & 1 \end{bmatrix}.$$

The bicodes are obtained by solving the equations.

$$Hx^T = H_1 x_1^T \cup H_2 \; x_2^T = (0) \cup (0).$$

There are $2^3 \cup 2^4$ bicode words given by

{(0 0 0 0 0 0), (0 1 1 0 1 1), (1 1 0 1 1 0), (0 0 1 1 1 1 0), (1 0 0 0 1 1), (1 1 1 0 0 0), (0 1 0 1 0 1), (1 0 1 1 0 1)} ∪ {(0 0 0 0 0 0 0), (1 0 0 0 1 0 1), (0 1 0 0 1 1 1), (0 0 1 0 1 1 0), (0 0 0 1 0 1 1), (1 1 0 0 0 1 0), (1 0 1 0 0 1 1), (1 0 0 1 1 1 1), (0 1 1 0 0 0 1), (0 1 0 1 1 0 0), (0 0 1 1 1 0 1), (1 1 1 0 1 0 0), (1 1 0 1 0 0 1), (1 0 1 1 0 0 0), (0 1 1 1 0 1 0), (1 1 1 1 1 1 1)}.

Clearly this is a bicode over $Z_2$. Now the main advantage of a bicode is that at a time two codes of same length or of different length are sent simultaneously and the resultant can be got.

As in case of codes if $x = x_1 \cup x_2$ is a sent message and $y = y_1 \cup y_2$ is a received message using the parity check bimatrix $H = H_1 \cup H_2$ one can easily verify whether the received message is a correct one or not.

For if we consider $Hx^T = H_1x_1^T \cup H_2x_2^T$ then $Hx^T = (0) \cup (0)$ for $x = x_1 \cup x_2$ is the bicode word which was sent. Let $y = y_1$ $\cup$ $y_2$ be the received word. Consider $Hy^T = H_1y_1^T \cup H_2y_2^T$, if $Hy^T = (0) \cup (0)$ then the received bicode word is a correct bicode if $Hy^T \neq (0)$, then we say there is error in the received bicode word. So we define bisyndrome $S^B(y) = Hy^T$ for any bicode $y = y_1 \cup y_2$.

Thus the bisyndrome

$$S^B(y) = S_1(y_1) \cup S_2(y_2) = H_1y_1^T \cup H_2y_2^T.$$

If the bisyndrome $S^B(y) = (0)$ then we say y is the bicode word. If $S^B(y) \neq (0)$ then the word y is not a bicode word. This is the fastest and the very simple means to check whether the received bicode word is a bicode word or not.

Now we proceed on to define the notion of bigenerator bimatrix or just the generator bimatrix of a bicode C(n, k) = $C_1(n_1, k_1) \cup C_2(n_2, k_2)$.

**DEFINITION 2.1.1:** *The generator bimatrix $G = (I_k, -A^T) = G_1 \cup G_2 = \left(I_{k_1}^1 \; -A_1^T\right) \cup \left(I_{k_2}^2 \; -A_2^T\right)$ is called the canonical generator bimatrix or canonical basic bimatrix or encoding bimatrix of a linear bicode, $C(n, k) = C(n_1 \cup n_2, k_1 \cup k_2)$ with parity check bimatrix; $H = H_1 \cup H_2 = (A_1, I_{n_1-k_1}^1) \cup (A_2, I_{n_2-k_2}^2)$. We have $GH^T = (0)$ i.e. $G_1H_1^T \cup G_2H_2^T = (0) \cup (0)$.*

We now illustrate by an example how a generator bimatrix of a bicode functions.

***Example 2.1.3:*** Consider $C(n, k) = C_1(n_1, k_1) \cup C_2(n_2, k_2)$ a bicode over $Z_2$ where $C_1(n_1, k_1)$ is a (7, 4) code and $C_2(n_2, k_2)$ is a (9, 3) code given by the generator bimatrix

$G = G_1 \cup G_2 =$

$$\begin{bmatrix} 1 & 0 & 0 & 0 & 1 & 0 & 1 \\ 0 & 1 & 0 & 0 & 1 & 1 & 1 \\ 0 & 0 & 1 & 0 & 1 & 1 & 0 \\ 0 & 0 & 0 & 1 & 0 & 1 & 1 \end{bmatrix} \cup \begin{bmatrix} 1 & 0 & 0 & 1 & 0 & 0 & 1 & 0 & 0 \\ 0 & 1 & 0 & 0 & 1 & 0 & 0 & 1 & 0 \\ 0 & 0 & 1 & 0 & 0 & 1 & 0 & 0 & 1 \end{bmatrix}.$$

One can obtain the bicode words by using the rule $x_1 = aG$ where a is the message symbol of the bicode. $x_1 \cup x_2 = a^1G_1 \cup a^2G_2$ with $a = a^1 \cup a^2$;
where

$$a^1 = a_1^1 a_2^1 a_3^1$$

and

$$a^2 = a_1^2\ a_2^2\ a_3^2\ a_4^2\ a_5^2\ a_6^2.$$

We give another example in which we calculate the bicode words.

***Example 2.1.4:*** Consider a $C(n, k) = C_1(6, 3) \cup C_2(4, 2)$ bicode over $Z_2$.

The $2^3 \cup 2^2$ code words $x^1$ and $x^2$ of the binary bicode can be found using the generator bimatrix $G = G_1 \cup G_2$.

$$G = \begin{bmatrix} 1 & 0 & 0 & 0 & 1 & 1 \\ 0 & 1 & 0 & 1 & 0 & 1 \\ 0 & 0 & 1 & 1 & 1 & 0 \end{bmatrix} \cup \begin{bmatrix} 1 & 0 & 1 & 1 \\ 0 & 1 & 0 & 1 \end{bmatrix}.$$

$x = aG$ i.e., $x^1 \cup x^2 = a^1G_1 \cup a^2G_2$ where $a^1 = a_1^1\ a_2^1\ a_3^1$ and $a^2 = a_1^2\ a_2^2$ where the message symbols are

$$\begin{Bmatrix} 000 & 001 & 100 & 010 \\ 110 & 011 & 101 & 111 \end{Bmatrix} \cup \begin{Bmatrix} 01 & 10 \\ 00 & 11 \end{Bmatrix}.$$

Thus $x = x^1 \cup x^2 = a^1G_1 \cup a^2G_2$.
We get the bicodes as follows:

$$\left\{\begin{matrix} 000000 & 011011 & 110110 & 001110 \\ 100011 & 111000 & 010101 & 101101 \end{matrix}\right\} \cup \left\{\begin{matrix} 0000 & 1010 \\ 0101 & 1110 \end{matrix}\right\}$$

Now we proceed on to define the notion of repetition bicode and parity check bicode.

**DEFINITION 2.1.2:** *If each bicode word of a bicode consists of only one message symbol* $a = a_1 \cup a_2 \in F_2$ *and the* $(n_1 - 1) \cup (n_2 - 1)$ *check symbols and* $x_2^1 = x_3^1 = \ldots = x_{n_1}^1 = a_1$ *and* $x_2^2 = x_3^2 = \ldots = x_{n_2}^2 = a_2$, $a_1$ *repeated* $n_1 - 1$ *times and* $a_2$ *is repeated* $n_2 - 1$ *times. We obtain the binary* $(n, 1) = (n_1, 1) \cup (n_2, 1)$ *repetition bicode with parity check bimatrix,* $H = H_1 \cup H_2$*;*

$$\text{i.e., H} = \begin{bmatrix} 110 & \cdots & 0 \\ 101 & \cdots & 0 \\ \vdots & & \vdots \\ 100 & & 1 \end{bmatrix}_{n_1-1\times n_1} \cup \begin{bmatrix} 110 & \ldots & 0 \\ 101 & & 0 \\ \vdots & \ldots & \vdots \\ 100 & & 1 \end{bmatrix}_{n_2-1\times n_2}$$

There are only 4 bicode words in a repetition bicode.
{(1 1 1 1…1), (0 0 0…0)} ∪ {(1 1 1 1 1 1…1), (0 0 … 0)}.

***Example 2.1.5:*** Consider the Repetition bicode (n, k) = (5, 1) ∪ (4, 1). The parity check bimatrix is given by

$$\begin{bmatrix} 1 & 1 & 0 & 0 & 0 \\ 1 & 0 & 1 & 0 & 0 \\ 1 & 0 & 0 & 1 & 0 \\ 1 & 0 & 0 & 0 & 1 \end{bmatrix}_{4\times 5} \cup \begin{bmatrix} 1 & 1 & 0 & 0 \\ 1 & 0 & 1 & 0 \\ 1 & 0 & 0 & 1 \end{bmatrix}_{3\times 4} .$$

The bicode words are
{(1 1 1 1 1), (0 0 0 0 0)} ∪ {(1 1 1 1), (0 0 0 0)}.

We now proceed on to define Parity-check bicode.

**DEFINITION 2.1.3:** *Parity check bicode is a (n, n-1) = ($n_1$, $n_1$-1) ∪ ($n_2$, $n_2$ – 1) bicode with parity check bimatrix H =* $(1\,1\,...\,1)_{1\times n_1} \cup (1\,1\,...\,1)_{1\times n_2}$.

*Each bicode word has one check symbol and all bicode words are given by all binary bivectors of length n = $n_1$ ∪ $n_2$ with an even number of ones. Thus if the sum of ones of a received bicode word is 1 at least an error must have occurred at the transmission.*

***Example 2.1.6:*** Let (n, k) = (4, 3) ∪ (5, 4) be a bicode with parity check bimatrix.

$$H = (1\ 1\ 1\ 1) \cup (1\ 1\ 1\ 1\ 1) = H_1 \cup H_2.$$

The bicodes related with the parity check bimatrix H is given by

$$\begin{Bmatrix} 0000 & 1001 \\ 0101 & 0011 \\ 1100 & 1010 \\ 0110 & 1111 \end{Bmatrix} \cup$$

$$\begin{Bmatrix} 00000 & 00011 & 11000 & 10100 \\ 10001 & 01100 & 00110 & 01010 \\ 01001 & 11011 & 10110 & 01111 \\ 00101 & 10010 & 11101 & 11110 \end{Bmatrix}.$$

Now the concept of linear bialgebra i.e. bimatrices are used when we define cyclic bicodes. So we define a cyclic bicode as the union of two cyclic codes for we know if C (n, k) is a cyclic code then if ($x_1 \ldots x_n$) ∈ C (n, k) it implies ($x_n\ x_1 \ldots x_{n-1}$) belongs to C (n, k).

We just give an example of a cyclic bicode using both generator bimatrix and the parity check bimatrix or to be more precise using a generator bipolynomial and a parity check bipolynomial.

Suppose $C(n, k) = C_1 (n_1, k_1) \cup C_2 (n_2, k_2)$ be a cyclic bicode then we need generator polynomials $g_1 (x) \mid x^{n_1} - 1$ and $g_2(x) \mid x^{n_2} - 1$ where $g(x) = g_1 (x) \cup g_2 (x)$ is a generator bipolynomial with degree of $g(x) = (m_1, m_2)$ i.e. degree of $g_1(x) = m_1$ and degree of $g_2(x) = m_2$. The $C(n, k)$ linear cyclic bicode is generated by the bimatrix; $G = G_1 \cup G_2$ which is given below

$$\begin{bmatrix} g_0^1 & g_1^1 & g_{m_1}^1 & 0 & 0 \\ 0 & g_0^1 & g_{m_1-1}^1 & g_{m_1}^1 & 0 \\ \vdots & & & & \\ 0 & 0 & g_0^1 \, g_1^1 & \cdots & g_{m_1}^1 \end{bmatrix}$$

$$\cup \begin{bmatrix} g_0^2 & g_1^2 & g_{m_2}^2 & 0 & 0 \\ 0 & g_0^2 & g_{m_2-1}^2 & g_{m_2}^2 & 0 \\ \vdots & & & & \\ 0 & 0 & g_0^2 & \cdots & g_{m_2}^2 \end{bmatrix}$$

$$= \begin{bmatrix} g_1 \\ xg_1 \\ \vdots \\ x^{k_1-1} g_1 \end{bmatrix} \cup \begin{bmatrix} g_2 \\ xg_2 \\ \vdots \\ x^{k_2-1} g_2 \end{bmatrix}.$$

Then C (n, k) is cyclic bicode.

***Example 2.1.7:*** Suppose $g = g_1 \cup g_2 = 1 + x^3 \cup 1 + x^2 + x^3$ be the generator bipolynomial of a $C (n, k) (= C_1 (6, 3) \cup C_2 (7, 3))$ bicode. The related generator bimatrix is given by

$$G = G_1 \cup G_2$$

i.e., $\begin{bmatrix} 1 & 0 & 0 & 1 & 0 & 0 \\ 0 & 1 & 0 & 0 & 1 & 0 \\ 0 & 0 & 1 & 0 & 0 & 1 \end{bmatrix} \cup \begin{bmatrix} 1 & 0 & 1 & 1 & 0 & 0 & 0 \\ 0 & 1 & 0 & 1 & 1 & 0 & 0 \\ 0 & 0 & 1 & 0 & 1 & 1 & 0 \\ 0 & 0 & 0 & 1 & 0 & 1 & 1 \end{bmatrix}$.

Clearly the cyclic bicode is given by

$$\begin{Bmatrix} 000000 & 100100 \\ 001001 & 101101 \\ 010010 & 110110 \\ 011011 & 111111 \end{Bmatrix} \cup$$

$$\begin{Bmatrix} 0000000 & 0001011 & 0110001 & 1101001 \\ 1000101 & 1100010 & 0101100 & 1011000 \\ 0100111 & 1010000 & 0011101 & 0111010 \\ 0010110 & 1001110 & 1110100 & 1111111 \end{Bmatrix}.$$

It is easily verified that the bicode is a cyclic bicode. Now we define how a cyclic bicode is generated by a generator bipolynomial. Now we proceed on to give how the check bipolynomial helps in getting the parity check bimatrix. Let g = $g_1 \cup g_2$ be the generator bipolynomial of a bicode C (n, k) = $C_1(n_1, k_1) \cup C_2 (n_2, k_2)$.

Then $h_1 = \dfrac{x^{n_1 - 1}}{g_1}$ and $h_2 = \dfrac{x^{n_2 - 1}}{g_2}$; $h = h_1 \cup h_2$ is called the check bipolynomial of the bicode C(n, k). The parity check bimatrix $H = H_1 \cup H_2$

$$= \begin{bmatrix} 0 & 0 & 0 & h^1_{k_1} & \cdots & h^1_1 & h^1_0 \\ 0 & 0 & h^1_{k_1} & h^1_{k_1 - 1} & \ldots & h^1_0 & 0 \\ \vdots & & & & \vdots & & \vdots \\ h^1_{k_1} & & \ldots & & h^1_0 & \ldots & 0 \end{bmatrix} \cup$$

$$\begin{bmatrix} 0 & 0 & 0 & h^2_{k_2} & \cdots & h^2_1 & h^2_0 \\ 0 & 0 & h^2_{k_2} & h^2_{k_2-1} & \ldots & h^2_0 & 0 \\ \vdots & & & & \vdots & & \vdots \\ h^2_{k_2} & & \ldots & & h^2_0 & \ldots & 0 \end{bmatrix}.$$

From the above example the parity check bipolynomial h = $h_1 \cup h_2$ is given by

$$(x^3 + 1) \cup (x^4 + x^3 + x^2 + 1).$$

The parity check bimatrix associated with this bipolynomial is given by $G = G_1 \cup G_2$ where,

$$G = \begin{bmatrix} 0 & 0 & 1 & 0 & 0 & 1 \\ 0 & 1 & 0 & 0 & 1 & 0 \\ 1 & 0 & 0 & 1 & 0 & 0 \end{bmatrix} \cup \begin{bmatrix} 0 & 0 & 1 & 0 & 1 & 1 & 1 \\ 0 & 1 & 0 & 1 & 1 & 1 & 0 \\ 1 & 0 & 1 & 1 & 1 & 0 & 0 \end{bmatrix}.$$

The linear bialgebra concept is used for defining biorthogonal bicodes or dual bicodes.

Let $C(n, k) = C_1(n_1, k_1) \cup C_2(n_2, k_2)$ be a linear bicode over $Z_2$. The dual (or orthogonal) bicode $C^{\perp}$ of C is defined by

$$C^{\perp} = \{u \mid u.v = 0, \forall\, v \in C\}$$
$$= \{u_1 \mid u_1 . v_1 = 0, \forall v_1 \in C_1\} \cup \{u_2 \mid u_2 . v_2 = 0, \forall v_2 \in C_2\}.$$

Thus in case of dual bicodes we have if for the bicode C (n, k), where $G = G_1 \cup G_2$ is the generator bimatrix and if its parity check bimatrix is $H = H_1 \cup H_2$ then for the orthogonal bicode or dual bicode, the generator bimatrix is $H_1 \cup H_2$ and the parity check bimatrix is $G = G_1 \cup G_2$.

Thus we have $HG^T = GH^T$.

In this section we define the notion of best biapproximation and pseudo inner biproduct. Here we give their applications to bicodes.

The applications of pseudo best approximation have been carried out in (2005) [39] to coding theory. We have recalled this notion in chapter I of this book.

Now we apply this pseudo best approximation, to get the most likely transmitted code word. Let C be code over $Z_q^n$. Clearly C is a subspace of $Z_q^n$. $Z_q^n$ is a vector space over $Z_p$ where $q = p^t$, ($p$ – a prime, $t > 1$).

We take C = W in the definition and apply the best biapproximation to $\beta$ where $\beta$ is the received code word for some transmitted word from C but $\beta \notin W = C$, for if $\beta \in W = C$ we accept it as the correct message. If $\beta \notin C$ then we apply the notion of pseudo best approximation to $\beta$ related to the subspace C in $Z_q^n$.

Let $\{c_1, \ldots , c_k\}$ be chosen as the basis of C then

$$\sum_{i=1}^{k} \langle \beta, c_i \rangle_p c_i$$

gives the best approximation to $\beta$ clearly

$$\sum_{i=1}^{k} \langle \beta, c_i \rangle_p c_i$$

belongs to C provided

$$\sum_{i=1}^{k} \langle \beta, c_i \rangle_p c_i \neq 0.$$

It is easily seen that this is the most likely received message. If

$$\sum_{i=1}^{k} \langle \beta, c_i \rangle_p c_i = 0,$$

we can choose another basis for C so that

$$\sum_{i=1}^{k} \langle \beta, c_i \rangle_p c_i \neq 0 .$$

Now we just adopt this argument in case of bicodes.

**DEFINITION 2.1.4:** *Let V = $V_1 \cup V_2$ be a bivector space over the finite field $Z_p$, with some pseudo inner biproduct $\langle , \rangle_p$ defined on V. Let W = $W_1 \cup W_2$ be the subbispace of V = $V_1 \cup V_2$. Let $\beta \in V_1 \cup V_2$ i.e. $\beta = \beta_1 \cup \beta_2$ related to W = $W_1 \cup W_2$ is defined as follows: Let*

$$\{\alpha_1, \ldots, \alpha_k\} = \{\alpha_1^1, \ldots, \alpha_{k_1}^1\} \cup \{\alpha_1^2, \ldots, \alpha_{k_2}^2\}$$

*be the chosen basis of the bisubspace W = $W_1 \cup W_2$. The pseudo best biapproximation to $\beta = \beta_1 \cup \beta_2$ if it exists is given by*

$$\sum_{i=1}^{k} \langle \beta, \alpha_i \rangle_p \alpha_i = \sum_{i=1}^{k_1} \langle \beta_1, \alpha_i^1 \rangle_{p_1} \alpha_i^1 \cup \sum_{i=1}^{k_2} \langle \beta_2, \alpha_i^2 \rangle_{p_2} \alpha_i^2 .$$

*If* $\sum_{i=1}^{k} \langle \beta, \alpha_i \rangle_p \alpha_i = 0$ *then we say that the pseudo best biapproximation does not exist for the set of basis*

$$\{\alpha_1, \alpha_2, \ldots, \alpha_k\} = \{\alpha_1^1, \alpha_2^1, \ldots, \alpha_{k_1}^1\} \cup \{\alpha_1^2, \alpha_2^2, \ldots, \alpha_{k_2}^2\} .$$

*In this case we choose another set of basis* $\{\alpha'_1, \ldots, \alpha'_{k_1}\}$ *and find* $\sum_{i=1}^{k} \langle \beta, \alpha'_i \rangle_p \alpha'_1$ *which is taken as the pseudo best biapproximation to $\beta = \beta_1 \cup \beta_2$.*

Now we apply it in case of bicodes in the following way to find the most likely bicode word. Let C = $C_1 \cup C_2$ be a bicode in $Z_q^n$. Clearly C is a bivector subspace of $Z_q^n$. Take in the definition C = W and apply the pseudo best biapproximation. If some bicode word x = $x_1 \cup x_2$ in C = $C_1 \cup C_2$ is transmitted and $\beta = \beta_1 \cup \beta_2$ is received bicode word then if $\beta \in$ C = $C_1 \cup C_2$ then it is accepted as the correct message if $\beta \notin C = C_1 \cup C_2$ then we apply pseudo best biapproximation to $\beta = \beta_1 \cup \beta_2$ related to the subbispace $C = C_1 \cup C_2$ in $Z_q^n$.

Here three cases occur if $\beta = \beta_1 \cup \beta_2 \notin C = C_1 \cup C_2$.

1) $\beta_1 \in C_1$ and $\beta_2 \notin C_2$ so that $\beta_1 \cup \beta_2 \notin C$
2) $\beta_1 \notin C_1$ and $\beta_2 \in C_2$ so that $\beta_1 \cup \beta_2 \notin C$
3) $\beta_1 \notin C_1$ and $\beta_2 \notin C_2$ so that $\beta_1 \cup \beta_2 \notin C$.

We first deal with (3) then indicate the working method in case of (1) and (2).

Given $\beta = \beta_1 \cup \beta_2$, with $\beta \notin C$ ; $\beta_1 \notin C_1$ and $\beta_2 \notin C_2$;choose a basis

$$(c_1, c_2, \ldots, c_k) = \left\{c_1^1, c_2^1, \ldots, c_{k_1}^1\right\} \cup \left\{c_1^2, c_2^2, \ldots, c_{k_2}^2\right\}$$

of the subbispace C. To find the pseudo best biapproximation to $\beta$ in C find

$$\sum_{i=1}^{k} \langle \beta / c_i \rangle_p \, c_i = \sum_{i=1}^{k_1} \langle \beta_1 / c_i^1 \rangle_1 \, c_i^1 \cup \sum_{i=1}^{k_2} \langle \beta_2 / c_i^2 \rangle_2 \, c_i^2 \, .$$

If both

$$\sum_{i=1}^{k_1} \langle \beta^1 / c_i^1 \rangle_1 \, c_i^1 \neq 0$$

and

$$\sum_{i=1}^{k_2} \langle \beta_2 / c_i^2 \rangle_2 \, c_i^2 \neq 0$$

then

$$\sum_{i=1}^{k} \langle \beta / c_i \rangle_p \, c_i \neq 0$$

is taken as the pseudo best biapproximation to $\beta$.

If one of

$$\sum_{i=1}^{k_t} \langle \beta_t \, c_i^t \rangle_t \, c_i^t = 0 \;\; (t = 1, 2)$$

say

$$\sum_{i=1}^{k_1} \langle \beta_i \, c_i^1 \rangle c_i^1 = 0$$

then choose a new basis, say $\{c'_1 \ldots c'_{k_t}\}$ and calculate

$$\sum_{i=1}^{k_t} \langle \beta_t \, c'_i \rangle_t \, c'_i \; ;$$

that is

$$\sum_{i=1}^{k_1} \langle \beta_1 c'_i \rangle_1 \, c'_i \cup \sum_{i=1}^{k_2} \langle \beta_2 c_i^2 \rangle_2 \, c_i^2$$

will be the pseudo best biapproximation to β in C. If both

$$\sum_{i=1}^{k_1} \langle \beta_1 c_i^1 \rangle_1 \, c_i^1 = 0 \text{ and } \sum_{i=1}^{k_2} \langle \beta_2 c_i^2 \rangle_2 \, c_i^2 = 0$$

then choose a new basis for C.

$$\{b_1, \ldots, b_k\} = \{b_1^1, \ldots, b_{k_1}^1\} \cup \{b_1^2, \ldots, b_{k_2}^2\}$$

and find

$$\sum_{i=1}^{k} \langle \beta, b_i \rangle_p \, b_i = \sum_{i=1}^{k_1} \langle \beta_1, b_i^1 \rangle_1 \, b_i^1 \cup \sum_{i=1}^{k_2} \langle \beta_2, b_i^2 \rangle_2 \, b_i^2 \; .$$

If this is not zero it will be the pseudo best biapproximation to β in C.

Now we proceed on to work for cases (1) and (2) if we work for one of (1) or (2) it is sufficient. Suppose we assume

$$\sum_{i=1}^{k} \langle \beta c_i \rangle_p \, c_i = \sum_{i=1}^{k_1} \langle \beta_1 \, c_i^1 \rangle_1 \, c_i^1 \cup \sum_{i=1}^{k_2} \langle \beta_2 \, c_i^2 \rangle_2 \, c_i^2$$

and say

$$\sum_{i=1}^{k_2} \langle \beta_2 \, c_i^2 \rangle_2 \, c_i^2 = 0$$

and

$$\sum_{i=1}^{k_1} \langle \beta_1 \, c_i^1 \rangle_1 \, c_i^1 \neq 0$$

then we choose only a new basis for $C_2$ and calculate the pseudo best approximation for $\beta_2$ relative to the new basis of $C_2$.

Now we illustrate this by the following example:

***Example 2.1.8:*** Let $C = C_1 (6, 3) \cup C_2 (8, 4)$ be a bicode over $Z_2$ generated by the parity check bimatrix $H = H_1 \cup H_2$ given by

$$\begin{bmatrix} 0 & 1 & 1 & 1 & 0 & 0 \\ 1 & 0 & 1 & 0 & 1 & 0 \\ 1 & 1 & 0 & 0 & 0 & 1 \end{bmatrix} \cup \begin{bmatrix} 1 & 1 & 0 & 1 & 1 & 0 & 0 & 0 \\ 0 & 0 & 1 & 1 & 0 & 1 & 0 & 0 \\ 1 & 0 & 1 & 0 & 0 & 0 & 1 & 0 \\ 1 & 1 & 1 & 1 & 0 & 0 & 0 & 1 \end{bmatrix}.$$

The bicode is given by $C = C_1 \cup C_2$

$$= \left\{ \begin{matrix} (000000) & (011011) & (110110) & (001110) \\ (100011) & (111000) & (100101) & (101101) \end{matrix} \right\} \cup$$

{(0 0 0 0 0 0 0 0), (1 0 0 0 1 0 1 1), (0 1 0 0 1 0 0 1), (0 0 1 0 0 1 1 1), (0 0 0 1 1 1 0 1), (1 1 0 0 0 0 1 0 ), (0 1 1 0 1 1 1 0), (0 0 1 1 1 0 1 0), (0 1 0 1 0 1 0 0), (1 0 1 0 1 1 0 0), (1 0 0 1 0 1 1 0), (1 1 1 0 0 1 0 1), (0 1 1 1 0 0 1 1), (1 1 0 1 1 1 1 1), (1 0 1 1 0 0 0 1), (1 1 1 1 1 0 0 0 )}.

Choose a basis $B = B_1 \cup B_2$ = {(0 0 1 1 1 0), (1 1 1 0 0 0), (0 1 0 1 0 1)} $\cup$ {(0 1 0 0 1 0 0 1), (1 1 0 0 0 0 1 0), (1 1 1 0 0 1 0 1), (1 1 1 1 1 0 0 0 )}. Let $\beta$ = (1 1 1 1 1 1) $\cup$ (1 1 1 1 1 1 1 1) be the received bicode. Clearly $\beta \notin C_1 \cup C_2 = C$. Now we find the best pseudo biapproximation to $\beta$ in C relative to the basis B = $B_1 \cup B_2$ under the pseudo inner biproduct.

Let $\alpha = \alpha_1 \cup \alpha_2$ = {$\langle$(1 1 1 1 1 1) / (0 0 1 1 1 0)$\rangle_1$ (0 0 1 1 1 0) + $\langle$(1 1 1 1 1 1) / (1 1 1 0 0 0)$\rangle_1$ (1 1 1 0 0 0) + $\langle$(1 1 1 1 1 1) / (0 1 0 1 0 1)$\rangle_1$ (0 1 0 1 0 1) } $\cup$ {$\langle$(1 1 1 1 1 1 1 1) / (0 1 0 0 1 0 0 1)$\rangle_2$ (0 1 0 0 1 0 0 1) + $\langle$(1 1 1 1 1 1 1 1) , (1 1 0 0 0 0 1 0)$\rangle_p$ (1 1 0 0 0 0 1 0) + $\langle$(1 1 1 1 1 1 1 1) / (1 1 1 0 0 1 0 1)$\rangle_2$ (1 1 1 0 0 1 0 1) + $\langle$(1 1 1 1 1 1 1 1) / (1 1 1 1 1 0 0 0)$\rangle_2$ (1 1 1 1 1 0 0 0)}

= {(0 0 1 1 1 0) + (1 1 1 0 0 0) + (0 1 0 1 0 1)} $\cup$ {(0 1 0 0 1 0 0 1) + (1 1 0 0 0 0 1 0) + (1 1 1 0 0 1 0 1) + (1 1 1 1 1 0 0 0)}

= (1 0 0 0 1 1) $\cup$ (1 0 0 1 0 1 1 0) $\in C_1 \cup C_2 = C$.

Thus this is the pseudo best biapproximation to the received bicode word {(1 1 1 1 1 1) $\cup$ (1 1 1 1 1 1 1 1)}.

Thus the method of pseudo best biapproximation to a received bicode word which is not a bicode word is always guaranteed. If one wants to get best of the best pseudo biapproximations one can vary the basis and find the resultant pseudo biapproximated bicodes, compare it with the received message, the bicode which gives the least Hamming bidistance from the received word is taken as the best of pseudo best biapproximated bicode.

*Note:* We just define the Hamming bidistance of two bicodes x = $x_1 \cup x_2$ and y = $y_1 \cup y_2$, which is given by $d^B$ (x, y) = d ($x_1$ $y_1$) $\cup$ d ($x_2$ $y_2$) where d ($x_1$ $y_1$) and d($x_2$, $y_2$) are the Hamming distance. The least of Hamming bidistances is taken as minimum of the sum of d ($x_1$ $y_1$) + d ($x_2$ $y_2$).

In this section we introduce the notion of tricodes and describe a few of its properties.

**DEFINITION 2.1.5:** *Let C = $C_1 \cup C_2 \cup C_3$ where $C_1$ $C_2$ and $C_3$ are distinct codes '$\cup$' just a symbol. C is then defined to be a tricode.*

*Note:* It is very important to know that each of the $C_i$ 's must be distinct for $1 \le i \le 3$.

***Example 2.1.9:*** Consider the tricode C = $C_1 \cup C_2 \cup C_3$ where $C_1$ is a (7, 4) code, $C_2$ is a (6, 3) code and $C_3$ is a (5, 4) code with the associated generator trimatrix

$$G = G_1 \cup G_2 \cup G_3$$

$$= \begin{bmatrix} 1 & 1 & 0 & 1 & 0 & 0 & 0 \\ 0 & 1 & 1 & 0 & 1 & 0 & 0 \\ 0 & 0 & 1 & 1 & 0 & 1 & 0 \\ 0 & 0 & 0 & 1 & 1 & 0 & 1 \end{bmatrix} \cup \begin{bmatrix} 1 & 0 & 0 & 1 & 0 & 0 \\ 0 & 1 & 0 & 0 & 1 & 0 \\ 0 & 0 & 1 & 0 & 0 & 1 \end{bmatrix} \cup$$

$$\begin{bmatrix} 1 & 1 & 0 & 0 & 0 \\ 0 & 1 & 1 & 0 & 0 \\ 0 & 0 & 1 & 1 & 0 \\ 0 & 0 & 0 & 1 & 1 \end{bmatrix}.$$

The tricode generated by G is given as follows. {(0 0 0 0 0 0 0) ∪ (0 0 0 0 0 0) ∪ (0 0 0 0 0), (1 1 0 1 0 0 0) ∪ (1 0 0 1 0 0) ∪ (1 1 0 0 0), (0 1 1 0 1 0 0) ∪ (0 1 0 0 1 0) ∪ (0 0 0 0 0), (1 1 0 1 0 0 0) ∪ (1 0 0 1 0 0) ∪ (1 0 0 0 1), (0 1 1 0 1 0 0) ∪ (1 0 0 1 0 0) ∪ (1 0 0 0 1), (0 1 1 0 1 0 0) ∪ (1 0 0 1 0 0) ∪ (1 0 0 0 1)} and so on}.

***Example 2.1.10:*** Let C = $C_1 \cup C_2 \cup C_3$ be a tricode where C is generated by the trimatrix G = $G_1 \cup G_2 \cup G_3$ with G =

$$\begin{bmatrix} 1 & 0 & 1 & 1 \\ 0 & 1 & 0 & 0 \end{bmatrix} \cup \begin{bmatrix} 1 & 0 & 0 & 1 & 1 \\ 0 & 1 & 1 & 0 & 1 \end{bmatrix} \cup \begin{bmatrix} 1 & 0 & 0 & 0 & 0 & 1 \\ 1 & 1 & 0 & 1 & 0 & 0 \end{bmatrix}.$$

The elements of the tricode C generated by G is as follows:

{(0 0 0 0) ∪ (0 0 0 0 0) ∪ (0 0 0 0 0 0), (1 0 1 1) ∪ (0 0 0 0 0) (0 0 0 0 0 0), (0 1 0 0) ∪ (0 0 0 0 0) ∪ (0 0 0 0 0 0), (1 1 1 1) ∪ (0 0 0 0 0) ∪ (0 0 0 0 0 0), (0 0 0 0) ∪ (1 0 0 1 1) ∪ (0 0 0 0 0 0), (0 0 0 0) ∪ (0 1 1 0 1) ∪ (0 0 0 0 0 0), (0 0 0 0) ∪ (1 1 1 1 0) ∪ (0 0 0 0 0 0), (0 0 0 0) ∪ (0 0 0 0 0) ∪ (1 0 0 0 0 1), (0 0 0 0) ∪ (0 0 0 0 0) ∪ (0 1 0 1 0 0), (0 0 0 0) ∪ (1 0 0 0 0 0) ∪ (0 1 0 1 0 1), (0 0 0 0) ∪ (1 1 1 1 0) ∪ (1 0 0 0 0 1), (0 0 0 0) ∪ (1 1 1 1 0) ∪ (0 1 0 1 0 1), (0 0 0 0) ∪ (1 1 1 1 0) ∪ (0 1 0 1 0 0), (0 0 0 0) ∪ (0 1 1 0 1) ∪ (1 0 0 0 0 1), (0 0 0 0) ∪ (0 1 1 0 1) ∪ (0 1 0 1 0 0), (0 0 0 0) ∪ (0 1 1 0 1) ∪ (0 1 0 1 0 1) and so on}.

Now having defined a tricode we proceed on to define the notion of n-code.

**DEFINITION 2.1.6:** *Let $C = C_1 \cup C_2 \cup \ldots \cup C_n$ ($n \geq 4$) is said to be n-code if each of the $C_i$ is a $(n_i, k_i)$ code, $1 \leq i \leq n$. Clearly C is generated by the n-matrix $G = G_1 \cup G_2 \cup \ldots \cup G_n$ where each $G_i$ generates a $(n_i, k_i)$ code. The parity check n-matrix of the n-code C is given by $H = H_1 \cup H_2 \cup \ldots \cup H_n$ with each $H_i$ in the standard form $1 \leq i \leq n$ where H is a n-matrix. Clearly $GH^T = (0) \cup (0) \cup (0) \cup \ldots \cup (0)$; i.e.,*

$$
\begin{aligned}
GH^T &= (G_1 \cup G_2 \cup \ldots \cup G_n)\,(H_1 \cup H_2 \cup \ldots \cup H_n)^T \\
&= (G_1 \cup G_2 \cup \ldots \cup G_n) \times (H_1^T \cup H_2^T \cup \ldots \cup H_n^T) \\
&= G_1 H_1^T \cup G_2 H_2^T \cup \ldots \cup G_n H_n^T \\
&= (0) \cup (0) \cup \ldots \cup (0)
\end{aligned}
$$

*is the zero n-matrix. Here each $C_i$ is a distinct code i.e., $C_i \neq C_j$ if $i \neq j$.*

Now we proceed give examples of n-codes.

***Example 2.1.11:*** Consider a 6-code C given by $C = C_1 \cup C_2 \cup \ldots \cup C_6$ where C is a 6-code generated by the 6-matrix given by

$$G = G_1 \cup G_2 \cup \ldots \cup G_6$$

$$= \begin{bmatrix} 1 & 0 & 0 & 1 \\ 0 & 1 & 1 & 1 \end{bmatrix} \cup \begin{bmatrix} 1 & 0 & 0 & 1 & 1 \\ 0 & 1 & 0 & 0 & 0 \\ 0 & 0 & 1 & 1 & 0 \end{bmatrix}$$

$$\cup \begin{bmatrix} 1 & 0 & 0 & 0 & 1 \\ 1 & 1 & 0 & 1 & 1 \end{bmatrix} \cup \begin{bmatrix} 1 & 0 & 1 & 1 \\ 0 & 1 & 0 & 0 \\ 0 & 1 & 1 & 0 \end{bmatrix}$$

$$\cup \begin{bmatrix} 1 & 0 & 0 & 1 & 1 & 1 \\ 0 & 1 & 1 & 0 & 0 & 1 \end{bmatrix} \cup \begin{bmatrix} 1 & 0 & 0 & 1 & 0 & 1 \\ 0 & 1 & 0 & 0 & 1 & 1 \\ 0 & 0 & 1 & 1 & 1 & 1 \end{bmatrix}.$$

The 6 codes generated by the generator 6-matrix is as follows:

{(0 0 0 0) ∪ (0 0 0 0 0) ∪ (0 0 0 0 0) ∪ (0 0 0 0) ∪ (0 0 0 0 0 0) ∪ (0 0 0 0 0 0), (1 0 0 1) ∪ (0 0 0 0 0) ∪ (0 0 0 0 0) ∪ (0 0 0 0) ∪ (0 0 0 0 0 0) ∪ (0 0 0 0 0 0), (1 0 0 1) ∪ (1 0 0 1 1) ∪ (1 0 0 0 1) ∪ (0 1 0 0) ∪ (1 0 0 1 1 1) ∪ (1 0 0 1 0 1), (1 0 0 1) ∪ (1 0 0 1 1) ∪ (1 0 0 0 1) ∪ (0 1 0 0) ∪ (1 0 0 1 1 1) ∪ (0 1 0 0 1 1), (1 0 0 1) ∪ (1 0 0 1 1) ∪ (1 0 0 0 1) ∪ (0 1 0 0) ∪ (1 0 0 1 1 1) ∪ (0 0 1 1 1 1) and so on}.

***Example 2.1.12:*** Let us consider the 5-code $C = C_1 \cup C_2 \cup C_3 \cup C_4 \cup C_5$ where the generator of the 5-code is given by the generator 5-matrix

$$G = G_1 \cup G_2 \cup G_3 \cup G_4 \cup G_5$$

$$= \begin{bmatrix} 1 & 0 & 0 & 0 & 0 \\ 0 & 1 & 1 & 1 & 0 \end{bmatrix} \cup \begin{bmatrix} 1 & 1 & 1 & 0 & 0 & 1 \\ 0 & 0 & 1 & 1 & 1 & 1 \end{bmatrix} \cup$$

$$\begin{bmatrix} 1 & 0 & 0 & 1 & 1 & 0 \\ 1 & 1 & 0 & 0 & 1 & 1 \end{bmatrix} \cup \begin{bmatrix} 1 & 0 & 1 & 1 & 0 \\ 1 & 0 & 0 & 1 & 0 \\ 1 & 1 & 0 & 0 & 1 \end{bmatrix} \cup \begin{bmatrix} 1 & 0 & 1 & 1 \\ 1 & 1 & 0 & 0 \\ 0 & 0 & 1 & 0 \end{bmatrix}.$$

The 5-code given by the G is as follows :

{(0 0 0 0 0) ∪ (0 0 0 0 0 0) ∪ (0 0 0 0 0 0) ∪ (0 0 0 0 0) ∪ (0 0 0 0), (1 0 0 0 0) ∪ (1 1 1 0 0 1) ∪ (1 0 0 1 1 0) ∪ (1 0 1 1 0) ∪ (1 0 1 1), (1 0 0 0 0) ∪ (0 0 0 0 0 0) ∪ (1 1 0 0 1 1) ∪ (1 1 0 0 1) ∪ (0 0 1 0), (0 1 1 1 0) ∪ (0 0 1 1 1 1) ∪ (1 1 0 0 1 1) ∪ (1 1 0 0 1) ∪ (0 0 1 0), (1 1 1 1 0) ∪ (1 1 0 1 1 0) ∪ (0 1 0 1 0 1) ∪ (1 1 1 0 1) ∪ (0 1 0 1), (1 1 1 1 0) ∪ (0 0 0 0 0 0) ∪ (0 0 0 0 0 0) ∪ (0 0 0 0 0) ∪ (0 0 0 0)), (1 1 1 1 0) ∪ (0 0 1 1 1 1) ∪

(0 0 0 0 0 0) ∪ (0 0 0 0 0) ∪ (0 0 0 0), (1 1 1 1 0) ∪ (0 0 1 1 1 1) ∪ (0 1 0 1 0 1) ∪ (0 0 0 0 0) ∪(0 0 0 0) and so on}.

We now illustrate and define the notion of repetition bicode and repetition n-code (n ≥ 3).

We can just non mathematically but technically define repetition bicode as ;

*Let us consider the bimatrix* $H = H_1 \cup H_2$ *where*

$$H_1 = \begin{bmatrix} 1 & 1 & 0 & \cdots & 0 \\ 1 & 0 & 1 & \cdots & 0 \\ \vdots & \vdots & \vdots & & \vdots \\ 1 & 0 & 0 & \cdots & 1 \end{bmatrix}$$

*be a* $n - 1 \times n$ *matrix with the first column in which all entries are 1 and the rest* $n - 1 \times n - 1$ *matrix is the identity matrix.* $H_2$ *is a* $m - 1 \times m$ *matrix* $m \neq n$ *with the first column in which all entries are 1 and the rest* $m - 1 \times m - 1$ *matrix is the identity matrix. The bicode* $C = C_1 \cup C_2$ *given by using this parity check bimatrix H is a repetition bicode.*
*There are only 4 bicode words namely* $\underbrace{(0\ 0 \cdots\ 0)}_{n-times} \cup \underbrace{(1\ 1\ 1\ \cdots 1)}_{m-times}$, *(0 0 … 0)* ∪ *(0 0 … 0), (1 1 … 1 1)* ∪ *(0 0 … 0) and (1 1 … 1)* ∪ *(1 1 1 … 1).*

Now we illustrate this by the following example.

***Example 2.1.13:*** Let $H = H_1 \cup H_2$ where

$$H_1 = \begin{bmatrix} 1 & 1 & 0 & 0 & 0 \\ 1 & 0 & 1 & 0 & 0 \\ 1 & 0 & 0 & 1 & 0 \\ 1 & 0 & 0 & 0 & 1 \end{bmatrix}$$

and

$$H_2 = \begin{bmatrix} 1 & 1 & 0 & 0 & 0 & 0 \\ 1 & 0 & 1 & 0 & 0 & 0 \\ 1 & 0 & 0 & 1 & 0 & 0 \\ 1 & 0 & 0 & 0 & 1 & 0 \\ 1 & 0 & 0 & 0 & 0 & 1 \end{bmatrix}.$$

The repetition bicode words given by H are (0 0 0 0 0) $\cup$ (0 0 0 0 0 0), (1 1 1 1 1) $\cup$ (0 0 0 0 0 0), (0 0 0 0 0) $\cup$ (1 1 1 1 1 1), (1 1 1 1 1) $\cup$ (1 1 1 1 1 1).

Now we give yet another example.

***Example 2.1.14:*** Let us consider the bicode $C = C_1 \cup C_2$ given by the parity check bimatrix $H = H_1 \cup H_2$

$$= \begin{bmatrix} 1 & 1 & 0 & 0 & 0 & 0 & 0 & 0 \\ 1 & 0 & 1 & 0 & 0 & 0 & 0 & 0 \\ 1 & 0 & 0 & 1 & 0 & 0 & 0 & 0 \\ 1 & 0 & 0 & 0 & 1 & 0 & 0 & 0 \\ 1 & 0 & 0 & 0 & 0 & 1 & 0 & 0 \\ 1 & 0 & 0 & 0 & 0 & 0 & 1 & 0 \\ 1 & 0 & 0 & 0 & 0 & 0 & 0 & 1 \end{bmatrix} \cup \begin{bmatrix} 1 & 1 & 0 & 0 & 0 & 0 \\ 1 & 0 & 1 & 0 & 0 & 0 \\ 1 & 0 & 0 & 1 & 0 & 0 \\ 1 & 0 & 0 & 0 & 1 & 0 \\ 1 & 0 & 0 & 0 & 0 & 1 \end{bmatrix}.$$

The repetition bicode words given by H are {(1 1 1 1 1 1 1 1) $\cup$ (1 1 1 1 1 1), (1 1 1 1 1 1 1 1) $\cup$ (0 0 0 0 0 0), (0 0 0 0 0 0 0 0) $\cup$ (1 1 1 1 1 1), (0 0 0 0 0 0 0 0) $\cup$ (0 0 0 0 0 0)}.

Clearly in the definition of repetition bicode we demand $m \neq n$.

We now proceed on to define the notion of repetition tricode.

**DEFINITION 2.1.7:** *Let H be parity check trimatrix i.e. $H = H_1 \cup H_2 \cup H_3$ where each $H_i$ is a $m_i - 1 \times m_i$ matrix with i = 1, 2, 3 and $m_i \neq m_j$ if $i \neq j$, j = 1, 2, 3, H is defined to be the repetition tricode parity check trimatrix of the tricode $C = C_1 \cup C_2 \cup C_3$.*

*The tricodes obtained by using this parity check trimatrix will be known as the repetition tricode.*

We demand $H_i \neq H_j$ if $i \neq j$, $i \leq i, j \leq 3$.

We illustrate this by the following example.

***Example 2.1.15:*** Let $H = H_1 \cup H_2 \cup H_3$ be a partly check trimatrix of a repetition tricode, $C = C_1 \cup C_2 \cup C_3$ where

$$H_1 = \begin{bmatrix} 1 & 1 & 0 & 0 & 0 \\ 1 & 0 & 1 & 0 & 0 \\ 1 & 0 & 0 & 1 & 0 \\ 1 & 0 & 0 & 0 & 1 \end{bmatrix},$$

$$H_2 = \begin{bmatrix} 1 & 1 & 0 & 0 & 0 & 0 & 0 \\ 1 & 0 & 1 & 0 & 0 & 0 & 0 \\ 1 & 0 & 0 & 1 & 0 & 0 & 0 \\ 1 & 0 & 0 & 0 & 1 & 0 & 0 \\ 1 & 0 & 0 & 0 & 0 & 1 & 0 \\ 1 & 0 & 0 & 0 & 0 & 0 & 1 \end{bmatrix}$$

and

$$H_3 = \begin{bmatrix} 1 & 1 & 0 & 0 & 0 & 0 & 0 & 0 \\ 1 & 0 & 1 & 0 & 0 & 0 & 0 & 0 \\ 1 & 0 & 0 & 1 & 0 & 0 & 0 & 0 \\ 1 & 0 & 0 & 0 & 1 & 0 & 0 & 0 \\ 1 & 0 & 0 & 0 & 0 & 1 & 0 & 0 \\ 1 & 0 & 0 & 0 & 0 & 0 & 1 & 0 \\ 1 & 0 & 0 & 0 & 0 & 0 & 0 & 1 \end{bmatrix}.$$

The repetition tricodes obtained by $H = H_1 \cup H_2 \cup H_3$ is

$\{(0\ 0\ 0\ 0\ 0) \cup (0\ 0\ 0\ 0\ 0\ 0\ 0) \cup (0\ 0\ 0\ 0\ 0\ 0\ 0\ 0)$,
$(0\ 0\ 0\ 0\ 0) \cup (0\ 0\ 0\ 0\ 0\ 0\ 0) \cup (1\ 1\ 1\ 1\ 1\ 1\ 1\ 1)$,

$(0\ 0\ 0\ 0\ 0) \cup (1\ 1\ 1\ 1\ 1\ 1\ 1) \cup (1\ 1\ 1\ 1\ 1\ 1\ 1\ 1)$,
$(0\ 0\ 0\ 0\ 0) \cup (1\ 1\ 1\ 1\ 1\ 1\ 1) \cup (0\ 0\ 0\ 0\ 0\ 0\ 0\ 0)$,
$(1\ 1\ 1\ 1\ 1) \cup (0\ 0\ 0\ 0\ 0\ 0\ 0) \cup (0\ 0\ 0\ 0\ 0\ 0\ 0\ 0)$,
$(1\ 1\ 1\ 1\ 1) \cup (1\ 1\ 1\ 1\ 1\ 1\ 1) \cup (1\ 1\ 1\ 1\ 1\ 1\ 1\ 1)$,
$(1\ 1\ 1\ 1\ 1) \cup (0\ 0\ 0\ 0\ 0\ 0\ 0) \cup (1\ 1\ 1\ 1\ 1\ 1\ 1\ 1)$ and
$(1\ 1\ 1\ 1\ 1) \cup (1\ 1\ 1\ 1\ 1\ 1\ 1) \cup (0\ 0\ 0\ 0\ 0\ 0\ 0\ 0)\}$.

Thus for the repetition tricode we have 8 choices to be operated upon by the sender and the receiver.

Now we proceed on to define the notion of repetition n-code using a n-matrix.

**DEFINITION 2.1.7:** *Let $H = H_1 \cup H_2 \cup \ldots \cup H_n$ be a n-matrix where each $H_i$ is a $m_i - 1 \times m_i$ parity check matrix related with a repetition code $1 \leq i \leq n$. i.e. each $H_i$ will have the first column to be ones and the remaining $m_i - 1 \times m_i - 1$ matrix will be a identity matrix with ones along the main diagonal and rest zeros. Clearly the n-matrix gives way to $2^n$ codes and H is called the n-parity check matrix and the codes are called repetition n-codes and has $2^n$ code words; when n = 1 we get the usual repetition code. When n = 2 we get the repetition bicode; when n = 3 the repetition 3-code or repetition tricode. When $n \geq 4$ we obtain the repetition n code which has $2^n$ code words.*

***Example 2.1.16:*** Let $H = H_1 \cup H_2 \cup H_3 \cup H_4 \cup H_5$ be a n-parity check matrix (n = 5); with values from $Z_2 = \{0, 1\}$. Here

$$H_1 = \begin{bmatrix} 1 & 1 & 0 & 0 \\ 1 & 0 & 1 & 0 \\ 1 & 0 & 0 & 1 \end{bmatrix},$$

$$H_2 = \begin{bmatrix} 1 & 1 & 0 & 0 & 0 & 0 \\ 1 & 0 & 1 & 0 & 0 & 0 \\ 1 & 0 & 0 & 1 & 0 & 0 \\ 1 & 0 & 0 & 0 & 1 & 0 \\ 1 & 0 & 0 & 0 & 0 & 1 \end{bmatrix},$$

$$H_3 = \begin{bmatrix} 1 & 1 & 0 & 0 & 0 \\ 1 & 0 & 1 & 0 & 0 \\ 1 & 0 & 0 & 1 & 0 \\ 1 & 0 & 0 & 0 & 1 \end{bmatrix},$$

$$H_4 = \begin{bmatrix} 1 & 1 & 0 & 0 & 0 & 0 & 0 \\ 1 & 0 & 1 & 0 & 0 & 0 & 0 \\ 1 & 0 & 0 & 1 & 0 & 0 & 0 \\ 1 & 0 & 0 & 0 & 1 & 0 & 0 \\ 1 & 0 & 0 & 0 & 0 & 1 & 0 \\ 1 & 0 & 0 & 0 & 0 & 0 & 1 \end{bmatrix}$$

and

$$H_5 = \begin{bmatrix} 1 & 1 & 0 \\ 1 & 0 & 1 \end{bmatrix}.$$

H is the n-parity check matrix for the repetition 5-code. The code words related with H are { (0 0 0 0) $\cup$ (0 0 0 0 0 0) $\cup$ (0 0 0 0 0) $\cup$ (0 0 0 0 0 0 0) $\cup$ (0 0 0), (0 0 0 0) $\cup$ (0 0 0 0 0 0) $\cup$ (0 0 0 0 0) $\cup$ (0 0 0 0 0 0 0) $\cup$ (1 1 1), (0 0 0 0) $\cup$ (0 0 0 0 0 0) $\cup$ (0 0 0 0 0) $\cup$ (1 1 1 1 1 1 1) $\cup$ (0 0 0), (0 0 0 0) $\cup$ (0 0 0 0 0 0) $\cup$ (0 0 0 0 0) $\cup$ (1 1 1 1 1 1 1) $\cup$ (1 1 1), …, (1 1 1 1) $\cup$ (1 1 1 1 1 1) $\cup$ (1 1 1 1 1) $\cup$ (1 1 1 1 1 1 1) $\cup$ (1 1 1) }.

We have $2^5$ i.e. 32 code words related with any repetition 5-code.

Now we proceed on to define parity check bicode differently to make one understand it better.

**DEFINITION 2.1.8:** *Let H = $H_1 \cup H_2$ where* $H_1 = \underbrace{(1\ 1\ \cdots\ 1)}_{n-times}$ *and* $H_2 = \underbrace{(1\ 1\ 1 \cdots 1)}_{m-times}$ *(m ≠ n). C = $C_1 \cup C_2$ is a parity check bicode related with the parity check bimatrix H = $H_1 \cup H_2$. This is a binary (n, n – 1) ∪ (m, m – 1) bicode with parity check bimatrix H = $H_1 \cup H_2$. Each of the bicode words are binary bivectors of length (n, m) with an even number of ones in both the vectors of the bivector. Thus if the sum of the ones in the bivectors is odd atmost two errors must have occurred during the transmission. The last digit is the control symbol or the control digit.*

***Example 2.1.17:*** Let H = (1 1 1 1) ∪ (1 1 1 1 1) be the parity check bimatrix of the bicode C = $C_1 \cup C_2$. The code words given by H are {(0 0 0 0) ∪ (0 0 0 0 0), (1 0 0 1) ∪ (0 0 0 0 0), (0 1 0 1) ∪ (0 0 0 0 0), (0 0 1 1) ∪ ( 0 0 0 0), (1 1 0 0) ∪ (0 0 0 0 0), (1 0 1 0) ∪ (0 0 0 0 0), (0 1 1 0) ∪ (0 0 0 0 0), (1 1 1 1) ∪ (0 0 0 0 0), (0 0 0 0) ∪ (1 0 0 0 1), (0 0 0 0) ∪ (1 1 0 0 0), (0 0 0 0) ∪ (0 0 0 1 1), (0 0 0 0) ∪ (0 1 1 0 0), (0 0 0 0) ∪ (1 0 1 0 0), (0 0 0 0) ∪ (1 0 0 1 0), (0 0 0 0) ∪ (0 1 1 0 0), (0 0 0 0) ∪ (0 1 0 1 0, (0 0 0 0) ∪ (0 1 0 0 1) so on (1 1 1 1) ∪ (0 1 1 1 1), (1 1 1 1) ∪ (1 0 1 1 1), (1 1 1 1) ∪ ( 1 1 0 1 1) (1 1 11 1) ∪ (1 1 1 0 1), (1 1 1 1) ∪ (1 1 1 1 0)}.

We see the sum of the ones in each of the vector of the bivector adds to an even number.

Now we give yet another simple example of a parity check bicode.

***Example 2.1.18:*** Let H = (1 1 1) ∪ (1 1 1 1) be the parity check bimatrix of the parity check bicode C = $C_1 \cup C_2$. The bicodes related with H are { (0 0 0) ∪ (0 0 0 0), (1 0 1) ∪ (0 0 0 0), (0 1 1) ∪ (0 0 0 0), (1 1 0) ∪ (0 0 0 0), (0 0 0) ∪ (1 1 0 0), (0 0 0) ∪ (1 0 0 1), (0 0 0) ∪ (1 0 1 0), (0 0 0) ∪ (0 1 1 0), (0 0 0) ∪ (0 0

1 1), (0 0 0) (0 1 0 1), so on (1 1 0) ∪ (1 1 1 1), (1 0 1) ∪ (1 1 1 1), (0 1 1) ∪ (1 1 1 1) }.

The reader is given the work of finding the number of parity check bicodes when H = $H_1 \cup H_2$ with $H_1$ having m coordinate and $H_2$ having n coordinates (m ≠ n).

Now we proceed on to define the notion of parity check tricodes.

**DEFINITION 2.1.9:** *Let H = $H_1 \cup H_2 \cup H_3$ be a parity check trimatrix where H is a row trivector, i.e. each $H_i$ is a row vector of different lengths with all its entries as ones, $1 \leq i \leq 3$.*

*The code C = $C_1 \cup C_2 \cup C_3$ related with H is a parity check tricode. Each code word in the row trivector C = $C_1 \cup C_2 \cup C_3$ is such that the sum of the elements in each of the row vector of the row trivector is always even.*

We denote this by the following example.

***Example 2.1.19:*** Let H = (1 1 1 1) ∪ (1 1 1) ∪ (1 1 1 1 1 1) = $H_1 \cup H_2 \cup H_3$ be the parity check trimatrix of the parity check tricode. Now the code words related with H are given by C = {(0 0 0 0) ∪ (0 0 0) ∪ (0 0 0 0 0 0), (0 0 0 0) ∪ (0 0 0) ∪ (1 1 0 0 0 0), (0 0 0 0) ∪ (0 0 0) ∪ (1 0 1 0 0 0), (0 0 0 0) ∪ (0 0 0) ∪ (1 0 0 1 0 0), (0 0 0 0) ∪ (0 0 0) ∪ (1 0 0 0 1 0), (0 0 0 0) ∪ (0 0 0) ∪ (1 0 0 0 0 1), (0 0 0 0) ∪ (0 0 0) ∪ (1 1 1 1 0 0), (0 0 0 0) ∪ (0 0 0) ∪ (1 1 1 0 1 0), (0 0 0 0) ∪ (0 0 0) ∪ (1 1 1 0 0 1), (0 0 0 0) ∪ (0 0 0) ∪ (1 1 0 0 1 1), (0 0 0 0) ∪ (0 0 0) ∪ (1 1 0 1 0 1), (0 0 0 0) ∪ (0 0 0) ∪ (1 0 1 1 1 0), (0 0 0 0) ∪ (0 0 0) ∪ (1 0 1 1 0 1), (0 0 0 0) ∪ (0 0 0) ∪ (1 0 0 1 1 1), (0 0 0 0) ∪ (0 0 0) ∪ (1 0 1 0 1 1), (0 0 0 0) ∪ (0 0 0) ∪ (0 1 1 1 1 0), (0 0 0 0) ∪ (0 0 0) ∪ (0 1 1 1 0 1), (0 0 0 0) ∪ (0 0 0) ∪ (0 1 1 0 1 1), (0 0 0 0) ∪ (0 0 0) ∪ (0 1 0 1 1 1), (0 0 0 0) ∪ (0 0 0) ∪ (1 1 1 1 1 1), (0 0 0 0) ∪ (0 0 0) ∪ (0 1 0 0 0 1), (0 0 0 0) ∪ (0 0 0) ∪ (0 1 0 1 0 0), (0 0 0 0) ∪ (0 0 0) ∪ (0 1 0 0 1 0), (0 0 0 0) ∪ (0 0 0) ∪ (0 1 1 0 0 0), (0 0 0 0) ∪ (0 0 0) ∪ (0 0 1 1 0 0), (0 0 0 0) ∪ (0 0 0) ∪ (0 0 1 0 1 0), (0 0 0 0) ∪ (0 0 0) ∪ (0 0 1 0 0 1), (0 0 0 0) ∪ (0 0 0) ∪ (0 0 0 110), (0 0 0 0) ∪ (0 0 0) ∪ (0 0 0 1 0 1), (0 0 0

0) ∪ (0 0 0) ∪ (0 0 0 0 1 1), (0 0 0 0) ∪ (0 0 0) ∪ (1 1 0 1 1 0), (1 1 0 0) ∪ (0 0 0) ∪ (0 0 0 0) and so on}.

As a simple exercise the reader is requested to find the number of elements in the parity check code given in the example 2.

Now we proceed on to define the notion of parity check n-code ($n \geq 4$).

**DEFINITION 2.1.10:** *Let us consider the parity check n-matrix ($n \geq 4$) $H_p$ = (1 1 1) ∪ (1 1 1 1) ∪ (1 1 1 1 1) ∪ … ∪ (1 11 1 1 1 1 1 1) = $H_1 \cup H_2 \cup \ldots \cup H_n$ where each $H_i$'s are distinct i.e. $H_i \neq H_j$, if $i \neq j$. The n-code obtained by $H_p$ the parity check n-matrix is a parity check n-code, $n \geq 4$.*

We just illustrate this by the following example.

***Example 2.1.20:*** Let us consider a parity check n-matrix (n = 5) given by H = $H_1 \cup H_2 \cup H_3 \cup H_4 \cup H_5$ = (1 1 1) ∪ (1 1 1 1 1) ∪ (1 1 1 1 1 1) ∪ (1 1 1 1) ∪ (1 1 1 1 1 1 1). We see H is associated with a parity check 5-code where C = { (0 0 0) ∪ (0 0 0 0 0) ∪ (0 0 0 0 0 0) ∪ (0 0 0 0) ∪ (0 0 0 0 0 0 0), (0 0 0) ∪ (0 0 0 0 0) ∪ … ∪ (1 1 0 0) ∪ (0 0 0 1 0 0 1), (1 1 0) ∪ (1 1 0 0 0) ∪ … ∪ (1 0 0 1) ∪ (1 1 0 0 0 0 0) and so on}. We see C is a parity check five code.

Interested reader can find the number of elements in C.

We just define the notion of binary Hamming bicode.

**DEFINITION 2.1.11:** *A binary Hamming bicode is a bicode C given by C = $C_{m_1} \cup C_{m_2}$ ($m_1 \neq m_2$) of length $n_1 = 2^{m_1} - 1$ and $n_2 = 2^{m_2} - 1$ with parity check bimatrix H = $H_1 \cup H_2$ where $H_i$ have all of its columns to consist of all non zero binary vectors of length $m_i$; i = 1, 2 (If $m_1 = m_2$ then it is important $H_1$ and $H_2$ generate different set of code words).*

Now we illustrate a binary Hamming bicode which is as follows.

***Example 2.1.21:*** Let C = $C_7 \cup C_{15}$, be a Hamming binary bicode C = $C_7 \cup C_{15}$ which has the following parity check bimatrix H = $H_1 \cup H_2$. $C_7$ has the following parity check matrix

$$H_1 = \begin{bmatrix} 0 & 0 & 1 & 1 & 1 & 0 & 1 \\ 0 & 1 & 0 & 1 & 0 & 1 & 1 \\ 1 & 0 & 0 & 0 & 1 & 1 & 1 \end{bmatrix}$$

and

$$H_2 = \begin{bmatrix} 0 & 0 & 1 & 0 & 1 & 1 & 1 & 0 & 0 & 0 & 0 & 1 & 1 & 1 & 1 \\ 0 & 1 & 0 & 0 & 1 & 0 & 0 & 1 & 0 & 1 & 1 & 1 & 0 & 1 & 1 \\ 0 & 0 & 0 & 1 & 0 & 1 & 0 & 1 & 1 & 0 & 1 & 1 & 1 & 0 & 1 \\ 1 & 0 & 0 & 0 & 0 & 0 & 1 & 0 & 1 & 1 & 1 & 0 & 1 & 1 & 1 \end{bmatrix};$$

is the parity check matrix of $C_{15}$. H = $H_1 \cup H_2$ is the parity check bimatrix associated with the Hamming binary bicode (Here $C_7$ is a C(7, 4) Hamming code and $C_{15}$ is a C (15, 11) Hamming code).

Now having seen an example of a binary Hamming bicode we proceed on to see some more properties of bicodes and various types of bicodes.

***Example 2.1.22:*** Now we can have also a Hamming bicode H = $H_1 \cup H_2$ where both are C(7, 4) codes but the parity check matrices are different; H = $H_1 \cup H_2$ where

$$H_1 = \begin{bmatrix} 0 & 0 & 1 & 1 & 1 & 0 & 1 \\ 0 & 1 & 0 & 1 & 0 & 1 & 1 \\ 1 & 0 & 0 & 0 & 1 & 1 & 1 \end{bmatrix}$$

and

$$H_2 = \begin{bmatrix} 1 & 1 & 0 & 0 & 1 & 1 & 0 \\ 1 & 1 & 1 & 0 & 0 & 0 & 1 \\ 0 & 1 & 0 & 1 & 0 & 1 & 1 \end{bmatrix}.$$

The set of bicodes are as follows: {(0 0 0 0 0 0 0) ∪ (0 0 0 0 0 0 0), (0 0 0 0 0 0 0) ∪ (1 0 0 0 1 1), (1 0 0 0 1 1 1) ∪ (0 0 0 0 0 0 0), (0 1 0 0 1 0 1) ∪ (0 0 0 0 0 0 0), (0 0 0 0 0 0 0) ∪ (0 1 0 0 1 1 1), (1 1 1 0 0 0 1) ∪ (1 1 1 0 0 0 1), (1 1 1 1 1 1 1) ∪ (1 1 1 1 1 1 1), …, (1 1 0 0 0 1 0) ∪ (1 1 0 0 0 1 1) and so on.

Next we proceed on to define the notion of binary Hamming tricode.

**DEFINITION 2.1.12:** *Let us consider a binary Hamming tricode given by a tricode* $C = C_{m_1} \cup C_{m_2} \cup C_{m_3}$ *where* $2^{m_i} - 1 = n_i$*; i = 1, 2, 3* $m_i \geq 2$*. The parity check trimatrix* $H = H_1 \cup H_2 \cup H_3$*, is such that we have each* $H_i$ *to have all of its columns to be non zero vectors of length* $m_i$*; i = 1, 2, 3.*

***Note:*** Even if $m_1 = m_2 = m_3$ still we can have a binary Hamming tricode, if each of the parity check matrices give way to a set of distinct codes. We can have a binary Hamming tricode provided $C_{m_1} \neq C_{m_2}$, $C_{m_1} \neq C_{m_3}$ and $C_{m_2} \neq C_{m_3}$. We illustrate this by a simple example.

***Example 2.1.23:*** Let $C = C_1 \cup C_7 \cup C_7$ be a binary Hamming tricode where $H = H_1 \cup H_2 \cup H_3$ is the associated parity check trimatrix of C which gives distinct sets of codes.

$$H = \begin{bmatrix} 1 & 1 & 1 & 0 & 1 & 0 & 0 \\ 1 & 0 & 0 & 0 & 1 & 1 & 1 \\ 1 & 0 & 1 & 1 & 0 & 0 & 1 \end{bmatrix} \cup$$

$$\begin{bmatrix} 0 & 0 & 0 & 1 & 1 & 1 & 1 \\ 0 & 1 & 1 & 1 & 0 & 1 & 0 \\ 1 & 1 & 0 & 1 & 0 & 0 & 1 \end{bmatrix} \cup$$

$$\begin{bmatrix} 1 & 0 & 1 & 0 & 1 & 1 & 0 \\ 1 & 0 & 1 & 1 & 0 & 0 & 1 \\ 0 & 1 & 1 & 1 & 1 & 0 & 0 \end{bmatrix}.$$

We see the parity check matrices associated with the trimatrix are distinct. Thus this gives a set of codes which are of same length and have the same collection of message symbols. Hence if one wishes to make use of the same set up to first 4 coordinates but with different elements in the 3 end coordinates in such cases we can make use of these Hamming trimatrix.

Next we proceed on to give yet another new Hamming trimatrix of different lengths.

***Example 2.1.24:*** Let us consider the associated Hamming trimatrix of a binary Hamming tricode C of different lengths; say $C = C_7 \cup C_{15} \cup C_7$ where C is given by the parity check Hamming trimatrix $H = H_1 \cup H_2 \cup H_3$ where

$$H_1 = \begin{bmatrix} 0 & 0 & 0 & 1 & 1 & 1 & 1 \\ 0 & 1 & 1 & 1 & 0 & 1 & 0 \\ 1 & 1 & 0 & 1 & 0 & 0 & 1 \end{bmatrix},$$

$$H_2 = \begin{bmatrix} 1 & 0 & 1 & 0 & 1 & 0 & 1 & 0 & 1 & 0 & 1 & 0 & 1 & 0 & 0 \\ 1 & 1 & 0 & 0 & 1 & 0 & 1 & 0 & 0 & 1 & 1 & 1 & 0 & 1 & 1 \\ 1 & 0 & 0 & 1 & 1 & 0 & 0 & 1 & 1 & 1 & 0 & 0 & 0 & 0 & 1 \\ 1 & 0 & 0 & 0 & 0 & 1 & 1 & 1 & 1 & 1 & 0 & 1 & 1 & 1 & 0 \end{bmatrix}$$

and

$$H_3 = \begin{bmatrix} 1 & 0 & 1 & 0 & 1 & 1 & 0 \\ 1 & 0 & 1 & 1 & 0 & 0 & 1 \\ 0 & 1 & 1 & 1 & 1 & 0 & 0 \end{bmatrix}.$$

We see this binary Hamming tricode is of different lengths.

It is left as an exercise for the interested reader to work with these Hamming tricodes.

Now we proceed on to define a binary Hamming n-code $n \geq 4$.

**DEFINITION 2.1.13:** *Let us consider a code* $C = C_{m_1} \cup C_{m_2} \cup \ldots \cup C_{m_n}$ *where* $m_i \geq 2$*, with parity check n-matrix* $H = H_1 \cup H_2 \cup \ldots \cup H_n$*; we call this code C to be a binary Hamming n-code;* $n \geq 4$*, it may so happen* $m_i = m_j$ $(i \neq j)$ *but* $H_i$ *and* $H_j$ *must be distinct matrices which generate distinct set of codes.*

Now we illustrate by an example a Hamming n-code, n = 5.

***Example 2.1.25:*** Let us consider the Hamming 5-code given by $C = C_7 \cup C_7 \cup C_{15} \cup C_7 \cup C_{15}$ generated by the parity check 5-matrix $H = H_1 \cup H_2 \cup H_3 \cup H_4 \cup H_5$ where

$$H_1 = \begin{bmatrix} 1 & 1 & 1 & 0 & 1 & 0 & 0 \\ 1 & 0 & 0 & 0 & 1 & 1 & 1 \\ 1 & 0 & 1 & 1 & 0 & 0 & 1 \end{bmatrix},$$

$$H_2 = \begin{bmatrix} 1 & 0 & 1 & 0 & 1 & 1 & 0 \\ 1 & 0 & 1 & 1 & 0 & 0 & 1 \\ 0 & 1 & 1 & 1 & 1 & 0 & 0 \end{bmatrix},$$

$$H_3 = \begin{bmatrix} 1 & 1 & 0 & 0 & 1 & 1 & 1 & 1 & 0 & 1 & 1 & 1 & 0 & 0 & 0 \\ 1 & 0 & 1 & 1 & 1 & 1 & 1 & 0 & 0 & 1 & 0 & 0 & 1 & 0 & 0 \\ 1 & 1 & 1 & 0 & 0 & 1 & 0 & 1 & 1 & 0 & 0 & 0 & 0 & 1 & 0 \\ 1 & 1 & 1 & 1 & 1 & 0 & 0 & 0 & 1 & 0 & 1 & 0 & 0 & 0 & 1 \end{bmatrix}$$

$$H_4 = \begin{bmatrix} 0 & 0 & 0 & 1 & 1 & 1 & 1 \\ 0 & 1 & 1 & 1 & 0 & 1 & 0 \\ 1 & 1 & 0 & 1 & 0 & 0 & 1 \end{bmatrix}$$

and

$$H_5 = \begin{bmatrix} 1 & 0 & 1 & 0 & 1 & 0 & 1 & 0 & 1 & 0 & 1 & 0 & 1 & 0 & 1 \\ 0 & 1 & 1 & 0 & 1 & 0 & 0 & 1 & 0 & 0 & 0 & 1 & 1 & 1 & 1 \\ 0 & 0 & 1 & 1 & 0 & 0 & 0 & 1 & 1 & 1 & 1 & 0 & 0 & 1 & 1 \\ 0 & 0 & 0 & 0 & 0 & 1 & 1 & 1 & 0 & 1 & 1 & 1 & 1 & 0 & 1 \end{bmatrix}.$$

Now we proceed on to define the notion of cyclic bicodes.

**DEFINITION 2.1.14:** *Let $C = C_1 \cup C_2$ be a linear $(n_1, k_1) \cup (n_2, k_2)$ bicode over $F_q$ (or $Z_2$). If each of $C_1$ and $C_2$ are cyclic linear codes then we call the bicode C to be a linear cyclic bicode.*

We illustrate a cyclic bicode by an example. We can also say a cyclic bicode will be generated by a bimatrix or a cyclic bicode can be associated with a parity check bimatrix. Before we give the definition of cyclic bicode using a bipolynomial we will give an example.

***Example 2.1.26:*** Let us consider the bicode generated by the generator bimatrix

$$G = G_1 \cup G_2$$

$$= \begin{bmatrix} 1 & 0 & 0 & 1 & 0 & 0 \\ 0 & 1 & 0 & 0 & 1 & 0 \\ 0 & 0 & 1 & 0 & 0 & 1 \end{bmatrix} \cup \begin{bmatrix} 1 & 1 & 1 & 0 & 1 & 0 & 0 \\ 0 & 1 & 1 & 1 & 0 & 1 & 0 \\ 0 & 0 & 1 & 1 & 1 & 0 & 1 \end{bmatrix}.$$

The message symbols associated with this cyclic bicode are (0 0 0) ∪ (0 0 0), (0 0 0) ∪ (1 0 0), (0 0 0) ∪ (0 1 0), (0 0 0) ∪ (0 0 1), (0 0 0) ∪ (1 1 0), (0 0 0) ∪ (0 1 1), (0 0 0) ∪ (1 0 1), (0 0 0) ∪ (1 1 1), (1 0 0) ∪ (0 0 0), (1 0 0) ∪ (1 0 0), … (1 1 1) ∪ (1 1 1). The cyclic bicodes generated by $G = G_1 \cup G_2$, is given by C = {(0 0 0 0 0 0) ∪ (0 0 0 0 0 0 0), (0 0 0 0 0 0) ∪ (1 1 1 0 1 0 0), (0 0 0 0 0 0) ∪ (0 1 1 1 0 1 0), (0 0 0 0 0 0) ∪ (0 0 1 1 1 0 1), …, (1 1 1 1 1 1) ∪ (1 0 1 0 0 1 1) and so on}.

Now one can define a cyclic bicode to be a code following a cyclic bishift.

The mapping $\eta : F_q^{n_1} \cup F_q^{n_2} \to F_q^{n_1} \cup F_q^{n_2}$ such that $(a_0, \ldots, a_{n_1-1}) \cup (b_0, \ldots, b_{n_2-1}) \to (a_{n_1-1}, a_0, \ldots, a_{n_1-2}) \cup (b_{n_2-1}, b_0, \ldots, b_{n_2-2})$ is a linear mapping called a cyclic bishift.

Now we proceed on to say a bicode $C = C_1 \cup C_2$ to be a cyclic bicode if and only if for all $\nu = \nu_1 \cup \nu_2$ in C such that if for all $\nu = (a_0, \ldots, a_{n_1-1}) \cup (b_0, b_1, \ldots b_{n_2-1})$ in $C = C_1 \cup C_2$ we have $\nu \in C$ implies $(a_{n_1-1}, a_0, \ldots, a_{n_1-2}) \cup (b_{n_2-1}, b_0, \ldots, b_{n_2-2})$ is in $C = C_1 \cup C_2$.

We now define the generator and parity check bimatrices of a cyclic bicode using generator bipolynomial.

Let us define $V = V_{n_1} \cup V_{n_2}$ to be a linear bispace where $V = \{v \in F_q[x] / \text{degree } v < n_1\} \cup \{u \in F_q[x] / \text{degree } u < n_2\} = \{v_0 + v_1x + \ldots + v_{n_1-1}x^{n_1-1} / v_i \in F_q, 0 \le i \le n_1-1\} \cup \{u_0 + u_1x + \ldots + u_{n_2-1}x^{n_2-1} / u_i \in F_q, 0 \le i \le n_2 - 1\}$.

We can define a biisomorphism of the two bispaces as

$T = T_1 \cup T_2$: $F_q^{n_1} \cup F_q^{n_2} \to V_{n_1} \cup V_{n_2}$;
$T : (v_0, v_1, \ldots, v_{n_1-1}) \cup (u_0, u_1, \ldots, u_{n_2-1}) \to (v_0 + v_1x + \ldots + v_{n_1-1}x^{n_1-1}) \cup (u_0 + u_1x + \ldots + u_{n_2-1}x^{n_2-1})$.

If $R = (F_q(x), +, .)$ is any polynomial ring with coefficients from $F_q$ then we can form $R / \langle x^{n_i-1} \rangle$ modulo the principal ideal generated by $x^{n_i-1}$ in R; i = 1, 2.

The bimapping $w = w_1 \cup w_2 : F_q^{n_1} \cup F_q^{n_2} \to R/\langle x^{n_1-1} \rangle \cup R/\langle x^{n_2-1} \rangle$ where

$$
\begin{aligned}
w(v, u) &= (w_1 \cup w_2)(v, u) \\
&= w_1(v) \cup w_2(u) \\
&= w_1(v_0\, v_1, \ldots, v_{n_1-1}) \\
&\quad \cup w_2(u_0, u_1, \ldots, u_{n_2-1}) \\
&= (v_0 + v_1x + \ldots + v_{n_1-1}x^{n_1-1}) \cup \\
&\quad (u_0 + u_1x + \ldots + u_{n_2-1}x^{n_2-1}).
\end{aligned}
$$

w can be checked to be a biisomorphism of the additive bigroup $F_q^{n_1} \cup F_q^{n_2}$ onto the factor bigroup of all bipolynomials of bidegree < $n_i$ over $F_q$, i = 1, 2 denoted by $V_{n_i}$ ; i = 1, 2. $V_{n_i}$ is also an algebra and $V_{n_1} \cup V_{n_2}$ is a bialgebra over $F_q$.

Let us define linear cyclic bicodes using generating bipolynomials. A polynomial is a bipolynomial of the form p(x) ∪ q(x) where p(x), q(x) ∈ $F_q$[x]. '∪' is just a symbol and nothing more. The only criteria that too in case of cyclic bicodes generated by bipolynomials is that we need them to be distinct i.e. p(x) ≠ q(x) for we need different sets of polynomials.

Using this concept of bipolynomial we generate a bicode which is cyclic.

**DEFINITION 2.1.15:** *Let g(x) = $g_1(x) \cup g_2(x)$ be a bipolynomial such that $g_i(x) \in V_{n_i}$ ; i = 1, 2 and $g_i(x) / (x^{n_i} - 1)$; i = 1, 2 with degree $g_i(x) = m_i < n_i$; i = 1, 2. Let C = $C_1 \cup C_2$ be a bicode with $C_i$ a $(n_i, k_i)$ code i.e. $k_i = n_i - m_i$ defined by the generator bimatrix G = $G_1 \cup G_2$; where*

$$G_1 = \begin{bmatrix} g_0^1 & g_1^1 & \cdots & g_{m_1}^1 & 0 & \cdots & 0 \\ 0 & g_0^1 & \cdots & g_{m_1-1}^1 & g_{m_1}^1 & \cdots & 0 \\ \vdots & \vdots & & \vdots & \vdots & & \vdots \\ 0 & 0 & & g_0^1 & g_1^1 & \cdots & g_{m_1}^1 \end{bmatrix}$$

$$= \begin{bmatrix} g_1 \\ xg_1 \\ \vdots \\ x^{k_1-1}g_1 \end{bmatrix}$$

*with* $g_1(x) = g_0^1 + g_0^1 x + \ldots + g_{m_1}^1 x^{m_1} \in F_{q_1}[x]$.

*Now*

$$G_2 = \begin{bmatrix} g_0^2 & g_1^2 & \cdots & g_{m2}^2 & 0 & \cdots & 0 \\ 0 & g_0^2 & \cdots & g_{m_2-1}^2 & g_{m_2}^2 & \cdots & 0 \\ \vdots & \vdots & & \vdots & \vdots & & \vdots \\ 0 & 0 & \cdots & g_0^2 & g_1^2 & \cdots & g_{m_2}^2 \end{bmatrix}$$

$$= \begin{bmatrix} g_2 \\ xg_2 \\ \vdots \\ x^{k_2-1}g_2 \end{bmatrix}$$

*with* $g_2(x) = g_0^2 + g_0^2 x + \ldots + g_{m_2}^2 x^{m_2} \in F_{q_2}[x]$.

*It is easily verified that the bimatrix* $G = G_1 \cup G_2$ *generates a cyclic bicode.*

Now we know if we have a generator matrix of a code we can obtain a parity check matrix; like wise a natural question would be can we find for a generator bimatrix obtained by bipolynomials a parity check bimatrix which also gives the same cyclic bicode.

The answer is yes and just we give a parity check bimatrix associated with the cyclic bicode whose generator matrix $G = G_1 \cup G_2$.

We see $h_i(x) = \dfrac{x^{n_i} - 1}{g_i(x)}$, i = 1, 2.

Let $h_i(x) = h^i_{k_i} x^{k_i} + h^i_{k_i-1} x^{k_i-1} + \cdots + h^i_1 x + h^i_0$; i = 1, 2. Now we obtain the parity check bimatrix $H = H_1 \cup H_2$ associated with the bipolynomial $h_1(x) \cup h_2(x)$.

$$H = H_1 \cup H_2 = \begin{bmatrix} 0 & \cdots & & 0 & h^1_{k_1} & \cdots & h^1_1 & h^1_0 \\ 0 & \cdots & 0 & h^1_{k_1} & h^1_{k_1-1} & \cdots & h^1_0 & 0 \\ \cdots & \cdots & \cdots & \cdots & \cdots & \cdots & \cdots & \cdots \\ \\ h^1_{k_1} & \cdots & \cdots & \cdots & h^1_0 & \cdots & 0 & 0 \end{bmatrix} \cup$$

$$\begin{bmatrix} 0 & \cdots & \cdots & 0 & h^2_{k_2} & \cdots & h^2_1 & h^2_0 \\ 0 & \cdots & 0 & h^2_{k_2} & h^2_{k_2-1} & \cdots & h^2_0 & 0 \\ & \cdots & & \cdots & \cdots & \cdots & \cdots & \cdots \\ \\ h^2_{k_2} & \cdots & \cdots & \cdots & h^2_0 & \cdots & 0 & 0 \end{bmatrix}$$

It is easily verified the parity check bimatrix for the cyclic bicode generated by $G = G_1 \cup G_2$ is $H = H_1 \cup H_2$.

Now we illustrate this by the following example.

***Example 2.1.27:*** Let us consider the cyclic bicode generated by the bipolynomial $(x^3 + x + 1) \cup (x^3 + x^2 + 1) = g_1(x) \cup g_2(x)$ with $g_1(x)/x^7 - 1$ and $g_2(x)/x^7 - 1$. The generator bimatrix associated with this bipolynomial is given by

$$G = G_1 \cup G_2$$

$$= \begin{bmatrix} 1 & 1 & 0 & 1 & 0 & 0 & 0 \\ 0 & 1 & 1 & 0 & 1 & 0 & 0 \\ 0 & 0 & 1 & 1 & 0 & 1 & 0 \\ 0 & 0 & 0 & 1 & 1 & 0 & 1 \end{bmatrix} \cup \begin{bmatrix} 1 & 0 & 1 & 1 & 0 & 0 & 0 \\ 0 & 1 & 0 & 1 & 1 & 0 & 0 \\ 0 & 0 & 1 & 0 & 1 & 1 & 0 \\ 0 & 0 & 0 & 1 & 0 & 1 & 1 \end{bmatrix}.$$

The parity check bipolynomial associated with the generator bipolynomial is given by $(x-1)\ (x^3 + x^2 + 1) \cup (x-1)\ (x^3 + x + 1) = (x^4 + x^2 + x + 1) \cup (x^4 + x^3 + x + 1)$.

The associated parity check bimatrix with the parity check bipolynomial is given by
$H = H_1 \cup H_2 =$

$$\begin{bmatrix} 0 & 0 & 1 & 0 & 1 & 1 & 1 \\ 0 & 1 & 0 & 1 & 1 & 1 & 0 \\ 1 & 0 & 1 & 1 & 1 & 0 & 0 \end{bmatrix} \cup \begin{bmatrix} 1 & 1 & 1 & 0 & 1 & 0 & 0 \\ 0 & 1 & 1 & 1 & 0 & 1 & 0 \\ 0 & 0 & 1 & 1 & 1 & 0 & 1 \end{bmatrix}.$$

The bicode is given by {(0 0 0 0 0 0 0) ∪ (0 0 0 0 0 0 0), (1 1 0 1 0 0 0) ∪ (1 0 0 0 1 0 1), (1 0 0 1 0 1 1) ∪ (1 1 1 1 1 1 1), (0 0 0 0 0 0 0) ∪ (1 1 1 1 1 1 1), (1 0 1 1 1 0 0) ∪ (1 1 0 0 0 1 0) and so on }.

We can also have cyclic bicodes of different lengths which is seen by the following example.

***Example 2.1.28:*** Let us consider a cyclic bicode generated by the bipolynomial

$$g = g_1(x) \cup g_2(x) = (x^3 + 1) \cup (x^4 + x^3 + x^2 + 1).$$

The generator bimatrix associated with the generator bipolynomial of the cyclic bicode $C = C_1 \cup C_2$ is given by $G = G_1 \cup G_2$ where

$$G_1 = \begin{bmatrix} 1 & 0 & 0 & 1 & 0 & 0 \\ 0 & 1 & 0 & 0 & 1 & 0 \\ 0 & 0 & 1 & 0 & 0 & 1 \end{bmatrix}$$

and

$$G_2 = \begin{bmatrix} 1 & 0 & 1 & 1 & 1 & 0 & 0 & 0 \\ 0 & 1 & 0 & 1 & 1 & 1 & 0 & 0 \\ 0 & 0 & 1 & 0 & 1 & 1 & 1 & 0 \\ 0 & 0 & 0 & 1 & 0 & 1 & 1 & 1 \end{bmatrix};$$

i.e. $G = \begin{bmatrix} 1 & 0 & 0 & 1 & 0 & 0 \\ 0 & 1 & 0 & 0 & 1 & 0 \\ 0 & 0 & 1 & 0 & 0 & 1 \end{bmatrix} \cup \begin{bmatrix} 1 & 0 & 1 & 1 & 1 & 0 & 0 & 0 \\ 0 & 1 & 0 & 1 & 1 & 1 & 0 & 0 \\ 0 & 0 & 1 & 0 & 1 & 1 & 1 & 0 \\ 0 & 0 & 0 & 1 & 0 & 1 & 1 & 1 \end{bmatrix}$.

The associated parity check bipolynomial of the given generator polynomial is $(x^3 + 1) \cup (x^3 + x^2 + 1)$. The corresponding parity check bimatrix is given by

$$H = H_1 \cup H_2 =$$

$$\begin{bmatrix} 0 & 0 & 1 & 0 & 0 & 1 \\ 0 & 1 & 0 & 0 & 1 & 0 \\ 1 & 0 & 0 & 1 & 0 & 0 \end{bmatrix} \cup \begin{bmatrix} 0 & 0 & 0 & 0 & 1 & 1 & 0 & 1 \\ 0 & 0 & 0 & 1 & 1 & 0 & 1 & 0 \\ 0 & 0 & 1 & 1 & 0 & 1 & 0 & 0 \\ 0 & 1 & 1 & 0 & 1 & 0 & 0 & 0 \\ 1 & 1 & 0 & 1 & 0 & 0 & 0 & 0 \end{bmatrix}.$$

The bicode which is cyclic is given by { (0 0 0 0 0 0) ∪ (0 0 0 0 0 0 0 0), (1 0 0 1 0 0) ∪ (1 0 1 1 1 0 0 0), (1 1 0 1 1 0) ∪ (1 0 1 1 1 0 0 0), (1 0 0 1 0 0) ∪ (1 1 1 0 0 1 0 0), (1 1 1 1 1 1) ∪ (1 1 0 1 1 1 0 1) and so on}, we see the lengths of the codes of the cyclic bicodes are 6 and 8. Thus we can also have cyclic bicodes of different lengths.

Now we proceed on to define the notion of cyclic tricodes before which we will just recall the new notion of tripolynomial.

**DEFINITION 2.1.16:** *Let R[x] be a polynomial ring. A tripolynomial $p(x) = p_1(x) \cup p_2(x) \cup p_3(x)$; where $p_i(x) \in R[x]$; i = 1, 2, 3. '∪' is only a symbol to denote a tripolynomial.*

For example $x^2 - 1 \cup (x^2 + 3)(x - 4) \cup x^7 + x^3 - 2x^2+1$ is a tripolynomial. If $p(x) = (x - 1) \cup (x^7 - 1) \cup (x^2 + 2x - 1)$ and $q(x) = x^2 - 1 \cup x^2+1 \cup x^8 - 3$ be two tripolynomials then

p(x) + q(x) = $\{(x - 1) \cup (x^7 - 1) \cup (x^2 + 2x - 1)\} + \{x^2 - 1 \cup x^2 + 1 \cup x^8 - 3\}$

= $\{(x - 1) + (x^2 - 1)\} \cup \{(x^7 - 1) + x^2 + 1\} \cup \{x^8 - 3 + x^2 + 2x - 1\}$

= $(x^2 + x - 2) \cup (x^7 + x^2) \cup (x^8 + x^2 + 2x - 4)$

is again a tripolynomial.

Likewise we can define the product

p(x) q(x) = $\{(x - 1) \cup (x^7 - 1) \cup (x^2 + 2x - 1)\} \times \{x^2 - 1 \cup x^2 + 1 \cup x^8 - 1\}$

= $(x - 1)(x^2 - 1) \cup (x^7 - 1)(x^2 + 1) \cup (x^2 + 2x - 1) \times (x^8 - 1)$

= $(x^3 - x^2 - x + 1) \cup (x^9 - x^7 + x^2 - 1) \cup (x^{10} + 2x^9 - x^8 - x^2 - 2x + 1)$

is again a tripolynomial.

Thus we can say $R_1[x] \cup R_2[x] \cup R_3[x]$ is a tripolynomial ring under polynomial addition and polynomial multiplication.

Now we define a cyclic tricode.

**DEFINITION 2.1.17:** *A tricode* $C = C(n_1, k_1) \cup C(n_2, k_2) \cup C(n_3, k_3)$ *where each* $C(n_i, k_i)$ *is a code i = 1, 2, 3. The only thing we demand is that each* $C(n_i, k_i)$ *is distinct i.e., the code* $C(n_i, k_i) \neq C(n_j, k_j)$*;* $i \neq j$*;* $1 \leq i, j \leq 3$*.*

We illustrate this by the following example.

***Example 2.1.29:*** C = C (5, 3) ∪ C (7, 3) ∪ C (6, 2) is a tricode.

**DEFINITION 2.1.18:** *Let* $C = C(n_1, k_1) \cup C(n_2, k_2) \cup C(n_3, k_3)$ *be a tricode. If each of* $C(n_i, k_i)$ *is a cyclic code; i = 1, 2, 3, then we call C to be a cyclic tricode.*

We illustrate with examples and define the generator and parity check trimatrices.

***Example 2.1.30:*** Let C = $C_1 \cup C_2 \cup C_3$ be a cyclic tricode where $C_i$ is generated by the generator matrix, $G_i$, i = 1, 2, 3; i.e., G = $G_1 \cup G_2 \cup G_3$ is the generator trimatrix given by

$$G_1 = \begin{bmatrix} 1 & 1 & 1 & 0 & 1 & 0 & 0 \\ 0 & 1 & 1 & 1 & 0 & 1 & 0 \\ 0 & 0 & 1 & 1 & 1 & 0 & 1 \end{bmatrix},$$

$$G_2 = \begin{bmatrix} 1 & 0 & 0 & 1 & 0 & 0 \\ 0 & 1 & 0 & 0 & 1 & 0 \\ 0 & 0 & 1 & 0 & 0 & 1 \end{bmatrix}$$

and

$$G_3 = \begin{bmatrix} 1 & 1 & 0 & 1 & 0 & 0 & 0 \\ 0 & 1 & 1 & 0 & 1 & 0 & 0 \\ 0 & 0 & 1 & 1 & 0 & 1 & 0 \\ 0 & 0 & 0 & 1 & 1 & 0 & 1 \end{bmatrix}.$$

The cyclic tricode given by G is as follows {(0 0 0 0 0 0 0) $\cup$ (0 0 0 0 0 0) $\cup$ (0 0 0 0 0 0 0), (1 1 1 0 1 0 0) $\cup$ (1 0 0 1 0 0) $\cup$ (1 1 0 1 0 0 0), (1 0 0 1 1 1 0) $\cup$ (1 0 0 1 0 0) $\cup$ (1 0 1 1 1 0 0), (1 0 1 0 0 1 1) $\cup$ (1 1 1 1 1 1) $\cup$ (1 0 0 1 0 1 1) and so on}.

Note the codes $C_1$ and $C_3$ are codes of same length 7 but are distinct. Thus C is a cyclic tricode generated by the generator trimatrix G = $G_1 \cup G_2 \cup G_3$.

Now we define the cyclic tricode in terms of the generator tripolynomial.

**DEFINITION 2.1.19:** *Let g(x) = $g^1(x) \cup g^2(x) \cup g^3(x)$ be a tripolynomial which generates a cyclic tricode, here each $g^i(x)$ divides $x^{n_i} - 1$, i = 1, 2, 3. Each $C_i$ is a cyclic ($n_i$, $k_i$) code i = 1, 2, 3 where C = $C_1 \cup C_2 \cup C_3$. $h^i(x) = \frac{x^{n_i} - 1}{g^i(x)}$ is the parity check polynomial of the cyclic code $C_i$, i = 1, 2, 3.*

*$h(x) = h^1(x) \cup h^2(x) \cup h^3(x)$ is the parity check tripolynomial related with the generator tripolynomial $g(x) = g^1(x) \cup g^2(x) \cup g^3(x)$.*

Now we illustrate a cyclic tricode together with its generator trimatrix and parity check trimatrix.

***Example 2.1.31:*** Let $C = C_1 \cup C_2 \cup C_3$ be a cyclic tricode generated by the tripolynomial $g(x) = g_1(x) \cup g_2(x) \cup g_3(x) = (x^3 + 1) \cup (x^3 + x^2 + 1) \cup (x^4 + x^3 + x^2 + 1)$.

The related parity check tripolynomial $h(x) = x^3 + 1 \cup x^4 + x^3 + x^2 + 1 \cup (x^3 + x^2 + 1)$. The generator trimatrix

$G = G_1 \cup G_2 \cup G_3$

$$= \begin{bmatrix} 1 & 0 & 0 & 1 & 0 & 0 \\ 0 & 1 & 0 & 0 & 1 & 0 \\ 0 & 0 & 1 & 0 & 0 & 1 \end{bmatrix} \cup$$

$$\begin{bmatrix} 1 & 0 & 1 & 1 & 0 & 0 & 0 \\ 0 & 1 & 0 & 1 & 1 & 0 & 0 \\ 0 & 0 & 1 & 0 & 1 & 1 & 0 \\ 0 & 0 & 0 & 1 & 0 & 1 & 1 \end{bmatrix} \cup$$

$$\begin{bmatrix} 1 & 0 & 1 & 1 & 1 & 0 & 0 \\ 0 & 1 & 0 & 1 & 1 & 1 & 0 \\ 0 & 0 & 1 & 0 & 1 & 1 & 1 \end{bmatrix}$$

generates the cyclic tricode. The related parity check trimatrix $H = H_1 \cup H_2 \cup H_3 =$

$$\begin{bmatrix} 0 & 0 & 1 & 0 & 0 & 1 \\ 0 & 1 & 0 & 0 & 1 & 0 \\ 1 & 0 & 0 & 1 & 0 & 0 \end{bmatrix} \cup$$

$$\begin{bmatrix} 0 & 0 & 1 & 1 & 1 & 0 & 1 \\ 0 & 1 & 1 & 1 & 0 & 1 & 0 \\ 1 & 1 & 1 & 0 & 1 & 0 & 0 \end{bmatrix} \cup$$

$$\begin{bmatrix} 0 & 0 & 0 & 1 & 1 & 0 & 1 \\ 0 & 0 & 1 & 1 & 0 & 1 & 0 \\ 0 & 1 & 1 & 0 & 1 & 0 & 0 \\ 1 & 1 & 0 & 1 & 0 & 0 & 0 \end{bmatrix},$$

is the associated parity check trimatrix of the cyclic tricode C = $C_1 \cup C_2 \cup C_3$. The cyclic tricode words are as follows: {(0 0 0 0 0 0) ∪ (0 0 0 0 0 0 0) ∪ (0 0 0 0 0 0 0), (0 0 0 0 0 0) ∪ (1 0 1 1 0 0 0) ∪ (0 0 0 0 0 0 0), (0 0 0 0 0 0) ∪ (0 0 0 0 0 0 0) ∪ (1 0 1 1 1 0 0), (1 0 0 1 0 0) ∪ ( 1 1 0 1 0 0) ∪ (1 1 0 0 1 0 1), (1 1 1 1 1 1) ∪ (1 1 0 1 0 0 1) ∪ (1 1 0 0 1 0 1), (1 1 1 1 1 1) ∪ (1 1 0 0 0 1 0) ∪ (1 1 1 0 0 1 0) and so on}.

Thus we have shown how a cyclic tricode is generated.

It is left as an exercise for the reader to determine the number of code words in the cyclic tricode given in the above example. If C = $C(n_1, k_1) \cup C(n_2, k_2) \cup C(n_3, k_3)$ is a cyclic tricode find the number of tricode words in C.

Now we show how many tricodes does the cyclic tricode C = C(4, 2) ∪ C(5, 2) ∪ C(6, 2), given the generator tripolynomial $g(x) = x^2 + 1 \cup x^4 + x^3 + x + 1 \cup x^4 + x^2 + 1$. The related generator trimatrix

$G = G_1 \cup G_2 \cup G_3 =$

$$\begin{bmatrix} 1 & 0 & 1 & 0 \\ 0 & 1 & 0 & 1 \end{bmatrix} \cup$$

$$\begin{bmatrix} 1 & 1 & 0 & 1 & 1 & 0 \\ 0 & 1 & 1 & 0 & 1 & 1 \end{bmatrix} \cup$$

$$\begin{bmatrix} 1 & 0 & 1 & 0 & 1 & 0 \\ 0 & 1 & 0 & 1 & 0 & 1 \end{bmatrix}.$$

The cyclic tricode is given by

{(0 0 0 0) ∪ (0 0 0 0 0 0) ∪ (0 0 0 0 0 0),
(0 0 0 0) ∪ (0 0 0 0 0 0) ∪ (1 0 1 0 1 0),
(0 0 0 0) ∪ (0 0 0 0 0 0) ∪ (0 1 0 1 0 1),
(0 0 0 0) ∪ (0 0 0 0 0 0) ∪ (1 1 1 1 1 1),
(0 0 0 0) ∪ (1 1 0 1 1 0) ∪ (0 0 0 0 0 0),
(0 0 0 0) ∪ (1 1 0 1 1 0) ∪ (1 0 1 0 1 0),
(0 0 0 0) ∪ (1 1 0 1 1 0) ∪ (0 1 0 1 0 1),
(0 0 0 0) ∪ (1 1 0 1 1 0) ∪ (1 1 1 1 1 1),
(0 0 0 0) ∪ (0 1 1 0 1 1) ∪ (0 0 0 0 0 0),
(0 0 0 0) ∪ (0 1 1 0 1 1) ∪ (1 0 1 0 1 0),
(0 0 0 0) ∪ (0 1 1 0 1 1) ∪ (1 0 1 0 1 0),
(0 0 0 0) ∪ (0 1 1 0 1 1) ∪ (1 1 1 1 1 1),
(0 0 0 0) ∪ (1 0 1 1 0 1) ∪ (0 0 0 0 0 0),
(0 0 0 0) ∪ (1 0 1 1 0 1) ∪ (1 0 1 0 1 0),
(0 0 0 0) ∪ (1 0 1 1 0 1) ∪ (0 1 0 1 0 1),
(0 0 0 0) ∪ (1 0 1 1 0 1) ∪ (1 1 1 1 1 1),
(1 0 1 0) ∪ (0 0 0 0 0 0) ∪ (0 0 0 0 0 0),
(1 0 1 0) ∪ (0 0 0 0 0 0) ∪ (1 0 1 0 1 0),
(1 0 1 0) ∪ (0 0 0 0 0 0) ∪ (0 1 0 1 0 1),
(1 0 1 0) ∪ (0 0 0 0 0 0) ∪ (1 1 1 1 1 1),
(1 0 1 0) ∪ (1 1 0 1 1 0) ∪ (0 0 0 0 0 0),
(1 0 1 0) ∪ (1 1 0 1 1 0) ∪ (1 0 1 0 1 0),
(1 0 1 0) ∪ (1 1 0 1 1 0) ∪ (0 1 0 1 0 1),
(1 0 1 0) ∪ (1 1 0 1 1 0) ∪ (1 1 1 1 1 1),
(1 0 1 0) ∪ (0 1 1 0 1 1) ∪ (0 0 0 0 0 0),
(1 0 1 0) ∪ (0 1 1 0 1 1) ∪ (1 0 1 0 1 0),
(1 0 1 0) ∪ (0 1 1 0 1 1) ∪ (0 1 0 1 0 1),
(1 0 1 0) ∪ (0 1 1 0 1 1) ∪ (1 1 1 1 1 1),

$(1\ 0\ 1\ 0) \cup (1\ 0\ 1\ 1\ 0\ 1) \cup (0\ 0\ 0\ 0\ 0\ 0)$,
$(1\ 0\ 1\ 0) \cup (1\ 0\ 1\ 1\ 0\ 1) \cup (1\ 0\ 1\ 0\ 1\ 0)$,
$(1\ 0\ 1\ 0) \cup (1\ 0\ 1\ 1\ 0\ 1) \cup (0\ 1\ 0\ 1\ 0\ 1)$,
$(1\ 0\ 1\ 0) \cup (1\ 0\ 1\ 1\ 0\ 1) \cup (1\ 1\ 1\ 1\ 1\ 1)$,
$(0\ 1\ 0\ 1) \cup (0\ 0\ 0\ 0\ 0\ 0) \cup (0\ 0\ 0\ 0\ 0\ 0)$,
$(0\ 1\ 0\ 1) \cup (0\ 0\ 0\ 0\ 0\ 0) \cup (1\ 0\ 1\ 0\ 1\ 0)$,
$(0\ 1\ 0\ 1) \cup (0\ 0\ 0\ 0\ 0\ 0) \cup (0\ 1\ 0\ 1\ 0\ 1)$,
$(0\ 1\ 0\ 1) \cup (0\ 0\ 0\ 0\ 0\ 0) \cup (1\ 1\ 1\ 1\ 1\ 1)$,
$(0\ 1\ 0\ 1) \cup (1\ 1\ 0\ 1\ 1\ 0) \cup (0\ 0\ 0\ 0\ 0\ 0)$,
$(0\ 1\ 0\ 1) \cup (1\ 1\ 0\ 1\ 1\ 0) \cup (1\ 0\ 1\ 0\ 1\ 0)$,
$(0\ 1\ 0\ 1) \cup (1\ 1\ 0\ 1\ 1\ 0) \cup (0\ 1\ 0\ 1\ 0\ 1)$,
$(0\ 1\ 0\ 1) \cup (1\ 1\ 0\ 1\ 1\ 0) \cup (1\ 1\ 1\ 1\ 1\ 1)$,
$(0\ 1\ 0\ 1) \cup (0\ 1\ 1\ 0\ 1\ 1) \cup (0\ 0\ 0\ 0\ 0\ 0)$,
$(0\ 1\ 0\ 1) \cup (0\ 1\ 1\ 0\ 1\ 1) \cup (1\ 0\ 1\ 0\ 1\ 0)$,
$(0\ 1\ 0\ 1) \cup (0\ 1\ 1\ 0\ 1\ 1) \cup (0\ 1\ 0\ 1\ 0\ 1)$,
$(0\ 1\ 0\ 1) \cup (0\ 1\ 1\ 0\ 1\ 1) \cup (1\ 1\ 1\ 1\ 1\ 1)$,
$(0\ 1\ 0\ 1) \cup (1\ 0\ 1\ 1\ 0\ 1) \cup (0\ 0\ 0\ 0\ 0\ 0)$,
$(0\ 1\ 0\ 1) \cup (1\ 0\ 1\ 1\ 0\ 1) \cup (1\ 0\ 1\ 0\ 1\ 0)$,
$(0\ 1\ 0\ 1) \cup (1\ 0\ 1\ 1\ 0\ 1) \cup (0\ 1\ 0\ 1\ 0\ 1)$,
$(0\ 1\ 0\ 1) \cup (1\ 0\ 1\ 1\ 0\ 1) \cup (1\ 1\ 1\ 1\ 1\ 1)$,
$(1\ 1\ 1\ 1) \cup (0\ 0\ 0\ 0\ 0\ 0) \cup (0\ 0\ 0\ 0\ 0\ 0)$,
$(1\ 1\ 1\ 1) \cup (0\ 0\ 0\ 0\ 0\ 0) \cup (1\ 1\ 1\ 1\ 1\ 1)$,
$(1\ 1\ 1\ 1) \cup (0\ 0\ 0\ 0\ 0\ 0) \cup (1\ 0\ 1\ 0\ 1\ 0)$,
$(1\ 1\ 1\ 1) \cup (0\ 0\ 0\ 0\ 0\ 0) \cup (0\ 1\ 0\ 1\ 0\ 1)$,
$(1\ 1\ 1\ 1) \cup (1\ 1\ 0\ 1\ 1\ 0) \cup (0\ 0\ 0\ 0\ 0\ 0)$,
$(1\ 1\ 1\ 1) \cup (1\ 1\ 0\ 1\ 1\ 0) \cup (1\ 1\ 1\ 1\ 1\ 1)$,
$(1\ 1\ 1\ 1) \cup (1\ 1\ 0\ 1\ 1\ 0) \cup (0\ 1\ 0\ 1\ 0\ 1)$,
$(1\ 1\ 1\ 1) \cup (1\ 1\ 0\ 1\ 1\ 0) \cup (1\ 0\ 1\ 0\ 1\ 0)$,
$(1\ 1\ 1\ 1) \cup (0\ 1\ 1\ 0\ 1\ 1) \cup (0\ 0\ 0\ 0\ 0\ 0)$,
$(1\ 1\ 1\ 1) \cup (0\ 1\ 1\ 0\ 1\ 1) \cup (1\ 1\ 1\ 1\ 1\ 1)$,
$(1\ 1\ 1\ 1) \cup (0\ 1\ 1\ 0\ 1\ 1) \cup (0\ 1\ 0\ 1\ 0\ 1)$,
$(1\ 1\ 1\ 1) \cup (0\ 1\ 1\ 0\ 1\ 1) \cup (1\ 0\ 1\ 0\ 1\ 0)$,
$(1\ 1\ 1\ 1) \cup (1\ 0\ 1\ 1\ 0\ 1) \cup (0\ 0\ 0\ 0\ 0\ 0)$,
$(1\ 1\ 1\ 1) \cup (1\ 0\ 1\ 1\ 0\ 1) \cup (0\ 1\ 0\ 1\ 0\ 1)$,
$(1\ 1\ 1\ 1) \cup (1\ 0\ 1\ 1\ 0\ 1) \cup (1\ 0\ 1\ 0\ 1\ 0)$,
and $(1\ 1\ 1\ 1) \cup (1\ 0\ 1\ 1\ 0\ 1) \cup (1\ 1\ 1\ 1\ 1\ 1)\}$,

The number of tricode words in this cyclic tricode is 64.

Thus if we have a cyclic tricode $C = C_1 \cup C_2 \cup C_3$ with $|C_1| = 2^{m_1}$, $|C_2| = 2^{m_2}$ and $|C_3| = 2^{m_3}$, then C has $2^{m_1+m_2+m_3}$ number of tricode words.

We now proceed on to define a cyclic n-code, $n \geq 4$.

**DEFINITION 2.1.20:** *A cyclic n-code $C = C(n_1, k_1) \cup C(n_2, k_2) \cup \ldots \cup C(n_n, k_n)$ is such that each $C(n_i, k_i)$ is a distinct cyclic $(n_i, k_i)$ code of length $n_i$ with $k_i$ message symbols, i = 1, 2, …, n. Thus every code word of the cyclic n code will be a n-tuple each tuple will be of length $n_i$; i = 1, 2, …, n.*

*The only demand made in this definition is that each of the cyclic codes $C(n_i, k_i)$ should be distinct, they can be of same length we do not bother about it but what we need is, $C(n_i, k_i)$ is a cyclic code different from every $C(n_j, k_j)$ if $i \neq j$; i = 1, 2, …, n and $1 \leq j, i \leq n$.*

Now having defined a cyclic n-code we give the definition of n-polynomials in a n-polynomial ring, $n \geq 4$.

**DEFINITION 2.1.21:** *Let R[x] be a polynomial ring. A polynomial $p(x) = p_1(x) \cup p_2(x) \cup \ldots \cup p_n(x)$ where $p_i(x) \in R[x]$, i = 1, 2, …, n is defined to be a n-polynomial. When n = 1 we just get the polynomial, n = 2 gives us the bipolynomial and n = 3 is the tripolynomial. When $n \geq 4$ we get the n-polynomial.*

***Example 2.1.32:*** Let $p(x) = p_1(x) + p_2(x) + p_3(x) + p_4(x) + p_5(x)$ be a 5-polynomial where $p(x) = (x^3 - 1) \cup (3x^2 + 7x + 1) \cup (8x^5 + 3x^3 + 2x + 1) \cup (7x^8 + 4x + 4) \cup (x^5 - 1)$ is in R[x].

**Note:** Even if the polynomial $p_i(x)$ in the n-polynomial are not distinct still we call them only as n-polynomials.

For instance $p(x) = x^7 \cup (x^8 - 1) \cup x^2 + 5x - 1 \cup (x^8 - 1) \cup x^7 + 5x + 1$ is a 5 polynomial. We see the polynomials $p_i(x)$ need not be distinct.

Now we define the notion of generator n-polynomial of a cyclic n-code.

**DEFINITION 2.1.22:** *Let us consider a cyclic n-code C = C($n_1$, $k_1$) ∪ C($n_2$, $k_2$) … ∪ C($n_n$, $k_n$) where each of the C($n_i$, $k_i$) are cyclic codes generated by the polynomial $g_i(x)$, i = 1, 2, …, n. Clearly each of the cyclic codes C($n_i$, $k_i$) must be distinct for otherwise the code C cannot be defined as the cyclic n-code. If each of the cyclic codes is to be distinct then we have the associated polynomials with them must also be distinct.*

*Let the n-polynomial which generates the cyclic n-code be denoted by g(x) = $g_1(x) \cup g_2(x) \cup \ldots \cup g_n(x)$. Each $g_i(x)$ is distinct, directly implies each of the generator matrices associated with the generator polynomial $g_i(x)$ are distinct. Thus let G denote the generator n-matrix, then G = $G_1 \cup G_2 \cup \ldots \cup G_n$ where each of the matrices $G_i$ are distinct i.e. $G_i = G_k$ if and only if i = k.*

*Now using each of the generator polynomials $g_i(x)$ we can obtain the parity check polynomial $h_i(x)$ given by $h_i(x) = \dfrac{x^{n_i} - 1}{g_i(x)}$; i = 1, 2, …, n. Thus the parity check n-polynomial h(x) is given by h(x) = $h_1(x) \cup h_2(x) \cup \ldots \cup h_n(x)$. Now the parity check n-matrix H associated with each of these parity check polynomial $h_i(x)$ i.e. with the parity check n-polynomial h(x) = $h_1(x) \cup h_2(x) \cup \ldots \cup h_n(x)$ is given by H = $H_1 \cup H_2 \cup \ldots \cup H_n$ where each matrix $H_i$ is the parity check matrix associated with the parity check polynomial $h_i(x)$; i = 1, 2, 3, …, n.*

We illustrate this definition by an example, when n = 4 i.e. the cyclic 4-code.

***Example 2.1.33:*** Let us consider a cyclic 4-code C = C(8, 4) ∪ C(3, 3) ∪ C(7, 3) ∪ C(7, 3), which is generated by the generator 4-polynomial g(x) = $x^4 + 1 \cup x^3 + 1 \cup x^3 + x + 1 \cup x^3 + x^2 + 1$ (we see each of the generator polynomials are distinct). The generator 4-matrix G associated with the cyclic 4 code is given by

$$G = G_1 \cup G_2 \cup G_3 \cup G_4 =$$

$$\begin{bmatrix} 1 & 0 & 0 & 0 & 1 & 0 & 0 & 0 \\ 0 & 1 & 0 & 0 & 0 & 1 & 0 & 0 \\ 0 & 0 & 1 & 0 & 0 & 0 & 1 & 0 \\ 0 & 0 & 0 & 1 & 0 & 0 & 0 & 1 \end{bmatrix} \cup$$

$$\begin{bmatrix} 1 & 0 & 0 & 1 & 0 & 0 \\ 0 & 1 & 0 & 0 & 1 & 0 \\ 0 & 0 & 1 & 0 & 0 & 1 \end{bmatrix} \cup$$

$$\begin{bmatrix} 1 & 0 & 1 & 1 & 0 & 0 & 0 \\ 0 & 1 & 0 & 1 & 1 & 0 & 0 \\ 0 & 0 & 1 & 0 & 1 & 1 & 0 \\ 0 & 0 & 0 & 1 & 0 & 1 & 1 \end{bmatrix} \cup$$

$$\begin{bmatrix} 1 & 1 & 0 & 1 & 0 & 0 & 0 \\ 0 & 1 & 1 & 0 & 1 & 0 & 0 \\ 0 & 0 & 1 & 1 & 0 & 1 & 0 \\ 0 & 0 & 0 & 1 & 1 & 0 & 1 \end{bmatrix}.$$

Now we give the parity check 4-polynomial h(x) and the associated parity check 4-matrix

$$\begin{aligned} h(x) &= h_1(x) \cup h_2(x) \cup h_3(x) \cup h_4(x) \\ &= (x^4 + 1) \cup (x^3 + 1) \cup (x^4 + x^2 + x + 1) \cup (x^4 + x^3 + x^2 + 1). \end{aligned}$$

The related parity check 4-matrix

$$H = H_1 \cup H_2 \cup H_3 \cup H_4 =$$

$$\begin{bmatrix} 0 & 0 & 0 & 1 & 0 & 0 & 0 & 1 \\ 0 & 0 & 1 & 0 & 0 & 0 & 1 & 0 \\ 0 & 1 & 0 & 0 & 0 & 1 & 0 & 0 \\ 1 & 0 & 0 & 0 & 1 & 0 & 0 & 0 \end{bmatrix} \cup$$

$$\begin{bmatrix} 0 & 0 & 1 & 0 & 0 & 1 \\ 0 & 1 & 0 & 0 & 1 & 0 \\ 1 & 0 & 0 & 1 & 0 & 0 \end{bmatrix} \cup$$

$$\begin{bmatrix} 0 & 0 & 1 & 0 & 1 & 1 & 1 \\ 0 & 1 & 0 & 1 & 1 & 1 & 0 \\ 1 & 0 & 1 & 1 & 1 & 0 & 0 \end{bmatrix} \cup$$

$$\begin{bmatrix} 0 & 0 & 1 & 1 & 1 & 0 & 1 \\ 0 & 1 & 1 & 1 & 0 & 1 & 0 \\ 1 & 1 & 1 & 0 & 1 & 0 & 0 \end{bmatrix}.$$

The following are the cyclic 4-codes;

{(0 0 0 0 0 0 0 0) ∪ (0 0 0 0 0 0) ∪ (0 0 0 0 0 0 0) ∪ (0 0 0 0 0 0 0), (0 0 0 0 0 0 0 0) ∪ (0 0 0 0 0 0) ∪ (0 0 0 0 0 0 0) ∪ (1 0 1 1 0 0 0), (0 0 0 0 0 0 0 0) ∪ (0 0 0 0 0 0) ∪ (1 1 0 1 0 0 0) ∪ (1 0 1 1 0 0 0), (0 0 0 0 0 0 0 0) ∪ (1 0 0 1 0 0) ∪ (1 1 0 1 0 0 0) ∪ (1 0 1 1 0 0 0), (0 0 0 0 0 0 0 0) ∪ (1 0 0 1 0 0) ∪ (1 1 0 1 0 0 0) ∪ (1 0 1 1 0 0 0), (1 0 0 0 1 0 0 0) ∪ (1 0 0 1 0 0) ∪ (0 0 0 0 0 0) ∪ (0 0 0 0 0 0 0), (0 0 1 1 0 0 1 1) ∪ (0 1 1 0 1 1) ∪ (0 1 1 1 0 0 1) ∪ (0 1 0 0 1 1 1) and so on}.

Interested reader is expected to find the number of code words in the above example.

Now we proceed on to define dual or orthogonal bicodes, orthogonal tricodes and the notion of orthogonal or dual n-codes.

**DEFINITION 2.1.23:** *Let C be a bicode given by $C = C_1 \cup C_2$, the dual bicode of C is defined to be $C^{\perp} = \{u_1 \mid u_1.v_1 = 0$ for all $v_1 \in C_1\} \cup \{u_2 \mid u_2.v_2 = 0$ for all $v_2 \in C_2\} = C_1^{\perp} \cup C_2^{\perp}$. If C is a $(k_1 \cup k_2)$ dimensional subbispace of the $n_1 \cup n_2$ dimensional vector space, the orthogonal vector space, i.e., the orthogonal complement is of dimension $(n_1 - k_1) \cup (n_2 - k_2)$. Thus if C is a $(n_1, k_1) \cup (n_2, k_2)$ bicode then $C^{\perp}$ is a $(n_1, n_1 - k_1) \cup (n_2, n_2 - k_2)$ bicode.*

*Further if $G = G_1 \cup G_2$ is the generator bimatrix and $H = H_1 \cup H_2$ is the parity check bimatrix then $C^{\perp}$ has generator bimatrix to be $H = H_1 \cup H_2$ and parity check bimatrix $G = G_1 \cup G_2$. Orthogonality of two bicodes can be expressed by $GH^T = HG^T = (0)$. That is $(G_1 \cup G_2) (H_1 \cup H_2)^T = (H_1 \cup H_2) (G_1 \cup G_2)^T = (G_1 \cup G_2) (H_1^T \cup H_2^T) = (H_1 \cup H_2) (G_1^T \cup G_2^T) = G_1H_1^T \cup G_2H_2^T = H_1G_1^T \cup H_2G_2^T = (0) \cup (0)$.*

Now using these properties we give an example of a bicode $C = C_1 \cup C_2$ and its orthogonal code $C^{\perp} = (C_1 \cup C_2)^{\perp} = C_1^{\perp} \cup C_2^{\perp}$.

***Example 2.1.34:*** Given the generator bimatrix of C a binary $(7, 3) \cup (4, 3)$ bicode; to find its dual or its orthogonal complement. Given

$$G = \begin{bmatrix} 0 & 0 & 0 & 1 & 1 & 1 & 1 \\ 0 & 1 & 1 & 0 & 0 & 1 & 1 \\ 1 & 0 & 1 & 0 & 1 & 0 & 1 \end{bmatrix} \cup \begin{bmatrix} 1 & 0 & 1 & 1 \\ 0 & 1 & 1 & 1 \\ 1 & 0 & 0 & 1 \end{bmatrix}.$$

The bicode given by the bimatrix G is as follows: $C = C_1 \cup C_2$ = { (0 0 0 0 0 0 0) ∪ (0 0 0 0), (0 0 0 0 0 0 0) ∪ (1 1 1 0), (0 0 0 0 0 0 0) ∪ (1 0 0 0), (0 0 0 0 0 0 0) ∪ (1 0 1 1), (0 0 0 0 0 0 0) ∪ (0 1 1 0), (0 0 0 0 0 0 0) ∪ (0 0 1 1), (0 0 0 0 0 0 0) ∪ (1 1 0 1), (0 0 0 0 0 0 0) ∪ (0 1 0 1), and so on} We can find $C^{\perp} = (C_1 \cup C_2)^{\perp} = C_1^{\perp} \cup C_2^{\perp}$.

The generator matrix of C serves as the parity check matrix of $C^{\perp}$. Thus using

$$G = \begin{bmatrix} 0 & 0 & 0 & 1 & 1 & 1 & 1 \\ 0 & 1 & 1 & 0 & 0 & 1 & 1 \\ 1 & 0 & 1 & 0 & 1 & 0 & 1 \end{bmatrix} \cup \begin{bmatrix} 1 & 0 & 1 & 1 \\ 0 & 1 & 1 & 1 \\ 1 & 0 & 0 & 1 \end{bmatrix}$$

as the parity check matrix i.e. G is taken as H. We find $C^{\perp}$ which is as follows {(0 0 0 0 0 0 0) ∪ (0 0 0 0), (0 0 0 0 0 0 0) ∪ (1 1 0 1), (1 0 0 0 0 1 1) ∪ (0 0 0 0), (1 0 0 0 0 1 1) ∪ (1 1 0 1), (0 1 0 0 1 0 1) ∪ (0 0 0 0), (0 1 0 0 1 0 1) ∪ (1 1 0 1), (0 0 1 0 1 1 0) ∪ (0 0 0 0), (0 0 1 0 1 1 0) ∪ (1 1 0 1), (0 0 0 1 1 1 1) ∪ (0 0 0 0), (0 0 0 1 1 1 1) ∪ (1 1 0 1), (1 1 0 0 0 0 1) ∪ (0 0 0 0), (1 1 0 0 0 0 1) ∪ (1 1 0 1), (1 0 1 0 1 0 1) ∪ (0 0 0 0), (1 0 1 0 1 0 1) ∪ (1 1 0 1), (1 0 0 1 1 0 0) ∪ (0 0 0 0), (1 0 0 1 1 0 0) ∪ (1 1 0 1), (0 1 1 0 0 1 1) ∪ (0 0 0 0), (0 1 1 0 0 1 1) ∪ (1 1 0 1), (0 1 0 1 0 1 0) ∪ (0 0 0 0), (0 1 0 1 0 1 0) ∪ (1 1 0 1), (0 0 1 1 1 0 0) ∪ (0 0 0 0), (0 0 1 1 1 0 0) ∪ (1 1 0 1), (1 1 1 0 0 0 0) ∪ (0 0 0 0), (1 1 1 0 0 0 0) ∪ (1 1 0 1), (1 1 0 1 0 0 1) ∪ (0 0 0 0), (1 1 0 1 0 0 1) ∪ (1 1 0 1), (1 0 1 1 0 1 0) ∪ (0 0 0 0), (1 0 1 1 0 1 0) ∪ (1 1 0 1), (0 1 1 1 0 1 0) ∪ (0 0 0 0), (0 1 1 1 0 1 0) ∪ (1 1 0 1), (1 1 1 1 1 1 1) ∪ (0 0 0 0), (1 1 1 1 1 1 1) ∪ (1 1 0 1) }

Now having defined the notion of the dual bicode of a bicode, C = $C_1 \cup C_2$ we now proceed on to define the notion of dual or orthogonal tricode of a tricode and then pass onto generalize it to n-codes (n ≥ 4).

**DEFINITION 2.1.24:** *Let us consider a tricode* $C = C_1 \cup C_2 \cup C_3 = C(n_1, k_1) \cup C(n_2, k_2) \cup C(n_3, k_3)$. *Let C be generated by the generator trimatrix* $G = G_1 \cup G_2 \cup G_3$. *The orthogonal complement of the tricode C or the dual tricode of C denoted by* $C^{\perp} = (C_1 \cup C_2 \cup C_3)^{\perp} = C_1^{\perp} \cup C_2^{\perp} \cup C_3^{\perp}$ *is associated with the trimatrix* $G = G_1 \cup G_2 \cup G_3$ *which acts as the parity check trimatrix i.e.,* $G = H = H_1 \cup H_2 \cup H_3$. *Thus the tricode having* $H = H_1 \cup H_2 \cup H_3$ *where G = H and* $H_i = G_i$; *i = 1, 2, 3 is a dual tricode of C and is got by using the parity check trimatrix H (= G) and the tricode* $C^{\perp} = C_1^{\perp} \cup C_2^{\perp} \cup C_3^{\perp}$ *is a* $C(n_1, n_1 - k_1) \cup C(n_2, n_2 - k_2) \cup C(n_3, n_3 - k_3)$ *tricode called the dual tricode of C.*

We illustrate this by the following example.

***Example 2.1.35:*** Let $C = C_1 \cup C_2 \cup C_3$ be a $C_1(7, 4) \cup C_2(7, 3) \cup C_3(6, 3)$ tricode having the generator trimatrix

$G = G_1 \cup G_2 \cup G_3 =$

$$\begin{bmatrix} 1 & 0 & 0 & 0 & 1 & 0 & 1 \\ 0 & 1 & 0 & 0 & 1 & 1 & 1 \\ 0 & 0 & 1 & 0 & 1 & 1 & 0 \\ 0 & 0 & 0 & 1 & 0 & 1 & 1 \end{bmatrix} \cup$$

$$\begin{bmatrix} 1 & 0 & 0 & 1 & 1 & 1 & \\ 0 & 1 & 1 & 0 & 0 & 1 & 1 \\ 1 & 0 & 1 & 0 & 1 & 0 & 1 \end{bmatrix} \cup$$

$$\begin{bmatrix} 1 & 0 & 0 & 0 & 1 & 1 \\ 0 & 1 & 0 & 1 & 0 & 1 \\ 0 & 0 & 1 & 1 & 1 & 0 \end{bmatrix}.$$

The tricode generated by the trimatrix $G = G_1 \cup G_2 \cup G_3$ is given by { (0 0 0 0 0 0 0) ∪ (0 0 0 0 0 0 0) ∪ (0 0 0 0 0 0), (0 0 0 0 0 0 0) ∪ (0 0 0 0 0 0 0) ∪ (1 0 0 0 1 1), (0 0 0 0 0 0 0) ∪ (0 0 0 0 0 0 0) ∪ (0 1 0 1 0 1), (0 0 0 0 0 0 0) ∪ (0 0 0 0 0 0 0) ∪ (0 0 1 1 1 0), (0 0 0 0 0 0 0) ∪ (0 0 0 0 0 0 0) ∪ (1 1 0 1 1 0), (0 0 0 0 0 0 0) ∪ (0 0 0 0 0 0 0) ∪ (0 1 1 0 1 1), (0 0 0 0 0 0 0) ∪ (0 0 0 0 0 0 0) ∪ (1 0 1 1 0 1), (0 0 0 0 0 0 0) ∪ (0 0 0 0 0 0 0) ∪ (1 1 1 0 0 0), (0 0 0 0 0 0 0) ∪ (1 0 0 1 1 1 1) ∪ (0 0 0 0 0 0), (0 0 0 0 0 0 0) ∪ (0 1 1 0 0 1 1) ∪ (0 0 0 0 0 0), (0 0 0 0 0 0 0) ∪ (1 0 1 0 1 0 1) ∪ (0 0 0 0 0 0), …, (1 0 0 0 1 0 1) ∪ (1 1 1 1 1 1 1) ∪ (0 1 0 1 0 0 1) ∪ (1 1 1 0 0 0)} now we find the dual of this tricode C using G the generator trimatrix of the tricode C as the parity check trimatrix of the orthogonal tricode $C^{\perp}$; i.e. G = H consequently we have $G_i = H_i$, i = 1, 2, 3. Thus

$G = H = H_1 \cup H_2 \cup H_3 =$

$$\begin{bmatrix} 1 & 0 & 0 & 0 & 1 & 0 & 1 \\ 0 & 1 & 0 & 0 & 1 & 1 & 1 \\ 0 & 0 & 1 & 0 & 1 & 1 & 0 \\ 0 & 0 & 0 & 1 & 0 & 1 & 1 \end{bmatrix} \cup$$

$$\begin{bmatrix} 1 & 0 & 0 & 1 & 1 & 1 & 1 \\ 0 & 1 & 1 & 0 & 0 & 1 & 1 \\ 1 & 0 & 1 & 0 & 1 & 0 & 1 \end{bmatrix} \cup$$

$$\begin{bmatrix} 1 & 0 & 0 & 0 & 1 & 1 \\ 0 & 1 & 0 & 1 & 0 & 1 \\ 0 & 0 & 1 & 1 & 1 & 0 \end{bmatrix};$$

is the parity trimatrix of the dual tricode $C^{\perp}$ of the tricode C. Clearly $C^{\perp}$ is a $C_1(7, 3) \cup C_1(7, 4) \cup C(6, 3)$ tricode associated with the parity check trimatrix $H = H_1 \cup H_2 \cup H_3$.

The dual tricode $C^{\perp}$ of C is given as follows: {(0 0 0 0 0 0 0) ∪ (0 0 0 0 0 0 0) ∪ (0 0 0 0 0 0), (0 0 0 0 0 0 0) ∪ (0 0 0 0 0 0 0) ∪ (0 1 1 1 0 0), (0 0 0 0 0 0 0) ∪ (0 0 0 0 0 0 0) ∪ (1 0 1 0 1 0), (0 0 0 0 0 0 0) ∪ (0 0 0 0 0 0 0) ∪ (1 1 0 0 0 1), (0 0 0 0 0 0 0) ∪ (0 0 0 0 0 0 0) ∪ (1 1 0 1 1 0), (0 0 0 0 0 0 0) ∪ (0 0 0 0 0 0 0) ∪ (0 1 1 0 1 1) ∪ (0 0 0 0 0 0 0) ∪ (0 0 0 0 0 0 0) ∪ (1 0 1 1 0 1), (0 0 0 0 0 0 0) ∪ (0 0 0 0 0 0 0) ∪ (0 0 0 1 1 1), and so on}.

Interested reader is left to the task of finding the orthogonal tricode $C^{\perp}$ of the given tricode C.

Now we proceed onto define the concept of dual n-code of an n-code C.

**DEFINITION 2.1.25:** *Let us consider the n-code* $C = C_1 (n_1, k_1) \cup C_2(n_2, k_2) \cup \ldots \cup C_n (n_n, k_n)$*; the dual code of C or the orthogonal code of C denoted by* $C^{\perp}$ *is defined as* $C^{\perp} = [(C_1 (n_1, k_1) \cup C_2 (n_2, k_2) \cup \ldots \cup C_n (n_n, k_n) ]^{\perp} = C_1^{\perp} (n_1, k_1) \cup C_2^{\perp} (n_2, k_2) \cup \ldots \cup C_n^{\perp} (n_n, k_n) = C_1 (n_1, n_1 - k_1) \cup C_2 (n_2, n_2 - k_2) \cup \ldots \cup C_n (n_n, n_n - k_n)$ *n-code. If* $G = G_1 \cup G_2 \cup \ldots \cup G_n$ *is the generator n-matrix of the n-code C then the orthogonal n-code* $C^{\perp}$ *has* $G = G_1 \cup G_2 \cup \ldots \cup G_n$ *to be its parity check n-matrix i.e. the generator n-matrix of G acts as the parity check n-matrix of* $C^{\perp}$ *and the parity check n-matrix of C acts as the generator n-matrix of* $C^{\perp}$*. We assume* $n \geq 4$*; when n = 1 we get the usual code C and its dual code* $C^{\perp}$ *when n = 2 we get the bicode* $C = C_1 \cup C_2$ *and the dual bicode* $C^{\perp} = C_1^{\perp} \cup C_2^{\perp}$*, when n = 3 we get the tricode* $C = C_1 \cup C_2 \cup C_3$ *and its dual* $C^{\perp} = C_1^{\perp} \cup C_2^{\perp} \cup C_3^{\perp}$*. When* $n \geq 4$ *we get the n-code and its dual n-code.*

Now are illustrate the dual n-code by an example when n = 5.

***Example 2.1.36:*** Consider the 5-code $C = C_1 \cup C_2 \cup C_3 \cup C_4 \cup C_5$ generated by the 5-matrix

$G = G_1 \cup G_2 \cup G_3 \cup G_4 \cup G_5 =$

$$\begin{bmatrix} 1 & 0 & 1 & 1 \\ 0 & 1 & 0 & 1 \end{bmatrix} \cup$$

$$\begin{bmatrix} 1 & 0 & 1 & 1 & 0 \\ 0 & 1 & 0 & 1 & 1 \end{bmatrix} \cup$$

$$\begin{bmatrix} 1 & 0 & 0 & 0 & 1 & 1 \\ 0 & 1 & 0 & 1 & 0 & 1 \\ 0 & 0 & 1 & 1 & 1 & 0 \end{bmatrix} \cup$$

$$\begin{bmatrix} 1 & 0 & 1 & 1 \\ 0 & 1 & 1 & 1 \\ 1 & 0 & 0 & 1 \end{bmatrix} \cup$$

$$\begin{bmatrix} 1 & 0 & 0 & 1 & 1 & 0 & 1 \\ 0 & 1 & 0 & 1 & 0 & 1 & 1 \\ 0 & 0 & 1 & 0 & 1 & 1 & 1 \end{bmatrix}.$$

Here $C_1$ = C(4, 2) code, $C_2$ = C(5, 2) code, $C_3$ = C(6, 3) code, $C_4$ = C(4, 3) code and $C_5$ = C(7, 3) code. The 5-codes generated by the 5-matrix G is as follows C = { (0 0 0 0) ∪ (0 0 0 0 0) ∪ (0 0 0 0 0 0) ∪ (0 0 0 0) ∪ (0 0 0 0 0 0 0), (1 0 1 1) ∪ (0 0 0 0 0) ∪ (0 0 0 0 0 0) ∪ (0 0 0 0) ∪ (0 0 0 0 0 0 0), (0 1 0 1) ∪ (0 0 0 0 0) ∪ (0 0 0 0 0 0) ∪ (0 0 0 0) ∪ (0 0 0 0 0 0 0), (1 1 1 0) ∪ (0 0 0 0 0) ∪ (0 0 0 0 0 0) ∪ (0 0 0 0) ∪ (0 0 0 0 0 0 0), (1 0 1 1) ∪ (1 0 1 1 0) ∪ (0 0 0 0 0 0) ∪ (0 0 0 0) ∪ (0 0 0 0 0 0 0), (1 0 1 1) ∪ (0 1 0 1 1) ∪ (0 0 0 0 0 0) ∪ (0 0 0 0) ∪ (0 0 0 0 0 0 0), (1 0 1 1) ∪ (1 1 1 0 1) ∪ (0 0 0 0 0 0) ∪ (0 0 0 0) ∪ (0 0 0 0 0 0 0), (1 0 1 1) ∪ (1 0 1 1 0) ∪ (1 0 0 0 1 1) ∪ (0 0 0 0) ∪ (0 0 0 0 0 0 0) and so on}.

Now $C^{\perp}$ the orthogonal 5-code of the 5-code works with G the generator 5-matrix of the 5-code as the parity check 5-matrix. Thus the parity check 5-matrix of the dual (orthogonal), 5-code $C^{\perp}$ is

$$H = \begin{bmatrix} 1 & 0 & 1 & 1 \\ 0 & 1 & 0 & 1 \end{bmatrix} \cup$$

$$\begin{bmatrix} 1 & 0 & 1 & 1 & 0 \\ 0 & 1 & 0 & 1 & 1 \end{bmatrix} \cup$$

$$\begin{bmatrix} 1 & 0 & 0 & 0 & 1 & 1 \\ 0 & 1 & 0 & 1 & 0 & 1 \\ 0 & 0 & 1 & 1 & 1 & 0 \end{bmatrix} \cup$$

$$\begin{bmatrix} 1 & 0 & 1 & 1 \\ 0 & 1 & 1 & 1 \\ 1 & 0 & 0 & 1 \end{bmatrix} \cup$$

$$\begin{bmatrix} 1 & 0 & 0 & 1 & 1 & 0 & 1 \\ 0 & 1 & 0 & 1 & 0 & 1 & 1 \\ 0 & 0 & 1 & 0 & 1 & 1 & 1 \end{bmatrix}.$$

The 5-code associated with the parity check 5-matrix H is $C^{\perp}$ given by

$C^{\perp}$ = {(0 0 0 0) ∪ (0 0 0 0 0) ∪ (0 0 0 0 0 0) ∪ (0 0 0 0) ∪ (0 0 0 0 0 0 0), (1 0 1 0) ∪ (0 0 0 0 0) ∪ (0 0 0 0 0 0) ∪ (0 0 0 0) ∪ (0 0 0 0 0 0 0), (0 1 1 1) ∪ (0 0 0 0 0) ∪ (0 0 0 0 0 0) ∪ (0 0 0 0) ∪ (0 0 0 0 0 0 0), (1 1 0 1) ∪ (0 0 0 0 0) ∪ (0 0 0 0 0 0) ∪ (0 0 0 0) ∪ (0 0 0 0 0 0 0), (1 0 1 0) ∪ (1 0 0 1 1) ∪ (0 0 0 0 0 0) ∪ (0 0 0 0) ∪ (0 0 0 0 0 0 0), (0 1 1 1) ∪ (0 1 0 0 1) ∪ (0 0 0 0 0 0) ∪ (0 0 0 0) ∪ (0 0 0 0 0 0 0), (1 1 0 1) ∪ (1 1 1 0 1) ∪ (0 0 0 0 0 0) ∪ (0 0 0 0) ∪ (0 0 0 0 0 0 0), (1 0 1 0) ∪ (1 0 0 1 1) ∪ (0 0 0 0 0 0) ∪ (1 1 0 1) ∪ (0 0 0 0 0 0 0) and so on}.

Thus we have shown how the orthogonal n-code of a given n code is obtained. It is also left for the reader to work with any n-code; n ≥ 4 and obtain the dual n-code.

A code is said to be reversible if $(a_0, a_1, \ldots, a_{n-1}) \in C$ holds if and only if $(a_{n-1}, a_{n-2}, \ldots, a_1\ a_0) \in C$. In general the reversible bicode and n-code can be defined.

**DEFINITION 2.1.26:** *A bicode $C = C_1 \cup C_2$ is said to be a reversible bicode , i.e., $(a_0\ a_1 \ldots a_{n-1}) \cup (b_0\ b_1 \ldots b_{m-1}) \in C = C_1 \cup C_2$ if and only if $(a_{n-1}\ a_{n-2}\ a_{n-3} \ldots a_1\ a_0) \cup (b_{m-1}\ b_{m-2}\ b_{m-3} \ldots b_1\ b_0) \in C = C_1 \cup C_2$.*

*Likewise if $C = C_1 \cup C_2 \cup C_3$ be a tricode then it is said to be reversible i.e., $(a_0\ a_1 \ldots a_{n-1}) \cup (b_0\ b_1 \ldots b_{m-1}) \cup (c_0\ c_1 \ldots c_{t-1}) \in C$ if and only if $(a_{n-1}\ a_{n-2} \ldots a_1\ a_0) \cup (b_{m-1}\ b_{m-2} \ldots b_1\ b_0) \cup (c_{t-1}\ c_{t-2} \ldots c_1\ c_0) \in C$.*

*The tricode generated by (1 1 1 1) $\cup$ (1 1 1 1 1 1) $\cup$ (1 1 1 1 1) is a reversible tricode.*

Is a cyclic n-code a reversible n-code? (n ≥ 2). To be more precise is a cyclic code a reversible code? Recall a code C is said to be self orthogonal if $C \subseteq C^{\perp}$.

We give the definition of self orthogonal n-code.

**DEFINITION 2.1.27:** *A bicode $C = C_1 \cup C_2$ is said to be self orthogonal bicode if each of $C_1$ and $C_2$ are self orthogonal codes and we denote it by $C \subseteq C^{\perp}$ i.e. $C = C_1 \cup C_2 \subseteq C^{\perp} = C_1^{\perp} \cup C_2^{\perp}$.*

Note: If only one of $C_1$ or $C_2$ alone is a self orthogonal code of the bicode $C = C_1 \cup C_2$ then we define the new notion of semi self orthogonal bicode.

**DEFINITION 2.1.28:** *Let $C = C_1 \cup C_2$ be a bicode we say C is said to be semi self orthogonal bicode if in C only one of $C_1$ or $C_2$ is self orthogonal. We denote a semi self orthogonal bicode by $C_S^{\perp} = C_1^{\perp} \cup C_2$, $C_1 \subseteq C_1^{\perp}$ or $C_S^{\perp} = C_1 \cup C_2^{\perp}$ and $C_2 \subseteq C_2^{\perp}$.*

Now we proceed onto define the notion of self orthogonal tricode and generalize the notion to n-codes.

**DEFINITION 2.1.29:** *Let $C = C_1 \cup C_2 \cup C_3$ be a tricode, we say C is a self orthogonal tricode if each of codes in C i.e. each $C_i$ (i = 1, 2, 3) are self orthogonal and this is denoted by $C \subseteq C^{\perp}$ i.e. $C_1 \cup C_2 \cup C_2 \subseteq (C_1 \cup C_2 \cup C_3)^{\perp} = C_1^{\perp} \cup C_2^{\perp} \cup C_3^{\perp}$.*

*If even one of the $C_i$'s in C is not self orthogonal (i = 1, 2, 3) then we say C is not self orthogonal tricode but C is only a semi self orthogonal tricode i.e. in case of semi self orthogonal tricode we denote it by $C_S^{\perp} = C_1^{\perp} \cup C_2^{\perp} \cup C_3$ where $C_1 \subseteq C_1^{\perp}$ and $C_2 \subseteq C_2^{\perp}$ but $C_3 \not\subseteq C_3^{\perp}$ or $C_S^{\perp} = C_1^{\perp} \cup C_2 \cup C_3^{\perp}$ where $C_1 \subseteq C_1^{\perp}$, $C_2 \not\subseteq C_2^{\perp}$ and $C_3 \subseteq C_3^{\perp}$ or $C_S^{\perp} = C_1^{\perp} \cup C_2 \cup C_3$ i.e. $C_1$ alone is self orthogonal but $C_2$ and $C_3$ are not self orthogonal.*

**DEFINITION 2.1.30:** *Let $C = C_1 \cup C_2 \cup \ldots \cup C_n$ be a n-code; C is said to be self orthogonal if each of the n codes in C is self orthogonal i.e. each $C_i^{\perp} \subseteq C_i$; i = 1, 2, 3, …, n and denoted by $C \subseteq C^{\perp}$ i.e. $C_1 \cup C_2 \cup \ldots \cup C_n \subseteq (C_1 \cup C_2 \cup \ldots \cup C_n)^{\perp} = C_1^{\perp} \cup C_2^{\perp} \cup \ldots \cup C_n^{\perp}$. Even if one of the $C_i$'s in C is not self orthogonal then we do not call the n-code to be a self orthogonal n-code.*

*However if in the n-code $C = C_1 \cup C_2 \cup \ldots \cup C_n$ some of the codes are self orthogonal and some of them are not self orthogonal then we call the n-code to be a semi self orthogonal n-code.*

Now we illustrate some examples of self orthogonal n-code, tricode and a bicode ($n \geq 4$).

***Example 2.1.37:*** Let us consider the bicode $C = C_1 \cup C_2$ where $C_1$ is a C(4, 2) code, generated by $x^2 + 1$ and $C_2$ is a C(6, 3) code generated by the polynomial $x^3 + 1$. Clearly $C = C_1 \cup C_2$ is a self orthogonal bicode.

***Example 2.1.38:*** Let $C = C_1 \cup C_2 \cup C_3$ be a tricode. $C_1$ the code generated by $x^3 + 1$, $C_2$ = C(8, 4) is the code generated by $x^4 + 1$ and $C_3$ = C(4, 2) is the code generated by $x^2 + 1$. Thus $C = C_1 \cup C_2 \cup C_3$ is self orthogonal tricode.

***Example 2.1.39:*** Let us consider the 6-code C generated by the 6-polynomial $\langle (x^4+1) \cup (x^3+1) \cup (x^5+1) \cup (x^6+1) \cup (x^2+1) \cup (x^7+1) \rangle$ where $C = C_1 \cup C_2 \cup C_3 \cup C_4 \cup C_5 \cup C_6$ with $C_1$ = C(8, 4), code generated by $x^4+1$, $C_2$ = C(6, 3) code generated by the polynomial $x^3+1$, $C_3$ is a code having the generator polynomial $(x^5+1)$, $C_4$ is a code with generator polynomial $x^6+1$, the code $C_5$ is generated by the polynomial $x^2+1$. Finally the polynomial $x^7 + 1$ generates the code $C_6$.

It is left as an exercise for the reader to verify that the 6-code is a self orthogonal code with the generator 6-matrix $G = G_1 \cup G_2 \cup G_3 \cup G_4 \cup G_5 \cup G_6$

$$= \begin{bmatrix} 1 & 0 & 0 & 0 & 1 & 0 & 0 & 0 \\ 0 & 1 & 0 & 0 & 0 & 1 & 0 & 0 \\ 0 & 0 & 1 & 0 & 0 & 0 & 1 & 0 \\ 0 & 0 & 0 & 1 & 0 & 0 & 0 & 1 \end{bmatrix} \cup \begin{bmatrix} 1 & 0 & 0 & 1 & 0 & 0 \\ 0 & 1 & 0 & 0 & 1 & 0 \\ 0 & 0 & 1 & 0 & 0 & 1 \end{bmatrix} \cup$$

$$\begin{bmatrix} 1 & 0 & 0 & 0 & 0 & 1 & 0 & 0 & 0 & 0 \\ 0 & 1 & 0 & 0 & 0 & 0 & 1 & 0 & 0 & 0 \\ 0 & 0 & 1 & 0 & 0 & 0 & 0 & 1 & 0 & 0 \\ 0 & 0 & 0 & 1 & 0 & 0 & 0 & 0 & 1 & 0 \\ 0 & 0 & 0 & 0 & 1 & 0 & 0 & 0 & 0 & 1 \end{bmatrix} \cup$$

$$\begin{bmatrix} 1 & 0 & 0 & 0 & 0 & 0 & 1 & 0 & 0 & 0 & 0 & 0 \\ 0 & 1 & 0 & 0 & 0 & 0 & 0 & 1 & 0 & 0 & 0 & 0 \\ 0 & 0 & 1 & 0 & 0 & 0 & 0 & 0 & 1 & 0 & 0 & 0 \\ 0 & 0 & 0 & 1 & 0 & 0 & 0 & 0 & 0 & 1 & 0 & 0 \\ 0 & 0 & 0 & 0 & 1 & 0 & 0 & 0 & 0 & 0 & 1 & 0 \\ 0 & 0 & 0 & 0 & 0 & 1 & 0 & 0 & 0 & 0 & 0 & 1 \end{bmatrix} \cup$$

$$\begin{bmatrix} 1 & 0 & 1 & 0 \\ 0 & 1 & 0 & 1 \end{bmatrix} \cup$$

$$\begin{bmatrix} 1 & 0 & 0 & 0 & 0 & 0 & 0 & 1 & 0 & 0 & 0 & 0 & 0 & 0 \\ 0 & 1 & 0 & 0 & 0 & 0 & 0 & 0 & 1 & 0 & 0 & 0 & 0 & 0 \\ 0 & 0 & 1 & 0 & 0 & 0 & 0 & 0 & 0 & 1 & 0 & 0 & 0 & 0 \\ 0 & 0 & 0 & 1 & 0 & 0 & 0 & 0 & 0 & 0 & 1 & 0 & 0 & 0 \\ 0 & 0 & 0 & 0 & 1 & 0 & 0 & 0 & 0 & 0 & 0 & 1 & 0 & 0 \\ 0 & 0 & 0 & 0 & 0 & 1 & 0 & 0 & 0 & 0 & 0 & 0 & 1 & 0 \\ 0 & 0 & 0 & 0 & 0 & 0 & 1 & 0 & 0 & 0 & 0 & 0 & 0 & 1 \end{bmatrix}.$$

## 2.2 Error Detection and Error Correction in n-codes ($n \geq 2$)

Now having defined the notion of n-codes ($n \geq 4$) we proceed on to define the methods of decoding and error detection. We know if we have any n-code $C = C_1 \cup C_2 \cup \ldots \cup C_n$ ($n \geq 4$), associated with C is the generator n-matrix $G = G_1 \cup G_2 \cup \ldots \cup G_n$ and the parity check n-matrix $H = H_1 \cup H_2 \cup \ldots \cup H_n$, then G and H are related in such a fashion that $GH^T = (0) \cup (0) \cup \ldots \cup (0)$. Now any n-code word in a n-code would be of the form $(x_1^1 x_2^1 \ldots x_{n_1}^1) \cup (x_1^2 x_2^2 \ldots x_{n_2}^2) \cup \ldots \cup (x_1^n x_2^n \ldots x_{n_n}^n)$.

Now the operation of any n-code word on the parity check n-matrix H is carried out as follows:

$$\begin{aligned}
Hx^T &= H[(x_1^1 x_2^1 \ldots x_{n_1}^1) \cup (x_1^2 x_2^2 \ldots x_{n_2}^2) \cup \ldots \cup (x_1^n x_2^n \ldots x_{n_n}^n)]^T \\
&= H[(x_1^1 x_2^1 \ldots x_{n_1}^1)^T \cup (x_1^2 x_2^2 \ldots x_{n_2}^2)^T \cup \ldots \\
&\quad \cup (x_1^n x_2^n \ldots x_{n_n}^n)^T] \\
&= [H_1 \cup H_2 \cup \ldots \cup H_n]
\end{aligned}$$

$$\begin{aligned}
&[(x_1^1 x_2^1 \ldots x_{n_1}^1)^T \cup (x_1^2 x_2^2 \ldots x_{n_2}^2)^T \cup \ldots \cup (x_1^n x_2^n \ldots x_{n_n}^n)^T] \\
&= H_1 (x_1^1 x_2^1 \ldots x_{n_1}^1)^T \cup H_2 (x_1^2 x_2^2 \ldots x_{n_2}^2)^T \cup \ldots \cup \\
&\quad H_n (x_1^n x_2^n \ldots x_{n_n}^n)^T. \\
&= (0) \cup (0) \cup \ldots \cup (0) \qquad \text{(I)}
\end{aligned}$$

Equation (I) implies x is a n-code word of C associated with the parity check n-matrix H. Even if one of $H_i x_i^T \neq (0)$, $1 \leq i \leq n$, then we conclude the code word is not from C. Thus we say $x \in C = C_1 \cup C_2 \cup \ldots \cup C_n$ if and only if $Hx^T = (0) \cup (0) \cup \ldots \cup (0)$ otherwise $x \notin C$ and error has occurred during transmission of the n-code.

Before we proceed further we define the notion of Hamming n-distance of two n-code words x and y from the n-code C. Let

$$x = (x_1^1 x_2^1 \ldots x_{n_1}^1) \cup (x_1^2 x_2^2 \ldots x_{n_2}^2) \cup \ldots \cup (x_1^n x_2^n \ldots x_{n_n}^n) \in C$$

where each $(x_1^i x_2^i \ldots x_{n_i}^i) \in C_i$ with $C_i$ a $(n_i, k_i)$ code.

The Hamming n-distance between x and y is defined to be a n-tuple of numbers $(t_1, \ldots, t_n)$ where

$$t_i \quad = \quad d\{(x_1^i x_2^i \ldots x_{n_i}^i),\ (y_1^i y_2^i \ldots y_{n_i}^i)\};\ 1 \le i \le n;$$

where $y = (y_1^1 y_2^1 \ldots y_{n_1}^1) \cup (y_1^2 y_2^2 \ldots y_{n_2}^2) \cup \ldots \cup (y_1^n y_2^n \ldots y_{n_n}^n)$. So

$$\begin{aligned} d(x, y) &= (d[(x_1^1 x_2^1 \ldots x_{n_1}^1), (y_1^1 y_2^1 \ldots y_{n_1}^1)], \\ &\quad d[(x_1^2 x_2^2 \ldots x_{n_2}^2), (y_1^2 y_2^2 \ldots y_{n_2}^2)], \ldots, \\ &\quad d[(x_1^n x_2^n \ldots x_{n_n}^n), (y_1^n y_2^n \ldots y_{n_n}^n)]) \\ &= (t_1, t_2, \ldots, t_n)\ . \\ d(x, x) &= (0, 0, \ldots, 0). \end{aligned}$$

when

$$\begin{aligned} x &= (1\ 1\ 0\ 1\ 0) \cup (1\ 1\ 1\ 1\ 1\ 1\ 1) \cup (0\ 0\ 0\ 0) \cup \\ &\quad (1\ 0\ 1\ 0\ 1\ 0) \cup (0\ 1\ 0\ 1\ 0\ 1\ 0\ 1\ 0\ 1) \cup \\ &\quad (1\ 1\ 1\ 0\ 0\ 0\ 1\ 1\ 0) \in C \end{aligned}$$

and

$$\begin{aligned} y &= (0\ 1\ 0\ 1\ 0) \cup (1\ 0\ 1\ 0\ 1\ 0\ 1) \cup (1\ 1\ 0\ 0) \cup \\ &\quad (1\ 1\ 1\ 1\ 0\ 0) \cup (1\ 0\ 1\ 0\ 1\ 1\ 1\ 1\ 0\ 1) \cup \\ &\quad (1\ 1\ 1\ 0\ 0\ 1\ 0\ 0\ 0) \in C. \end{aligned}$$

Now

$$\begin{aligned} d(x, y) &= \{d[(1\ 1\ 0\ 1\ 0\ ), (0\ 1\ 0\ 1\ 0)], \\ &\quad d[(1\ 1\ 1\ 1\ 1\ 1\ 1), (1\ 0\ 1\ 0\ 1\ 0\ 1)], \\ &\quad d[(0\ 0\ 0\ 0), (1\ 1\ 0\ 0)], \\ &\quad d[(1\ 0\ 1\ 0\ 1\ 0), (1\ 1\ 1\ 1\ 0\ 0)], \\ &\quad d[(0\ 1\ 0\ 1\ 0\ 1\ 0\ 1\ 0\ 1), (1\ 0\ 1\ 0\ 1\ 1\ 1\ 1\ 0\ 1)], \\ &\quad d[(1\ 1\ 1\ 0\ 0\ 0\ 1\ 1\ 0), (1\ 1\ 1\ 0\ 0\ 1\ 0\ 0\ 0)]\} \\ &= (1, 3, 2, 3, 6, 3). \end{aligned}$$

Here d is a Hamming 6-distance i.e., n = 6. The Hamming n-weight of the a n-code x is defined as w(x) = d(x, 0). In case

$$\begin{aligned} x &= (1\ 1\ 0\ 1\ 0) \cup (1\ 1\ 1\ 1\ 1\ 1\ 1) \cup (0\ 0\ 0\ 0) \cup \\ &\quad (1\ 0\ 1\ 0\ 1\ 0) \cup (0\ 1\ 0\ 1\ 0\ 1\ 0\ 1\ 0\ 1) \cup \\ &\quad (1\ 1\ 1\ 0\ 0\ 0\ 1\ 1\ 0) \in C_6 \end{aligned}$$

we have w(x) = d(x, 0) = (3, 7, 0, 3, 5, 5). When n = 1 we get the usual Hamming distance and Hamming weight. When n = 2 we get the Hamming bidistance or Hamming 2 distance and Hamming biweight or Hamming 2-weight. When n = 3 we get the Hamming tridistance or Hamming 3-distance and Hamming 3-weight or Hamming triweight and when $n \geq 4$ we have the Hamming n-distance and Hamming n-weight. So it is analogously defined for linear n-codes $d(x, y) = d(x - y, 0) = w(x - y)$. Thus the minimum n-distance of a n-code C is equal to the least of n-weight of all non zero n-code words. If $(d_1, d_2, \ldots, d_n)$, denotes the minimum distance of a linear (n, k) code then we call this n-code C as $((n_1, n_2, \ldots, n_n), (k_1, k_2, \ldots, k_n), (d_1, d_2, \ldots, d_n))$ code.

Now this n-code equivalently be denoted by $(n_1, k_1, d_1) \cup (n_2, k_2, d_2) \cup \ldots \cup (n_n, k_n, d_n)$.

Now we proceed on to define the notion of n-sphere of $(r_1, \ldots, r_n)$ radius of a n-code word $x \in C_1 \cup C_2 \cup \ldots \cup C_n$. Let $x = x_1 \cup x_2 \cup \ldots \cup x_n \in C$; $x_i \in C_i$, $1 \leq i \leq n$.

$$\begin{aligned} S_{(r_1, r_2, \ldots, r_n)}(x) &= S_{(r_1, r_2, \ldots, r_n)}(x_1 \cup x_2 \cup \ldots \cup x_n) \\ &= S_{r_1}(x_1) \cup S_{r_2}(x_2) \cup \ldots \cup S_{r_n}(x_n) \\ &= \{y_1 \in Z_q^{n_1} / d(x_1, y_1) \leq r_1\} \cup \{y_2 \in Z_q^{n_2} / d(x_2, \\ &\quad y_2) \leq r_2\} \cup \ldots \cup \{y_n \in Z_q^{n_n} / d(x_n, y_n) \leq r_n\} \end{aligned}$$

where $x_i = (x_1^i x_2^i \ldots x_{n_i}^i) \in C_i$; $1 \leq i \leq n$ where $C_i$ is a $(n_i, k_i)$ - code.

In decoding we distinguish between the detection and the correction of errors. We say a n-code can correct $(t_1, \ldots, t_n) + (s_1, \ldots, s_n) \geq 0$ errors, if the structure of the n-code makes it possible to correct upto $(t_1, \ldots, t_n)$ tuple of errors to detect $(t_1, \ldots, t_n) + (j_1, \ldots, j_n)$, $0 \leq j_i \leq s_i$, $i = 1, 2, \ldots, n$, tuple errors, which has occurred during the transmission over a channel. Now we have the following result.

*Result:* A linear n-code $C = C_1 \cup C_2 \cup \ldots \cup C_n$ with minimum distance $(d_1, \ldots, d_n)$ can correct up to $(t_1, \ldots, t_n)$ tuple errors and can detect $(t_1 + j_1, t_2 + j_2, \ldots, t_n + j_n)$ tuple errors with $0 \leq j \leq s_i$, $i = 1, 2, \ldots, n$ if and only if each $2t_i + s_i \leq d_i$; $1 \leq i \leq n$. Here $C_i$ is a $(n_i, k_i, d_i)$ code, $i = 1, 2, \ldots, n$.

Now we proceed on to define the notion of Hamming n-bound. The parameter tuples $(q_1, \ldots, q_n)$, $(n_1, \ldots, n_n)$, $(t_1, \ldots, t_n)$, $(M_1, \ldots, M_n)$ of a $(t_1, \ldots, t_n)$ tuple error correcting n-code, $C = C(n_1, k_1, d_1) \cup C(n_2, k_2, d_2) \cup \ldots \cup C(n_n, k_n, d_n)$ over $F_{q_1} \cup F_{q_2} \cup \ldots \cup F_{q_n}$ of length $(n_1, \ldots, n_n)$ with $(m_1, \ldots, m_n)$, n-code words satisfy the n tuple inequality

$$M_1\left(1+(q_1-1)\binom{n_1}{1}+\ldots+(q_1-1)^{t_1}\binom{n_1}{t_1}\right)\cup$$
$$M_2\left(1+(q_2-1)\binom{n_2}{1}+\ldots+(q_2-1)^{t_2}\binom{n_2}{t_2}\right)\cup\ldots\cup$$
$$M_n\left(1+(q_n-1)\binom{n_n}{1}+\ldots+(q_n-1)^{t_n}\binom{n_n}{t_n}\right)$$
$$\leq (q_1^{n_1}, q_1^{n_2}, \ldots, q_n^{n_n})$$

i.e., each

$$M_i\left[(1+(q_i-1)\binom{n_i}{1}+\ldots+(q_i-1)^{t_i}\binom{n_i}{t_i}\right]\leq q^{n_i}\text{ ; } i = 1, 2, \ldots, n.$$

As in case of linear codes we define in case of linear n-codes the notion of a special class of n-codes where the n-tuple of vectors are within or on n-spheres of radius tuple $(t_1, \ldots, t_n)$ about the n-code words of a linear $(((n_1, \ldots, n_n), (k_1, \ldots, k_n)))$ code.

Now we illustrate a method by which error first is detected, this is described for a n-code; n = 5.

***Example 2.2.1:*** Let us consider a 5-code $C = C_1 \cup C_2 \cup \ldots \cup C_5$ generated by the generator 5-matrix $G = G_1 \cup G_2 \cup \ldots \cup G_5$ where

$$G = \begin{bmatrix} 1 & 0 & 0 & 0 & 1 & 1 \\ 0 & 1 & 0 & 1 & 0 & 1 \\ 0 & 0 & 1 & 1 & 1 & 0 \end{bmatrix} \cup$$

$$\begin{bmatrix} 1 & 0 & 1 & 1 \\ 0 & 1 & 0 & 1 \end{bmatrix} \cup$$

$$\begin{bmatrix} 1 & 0 & 0 & 0 & 1 & 0 & 0 \\ 0 & 1 & 0 & 0 & 1 & 0 & 1 \\ 0 & 0 & 1 & 0 & 0 & 1 & 1 \\ 0 & 0 & 0 & 1 & 1 & 1 & 0 \end{bmatrix} \cup$$

$$\begin{bmatrix} 1 & 0 & 1 & 0 & 1 \\ 0 & 1 & 0 & 1 & 0 \end{bmatrix} \cup$$

$$\begin{bmatrix} 1 & 0 & 0 & 1 & 1 & 0 & 1 \\ 0 & 1 & 0 & 0 & 1 & 1 & 0 \\ 0 & 0 & 1 & 0 & 0 & 1 & 1 \end{bmatrix}.$$

is a 5-matrix.

The corresponding parity check 5-matrix;

$$H = H_1 \cup H_2 \cup H_3 \cup H_4 \cup H_5.$$

$$= \begin{bmatrix} 0 & 1 & 1 & 1 & 0 & 0 \\ 1 & 0 & 1 & 0 & 1 & 0 \\ 1 & 1 & 0 & 0 & 0 & 1 \end{bmatrix} \cup$$

$$\begin{bmatrix} 1 & 0 & 1 & 0 \\ 1 & 1 & 0 & 1 \end{bmatrix} \cup$$

$$\begin{bmatrix} 1 & 1 & 0 & 1 & 1 & 0 & 0 \\ 0 & 0 & 1 & 1 & 0 & 1 & 0 \\ 0 & 1 & 1 & 0 & 0 & 0 & 1 \end{bmatrix} \cup$$

$$\begin{bmatrix} 1 & 0 & 1 & 0 & 0 \\ 0 & 1 & 0 & 1 & 0 \\ 1 & 0 & 0 & 0 & 1 \end{bmatrix} \cup$$

$$\begin{bmatrix} 1 & 0 & 0 & 1 & 0 & 0 & 0 \\ 1 & 1 & 0 & 0 & 1 & 0 & 0 \\ 0 & 1 & 1 & 0 & 0 & 1 & 0 \\ 1 & 0 & 1 & 0 & 0 & 0 & 1 \end{bmatrix}.$$

Now we give some of the 5-codes generated by G in the following;

{(0 0 0 0 0 0) ∪ (0 0 0 0) ∪ (0 0 0 0 0 0 0) ∪ (0 0 0 0 0) ∪ (0 0 0 0 0 0 0), (0 1 1 1 0 0) ∪ (1 0 1 0) ∪ (1 1 0 1 1 0 0) ∪ (1 0 1 0 0) ∪ (1 0 0 1 0 0 0), (0 0 0 1 1 1) ∪ (0 1 1 1) ∪ (1 0 0 0 1 1 1) ∪ (0 1 1 1 1) ∪ (1 0 0 1 1 1 1), (0 0 0 0 0 0) ∪ (1 1 0 1) ∪ (0 1 1 0 0 0 1) ∪ (1 0 0 0 1) ∪ (1 0 1 0 0 0 1) and so on}.

Suppose we are sending a code x ∈ C = C(6, 3) ∪ C(4, 2) ∪ C(7, 3) ∪ C(5, 3) ∪ C(7, 4) the 5-code and we have received a message

y = (1 1 1 1 0 0) ∪ (1 0 1 1) ∪ (1 1 0 1 1 0 1) ∪ (1 1 1 0 0) ∪ (0 1 1 1 0 1 1);

$$Hy^T = \left[\begin{bmatrix} 0 & 1 & 1 & 1 & 0 & 0 \\ 1 & 0 & 1 & 0 & 1 & 0 \\ 1 & 1 & 0 & 0 & 0 & 1 \end{bmatrix} \cup \begin{bmatrix} 1 & 0 & 1 & 0 \\ 1 & 1 & 0 & 1 \end{bmatrix} \cup\right.$$

$$\begin{bmatrix} 1 & 1 & 0 & 1 & 1 & 0 & 0 \\ 0 & 0 & 1 & 1 & 0 & 1 & 0 \\ 0 & 1 & 1 & 0 & 0 & 0 & 1 \end{bmatrix} \cup \begin{bmatrix} 1 & 0 & 1 & 0 & 0 \\ 0 & 1 & 0 & 1 & 0 \\ 1 & 0 & 0 & 0 & 1 \end{bmatrix} \cup$$

$$\begin{bmatrix} 1 & 0 & 0 & 1 & 0 & 0 & 0 \\ 1 & 1 & 0 & 0 & 1 & 0 & 0 \\ 0 & 1 & 1 & 0 & 0 & 1 & 0 \\ 1 & 0 & 1 & 0 & 0 & 0 & 1 \end{bmatrix}\Bigg] \times$$

[(1 1 1 1 0 0) ∪ (1 0 1 1) ∪ (1 1 0 1 1 0 1) ∪ (1 1 1 0 0) ∪ (0 1 1 1 0 1 1)]$^T$

$$= \begin{bmatrix} 0 & 1 & 1 & 1 & 0 & 0 \\ 1 & 0 & 1 & 0 & 1 & 0 \\ 1 & 1 & 0 & 0 & 0 & 1 \end{bmatrix} \begin{bmatrix} 1 \\ 1 \\ 1 \\ 1 \\ 0 \\ 0 \end{bmatrix} \cup \begin{bmatrix} 1 & 0 & 1 & 0 \\ 1 & 1 & 0 & 1 \end{bmatrix} \begin{bmatrix} 1 \\ 0 \\ 1 \\ 1 \end{bmatrix} \cup$$

$$\begin{bmatrix} 1 & 1 & 0 & 1 & 1 & 0 & 0 \\ 0 & 0 & 1 & 1 & 0 & 1 & 0 \\ 0 & 1 & 1 & 0 & 0 & 0 & 1 \end{bmatrix} \begin{bmatrix} 1 \\ 1 \\ 0 \\ 1 \\ 1 \\ 0 \\ 1 \end{bmatrix} \cup$$

$$\begin{bmatrix} 1 & 0 & 1 & 0 & 0 \\ 0 & 1 & 0 & 1 & 0 \\ 1 & 0 & 0 & 0 & 1 \end{bmatrix} \begin{bmatrix} 1 \\ 1 \\ 1 \\ 0 \\ 0 \end{bmatrix} \cup$$

$$\begin{bmatrix} 1 & 0 & 0 & 1 & 0 & 0 & 0 \\ 1 & 1 & 0 & 0 & 1 & 0 & 0 \\ 0 & 1 & 1 & 0 & 0 & 1 & 0 \\ 1 & 0 & 1 & 0 & 0 & 0 & 1 \end{bmatrix} \begin{bmatrix} 0 \\ 1 \\ 1 \\ 1 \\ 0 \\ 1 \\ 1 \end{bmatrix}.$$

$= [1\ 0\ 0]^T \cup [0\ 0]^T \cup [0\ 1\ 0]^T \cup [0\ 1\ 1]^T \cup [1\ 1\ 1\ 0]^T$.

This clearly shows the received message y is not a 5-code word from C as $Hy^T \neq (0) \cup (0) \cup (0) \cup (0) \cup (0)$. Thus we have to formulate a method to obtain the correct code word.

We make use of coset n-leader which we would be defining now.

**DEFINITION 2.2.1:** *Let $G = G_1 \cup G_2 \cup \ldots \cup G_n$ be a n-group. H is a n-subgroup of G i.e., $H = H_1 \cup H_2 \cup \ldots \cup H_n$ where each $H_i$ is a subgroup of $G_i$; i = 1, 2, …, n. The coset of the n-subgroup H of the n-group G is defined as,*

$$gH = (g_1 \cup g_2 \cup \ldots \cup g_n)\, H$$

*where*

$$g = g_1 \cup g_2 \cup \ldots \cup g_n \in G$$

*with $g_i \in G_i$, $1 \leq i \leq n$; now*

$$\begin{aligned} gH &= (g_1 \cup g_2 \cup \ldots \cup g_n)\, H \\ &= (g_1 \cup g_2 \cup \ldots \cup g_n)(H_1 \cup \cup \ldots \cup H_n) \\ &= (g_1H_1 \cup g_2H_2 \cup \ldots \cup g_nH_n) \end{aligned}$$

*where $g_iH_i$ is the coset of $H_i$ in $G_i$.*

$$gH = \{g_1h^1 \,/\, h^1 \in H_1\} \cup \{g_2h^2 \,/\, h^2 \in H_2\} \cup \ldots \cup \{g_nh^n \,/\, h^n \in H_n\}.$$

Now we can have a n-field $F_{q_1} \cup F_{q_2} \cup \ldots \cup F_{q_n}$ or $Z_2 \cup Z_2 \cup \ldots \cup Z_2$ also for when we have to define n codes we need the notion of n-vector spaces and n-subspaces so we accept the structure of $F_{q_1} \cup F_{q_2} \cup \ldots \cup F_{q_n}$ even if each of $q_i$ are the same prime number as special n-field which we can call as a pseudo n-field. Thus $Z_5 \cup Z_7 \cup Z_3 \cup Z_{11} \cup Z_{19}$ is a n-field or to be more precise a 5-field now, $Z_3 \cup Z_3 \cup Z_3 \cup Z_3 \cup Z_7$ is a pseudo n field; likewise $Z_2 \cup Z_2 \cup Z_2 \cup Z_2 \cup Z_2 \cup Z_2 \cup Z_2$ is a pseudo 6-field. It is pertinent to mention here a n-field or a pseudo n-field is not a field. Now we can define n-vector space over the n-field or pseudo n-field. We say $V = V_1 \cup V_2 \cup \ldots \cup V_n$ to be n-vector space over the n-field, $F = F_{q_1} \cup F_{q_2} \cup \ldots \cup F_{q_n}$ ($q_i$ – prime number i = 1,2, ..., n) or over a n-pseudo field $F = Z_2 \cup Z_2 \cup \ldots \cup Z_2$ if each $V_i$ is a vector space over $F_{q_i}$ (or $Z_2$); true for i = 1, 2, ..., n. So $V = Z_2^{n_1} \cup Z_2^{n_2} \cup \ldots \cup Z_2^{n_n}$ is a n-vector space over the n-pseudo field F.

Now if each $n_i$ is distinct i.e. $n_i \neq n_j$, if $i \neq j$ for $1 \leq i, j \leq n$, then we say V is a n-vector space even if some of the $n_i$'s are not distinct then we call V to be a pseudo n-vector space over F a n-field or over a pseudo n-field. Clearly if $V = V_1 \cup V_2 \cup \ldots \cup V_n$ is a n-vector space over the n-field $F_{q_1} \cup F_{q_2} \cup \ldots \cup F_{q_n}$ then each $V_i$ is a vector space over $F_{q_i}$, i = 1, 2, ..., n. Further if each of the vector space $V_i$ over $F_{q_i}$ is finite dimensional say of dimension $n_i$ then each $V_i$ is isomorphic to $F_{q_i}^{n_i}$; i = 1, 2, ..., n. Thus

$$V = F_{q_1}^{n_1} \cup F_{q_2}^{n_2} \cup \ldots \cup F_{q_n}^{n_n}$$

is a finite dimensional n-vector space over the n-field, $F = F_{q_1} \cup F_{q_2} \cup \ldots \cup F_{q_n}$ and we say the n-vector space of finite dimension i.e., V is $(n_1, n_2, \ldots, n_n)$-dimension over F.

Suppose $W = W_1 \cup W_2 \cup \ldots \cup W_n$ be a n-subspace of $V = V_1 \cup V_2 \cup \ldots \cup V_n$ i.e. each $W_i$ is a subspace of $V_i$, i = 1, 2, ..., n and if each $W_i$ is of dimension $k_i$ over $F_{q_i}$, $1 \leq i \leq n$ then we say W is a n-subspace of the n-vector space over $F = F_{q_1} \cup F_{q_2} \cup \ldots \cup F_{q_n}$ and W is a $(k_1, k_2, \ldots, k_n)$ dimensional subspace of V and $W = W_1 \cup W_2 \cup \ldots \cup W_n \cong F_{q_1}^{k_1} \cup F_{q_2}^{k_2} \cup \ldots \cup F_{q_n}^{k_n}$;

where $1 \leq k_i \leq n_i$; i = 1, 2, ..., n. When we have the n-code these concepts will be used for error correction as well as error detection.

Now we just recall when we have a n-code $C = C_1 \cup C_2 \cup \ldots C_n = C(n_1, k_1) \cup C(n_2, k_2) \cup \ldots \cup C(n_n, k_n)$, then each $C(n_i, k_i)$ can be realized as a subspace of the vector space $Z_2^{n_i}$ over $Z_2$ or to be more precise each $C(n_i, k_i)$ is the subspace $Z_2^{k_i} \subset Z_2^{n_i}$ over $Z_2$; for i = 1, 2, ..., n. Thus we can say

$$\begin{aligned} C &= C(n_1, k_1) \cup C(n_2, k_2) \cup \ldots \cup C(n_n, k_n) \\ &\cong Z_2^{k_1} \cup Z_2^{k_2} \cup \ldots \cup Z_2^{k_n} \subseteq Z_2^{n_1} \cup Z_2^{n_2} \cup \ldots \cup Z_2^{n_n} \end{aligned}$$

over $Z_2$. This clearly shows each code in C is of length $n_i$ having $k_i$ number of messages and $n_i - k_i$ check symbols; $1 \leq i \leq n$.

Now we proceed onto define the coset n-leader.

We see $V = V_1 \cup V_2 \cup \ldots \cup V_n$ is a n-vector space over the pseudo n-field $F = Z_2 \cup Z_2 \cup \ldots \cup Z_2$ of finite dimension $(n_1, n_2, \ldots, n_n)$. So

$$\begin{aligned} V &= V_1 \cup V_2 \cup \ldots \cup V_n \\ &\cong Z_2^{n_1} \cup Z_2^{n_2} \cup \ldots \cup Z_2^{n_n} . \end{aligned}$$

Let C be a n code over $Z_2$ where C is a $(k_1, k_2, \ldots, k_n)$ dimensional code. Clearly C is a n-subspace of V and $V \cong Z_2^{n_1} \cup Z_2^{n_2} \cup \ldots \cup Z_2^{n_n}$ i.e. $C = Z_2^{k_1} \cup Z_2^{k_2} \cup \ldots \cup Z_2^{k_n}$ is contained in $Z_2^{n_1} \cup Z_2^{n_2} \cup \ldots \cup Z_2^{n_n}$ .

Now we see C is a n-subgroup of the n-group V, so that we can define the n-coset of C in V to be

$$\begin{aligned} x + C &= \{(x = x_1 \cup x_2 \cup \ldots \cup x_n) + C \mid x \in V\} \\ &= \{(x_1 \cup x_2 \cup \ldots \cup x_n) + (Z_2^{k_1} \cup Z_2^{k_2} \cup \ldots \cup Z_2^{k_n}) \end{aligned}$$

(where $x = (x_1 \cup x_2 \cup \ldots \cup x_n) \in Z_2^{n_1} \cup Z_2^{n_2} \cup \ldots \cup Z_2^{n_n}$ ) = $\{(x_1 + Z_2^{k_1}) \cup (x_2 + Z_2^{k_2}) \cup \ldots \cup (x_n + Z_2^{k_n}) \mid x \in Z_2^{n_1} \cup Z_2^{n_2} \cup \ldots \cup Z_2^{n_n}$ ) i.e. $x_i \in Z_2^{n_i}$ ; $1 \leq i \leq n$ }.

We know the coset of a subgroup H of a group G partitions G; likewise the n-coset of a n-subgroup H of a n-group G will partition G. In fact here we can realize V = $Z_2^{n_1} \cup Z_2^{n_2} \cup \ldots \cup Z_2^{n_n}$, the n-group will be partitioned by its n-subgroup C = $Z_2^{k_1} \cup Z_2^{k_2} \cup \ldots \cup Z_2^{k_n}$. The n-coset contains $(2^{k_1}, 2^{k_2}, \ldots, 2^{k_n})$ n-vectors. Thus V is partitioned in the form

$$\begin{aligned} V &= Z_2^{n_1} \cup Z_2^{n_2} \cup \ldots \cup Z_2^{n_n} \\ &= \{Z_2^{k_1} \cup Z_2^{k_2} \cup \ldots \cup Z_2^{k_n}\} \cup \{(x_1^{(1)} + Z_2^{n_1}) \cup (x_1^{(2)} + Z_2^{n_1}) \cup \ldots \cup (x_1^{(t_1)} + Z_2^{n_1})\} \cup \{(x_2^{(1)} + Z_2^{n_2}) \cup \\ &\quad (x_2^{(2)} + Z_2^{n_2}) \cup \ldots \cup (x_2^{(t_2)} + Z_2^{n_2}\} \cup \ldots \cup \{(x_n^{(1)} + Z_2^{n_n}) \cup (x_n^{(2)} + Z_2^{n_n}) \cup \ldots \cup (x_n^{(t_n)} + Z_2^{n_n})\} \end{aligned}$$

where $t_i = q^{n_i - k_i} - 1$; i = 1, 2, ..., n. Thus we see a partition of the n-group V = $Z_2^{n_1} \cup Z_2^{n_2} \cup \ldots \cup Z_2^{n_n}$ by C can be carried out in the following way.

If

$$\begin{aligned} C &= C_1 \cup C_2 \cup \ldots \cup C_n \\ &= C(n_1, k_1) \cup C(n_2, k_2) \cup \ldots \cup C(n_n, k_n) \end{aligned}$$

be the n-code, then we see

$$\begin{aligned} V &= Z_2^{n_1} \cup Z_2^{n_2} \cup \ldots \cup Z_2^{n_n} \\ &= \{x + C \mid x \in V\} \\ &= \{(x_1 \cup x_2 \cup \ldots \cup x_n) + (C_1 \cup C_2 \cup \ldots \cup C_n) \\ &= \{\{x_1 + C_1\} \cup \{x_2 + C_2\} \cup \ldots \cup \{x_n + C_n\} \mid x_i \in Z_2^{n_i};\ i = 1, 2, \ldots, n\}. \end{aligned}$$

Now if y is the received n-code word then certainly $y = y_1 \cup y_2 \cup \ldots \cup y_n$ is such that y is in one of the n-cosets. If the n-code word x has been transmitted then the error n-vector

$$\begin{aligned} e &= (e_1 \cup e_2 \ldots \cup e_n) \\ &= (y - x) \\ &= (y_1 - x_1) \cup \ldots \cup (y_n - x_n) \end{aligned}$$

where each $(y_1 - x_i) \in x_i^{(j)} + C_i$, i = 1, 2, ..., n and $1 \leq j \leq t_i$.

We will illustrate this by an example.

***Example 2.2.2:*** Let $C = C_1 \cup C_2 = C(4, 2) \cup C(6, 3)$ be a bicode associated with the parity check bimatrix,

$$H = H_1 \cup H_2$$

$$= \begin{bmatrix} 1 & 0 & 1 & 0 \\ 1 & 1 & 0 & 1 \end{bmatrix} \cup \begin{bmatrix} 0 & 1 & 1 & 1 & 0 & 0 \\ 1 & 0 & 1 & 0 & 1 & 0 \\ 1 & 1 & 0 & 0 & 0 & 1 \end{bmatrix}.$$

and generator bimatrix

$$G = G_1 \cup G_2$$

$$= \begin{bmatrix} 1 & 0 & 1 & 1 \\ 0 & 1 & 0 & 1 \end{bmatrix} \cup \begin{bmatrix} 1 & 0 & 0 & 0 & 1 & 1 \\ 0 & 1 & 0 & 1 & 0 & 1 \\ 0 & 0 & 1 & 1 & 1 & 0 \end{bmatrix}.$$

The set of bicodes are given by
{(0 0 0 0) ∪ (0 0 0 0 0 0), (1 0 1 1) ∪ (0 0 0 0 0 0), (0 1 0 1) ∪ (0 0 0 0 0 0), (1 1 1 0) ∪ (0 0 0 0 0 0), (0 0 0 0) ∪ (1 0 0 0 1 1), (1 0 1 1) ∪ (1 0 0 0 1 1), (0 1 0 1) ∪ (1 0 0 0 1 1), (1 1 1 0) ∪ (1 0 0 0 1 1), (0 0 0 0) ∪ (0 0 1 1 1 0), (1 1 1 0) ∪ (0 0 1 1 1 0), (0 1 0 1) ∪ ( 0 0 1 1 1 0), (1 0 1 1) ∪ (0 0 1 1 1 0) and so on}.

Now we give the partition of $Z_2^4 \cup Z_2^6$ using the bicode $C = C_1 \cup C_2$.

| | | | | |
|---|---|---|---|---|
| 0 0 0 0 | 1 0 1 1 | 0 1 0 1 | 1 1 1 0 | code words |
| 1 0 0 0 | 0 0 1 0 | 1 1 0 1 | 0 1 1 0 | other cosets |
| 0 1 0 0 | 1 1 1 1 | 0 0 0 1 | 1 0 1 0 | |
| 0 0 1 0 | 1 0 0 1 | 0 1 1 1 | 1 1 0 0 | |

coset
leader

| | | | | |
|---|---|---|---|---|
| 0 0 0 0 0 0 | 1 0 0 0 1 1 | 0 0 1 1 1 0 | 0 1 1 0 1 1 | code words |
| 1 0 0 0 0 0 | 0 0 0 0 1 1 | 1 0 1 1 1 0 | 1 1 1 0 1 1 | other cosets |
| 0 1 0 0 0 0 | 1 1 0 0 1 1 | 0 1 1 1 1 0 | 0 0 1 0 1 1 | |
| 0 0 1 0 0 0 | 1 0 1 0 0 0 | 0 0 0 1 1 0 | 0 1 0 0 1 1 | |
| 0 0 0 1 0 0 | 1 0 0 1 1 1 | 0 0 1 0 1 0 | 0 1 1 1 1 1 | |
| 0 0 0 0 1 0 | 1 0 0 0 0 1 | 0 0 1 1 0 0 | 0 1 1 0 0 1 | |
| 0 0 0 0 0 1 | 1 0 0 0 1 0 | 0 0 1 1 1 1 | 0 1 1 0 1 0 | |
| 0 0 1 0 0 1 | 1 0 1 0 0 1 | 0 0 0 1 1 1 | 0 1 0 0 1 0 | |

coset leader

| | | | | |
|---|---|---|---|---|
| 1 1 0 1 1 0 | 1 1 1 0 0 0 | 0 1 0 1 0 1 | 1 0 1 1 0 1 | code words |
| 0 1 0 1 1 0 | 0 1 1 0 0 0 | 1 1 0 1 0 1 | 0 0 1 1 0 1 | other cosets |
| 1 0 0 1 1 0 | 1 0 1 0 0 0 | 0 0 0 1 0 1 | 1 1 1 1 0 1 | |
| 1 1 1 1 1 0 | 1 1 0 0 0 0 | 0 1 1 1 0 1 | 1 0 0 1 0 1 | |
| 1 1 0 0 1 0 | 1 1 1 1 0 0 | 0 1 0 0 0 1 | 1 0 1 0 0 1 | |
| 1 1 0 1 0 0 | 1 1 1 0 1 0 | 0 1 0 1 1 1 | 1 0 1 1 1 0 | |
| 1 1 0 1 1 1 | 1 1 1 0 0 1 | 0 1 0 1 0 0 | 1 0 1 1 0 1 | |
| 1 1 1 1 1 1 | 1 1 0 0 0 1 | 0 1 1 1 0 0 | 1 0 0 1 0 0 | |

coset leader

is the bicoset of C. Suppose y = (1 1 1 1) ∪ (1 0 0 1 1 1) be the received bicode. We find

$$Hy^T = \left[\begin{bmatrix} 1 & 0 & 1 & 0 \\ 1 & 1 & 0 & 1 \end{bmatrix} \cup \begin{bmatrix} 0 & 1 & 1 & 1 & 0 & 0 \\ 1 & 0 & 1 & 0 & 1 & 0 \\ 1 & 1 & 0 & 0 & 0 & 1 \end{bmatrix}\right]$$

$$\times [(1\ 1\ 1\ 1) \cup (1\ 0\ 0\ 1\ 1\ 1)]^T$$

$$= \left[\begin{bmatrix} 1 & 0 & 1 & 0 \\ 1 & 1 & 0 & 1 \end{bmatrix} \cup \begin{bmatrix} 0 & 1 & 1 & 1 & 0 & 0 \\ 1 & 0 & 1 & 0 & 1 & 0 \\ 1 & 1 & 0 & 0 & 0 & 1 \end{bmatrix}\right] \left[\begin{bmatrix} 1 \\ 1 \\ 1 \\ 1 \end{bmatrix} \cup \begin{bmatrix} 1 \\ 0 \\ 0 \\ 1 \\ 1 \\ 1 \end{bmatrix}\right]$$

$$= \begin{bmatrix} 1 & 0 & 1 & 0 \\ 1 & 1 & 0 & 1 \end{bmatrix} \begin{bmatrix} 1 \\ 1 \\ 1 \\ 1 \end{bmatrix} \cup \begin{bmatrix} 0 & 1 & 1 & 1 & 0 & 0 \\ 1 & 0 & 1 & 0 & 1 & 0 \\ 1 & 1 & 0 & 0 & 0 & 1 \end{bmatrix} \begin{bmatrix} 1 \\ 0 \\ 0 \\ 1 \\ 1 \\ 1 \end{bmatrix}$$

$$= (0\ 1) \cup (1\ 0\ 0)$$
$$\neq (0) \cup (0).$$

So y $\notin$ C. Hence we now inspect from the coset bileader from the table and we see the coset bileader where (1 1 1 1) $\cup$ (1 0 0 1 1 1) occurs is (0 1 0 0) $\cup$ (0 0 0 1 0 0) = e = $e_1 \cup e_2$.

Now y = $y_1 \cup y_2$. Find

x = e + y
= [(0 1 0 0) $\cup$ (0 0 0 1 0 0)] + [(1 1 1 1) $\cup$ (1 0 0 1 1 1)]
= [(0 1 0 0) + (1 1 1 1)] $\cup$ [(0 0 0 1 0 0) + (1 0 0 1 1 1)]
= (1 0 1 1) $\cup$ (1 0 0 0 1 1) $\in$ C.

Thus the sent message is (1 0 1 1) $\cup$ (1 0 0 0 1 1).

Now we have the following set to be the bicoset bileader

{(0 0 0 0) $\cup$ (0 0 0 0 0 0), (1 0 0 0) $\cup$ (1 0 0 0 0 0), (1 0 0 0) $\cup$ (0 1 0 0 0 0), (1 0 0 0) $\cup$ (0 0 1 0 0 0), (1 0 0 0) $\cup$ (0 0 0 1 0 0), (1 0 0 0) $\cup$ (0 0 0 0 0 0 0), (1 0 0 0) $\cup$ (0 0 0 0 1 0), (1 0 0 0) $\cup$ (0 0 0 0 0 1), (1 0 0 0) $\cup$ (0 0 1 0 0 1), (0 0 0 0) $\cup$ (1 0 0 0 0 0), (0 0 0 0) $\cup$ (0 1 0 0 0 0), (0 0 0 0) $\cup$ (0 0 1 0 0 0), (0 0 0 0) $\cup$ (0 0 0 1 0 0), (0 0 0 0) $\cup$ (0 0 0 0 1 0), (0 0 0 0) $\cup$ (0 0 0 0 0 1), (0 0 0 0) $\cup$ (0 0 1 0 0 1), (0 1 0 0) $\cup$ (0 0 0 0 0 0), (0 1 0 0) $\cup$ (1 0 0 0 0), (0 1 0 0) $\cup$ (0 1 0 0 0 0), (0 1 0 0) $\cup$ (0 0 1 0 0 0), (0 1 0 0) $\cup$ (0 0 0 1 0 0), (0 1 0 0) $\cup$ (0 0 0 0 1 0), (0 1 0 0) $\cup$ (0 0 0 0 0 1), (0 1 0 0) $\cup$ (0 0 1 0 0 1, (0 0 1 0) $\cup$ (0 0 0 0 0 0), (0 0 1 0 ) $\cup$ (1 0 0 0 0 0), (0 0 1 0) $\cup$ (0 1 0 0 0 0), (0 0 1 0) $\cup$ (0 0 1 0 0 0), (0 0 1 0) $\cup$ (0 0 0 1 0 0), (0 0 1 0) $\cup$ (0 0 0 0 1 0), (0 0 1 0) $\cup$ (0 0 0 0 0 1), (0 0 1 0) $\cup$ (0 0 1 0 0 1)}

for the code C = C(4, 2) ∪ C(6, 3).

Thus finding the sent message which is the bicode x and if y is the received bicode and y ∉ C, y is determined by finding

$$
\begin{aligned}
Hy^T &= (H_1 \cup H_2)y^T \\
&= (H_1 \cup H_2)\,(y_1 \cup y_2)^T \\
&= (H_1 \cup H_2)\,(y_1^T \cup y_2^T) \\
&= H_1y_1^T \cup H_2y_2^T.
\end{aligned}
$$

If $H_1y^T \cup H_2y_2^T = (0) \cup (0)$ then $y \in C$ if $H_1y_1^T \cup H_2y_2^T \neq (0) \cup (0)$, we find the coset leader to which the code $y_1$ belongs to and if $e_1$ is the coset leader to which $y_1$ belongs to, then we find the coset leader to which $y_2$ belongs to. Let $e_2$ be the coset leader to which $y_2$ belongs to. Then $e = e_1 \cup e_2$ is the bierror or the bierror bivector. Now y + e = x gives the sent bicode word. $e = e_1 \cup e_2$ is the bicoset bileader which will be found in the class of bicoset leaders.

Now we proceed on to describe how in the tricode the error is both detected and corrected using the tricoset leaders.

***Example 2.2.3:*** Let $C = C_1 \cup C_2 \cup C_3$ be a tricode with the associated parity check trimatrix

$$H = H_1 \cup H_2 \cup H_3 =$$

$$\begin{bmatrix} 1 & 1 & 1 & 0 & 1 & 0 & 0 \\ 0 & 1 & 1 & 1 & 0 & 1 & 0 \\ 1 & 1 & 0 & 1 & 0 & 0 & 1 \end{bmatrix} \cup$$

$$\begin{bmatrix} 1 & 0 & 0 & 1 & 0 & 0 \\ 0 & 1 & 0 & 0 & 1 & 0 \\ 0 & 0 & 1 & 0 & 0 & 1 \end{bmatrix} \cup \begin{bmatrix} 1 & 0 & 1 & 0 \\ 1 & 1 & 0 & 1 \end{bmatrix}.$$

The related generator trimatrix of the tricode is given by

$$G = G_1 \cup G_2 \cup G_3 =$$

$$\begin{bmatrix} 1 & 0 & 0 & 0 & 1 & 0 & 1 \\ 0 & 1 & 0 & 0 & 1 & 1 & 1 \\ 0 & 0 & 1 & 0 & 1 & 1 & 0 \\ 0 & 0 & 0 & 1 & 0 & 1 & 1 \end{bmatrix} \cup$$

$$\begin{bmatrix} 1 & 0 & 0 & 1 & 0 & 0 \\ 0 & 1 & 0 & 0 & 1 & 0 \\ 0 & 0 & 1 & 0 & 0 & 1 \end{bmatrix} \cup \begin{bmatrix} 1 & 0 & 1 & 1 \\ 0 & 1 & 0 & 1 \end{bmatrix}.$$

Clearly C = C (7, 4) ∪ C (6, 3) ∪ C (4, 2) is the tricode. The elements of the tricode are as follows:

{(0 0 0 0 0 0 0) ∪ (0 0 0 0 0 0) ∪ (0 0 0 0), (1 0 0 0 1 0 1) ∪ (1 0 0 1 0 0) ∪ (1 0 1 1), (0 1 0 1 1 1) ∪ (0 1 0 0 1 0) ∪ (0 1 0 1), (0 0 1 0 1 1 0) ∪ (0 0 1 1) ∪ (1 1 1 0), ( 0 0 1 0 1 1) ∪ (1 1 1 1 1 1) ∪ (0 0 0 0), (1 1 1 1 1 1 1) ∪ (1 1 1 1 1 1) ∪ (1 1 1 0), (1 1 0 0 0 1 0) ∪ (1 1 0 1 1 0) ∪ (1 1 1 0) and so on}.

Suppose y = (0 1 1 1 1 1 1) ∪ (0 1 1 1 1 1) ∪ (1 1 1 1) be the received tricode word. Find whether the received tricode word is the correct message. If y is to be a correctly received message than we must have $Hy^T$ = (0) ∪ (0) ∪ (0). Now consider

$$Hy^T = \left[ \begin{bmatrix} 1 & 1 & 1 & 0 & 1 & 0 & 0 \\ 0 & 1 & 1 & 1 & 0 & 1 & 0 \\ 1 & 1 & 0 & 1 & 0 & 0 & 1 \end{bmatrix} \cup \right.$$

$$\left. \begin{bmatrix} 1 & 0 & 0 & 1 & 0 & 0 \\ 0 & 1 & 0 & 0 & 1 & 0 \\ 0 & 0 & 1 & 0 & 0 & 1 \end{bmatrix} \cup \begin{bmatrix} 1 & 0 & 1 & 0 \\ 1 & 1 & 0 & 1 \end{bmatrix} \right]$$

$$\left[ \begin{bmatrix} 0 & 1 & 1 & 1 & 1 & 1 & 1 \end{bmatrix} \cup \begin{bmatrix} 0 & 1 & 1 & 1 & 1 & 1 \end{bmatrix} \cup \begin{bmatrix} 1 & 1 & 1 & 1 \end{bmatrix} \right]^T$$

$$= \begin{bmatrix} 1 & 1 & 1 & 0 & 1 & 0 & 0 \\ 0 & 1 & 1 & 1 & 0 & 1 & 0 \\ 1 & 1 & 0 & 1 & 0 & 0 & 1 \end{bmatrix} \begin{bmatrix} 0 \\ 1 \\ 1 \\ 1 \\ 1 \\ 1 \\ 1 \end{bmatrix} \cup \begin{bmatrix} 1 & 0 & 0 & 1 & 0 & 0 \\ 0 & 1 & 0 & 0 & 1 & 0 \\ 0 & 0 & 1 & 0 & 0 & 1 \end{bmatrix} \begin{bmatrix} 0 \\ 1 \\ 1 \\ 1 \\ 1 \\ 1 \end{bmatrix} \cup$$

$$\begin{bmatrix} 1 & 0 & 1 & 0 \\ 1 & 1 & 0 & 1 \end{bmatrix} \begin{bmatrix} 1 \\ 1 \\ 1 \\ 1 \end{bmatrix}$$

$$= \left[ \begin{bmatrix} 1 & 0 & 1 \end{bmatrix}^T \cup \begin{bmatrix} 1 & 0 & 0 \end{bmatrix}^T \cup \begin{bmatrix} 0 & 1 \end{bmatrix}^T \right]$$

$$\neq (0) \cup (0) \cup (0) .$$

Thus y has error, i.e. the received message has some error. Now we have to detect the error and obtain the correct message. Consider the coset trileader e = (1 0 0 0 0 0 0) ∪ (1 0 0 0 0 0) ∪ (0 0 0 1); for the vector (0 1 1 1 1 1 1) of $Z_2^7$ falls in the coset with coset leader (1 0 0 0 0 0 0); the vector (0 1 1 1 1 1) of $Z_2^6$ falls in the coset with coset leader (1 0 0 0 0 0) and the vector (1 1 1 1) of $Z_2^4$ falls in the coset with coset leader (0 0 0 1). Thus the trivector (0 1 1 1 1 1 1) ∪ (0 1 1 1 1 1) ∪ (1 1 1 1) falls in the tricoset with the tricoset leader (10 0 0 0 0 0) ∪ (1 0 0 0 0 0) ∪ (0 0 0 1) = e = $e_1 \cup e_2 \cup e_3$. Thus y + e = x is the correct tricode word which must be the sent message. Thus

y + e = [(0 1 1 1 1 1 1) ∪ ( 0 1 1 1 1 1) ∪ (1 1 1 1) ] + [ (1 0 0 0 0 0 0) ∪ (1 0 0 0 0 0) ∪ (0 0 0 1) ]
= [((0 1 1 1 1 1 1) + (1 0 0 0 0 0 0)] ∪ [(0 1 1 1 1 1) + (1 0 0 0 0 0) ] ∪ ( ( 1 1 1 1) + (0 0 0 1) ) ∈ C

$= \quad C_1 \cup C_2 \cup C_3.$

Now the same method can be used for any n-code; n ≥ 4. However we illustrate for a 5-code by an explicit example before we proceed to dewell about other methods.

***Example 2.2.4:*** Consider the 5-code $C = C_1 \cup C_2 \cup C_3 \cup C_4 \cup C_5$ associated with the parity check 5-matrix

$$H = H_1 \cup H_2 \cup H_3 \cup H_4 \cup H_5 =$$

$$\begin{bmatrix} 1 & 0 & 0 & 1 & 0 & 0 \\ 0 & 1 & 0 & 0 & 1 & 0 \\ 0 & 0 & 1 & 0 & 0 & 1 \end{bmatrix} \cup \begin{bmatrix} 1 & 0 & 1 & 0 \\ 1 & 1 & 0 & 1 \end{bmatrix} \cup$$

$$\begin{bmatrix} 1 & 1 & 1 & 1 & 1 & 0 & 0 \\ 0 & 1 & 0 & 1 & 0 & 1 & 0 \\ 1 & 0 & 1 & 0 & 0 & 0 & 1 \end{bmatrix} \cup$$

$$\begin{bmatrix} 1 & 1 & 1 & 1 & 0 & 0 & 0 \\ 0 & 0 & 1 & 0 & 1 & 0 & 0 \\ 1 & 0 & 1 & 0 & 0 & 1 & 0 \\ 0 & 1 & 1 & 0 & 0 & 0 & 1 \end{bmatrix} \cup \begin{bmatrix} 0 & 1 & 0 & 0 & 1 & 0 \\ 0 & 1 & 1 & 1 & 0 & 1 \end{bmatrix}.$$

The corresponding generator 5 matrix

$$G = G_1 \cup G_2 \cup G_3 \cup G_4 \cup G_5$$

$$= \begin{bmatrix} 1 & 0 & 0 & 1 & 0 & 0 \\ 0 & 1 & 0 & 0 & 1 & 0 \\ 0 & 0 & 1 & 0 & 0 & 1 \end{bmatrix} \cup \begin{bmatrix} 1 & 0 & 1 & 1 \\ 0 & 1 & 0 & 1 \end{bmatrix}$$

$$\begin{bmatrix} 1 & 0 & 0 & 0 & 1 & 0 & 1 \\ 0 & 1 & 0 & 0 & 1 & 1 & 0 \\ 0 & 0 & 1 & 0 & 1 & 0 & 1 \\ 0 & 0 & 0 & 1 & 1 & 1 & 0 \end{bmatrix} \cup$$

$$\cup \begin{bmatrix} 1 & 0 & 0 & 1 & 0 & 1 & 0 \\ 0 & 1 & 0 & 1 & 0 & 0 & 1 \\ 0 & 0 & 1 & 1 & 1 & 1 & 1 \end{bmatrix} \cup \begin{bmatrix} 1 & 0 & 0 & 0 & 0 & 0 \\ 0 & 1 & 0 & 0 & 1 & 1 \\ 0 & 0 & 1 & 0 & 0 & 1 \\ 0 & 0 & 0 & 1 & 0 & 1 \end{bmatrix}.$$

The related 5-code words are as follows:

{(0 0 0 0 0 0) ∪ (0 0 0 0) ∪ (0 0 0 0 0 0 0) ∪ (0 0 0 0 0 0 0) ∪ (0 0 0 0 0 0), (1 0 0 1 0 0) ∪ (1 0 1 1) ∪ (1 0 0 0 1 0 1) ∪ (1 0 0 1 0 1 0) ∪ (1 0 0 0 0 0), (0 1 0 0 1 0) ∪ (1 1 1 0) ∪ (0 1 0 0 1 1 0) ∪ (0 1 0 1 0 0 1) ∪ (0 1 0 0 1 1), (0 0 1 0 0 1) ∪ (1 1 1 0) ∪ (0 0 1 0 1 0 1) ∪ (0 0 1 1 1 1 1) ∪ (0 0 1 0 0 1) ∪ (0 0 1 0 0 1) ∪ (1 1 1 0) ∪ (0 0 1 0 1 0 1) ∪ (0 0 1 1 1 1 1) ∪ (0 0 0 1 0 1) and so on}.

Suppose $y = y_1 \cup y_2 \cup y_3 \cup y_4 \cup y_5$ = ( (1 1 1 1 0 0) ∪ (1 1 1 1) ∪ (1 1 1 1 1 0 0) ∪ (0 1 1 0 0 0 1) ∪ (0 1 0 1 0 1) ) be the received word, we have to find whether y is a correct received message. The simple test for this is to find $Hy^T$, if $Hy^T$ = (0) ∪ (0) ∪ (0) ∪ (0) ∪ (0) then $y \in C$ i.e. y is the correct received message. Now consider

$$Hy^T = \left\{ \begin{bmatrix} 1 & 0 & 0 & 1 & 0 & 0 \\ 0 & 1 & 0 & 0 & 1 & 0 \\ 0 & 0 & 1 & 0 & 0 & 1 \end{bmatrix} \cup \begin{bmatrix} 1 & 0 & 1 & 0 \\ 1 & 1 & 0 & 1 \end{bmatrix} \cup \right.$$

$$\begin{bmatrix} 1 & 1 & 1 & 1 & 1 & 0 & 0 \\ 0 & 1 & 0 & 1 & 0 & 1 & 0 \\ 1 & 0 & 1 & 0 & 0 & 0 & 1 \end{bmatrix}$$

$$\cup \begin{bmatrix} 1 & 1 & 1 & 1 & 0 & 0 & 0 \\ 0 & 0 & 1 & 0 & 1 & 0 & 0 \\ 1 & 0 & 1 & 0 & 0 & 1 & 0 \\ 0 & 1 & 1 & 0 & 0 & 0 & 1 \end{bmatrix} \cup \begin{bmatrix} 0 & 1 & 0 & 0 & 1 & 0 \\ 0 & 1 & 1 & 1 & 0 & 1 \end{bmatrix} \Bigg\}$$

× [(1 1 1 1 0 0)$^T$ ∪ (1 1 1 1)$^T$ ∪ (1 1 1 1 1 0 0)$^T$ ∪ (0 1 1 0 0 0 1)$^T$ ∪ (0 1 0 1 0 1)$^T$ ]

$$= \begin{bmatrix} 1 & 0 & 0 & 1 & 0 & 0 \\ 0 & 1 & 0 & 0 & 1 & 0 \\ 0 & 0 & 1 & 0 & 0 & 1 \end{bmatrix} \begin{bmatrix} 1 \\ 1 \\ 1 \\ 1 \\ 0 \\ 0 \end{bmatrix} \cup \begin{bmatrix} 1 & 0 & 1 & 0 \\ 1 & 1 & 0 & 1 \end{bmatrix} \begin{bmatrix} 1 \\ 1 \\ 1 \\ 1 \end{bmatrix}$$

$$\cup \begin{bmatrix} 1 & 1 & 1 & 1 & 1 & 0 & 0 \\ 0 & 1 & 0 & 1 & 0 & 1 & 0 \\ 1 & 0 & 1 & 0 & 0 & 0 & 1 \end{bmatrix} \begin{bmatrix} 1 \\ 1 \\ 1 \\ 1 \\ 1 \\ 0 \\ 0 \end{bmatrix} \cup \begin{bmatrix} 1 & 1 & 1 & 1 & 0 & 0 & 0 \\ 0 & 0 & 1 & 0 & 1 & 0 & 0 \\ 1 & 0 & 1 & 0 & 0 & 1 & 0 \\ 0 & 1 & 1 & 0 & 0 & 0 & 1 \end{bmatrix} \begin{bmatrix} 0 \\ 1 \\ 1 \\ 0 \\ 0 \\ 0 \\ 1 \end{bmatrix}$$

$$\begin{bmatrix} 0 & 1 & 0 & 0 & 1 & 0 \\ 0 & 1 & 1 & 1 & 0 & 1 \end{bmatrix} \begin{bmatrix} 0 \\ 1 \\ 0 \\ 1 \\ 0 \\ 1 \end{bmatrix}$$

$$= \begin{bmatrix} 0 \\ 1 \\ 1 \end{bmatrix}^T \cup \begin{bmatrix} 0 \\ 1 \end{bmatrix}^T \cup \begin{bmatrix} 1 \\ 1 \\ 0 \end{bmatrix}^T \cup \begin{bmatrix} 0 \\ 1 \\ 1 \\ 1 \end{bmatrix}^T \cup \begin{bmatrix} 1 \\ 1 \end{bmatrix}^T$$

$$\neq (0) \cup (0) \cup (0) \cup (0) \cup (0)\ .$$

Thus the received vector $y \notin C$ i.e. y is not the correct message; now we have detected that this message has error, now we proceed on to correct the error. This is done by finding the coset-5-leader in which the 5-vector $y = y_1 \cup y_2 \cup y_3 \cup y_4 \cup y_5$ occurs. Suppose, $e = e_1 \cup e_2 \cup e_3 \cup e_4 \cup e_5$ then $y + e = x$ would be the corrected message and $x \in C = C_1 \cup C_2 \cup C_3 \cup C_4 \cup C_5$.

Now the coset 5-leader in which y occurs is given by

$$\begin{aligned}
e &= (0\ 1\ 1\ 0\ 0\ 0) \cup (0\ 0\ 0\ 1) \cup (0\ 0\ 0\ 0\ 1\ 0\ 0) \cup (0\ 0\ 1\ 1\ 0\ 0\ 0) \cup (0\ 1\ 0\ 0\ 0\ 0). \\
y + e &= [(1\ 1\ 1\ 1\ 0\ 0) \cup (1\ 1\ 1\ 1) \cup (1\ 1\ 1\ 1\ 1\ 0\ 0) \cup (0\ 1\ 1\ 0\ 0\ 0\ 1) \cup (0\ 1\ 0\ 1\ 0\ 1)] + [(0\ 1\ 1\ 0\ 0\ 0) \cup (0\ 0\ 0\ 1) \cup (0\ 0\ 0\ 0\ 1\ 0\ 0) \cup (0\ 0\ 1\ 1\ 0\ 0\ 0) \cup (0\ 1\ 0\ 0\ 0\ 0)] \\
&= [(1\ 1\ 1\ 1\ 0\ 0) + (0\ 1\ 1\ 0\ 0\ 0)] \cup [(1\ 1\ 1\ 1) + (0\ 0\ 0\ 1)] \cup [(1\ 1\ 1\ 1\ 1\ 0\ 0) + (0\ 0\ 0\ 0\ 1\ 0\ 0)\ ] \cup [(0\ 1\ 1\ 0\ 0\ 0\ 1) + (0\ 0\ 1\ 1\ 0\ 0\ 0)] \cup [(0\ 1\ 0\ 1\ 0\ 1) + (0\ 1\ 0\ 0\ 0\ 0)] \\
&= (1\ 0\ 0\ 1\ 0\ 0) \cup (1\ 1\ 1\ 0) \cup (1\ 1\ 1\ 1\ 0\ 0\ 0) \cup (0\ 1\ 0\ 1\ 0\ 0\ 1) \cup (0\ 0\ 0\ 1\ 0\ 1) \\
&= x.
\end{aligned}$$

Now find

$$Hx^T = \left\{ \begin{bmatrix} 1 & 0 & 0 & 1 & 0 & 0 \\ 0 & 1 & 0 & 0 & 1 & 0 \\ 0 & 0 & 1 & 0 & 0 & 1 \end{bmatrix} \cup \begin{bmatrix} 1 & 0 & 1 & 0 \\ 1 & 1 & 0 & 1 \end{bmatrix} \cup \right.$$

$$\cup \begin{bmatrix} 1 & 1 & 1 & 1 & 1 & 0 & 0 \\ 0 & 1 & 0 & 1 & 0 & 1 & 0 \\ 1 & 0 & 1 & 0 & 0 & 0 & 1 \end{bmatrix} \cup \begin{bmatrix} 1 & 1 & 1 & 1 & 0 & 0 & 0 \\ 0 & 0 & 1 & 0 & 1 & 0 & 0 \\ 1 & 0 & 1 & 0 & 0 & 1 & 0 \\ 0 & 1 & 1 & 0 & 0 & 0 & 1 \end{bmatrix} \cup$$

$$\left.\begin{bmatrix} 0 & 1 & 0 & 0 & 1 & 0 \\ 0 & 1 & 1 & 1 & 0 & 1 \end{bmatrix}\right\} \times$$

$$\left( [1 \quad 0 \quad 0 \quad 1 \quad 0 \quad 0]^T \cup [1 \quad 1 \quad 1 \quad 0]^T \cup \right.$$

$$[1 \quad 1 \quad 1 \quad 1 \quad 0 \quad 0 \quad 0]^T \cup [0 \quad 1 \quad 0 \quad 1 \quad 0 \quad 0 \quad 1]^T \cup$$

$$\left. [0 \quad 0 \quad 0 \quad 1 \quad 0 \quad 1]^T \right)$$

$$= \begin{bmatrix} 1 & 0 & 0 & 1 & 0 & 0 \\ 0 & 1 & 0 & 0 & 1 & 0 \\ 0 & 0 & 1 & 0 & 0 & 1 \end{bmatrix} \begin{bmatrix} 1 \\ 0 \\ 0 \\ 1 \\ 0 \\ 0 \end{bmatrix}$$

$$\cup \begin{bmatrix} 1 & 0 & 1 & 0 \\ 1 & 1 & 0 & 1 \end{bmatrix} \begin{bmatrix} 1 \\ 1 \\ 1 \\ 0 \end{bmatrix} \cup \begin{bmatrix} 1 & 1 & 1 & 1 & 1 & 0 & 0 \\ 0 & 1 & 0 & 1 & 0 & 1 & 0 \\ 1 & 0 & 1 & 0 & 0 & 0 & 1 \end{bmatrix} \begin{bmatrix} 1 \\ 1 \\ 1 \\ 1 \\ 0 \\ 0 \\ 0 \end{bmatrix}$$

$$\cup \begin{bmatrix} 1 & 1 & 1 & 1 & 0 & 0 & 0 \\ 0 & 0 & 1 & 0 & 1 & 0 & 0 \\ 1 & 0 & 1 & 0 & 0 & 1 & 0 \\ 0 & 1 & 1 & 0 & 0 & 0 & 1 \end{bmatrix} \begin{bmatrix} 0 \\ 1 \\ 0 \\ 1 \\ 0 \\ 0 \\ 1 \end{bmatrix} \cup \begin{bmatrix} 0 & 1 & 0 & 0 & 1 & 0 \\ 0 & 1 & 1 & 1 & 0 & 1 \end{bmatrix} \begin{bmatrix} 0 \\ 0 \\ 0 \\ 1 \\ 0 \\ 1 \end{bmatrix}$$

$$= \quad (0\ 0\ 0)^T \cup (0\ 0\ 0\ 0)^T \cup (0\ 0\ 0)^T \cup (0\ 0\ 0\ 0)^T \cup (0\ 0)^T.$$

Thus x is a 5-code word of C which is the corrected message i.e. the sent message.

Now we proceed on to describe this more technically by defining the term n-syndrome n ≥ 4, bisyndrome and trisyndrome.

**DEFINITION 2.2.2:** *Let $H = H_1 \cup H_2$ be a parity check bimatrix of a linear $(n_1, k_1) \cup (n_2, k_2)$ bicode. Then the bivector*

$$\begin{aligned} S(y) &= S(y_1) \cup S(y_2) \\ &= Hy^T \\ &= (H_1 \cup H_2)(y_1 \cup y_2)^T \\ &= (H_1 \cup H_2)(y_1^T \cup y_2^T) \\ &= H_1 y_1^T \cup H_2 y_2^T \end{aligned}$$

*of length $(n_1 - k_1) \cup (n_2 - k_2)$, is called the bisyndrome of the bivector $y = y_1 \cup y_2$.*

$$S(y) = S(y_1) \cup S(y_2) = (0) \cup (0)$$

*if and only if $y \in C_1 \cup C_2 = C$.*

*$S(y^{(1)}) = S(y^{(2)}) = S(y_{11}^{(1)} \cup y_{21}^{(1)}) = S(y_{12}^{(2)} \cup y_{22}^{(2)})$ if and only if $y^{(1)} + C = y^{(2)} + C$ i.e. if and only if $(y_{11}^{(1)} \cup y_{21}^{(1)}) + (C_1 \cup C_2) = (y_{12}^{(2)} \cup y_{22}^{(2)}) + (C_1 \cup C_2)$ i.e. $(y_{11}^{(1)} + C_1) \cup (y_{21}^{(1)} + C_2) = (y_{12}^{(2)} + C_1) \cup (y_{22}^{(2)} + C_2)$.*

*If $y = y_1 \cup y_2$ is the received message and $e = e_1 \cup e_2$ is the error with $y = x + e$ then*

$$S(y) = S(y_1) \cup S(y_2)$$
$$= S(x+e)$$
$$= S((x_1 \cup x_2) + (e_1 + e_2))$$
$$= S((x_1 + e_1) \cup (x_2 + e_2))$$
$$= S(x_1 + e_1) \cup S(x_2 + e_2)$$
$$= \{S(x_1) + S(e_1)\} \cup \{S(x_2) + S(e_2)\}$$
$$= S(e_1) \cup S(e_2)$$

*as* $x = x_1 \cup x_2 \in C = C_1 \cup C_2$ *so* $S(x) = S(x_1 \cup x_2) = S(x_1) \cup S(x_2) = (0) \cup (0)$ *i.e. y and e are in the same bicoset.*

*Let us consider how a trisyndrome of a tricode functions and the tricoset leader. Let* $H = H_1 \cup H_2 \cup H_3$ *be a parity check trimatrix associated with the tricode* $C = C_1 \cup C_2 \cup C_3$*. Any tricode word in C would be of the form* $y = y_1 \cup y_2 \cup y_3$ *where* $y_1$ *is the code word in* $C_1$*,* $y_2$ *a code word in* $C_2$ *and* $y_3$ *the code word in* $C_3$*.* $Hy^T = S(y)$ *i.e.*

$$S(y_1 \cup y_2 \cup y_3) = (H_1 \cup H_2 \cup H_3)(y_1^T \cup y_2^T \cup y_3^T)$$
$$= H_1y_1^T \cup H_2y_2^T \cup H_3y_3^T$$
$$= S(y_1) \cup S(y_2) \cup S(y_3).$$

*This* $S(y) = (0) \cup (0) \cup (0)$ *if and only if* $y = y_1 \cup y_2 \cup y_3 \in C = C_1 \cup C_2 \cup C_3$*. If* $S(y) \neq (0)$ *then* $y \notin C$*.* $S(y) = S(y_1) \cup S(y_2) \cup S(y_3)$ *is of length* $(n_1 - k_1) \cup (n_2 - k_2) \cup (n_3 - k_3)$ *is called the trisyndrome of the trivector* $y = y_1 \cup y_2 \cup y_3$*. Thus* $S(y) = S(y_1) \cup S(y_2) \cup S(y_3) = (0) \cup (0) \cup (0)$ *if and only if* $y \in C_1 \cup C_2 \cup C_3 = C$*.*

$S(y^{(1)}) = S(y^{(2)}) = S(y_{11}^{(1)} \cup y_{21}^{(1)} \cup y_{31}^{(1)}) = S(y_{12}^{(2)} \cup y_{22}^{(2)} \cup y_{32}^{(2)})$ *if and only if*

$$[y_{11}^{(1)} \cup y_{21}^{(1)} \cup y_{31}^{(1)}] + (C_1 \cup C_2 \cup C_3)$$
$$= (y_{12}^{(2)} \cup y_{22}^{(2)} \cup y_{32}^{(2)}) + (C_1 \cup C_2 \cup C_3)$$
$$= (y_{11}^{(1)} + C_1) \cup (y_{21}^{(1)} + C_2) \cup (y_{31}^{(1)} + C_3)$$
$$= (y_{12}^{(2)} + C_1) \cup (y_{22}^{(2)} + C_2) \cup (y_{32}^{(2)} + C_3).$$

*If* $y = y_1 \cup y_2 \cup y_3$ *is the received message and* $e = e_1 \cup e_2 \cup e_3$ *is the error with* $y = x + e$ *then*

$$S(y) = S(y_1) \cup S(y_2) \cup S(y_3)$$
$$= S(x + e)$$
$$= S((x_1 \cup x_2 \cup x_3) + (e_1 \cup e_2 \cup e_3))$$

$$= S(x_1 + e_1) \cup S(x_2 + e_2) \cup S(x_3 + e_3)$$
$$= (S(x_1) + S(e_1)) \cup (S(x_2)+S(e_2)) \cup (S(x_3) + S(e_3))$$
$$= S(e_1) \cup S(e_2) \cup S(e_3)$$

*as $S(x_i) = (0)$ for $i = 1, 2, 3$. i.e., y and e are in the same tricoset and e is taken as the leader of the tricoset to which the trivector $y = y_1 \cup y_2 \cup y_3$ belongs.*

*Now in a similar manner we have the notion of n-syndrome and n coset leader and n-error vector associated with the n-code.*

*Let $C = C_1 \cup C_2 \cup \ldots \cup C_n$ be a n code where $C = C(n_1, k_1) \cup C(n_2, k_2) \cup \ldots \cup C(n_n, k_n)$ a n-code ($n \geq 4$). Any n vector y would be of the form $y = y_1 \cup y_2 \cup \ldots \cup y_n$. Any n-code $x = x_1 \cup x_2 \cup \ldots x_n$ would be such that each $x_i \in C(n_i, k_i)$ is a code of length $n_i$ with $k_i$ number of message symbols; $i = 1, 2, \ldots, n$.*

*Let $H = H_1 \cup H_2 \cup \ldots \cup H_n$ be the parity check n-matrix associated with $C = C_1 \cup C_2 \cup \ldots \cup C_n$, the n-code. Let $x = x_1 \cup x_2 \cup \ldots \cup x_n \in C = C_1 \cup C_2 \cup \ldots \cup C_n$. $Hx^T = (H_1 \cup H_2 \cup \ldots \cup H_n)(x^T) = (H_1 \cup H_2 \cup \ldots \cup H_n)(x_1^T \cup x_2^T \cup \ldots \cup x_n^T) = H_1x_1^T \cup H_2x_2^T \cup \ldots \cup H_nx_n^T = (0) \cup (0) \cup \ldots \cup (0)$.*

*If y is any n-vector if $Hy^T \neq (0) \cup (0) \cup \ldots \cup (0)$ then $y \notin C = C_1 \cup C_2 \cup \ldots \cup C_n$.*

*Now*

$$S(y) = Hy^T$$
$$= (H_1 \cup H_2 \cup \ldots \cup H_n)(y_1 \cup y_2 \cup \ldots \cup y_n)^T$$
$$= H_1y_1^T \cup H_2y_2^T \cup \ldots \cup H_ny_n^T$$

*is a $(n_1 - k_1) \cup (n_2 - k_2) \cup \ldots \cup (n_n - k_n)$ length n-vector defined to be the n-syndrome of y. $S(y) = (0) \cup (0) \cup \ldots \cup (0)$ if and only if $y \in C = C_1 \cup C_2 \cup \ldots \cup C_n$.*

*Now if y is a received message to find the correct sent message. In the first stage we detect whether the received message is correct or not. This is done by calculating S(y) i.e. by finding the n-syndrome of y. i.e.*

$$S(y) = Hy^T$$
$$= (H_1 \cup H_2 \cup \ldots \cup H_n)y^T$$
$$= (H_1 \cup H_2 \cup \ldots \cup H_n)(y_1^T \cup y_2^T \cup \ldots \cup y_n^T)$$
$$= H_1y_1^T \cup H_2y_2^T \cup \ldots \cup H_ny_n^T.$$

$$\textit{If } S(y) = H_1y_1^T \cup H_2y_2^T \cup \ldots \cup H_ny_n^T$$
$$= (0) \cup (0) \cup \ldots \cup (0);$$

*then we say the received message is a correct message and we accept y as the sent message. If on the other hand* $S(y) \neq (0) \cup (0) \cup \ldots \cup (0)$ *then we say the received message y has error. The error has to be detected as* $S(y) \neq (0) \cup (0) \cup \ldots \cup (0)$.

*Now we find the error n-vector* $e = e_1 \cup e_2 \cup \ldots \cup e_n$ *and add e with y which is the corrected message. Suppose* $y = y_1 \cup y_2 \cup \ldots \cup y_n$ *where each* $y_i$ *is a* $n_i$*-length vector with* $k_i$ *message symbols and* $n_i - k_i$ *check symbols; true for* $i = 1, 2, \ldots, n$.

*Now we find the n-coset leader for the n vector* $y = y_1 \cup y_2 \cup \ldots \cup y_n$. *Suppose* $e = e_1 \cup e_2 \cup \ldots \cup e_n$ *is the coset leader then we see* $y + e = x = x_1 \cup x_2 \cup \ldots \cup x_n$. *We have each* $e_i$ *is the coset leader of the vector* $y_i$ *in y for* $i = 1, 2, \ldots, n$.

*Now it is easily verified* $x \in C = C_1 \cup C_2 \cup \ldots \cup C_n$ *i.e.* $Hx^T = S(x) = (0) \cup (0) \cup \ldots \cup (0)$.

*Note:* It may so happen that at times when $y = y_1 \cup y_2 \cup \ldots \cup y_n$ is the received n-vector and $H = H_1 \cup H_2 \cup \ldots \cup H_n$ the parity check matrix of the n-code we have

$$\begin{aligned} Hy^T &= (H_1 \cup H_2 \cup \ldots \cup H_n)(y_1^T \cup y_2^T \cup \ldots \cup y_n^T) \\ &= H_1y_1^T \cup H_2y_2^T \cup \ldots \cup H_ny_n^T \\ &\neq (0) \cup (0) \cup \ldots \cup (0) \end{aligned}$$

but some of the $H_iy_i^T = (0)$, and some of the $H_jy_j^T \neq (0)$; $i \neq j$ so that resulting in $Hy^T \neq (0) \cup (0) \cup \ldots \cup (0)$ in which case we take the coset leaders of those i for which $H_iy_i^T = (0)$ to be (0 0 ... 0) i.e. the coset leader of the n-code, $C = C_1 \cup C_2 \cup \ldots \cup C_n$. Thus we may not or need not in general have each $H_ty_t^T = (0)$; t = 1, 2, ..., n; some may be zero and some may not be zero, at times it may so happen each $H_ty_t^T \neq (0)$, t = 1, 2, ..., n in which case we will have each coset leader to be non zero. We have illustrated this case in the 5-code given in the earlier example.

Now we proceed onto give the pseudo best n-approximation to the received message when the received message has errors.

This method gives us approximately the best n-code word which is very close to the received message. To this end we define the new notion of pseudo best n-approximation and pseudo n inner product of a finite dimensional n-vector space V = $V_1 \cup V_2 \cup \ldots \cup V_n$ over n-field $Z_2 \cup Z_2 \cup \ldots \cup Z_2$. For the sake of completeness we recall the definition of n-vector space, n-subvector space and a n-code.

Just for the sake of clarity and simplicity of understanding we use the prime field $Z_2$ of characteristic 2. It is stated that the results and definitions are true for any prime field $Z_p$ of characteristic p (p a prime) or even over any non prime field of characteristic p.

**DEFINITION 2.2.3:** *Let $V = V_1 \cup V_2 \cup \ldots \cup V_n$ be a n-vector space of finite dimension $(n_1, n_2, \ldots, n_n)$ over the n-field $F = Z_2 \cup Z_2 \cup \ldots \cup Z_2$ i.e. each $V_i$ is a $n_i$ dimensional vector space over $Z_2$; $i = 1, 2, \ldots, n$. The pseudo n-inner product on the n-vector space $V = V_1 \cup V_2 \cup \ldots \cup V_n$ is a n-map $\langle \, , \rangle_p$ (Here p = 2) from $V \times V \to F = Z_2 \cup Z_2 \cup \ldots \cup Z_2$ satisfying the following conditions.*

(1) *The n-map $\langle \, , \rangle_p = \langle \, \rangle^1_p \cup \langle \, \rangle^2_p \cup \ldots \cup \langle \, \rangle^n_p : (V_1 \cup V_2 \cup \ldots \cup V_n)$ x $(V_1 \cup V_2 \cup \ldots \cup V_n) \to Z_2 \cup Z_2 \cup \ldots \cup Z_2$ such that $\langle x, x \rangle_p \geq 0$ for all $x \in V$ i.e. $(x_1 \cup x_2 \cup \ldots \cup x_n) \in V_1 \cup V_2 \cup \ldots \cup V_n$ i.e. $\langle x_1, x_1 \rangle^1_p \cup \langle x_2, x_2 \rangle^2_p \cup \ldots \cup \langle x_n, x_n \rangle^n_p \geq 0 \cup 0 \cup \ldots \cup 0$ i.e. each $\langle x_i, x_i \rangle^i_p \geq 0$ for $i = 1, 2, \ldots, n$ with $x_i \in V_i$ i.e. $\langle \, , \rangle^i : V_i$ x $V_i \to Z_2$ and $\langle x, x \rangle_p = 0$ does not imply each $x = x_1 \cup x_2 \cup \ldots \cup x_n = (0) \cup (0) \cup \ldots \cup (0)$.*

(2) *$\langle x, y \rangle_p = \langle y, x \rangle_p$ for all $x, y \in V = V_1 \cup V_2 \cup \ldots \cup V_n$ where $x = x_1 \cup x_2 \cup \ldots \cup x_n$ and $y = y_1 \cup y_2 \cup \ldots \cup y_n$ i.e. $\langle x_i, y_i \rangle^i_p = \langle y_i, x_i \rangle^i_p$; $i = 1, 2, \ldots, n$.*

(3) *$\langle x + z, y \rangle_p = \langle x, y \rangle_p + \langle z, y \rangle_p$ for all $x, y, z \in V$ i.e. $\langle (x_1 \cup x_2 \cup \ldots \cup x_n) + (z_1 \cup z_2 \cup \ldots \cup z_n) \, , \, (y_1 \cup y_2 \cup \ldots \cup y_n) \rangle = \{(\langle x_1, y_1 \rangle^1_p \cup \langle x_2, y_2 \rangle^2_p \cup \ldots \cup \langle x_n, y_n \rangle^n_p)\} + \{ \langle z_1, y_1 \rangle^1_p \cup \langle z_2, y_2 \rangle^2_p \cup \ldots \cup \langle z_n, y_n \rangle^n_p \}$ for all $x, y, z \in V$.*

(4) $\langle x, y + z\rangle_p = \langle x_1 \cup x_2 \cup ... \cup x_n, (y_1 \cup y_2 \cup ... \cup y_n) + (z_1 \cup z_2 \cup ... \cup z_n)\rangle = \langle x_1 \cup x_2 \cup ... \cup x_n, (y_1 + z_1) \cup (y_2 + z_2) \cup ... \cup (y_n + z_n)\rangle = (\langle x_1, y_1\rangle^1_p \cup \langle x_2, y_2\rangle^2_p \cup ... \cup \langle x_n, y_n\rangle^n_p) + (\langle x_1, z_1\rangle^1_p \cup \langle x_2, z_2\rangle^2_p \cup ... \cup \langle x_n, z_n\rangle^n_p)$ *for all* $x, y, z \in V = V_1 \cup V_2 \cup ... \cup V_n$.

(5) $\langle \alpha x, y\rangle_p = \alpha \langle x, y\rangle_p$ *where*

$$\alpha = \alpha_1 \cup \alpha_2 \cup ... \cup \alpha_n \in Z_2 \cup Z_2 \cup ... \cup Z_2$$
$$x = x_1 \cup x_2 \cup ... \cup x_n \in V = V_1 \cup V_2 \cup ... \cup V_n$$

*and*

$$y = y_1 \cup y_2 \cup ... \cup y_n \in V = V_1 \cup V_2 \cup ... \cup V_n.$$

$\langle \alpha x, y\rangle_p$

$= \langle(\alpha_1 \cup \alpha_2 \cup ... \cup \alpha_n)(x_1 \cup x_2 \cup ... \cup x_n), (y_1 \cup y_2 \cup ... \cup y_n)\rangle_p$

$= \langle(\alpha_1 x_1 \cup \alpha_2 x_2 \cup ... \alpha_n x_n), (y_1 \cup y_2 \cup ... \cup y_n)\rangle_p$

$= \langle \alpha_1 x_1, y_1\rangle^1_p \cup \langle \alpha_2 x_2, y_2\rangle^2_p \cup ... \cup \langle \alpha_n x_n, y_n\rangle^n_p$

$= \alpha_1 \langle x_1, y_1\rangle^1_p \cup \alpha_2 \langle x_2, y_2\rangle^2_p \cup ... \cup \alpha_n \langle x_n, y_n\rangle^n_p$

$= (\alpha_1 \cup \alpha_2 \cup ... \cup \alpha_n) [\langle x_1, y_1\rangle^1_p \cup \langle x_2, y_2\rangle^2_p \cup ... \cup \langle x_n, y_n\rangle^n_p]$

$= \alpha [\langle x, y\rangle_p]$.

*Similarly* $\langle x, \beta y\rangle_p = \beta \langle x, y\rangle_p$ *for all* $\alpha, \beta \in F = Z_2 \cup Z_2 \cup ... \cup Z_2$ *and* $x, y \in V = V_1 \cup V_2 \cup ... \cup V_n$. *We say the n-vector space* $V = V_1 \cup V_2 \cup ... \cup V_n$ *over the n-field* $F = Z_2 \cup Z_2 \cup ... \cup Z_2$ *of characteristic 2 to be a pseudo n-inner product space if the n-map* $\langle\, ,\rangle_p = \langle\,\rangle^1_p \cup \langle\,\rangle^2_p \cup ... \cup \langle\,\rangle^n_p$ *which is a pseudo n-inner product defined on the n-vector space* $V = V_1 \cup V_2 \cup ... \cup V_n$ *over the pseudo n-field* $F = Z_2 \cup Z_2 \cup ... \cup Z_2$.

Now we define the notion of n-subspace of the n-vector space V over the n-field F.

**DEFINITION 2.2.4:** *Let* $V = V_1 \cup V_2 \cup ... \cup V_n$ *be a n-vector space over the n-field* $F = Z_2 \cup Z_2 \cup ... \cup Z_2$ *of dimension* $(n_1 \cup n_2 \cup \ldots \cup n_n)$ *over F or of dimension* $(n_1, n_2, \ldots, n_n)$ *over* $F = Z_2 \cup Z_2 \cup \ldots \cup Z_2$, *then we have* $V \cong Z_2^{n_1} \cup Z_2^{n_2} \cup ... \cup Z_2^{n_n}$. *Let W be a proper subset of V i.e.* $W = W_1 \cup W_2 \cup ... \cup W_n$ *where* *each $W_i$ is a subspace of dimension $k_i$ over $Z_2$, for i = 1, 2, ..., n, then W is a n-subspace of dimension ($k_1 \cup k_2 \cup \ldots \cup k_n$) or ($k_1$, $k_2$, ..., $k_n$) dimension over F = $Z_2 \cup Z_2 \cup \ldots \cup Z_2$ i.e. W = $W_1 \cup W_2 \cup \ldots \cup W_n$ where each $W_i$ is a subspace of dimension $k_i$ over $Z_2$ for i = 1, 2, ..., n then W is a n-subspace of dimension ($k_1 \cup k_2 \cup \ldots \cup k_n$) or ($k_1$, $k_2$, ..., $k_n$) dimension over F = $Z_2 \cup Z_2 \cup \ldots \cup Z_2$ i.e.*

$$\begin{aligned} W &= W_1 \cup W_2 \cup \ldots \cup W_n \\ &\cong Z_2^{k_1} \cup Z_2^{k_2} \cup \ldots \cup Z_2^{k_n} \subseteq Z_2^{n_1} \cup Z_2^{n_2} \cup \ldots \cup Z_2^{n_n} \end{aligned}$$

*i.e. $k_i \leq n_i$, for i = 1, 2, ..., n.*

*Now we proceed on to define the notion of pseudo best n-approximation of any n-vector $\beta = \beta_1 \cup \beta_2 \cup \ldots \cup \beta_n \in V$ relative to the ($k_1$, $k_2$, ..., $k_n$) subspace W of V; i.e. W = $W_1 \cup W_2 \cup \ldots \cup W_n \subset V = V_1 \cup V_2 \cup \ldots \cup V_n$ is a n-subspace of dimension ($k_1$, $k_2$, ... $k_n$) i.e. $W \cong Z_2^{k_1} \cup Z_2^{k_2} \cup \ldots \cup Z_2^{k_n} \subseteq Z_2^{n_1} \cup Z_2^{n_2} \cup \ldots \cup Z_2^{n_n}$. Let the n-vector space V = $V_1 \cup V_2 \cup \ldots \cup V_n$ over F = $Z_2 \cup Z_2 \cup \ldots \cup Z_n$ be a pseudo n-inner product space with the pseudo n-inner product $\langle\ \rangle_p = \langle\ \rangle_p^1 \cup \langle\ \rangle_p^2 \cup \ldots \cup \langle\ \rangle_p^n$ defined on V. Let $\beta = \beta_1 \cup \beta_2 \cup \ldots \cup \beta_n \in V = V_1 \cup V_2 \cup \ldots \cup V_n$. Let W = $W_1 \cup W_2 \cup \ldots \cup W_n \subseteq V = V_1 \cup V_2 \cup \ldots \cup V_n$ be a n-subspace of dimension ($k_1 \cup k_2 \cup \ldots \cup k_n$) of V. The pseudo best n-approximation to $\beta$ related to W is defined as follows. Let $\alpha$ = { $\alpha_1$, $\alpha_2$, ..., $\alpha_k$ } be a n-basis of W i.e. $\alpha_1 = \alpha_1^1 \cup \alpha_2^1 \cup \ldots \cup \alpha_n^1$, $\alpha_2 = \alpha_1^2 \cup \alpha_2^2 \cup \ldots \cup \alpha_n^2$, ..., $\alpha_k = \alpha_1^k \cup \alpha_2^k \cup \ldots \cup \alpha_n^k$ is the n-basis of W = $W_1 \cup W_2 \cup \ldots \cup W_n$ i.e., $\alpha$ = { $\alpha_1^1 \cup \alpha_2^1 \cup \ldots \cup \alpha_n^1$, $\alpha_1^2 \cup \alpha_2^2 \cup \ldots \cup \alpha_n^2$, ..., $\alpha_1^k \cup \alpha_2^k \cup \ldots \cup \alpha_n^k$ } is the chosen n-basis of the k = ($k_1 \cup k_2 \cup \ldots \cup k_n$) dimensional subspace W of V. The pseudo best n-approximation to $\beta$ related to W if it exists is defined as follows:*

$$\sum_{i=1}^{k} \langle \beta, \alpha_i \rangle_p \alpha_i = \bigcup_{i=1}^{k} \sum_{t=1}^{n} \langle \beta_t, \alpha_i^t \rangle_p^t \alpha_i^t .$$

*If* $\sum_{i=1}^{k} \langle \beta, \alpha_i \rangle_p \alpha_i = 0$ *we say the pseudo best n-approximation does not exist for this set of k-basis* $\{\alpha_1, ..., \alpha_k\} = \{(\alpha_1^1\ \alpha_2^1\ ...\ \alpha_n^1) \cup (\alpha_1^2\ \alpha_2^2\ ...\ \alpha_n^2) \cup ... \cup (\alpha_1^k\ \alpha_2^k\ ...\ \alpha_n^k)\}$*. In this case we choose yet another set of k-basis for W say* $(\gamma_1, \gamma_2, ..., \gamma_k)$ *and calculate* $\sum_{i=1}^{k} \langle \beta, \gamma_i \rangle_p \gamma_i$ *and take this as the pseudo best n-approximation to* $\beta$*.*

Now we apply this pseudo best n-approximation to get the most likely n-code word. Let C be the n-code over the n field $Z_2 \cup ... \cup Z_2 = F$; i.e. $C = C_1 \cup C_2 \cup ... \cup C_n$ be a $(k_1, n_1) \cup (k_2, n_2) \cup ... \cup (k_n, n_n))$ n-code i.e. C is a n-subspace of $V = V_1 \cup V_2 \cup ... \cup V_n \cong Z_2^{n_1} \cup Z_2^{n_2} \cup ... \cup Z_2^{n_n} = V$.

$C \cong Z_2^{k_1} \cup Z_2^{k_2} \cup ... \cup Z_2^{k_n}$. Now in the above definition we take C = W in the definition of n-pseudo best approximation. If some n-code word is transmitted and β is the received n-code word then

(1) If $\beta \in C$, then β is accepted as a correct message.
(2) If $\beta \notin C$ then we apply the notion of pseudo n-best approximation to β related to the n-subspace C in $Z_2^{n_1} \cup Z_2^{n_2} \cup ... \cup Z_2^{n_n}$.

Let us consider a k-basis for C say $\alpha = \{c_1, c_2, ..., c_k\}$ of C. Now $\beta \notin C$ so relative to this basis we have the pseudo best n-approximation as $\sum \langle \beta, c_i \rangle_p c_i \neq 0$ each $c_i$ is a basis of the code $C_i$ which is a $(n_i, k_i)$ code i.e., $c_i = (c_1^i, c_2^i, ..., c_{k_i}^i)$, i = 1, 2, ..., k. Thus $\alpha = \{(c_1^1, c_2^1, ..., c_{k_1}^1) \cup (c_1^2, c_2^2, ..., c_{k_2}^2) \cup ... \cup (c_1^k, c_2^k, ..., c_{k_k}^k)\}$.

If $\Sigma \langle \beta_i\ c_i \rangle_p\ c_i = 0$ then choose another set of k-basis $D = \{D^1, D^2, ..., D^k\}$ of C and find $x = \sum_{i=1}^{k} \langle \beta, D^i \rangle_p D^i$.

Then x is the pseudo best n-approximated message to $\beta$ and $x \in C = C_1 \cup C_2 \cup ... \cup C_n$.

*Note:* With the advent of computers it is left for the computer scientist to find programs to obtain the sets of basis for any C code and also a program to find the pseudo best n-approximation of $\beta$ say $\alpha$ in C, $\beta \notin C$.

It is still a challenge for them to find sets of n-basis for the n-code $C = C_1 \cup C_2 \cup ... \cup C_n$ and find the pseudo best n-approximation to $\beta$ in order to use the formula $\alpha = \sum_{i=1}^{k} \langle \beta, c_i \rangle_p c_i$.

Here it is pertinent to mention that a nice programming can be made by the computer scientist / engineer to find a method of obtaining different sets of basis. Once this is made it is easy for any coding theorist to obtain the pseudo best n-approximation to any received message which has error.

Still it is interesting to note that one can find sets of n basis for the given n-code C which is the n-subspace of the space C over $Z_2$. Using these sets of n-basis the approximately correct set of codes can be obtained.

Now using the Hamming n-distance i.e. if $x = (x_1^1, x_2^1, ..., x_{n_1}^1) \cup (x_1^2, x_2^2, ..., x_{n_2}^2) \cup ... \cup (x_1^n, x_2^n, ..., x_{n_n}^n)$ and $y = (y_1^1, y_2^1, ..., y_{n_1}^1) \cup (y_1^2, y_2^2, ..., y_{n_2}^2) \cup ... \cup (y_1^n, y_2^n, ..., y_{n_n}^n)$ are in $C = C_1 \cup C_2 \cup ... \cup C_n$, then the Hamming n-distance between x and y is the Hamming distance between each of $(x_1^i\ x_2^i\ \ldots\ x_{n_i}^i)$ and $(y_1^i\ y_2^i\ \ldots\ y_{n_i}^i)$ say $t_i$, $1 \leq i \leq n$ then the Hamming n distance between x and y is the n-tuple given by $(t_1, t_2, ..., t_n)$ where any $t_j$ is the usual Hamming distance between $(x_1^j\ x_2^j\ \ldots\ x_n^j)$ and $(y_1^j\ y_2^j\ \ldots\ y_n^j)$, j = 1, 2, ..., n. We will choose that best approximated message y′ which has least number of differences from the sent and received message y.

We illustrate this by the following example.

***Example 2.2.5:*** Let us consider the 4-code $C = C_1 \cup C_2 \cup C_3 \cup C_4$ where $C_1$ = C (4, 2), $C_2$ = C(6, 3), $C_3$ = C (6,3) and $C_4$ = C (5,2) with generator 4-matrix

$G = G_1 \cup G_2 \cup G_3 \cup G_4$ with

$$G_1 = \begin{bmatrix} 1 & 0 & 1 & 1 \\ 0 & 1 & 0 & 1 \end{bmatrix},$$

$$G_2 = \begin{bmatrix} 1 & 0 & 0 & 0 & 0 & 1 \\ 0 & 1 & 0 & 0 & 1 & 0 \\ 0 & 0 & 1 & 1 & 0 & 0 \end{bmatrix},$$

$$G_3 = \begin{bmatrix} 1 & 0 & 0 & 1 & 0 & 0 \\ 0 & 1 & 0 & 0 & 1 & 0 \\ 0 & 0 & 1 & 0 & 0 & 1 \end{bmatrix}$$

and

$$G_4 = \begin{bmatrix} 1 & 0 & 1 & 1 & 0 \\ 0 & 1 & 0 & 1 & 1 \end{bmatrix}.$$

Now the 4-codes are the collection of 4-code words

$C = C_1 \cup C_2 \cup C_3 \cup C_4$
= {(0 0 0 0), (1 0 1 1), (0 1 0 1), (1 1 1 0)} $\cup$ {(0 0 0 0 0 0), (1 0 0 0 0 1), (0 1 0 0 1 0), (0 0 1 1 0 0), (1 1 0 0 1 1), (1 0 1 1 0 1), (0 1 1 1 1 0), (1 1 1 1 1 1)} $\cup$ {(0 0 0 0 0 0), (1 0 0 1 0 0), (0 1 0 0 1 0), (0 0 1 0 0 1), (1 1 0 1 1 0), (1 0 1 1 0 1), (1 1 1 1 1 1) (0 1 1 0 1 1)} $\cup$ {(0 0 0 0 0), (1 0 1 1 0), (0 1 0 1 1), (1 1 1 0 1)}.

Any 4-code word $x = x_1 \cup x_2 \cup x_3 \cup x_4$, where $x_i \in C_i$, i = 1, 2, 3, 4. Now $C = C_1 \cup C_2 \cup C_3 \cup C_4$ is a 4-subspace of the 4-vector space $V = Z_2^4 \cup Z_2^6 \cup Z_2^6 \cup Z_2^5$ i.e., $C \cong Z_2^2 \cup Z_2^3 \cup Z_2^3 \cup Z_2^2$. Clearly all the four subspaces $Z_2^2$, $Z_2^3$, $Z_2^3$ and $Z_2^2$ are distinct as $C_i \neq C_j$; $i \neq j$; $1 \leq i, j \leq 4$.

Suppose we have received the message y = (1 1 1 1) ∪ (1 1 1 1 0 0) ∪ (0 1 1 1 1 1) ∪ (1 1 1 1 1) = $y_1 \cup y_2 \cup y_3 \cup y_4$. Clearly $y \notin C$ as the 4-syndrome $S(y) = S(y_1) \cup S(y_2) \cup S(y_3) \cup S(y_4) \neq (0) \cup (0) \cup (0) \cup (0)$. This is attained by using

$H = H_1 \cup H_2 \cup H_3 \cup H_4$

$$= \begin{bmatrix} 1 & 1 & 1 & 0 \\ 0 & 1 & 0 & 1 \end{bmatrix}$$

$$\cup \begin{bmatrix} 0 & 0 & 1 & 1 & 0 & 0 \\ 0 & 1 & 0 & 0 & 1 & 0 \\ 1 & 0 & 0 & 0 & 0 & 1 \end{bmatrix}$$

$$\cup \begin{bmatrix} 1 & 0 & 0 & 1 & 0 & 0 \\ 0 & 1 & 0 & 0 & 1 & 0 \\ 0 & 0 & 1 & 0 & 0 & 1 \end{bmatrix}$$

$$\cup \begin{bmatrix} 1 & 0 & 1 & 0 & 0 \\ 1 & 1 & 0 & 1 & 0 \\ 0 & 1 & 0 & 0 & 1 \end{bmatrix};$$

as

$$\begin{aligned} S(y) &= Hy^T \\ &= (H_1 \cup H_2 \cup H_3 \cup H_4)\,(y_1 \cup y_2 \cup y_3 \cup y_4)^T \\ &= H_1y_1^T \cup H_2y_2^T \cup H_3y_3^T \cup H_4y_4^T \\ &\neq (0) \cup (0) \cup (0) \cup (0). \end{aligned}$$

We get the basis of the 4-subspace C of V. Let

B = {(1 0 1 1), (0 1 0 1)} ∪ {(1 0 0 0 0 1), (1 0 1 1 0 1), (0 1 0 0 1 0)} ∪ {(1 0 0 1 0 0), (0 1 0 0 1 0), (0 0 1 0 0 1)} ∪ {(1 0 1 1 0) ∪ (0 1 0 1 1)}

be a 4-basis of C ⊆ V. To find the pseudo best 4-approximation of y relative to the 4-subspace C. Clearly y ∉ C. Let x be the pseudo best 4-approximation to y. Then

x = {⟨(1 1 1 1), (1 0 1 1)⟩$_p$ (1 0 1 1) + ⟨(1 1 1 1), (0 1 0 1)⟩$_p$ (0 1 0 1)} ∪ {⟨(1 1 1 1 0 0), (1 0 0 0 0 1)⟩$_p$ (1 0 0 0 0 1) + ⟨(1 1 1 1 0 0⟩, (1 0 1 1 0 1)⟩$_p$ (1 0 1 1 0 1) + ⟨(0 1 0 0 1 0), (1 1 1 1 0 0)⟩$_p$ (0 1 0 0 1 0)} ∪ {⟨(0 1 1 1 1 1), (1 0 0 1 0 0)⟩$_p$ (1 0 0 1 0 0) + ⟨(0 1 1 1 1 1), (0 1 0 0 1 0)⟩$_p$ (0 1 0 0 1 0) + ⟨(0 1 1 1 1 1), (0 0 1 0 0 1)⟩$_p$ (0 0 1 0 0 1)} ∪ ⟨(1 1 1 1 1), (1 0 1 1 0)⟩$_p$ (1 0 1 1 0) + ⟨(1 1 1 1 1), (0 1 0 1 1)⟩$_p$ (0 1 0 1 1)⟩}

= {(1 0 1 1) + 0} ∪ {(1 0 0 0 0 1) + (1 0 1 1 0 1) + (0 1 0 0 1 0)} ∪ {(1 0 0 1 0 0) + 0 + 0} ∪ {(1 0 1 1 0) + (0 1 0 1 1)}

= (1 0 1 1) ∪ (0 1 1 1 1 0) ∪ (1 0 0 1 0 0) ∪ (1 1 1 0 1).

This x belong to C. Now the 4-Hamming distance is given by (1, 3, 5, 1). Now if we {ϕ} is put small {ϕ} wish we can take this as our sent code. Otherwise we now use a different set of basis for the subcodes which gave the distance as 3 and 5 and work for a better approximation. Let the set of new-basis for C be {ϕ} ∪ {(1 0 0 0 0 1), (0 1 0 0 1 0), (0 0 1 1 0 0)} ∪ {(1 0 0 1 0 0), (0 0 1 0 0 1), (1 1 0 1 1 0)} ∪ {ϕ}. As the differences hamming distance between the $x_1$ and $y_1$ is 1 and that of $x_4$ and $y_4$ be 1. The new pseudo best 4-approximated code is

z = (1 0 1 1) ∪ {⟨(1 1 1 1 0 0) (1 0 0 0 0 1)⟩$_p$ (1 0 0 0 0 1) + ⟨(1 1 1 1 0 0), (0 1 0 0 1 0)⟩$_p$ (0 1 0 0 1 0) + ⟨(1 1 1 1 0 0), (0 0 1 1 0 0)⟩$_p$ (0 0 1 1 0 0)} ∪ {⟨(0 1 1 1 1 1) (1 0 0 1 0 0)⟩$_p$ (1 0 0 1 0 0) + ⟨(0 1 1 1 1 1), (0 0 1 0 0 1)⟩$_p$ (0 0 1 0 0 1) + ⟨(0 1 1 1 1 1), (1 1 0 1 1 0)⟩$_p$ (1 1 0 1 1 0)⟩ ∪ {(1 1 1 0 1)}

= (1 0 1 1) $\cup$ {(1 0 0 0 0 1) + (0 1 0 0 1 0) + 0} $\cup$ {(1 0 0 1 0 0) + 0 + (1 1 0 1 1 0)} $\cup$ (1 1 1 0 1)

= (1 0 1 1) $\cup$ (1 1 0 0 1 1) $\cup$ (0 1 0 0 1 0) $\cup$ (1 1 1 0 1).

Now the Hamming 4-distance between the received vector and z is (1, 4, 3, 1). Clearly z $\in$ C. Now one can take this as the received vector z or y, some may prefer z to y and others may prefer y to z.

Now we have seen the methods of finding, detecting error and correcting it by finding the pseudo best n-approximations.

The main use of pseudo best n-approximation than finding the correct code by the n-coset method is by the n-coset leader method we would get only one solution but in case of pseudo best n-approximations we can find for the closest solution as close to the received message by both varying the n-basis or the very pseudo dot product used in the construction of the pseudo vector space. So we have choice to choose the pseudo best approximated message when the received message is erroneous.

We have mainly used these n-codes to construct the new class Periyar linear codes. The n-codes will find their applications in computers, in successfully sending the message to satellites in cryptography and so on. These n-codes can also used as the data storage matrices. We as a last class of n-codes define the new notion of false n-code (n $\geq$ 2).

## 2.3 False n-matrix and Pseudo False n-matrix

We define the new notion of false n-matrix and m-pseudo false n-matrix (m < n). This is mainly introduced for the sake of both defining the 1-pseudo false generator n-matrix, m-pseudo false generator n-matrix in case of generator matrices of these special codes. Also for defining the 1-pseudo false parity check n-matrix, m-pseudo false parity check n-matrix (m < n) and the notion of false parity check matrix.

**DEFINITION 2.3.1:** *Let* $M = M_1 \cup M_2 \cup \ldots \cup M_n$ *if all of them are equal i.e.* $M_1 = M_2 = M_3 = \ldots = M_n$ *then we call M to be a false n-matrix.* $M_i$ *can be a square matrix or a rectangular* *matrix. If each of the matrix $M_i$ is a square matrix then we call M a false square n-matrix. If each of the matrix $M_i$ is a rectangular matrix then we call M to be a false rectangular n-matrix. Now when n = 2 we get the false bimatrix, when n = 3 we get the false trimatrix for $n \geq 4$ we get the false n-matrix.*

Now we proceed on to define the notion of 1-pseudo false n-matrix.

**DEFINITION 2.3.2:** *Let $M = M_1 \cup M_2 \cup \ldots \cup M_i \cup \ldots \cup M_n$ be a union of $n (n \geq 3)$ matrices if $M_1 = M_2 = \ldots = M_{i-1} = M_{i+1} = \ldots = M_n$ and $M_i$ alone distinct in this collection then we call M to be a 1-pseudo false n-matrix. If n = 2 we get the usual bimatrix. If n = 3 we get the 1-pseudo false trimatrix, for $n \geq 4$ we have 1-pseudo false n-matrix. If all the n matrices $M_1, M_2, \ldots, M_{i-1}, M_i, M_i+1, \ldots, M_n$ are $m \times m$ square matrices then we call M to be a 1-pseudo false n, $m \times m$ square matrix. If all the n-matrices $M_1, M_2, \ldots, M_{i-1}, M_i, M_{i+1}, \ldots, M_n$ are $t \times s$ rectangular matrices then we call M to be a 1-pseudo false $t \times s$ rectangular n-matrix.*

*If all the matrices $M_1, \ldots, M_{i-1}, M_{i+1}, \ldots, M_n$ are $m \times m$ square matrices and $M_i$ a $t \times t$ square matrix $(t \neq s)$ then we call M to be a 1-pseudo false mixed square n-matrix.*

*If all the matrices $M_1, M_2, \ldots, M_{i-1}, M_{i+1}, \ldots, M_n$ are rectangular $s \times t$ matrices and $M_i$ is a $p \times q$ rectangular matrix then we call M to be a 1-pseudo false mixed rectangular n-matrices.*

Now we illustrate them by simple examples.

***Example 2.3.1:*** Let $M = M_1 \cup M_2 \cup M_3 \cup M_4 \cup M_5$ where

$$M_1 = \begin{bmatrix} 3 & 0 & 1 & 2 \\ 2 & 1 & 8 & 9 \\ 1 & 6 & 5 & -1 \\ 0 & 7 & 2 & 5 \end{bmatrix} = M_2 = M_3 = M_5$$

and

$$M_4 = \begin{bmatrix} 5 & 1 & 2 & 3 \\ 4 & 5 & 6 & 7 \\ 8 & 9 & 0 & 1 \\ 2 & 3 & 4 & 5 \end{bmatrix}.$$

Clearly M is a 1-pseudo false square 5-matrix.

***Example 2.3.2:*** Let $M = M_1 \cup M_2 \cup M_3 \cup M_4$ where

$$M_1 = \begin{bmatrix} 5 & 7 & 9 & 11 & 8 \\ 1 & 2 & 3 & 4 & 5 \\ 6 & 7 & 8 & 9 & 10 \\ 11 & 12 & 13 & 4 & 3 \\ 2 & 1 & 0 & 9 & 6 \end{bmatrix} = M_2 = M_4$$

and

$$M_3 = \begin{bmatrix} 1 & 2 & 3 & 4 \\ 5 & 6 & 7 & 8 \\ 9 & 0 & 1 & 2 \\ 3 & 4 & 5 & 6 \end{bmatrix},$$

clearly M is a 1-pseudo false mixed square 4-matrix.

***Example 2.3.3:*** Consider $M = M_1 \cup M_2 \cup M_3 \cup M_4 \cup M_5 \cup M_6$ where

$$M_1 = \begin{bmatrix} 1 & 2 & 3 & 4 & 5 & 6 & 7 \\ 8 & 9 & 0 & 1 & 2 & 3 & 4 \\ 6 & 7 & 8 & 9 & 2 & 1 & 2 \\ 1 & 9 & 0 & 0 & 8 & 5 & 9 \end{bmatrix} = M_2 = M_3 = M_4 = M_5$$

and

$$M_6 = \begin{bmatrix} 0 & 1 & 0 & 1 & 9 & 6 & 5 \\ 9 & 0 & 1 & 8 & 3 & 2 & 1 \\ 1 & 3 & 5 & 7 & 6 & 1 & 5 \\ 2 & 4 & 6 & 9 & 0 & 5 & 1 \end{bmatrix}.$$

We see $M_6$ is also a rectangular 4 × 7 matrix. M is a 1-pseudo 4 × 7 rectangular false 6-matrix.

***Example 2.3.4:*** Consider the 1-pseudo false 5-matrix $N = N_1 \cup N_2 \cup N_3 \cup N_4 \cup N_5$ where

$$N_1 = \begin{bmatrix} 3 & 1 & 1 & 2 & 3 & 1 \\ 0 & 5 & 1 & 0 & 5 & 0 \\ 1 & 0 & 1 & 2 & 6 & 3 \\ 2 & 1 & 1 & 4 & 0 & 0 \end{bmatrix} = N_3 = N_4 = N_5$$

and

$$N_2 = \begin{bmatrix} 3 & 1 & 2 \\ 9 & 8 & 7 \\ 6 & 7 & 8 \\ 4 & 3 & 2 \\ 1 & 0 & 5 \\ 1 & 1 & 3 \end{bmatrix}.$$

We see $N_i$, i = 1, 3, 4, 5 is the same 4 × 6 rectangular matrix where as $N_2$ is a 6 × 3 rectangular matrix. Hence N is a 1-pseudo false mixed rectangular 5-matrix.

Now we give examples of 1-pseudo false mixed n-matrices.

***Example 2.3.5:*** Consider the 1-pseudo false 5-matrix $M = M_1 \cup M_2 \cup M_3 \cup M_4 \cup M_5$ where

$$M_1 = M_2 = M_3 = M_4 = \begin{bmatrix} 3 & 9 & 0 & 1 & 6 \\ 0 & 8 & 1 & 2 & 7 \\ 1 & 7 & 2 & 3 & 8 \\ 1 & 6 & 3 & 4 & 9 \\ 4 & 5 & 4 & 5 & 0 \end{bmatrix},$$

a 5 × 5 square matrix and

$$M_5 = \begin{bmatrix} 3 & 4 & 7 & 3 & 7 & 2 \\ 1 & 5 & 8 & 1 & 5 & 6 \\ 2 & 6 & 9 & 8 & 4 & 9 \end{bmatrix}$$

a 3 × 6 rectangular matrix.

We call M to be 1-pseudo mixed false 5-matrix.

***Example 2.3.6:*** Let $M = M_1 \cup M_2 \cup M_3 \cup M_4 \cup M_5 \cup M_6$ be a 1-pseudo false 6-matrix where

$$M_2 = M_3 = M_4 = M_5 = M_6 = \begin{bmatrix} 3 & 4 & 1 & 2 \\ 9 & 8 & 7 & 6 \\ 1 & 2 & 3 & 4 \\ 5 & 0 & 1 & 0 \\ 7 & 8 & 1 & 2 \\ 9 & 3 & 0 & 9 \end{bmatrix}$$

and

$$M_1 = \begin{bmatrix} 2 & 1 & 4 & 3 \\ 9 & 8 & 6 & 7 \\ 3 & 9 & 7 & 5 \\ 4 & 0 & 1 & 9 \end{bmatrix}$$

the 4 × 4 square matrix. Clearly M is a 1-pseudo mixed false 6-matrix.

Now having seen all the 5 types of 1-pseudo false n-matrices we now proceed on to define the notion of 2-pseudo false n-matrix ($n \geq 4$). We see when n = 3 we will have only false trimatrix and the notion of 1-pseudo false trimatrix, cannot exist i.e., we do not have the notion of 2-pseudo false trimatrix. For 2-pseudo false trimatrix to be defined we need $n \geq 4$.

**DEFINITION 2.3.3:** *Let $M = M_1 \cup M_2 \cup \ldots \cup M_n$ ($n \geq 4$) be a false matrix in which $M_1 = M_2 = \ldots = M_{n-2}$ and $M_{n-1} = M_n$. We call M to be a 2-pseudo false n-matrix.*

Before we proceed on to discuss more properties about these class of matrices we give an example of it.

***Example 2.3.7:*** Let $M = M_1 \cup M_2 \cup M_3 \cup M_4 \cup M_5 \cup M_6 \cup M_7$ be a false 7-matrix.

Suppose

$$M_1 = \begin{bmatrix} 3 & 4 & 5 & 6 & 7 & 8 & 1 \\ 1 & 0 & 2 & 3 & 4 & 5 & 6 \\ 7 & 8 & 9 & 0 & 1 & 2 & 3 \end{bmatrix} = M_3 = M_4 = M_5 = M_6$$

and

$$M_2 = \begin{bmatrix} 3 & 0 & 1 & 2 & 3 \\ 4 & 4 & 5 & 6 & 7 \\ 5 & 8 & 9 & 0 & 1 \\ 7 & 2 & 3 & 4 & 5 \\ 1 & 6 & 7 & 8 & 9 \end{bmatrix} = M_7$$

then we call M to be a 2-pseudo false 7-matrix.

It is very important to note that in case of 2-pseudo false n matrix $n \geq 4$. We give a few subclasses of 2-pseudo false n-matrices.

**DEFINITION 2.3.4:** *Let $M = M_1 \cup M_2 \cup \ldots \cup M_n$ be a 2-pseudo false n-matrix, suppose $M_1 = M_2 = \ldots = M_{n-2}$ be the same square m × m matrix*

$$\begin{bmatrix} a_{11} & a_{12} & \cdots & a_{1m} \\ a_{21} & a_{22} & \cdots & a_{2m} \\ \vdots & \vdots & \vdots & \vdots \\ a_{m1} & a_{m2} & \cdots & a_{mm} \end{bmatrix}$$

*and*

$$M_{n-1} = M_n = \begin{bmatrix} b_{11} & b_{12} & \cdots & b_{1m} \\ b_{21} & b_{22} & \cdots & b_{2m} \\ \vdots & \vdots & \vdots & \vdots \\ b_{m1} & b_{m2} & \cdots & b_{mm} \end{bmatrix}$$

*a square m × m matrix where $M_i \neq M_{n-1}$; $i = 1, 2, \ldots, n-2$, then we call M to be a 2-pseudo m × m square false n-matrix.*

We give an example of a 2-pseudo m × m square false n-matrix.

***Example 2.3.8:*** Let $M = M_1 \cup M_2 \cup M_3 \cup M_4 \cup M_5 \cup M_6$ where

$$M_1 = M_3 = M_4 = M_5 = \begin{bmatrix} 3 & 2 & 0 & 1 & 1 \\ 5 & 1 & 9 & 8 & 3 \\ 7 & 3 & 2 & 9 & 0 \\ 8 & 5 & 1 & 8 & 9 \\ 6 & 6 & 7 & 1 & 0 \end{bmatrix}$$

and

$$M_2 = M_6 = \begin{bmatrix} 1 & 2 & 3 & 4 & 5 \\ 6 & 7 & 8 & 9 & 0 \\ 1 & 3 & 5 & 7 & 9 \\ 2 & 4 & 6 & 8 & 0 \\ 8 & 7 & 9 & 2 & 1 \end{bmatrix}$$

be 5 × 5 matrices. M is a 2-pseudo false square 6-matrix.

**DEFINITION 2.3.5:** *Let $G = G_1 \cup G_2 \cup G_3 \cup G_4 \cup \ldots \cup G_m$ be a 2-pseudo false M matrix if $G_1 = G_2 = G_3 = \ldots = G_{m-2} = (a_{ij})$ be a square n × n matrix and $G_{m-1} = G_m$ be a square t × t matrix ($t \neq n$) then we call G to be the 2-pseudo false m-matrix to be a 2-pseudo false mixed square m-matrix.*

We illustrate this by the following example.

***Example 2.3.9:*** Consider the 2-pseudo false 7-matrix given by $M = M_1 \cup M_2 \cup M_3 \cup M_4 \cup M_5 \cup M_6 \cup M_7$ where

$$M_1 = M_3 = M_4 = M_6 = M_7 = \begin{bmatrix} 3 & 8 & 11 & 4 & 5 \\ 7 & 12 & 1 & 2 & 3 \\ 4 & 5 & 6 & 7 & 8 \\ 9 & 10 & 11 & 12 & 13 \\ 1 & 2 & 12 & 5 & 15 \end{bmatrix}$$

a 5 × 5 square matrix.

$$M_2 = M_5 = \begin{bmatrix} 3 & 1 & 7 \\ 1 & 2 & 3 \\ 4 & 5 & 6 \end{bmatrix}$$

be a 3 × 3 square matrix. We call M to be a 2-pseudo false mixed square 7-matrix.

Now we proceed on to define the notion of 2-pseudo false rectangular n-matrix.

**DEFINITION 2.3.6**: *Let $M = M_1 \cup M_2 \cup \ldots \cup M_n$ be a 2-pseudo false n-matrix. We call M to be a 2-pseudo false s × t rectangular n-matrix if $M_1 = M_2 = \ldots = M_{n-2} = (a_{ij})$, $1 \leq i \leq s$ and $1 \leq j \leq t$ and $M_{n-1} = M_{n-2} = (b_{ij})$; $1 \leq i \leq s$ and $1 \leq j \leq t$; $a_{ij} \neq b_{ij}$.*

We now illustrate this by a simple example.

***Example 2.3.10:*** Let $M = M_1 \cup M_2 \cup M_3 \cup M_4 \cup M_5$ where

$$M_1 = M_3 = M_5 = \begin{bmatrix} 3 & 0 & 1 & 2 & 3 & 4 \\ 1 & 5 & 6 & 7 & 8 & 9 \\ 2 & 1 & 1 & 0 & 2 & 2 \\ 5 & 3 & 0 & 4 & 0 & 5 \end{bmatrix}$$

be a 4 × 6 rectangular matrix and

$$M_2 = M_4 = \begin{bmatrix} 1 & 5 & 9 & 3 & 7 & 1 \\ 2 & 6 & 0 & 4 & 8 & 2 \\ 3 & 7 & 1 & 5 & 9 & 3 \\ 4 & 8 & 2 & 6 & 0 & 5 \end{bmatrix}$$

a 4 × 6 rectangular matrix. We call M to be a 2-pseudo false 4 × 6 rectangular 5- matrix.

Now we proceed on to define the notion of 2-pseudo false mixed rectangular 5-matrix.

**DEFINITION 2.3.7:** *Consider the 2-pseudo false m-matrix $N = N_1 \cup N_2 \cup N_3 \cup \ldots \cup N_m$ ($m \geq 4$). We say N is a 2-pseudo false mixed rectangular m-matrix if $M_1 = M_2 = \ldots = M_{m-2} = (a_{ij})$ is a s × t rectangular matrix $s \neq t$ and $1 \leq i \leq s$ and $1 \leq j \leq t$ and $M_{m-1} = M_m = (b_{ij})$, a p × q rectangular matrix $p \neq q$ and $p \neq s$ (or $t \neq q$).*

Now we illustrate this by the following example.

***Example 2.3.11:*** Let $M = M_1 \cup M_2 \cup M_3 \cup M_4 \cup M_5 \cup M_6 \cup M_7 \cup M_8$ be a 2-pseudo false 8-matrix.

$$M_1 = M_2 = M_4 = M_5 = M_6 = M_8 = \begin{bmatrix} 2 & 6 & 3 & 4 & 1 & 2 & 3 \\ 0 & 1 & 2 & 3 & 4 & 5 & 6 \\ 7 & 8 & 9 & 0 & 1 & 3 & 2 \end{bmatrix}$$

be a 3 × 7 matrix. Let

$$M_3 = M_7 = \begin{bmatrix} 3 & 1 & 3 & 1 & 5 & 9 \\ 4 & 2 & 5 & 2 & 6 & 0 \\ 5 & 0 & 2 & 3 & 7 & 1 \\ 6 & 1 & 8 & 4 & 8 & 5 \end{bmatrix}$$

a 4 × 6 rectangular matrix. M is a 2-pseudo mixed rectangular false 8-matrix.

Now we proceed on to define the notion of 2-pseudo mixed false n-matrix.

**DEFINITION 2.3.8:** *Let $M = M_1 \cup M_2 \cup M_3 \cup \ldots \cup M_n$ be a 2-pseudo false n-matrix. If $M_1 = M_2 = \ldots = M_{n-2}$ be a rectangular (or square) p × q (t × t) matrix and $M_{n-1} = M_n$ be the square (or rectangular) t × t (or p × q) matrix, then we define M to be a 2-pseudo false mixed n-matrix. Here 'or' is used in the mutually exclusive sense.*

We now illustrate this by the following example.

***Example 2.3.12:*** Let $G = G_1 \cup G_2 \cup G_3 \cup G_4 \cup G_5$ be a 2-pseudo false 5-matrix where

$$G_1 = G_3 = G_5 = \begin{bmatrix} 1 & 3 & 4 & 5 & 6 & 7 & 8 \\ 2 & 9 & 8 & 7 & 6 & 5 & 4 \\ 3 & 2 & 3 & 2 & 1 & 0 & 6 \\ 4 & 1 & 0 & 1 & 5 & 2 & 1 \end{bmatrix}$$

is a 4 × 7 matrix and

$$G_2 = G_n = \begin{bmatrix} 3 & 4 & 5 & 6 \\ 7 & 8 & 9 & 0 \\ 1 & 2 & 3 & 4 \\ 0 & 1 & 9 & 3 \end{bmatrix}$$

is a 4 × 4 square matrix. Hence G is a 2-pseudo false mixed 5-matrix.

***Example 2.3.13:*** Let $P = P_1 \cup P_2 \cup P_3 \cup P_4 \cup P_5 \cup P_6 \cup P_7$ be a 2-pseudo false 7-matrix. Here

$$P_1 = P_3 = P_4 = P_5 = P_6 = \begin{bmatrix} 1 & 2 & 3 & 4 & 5 \\ 6 & 7 & 8 & 9 & 10 \\ 11 & 12 & 1 & 3 & 14 \\ 15 & 7 & 14 & 10 & 1 \\ 6 & 5 & 4 & 3 & 2 \end{bmatrix}$$

is a 5 × 5 square matrix and

$$P_2 = P_7 = \begin{bmatrix} 1 & 4 & 7 & 0 & 15 & 3 \\ 2 & 5 & 8 & 1 & 7 & 18 \\ 3 & 6 & 9 & 10 & 8 & 2 \end{bmatrix}$$

is a 3 × 6 matrix. Thus P is a 2-pseudo false mixed 7-matrix.

Now we proceed on to define the notion of 3-pseudo false n-matrix $n \geq 6$. We see when n = 3 (or 4 and 5) we can get only a 1-pseudo false matrix (or a two pseudo false matrix) respectively. For n = 4 we cannot have a 3-pseudo false matrix for it will turn out to be a 1-pseudo false matrix.

Thus for one to define even the notion of 3-pseudo false n-matrix we need $n \geq 6$.

**DEFINITION 2.3.9:** *Let $M = M_1 \cup M_2 \cup \ldots \cup M_n$ be a pseudo false n-matrix ($n \geq 6$) where $M_1 = M_2 = \ldots = M_{n-3} = (a_{ij}) = A$, a matrix and $M_{n-2} = M_{n-1} = M_n = B$ with $A \neq B$. Then we call M to be a 3-pseudo false n-matrix. If in a 3-pseudo false n-matrix M both A and B are $m \times m$ square matrices then we call M to be a 3-pseudo false square n-matrix. If in the 3-pseudo false n-matrix M both A and B happen to be a $p \times q$ rectangular matrices then we call M to be a 3-pseudo false rectangular n-matrix.*

*If in the 3-pseudo false n-matrix $M = M_1 \cup M_2 \cup \ldots \cup M_n$, $A = M_1 = M_2 = \ldots = M_{n-3}$ happen to be a $m \times m$ square matrix and $B = M_{n-2} = M_{n-1} = M_n$ be a $t \times t$ square matrix ($m \neq t$) then we call M to be a 3-pseudo mixed square false n-matrix. If in the 3-pseudo false n-matrix $M = M_1 \cup M_2 \cup \ldots \cup M_n$, $M_1 = M_2 = \ldots = M_{n-3} = A$ be a $m \times m$ square matrix (or a $t \times s$ rectangular matrix) and $M_{n-2} = M_{n-1} = M_n = B$ be a $p \times q$ rectangular matrix (or a $s \times s$ square matrix) then we call M to be a 3-pseudo false mixed n-matrix.*

Now we will illustrate this by the following examples.

***Example 2.3.14:*** Let $M = M_1 \cup M_2 \cup \ldots \cup M_8$ be a 3-pseudo false 8-matrix where

$$M_1 = M_3 = M_5 = M_7 = M_8 = \begin{bmatrix} 3 & 7 & 1 & 5 & 8 & 10 \\ 1 & 2 & 3 & 4 & 5 & 6 \\ 7 & 8 & 9 & 10 & 11 & 1 \end{bmatrix}$$

and

$$M_2 = M_4 = M_6 = \begin{bmatrix} 1 & 2 & 3 & 4 & 5 & 6 \\ 7 & 8 & 9 & 0 & 1 & 2 \\ 3 & 4 & 5 & 6 & 7 & 8 \end{bmatrix}$$

where both the matrices are 3 × 6 rectangular matrices. Clearly M is a 3-pseudo rectangular false 8-matrix.

Now we proceed on to give an example of a 3-pseudo square false n-matrix.

***Example 2.3.15:*** Let $M = M_1 \cup M_2 \cup M_3 \cup \ldots \cup M_7$ be a 3-pseudo false 7-matrix, where

$$M_2 = M_3 = M_4 = M_7 = \begin{bmatrix} 3 & 1 & 5 & 9 \\ 1 & 2 & 6 & 0 \\ 5 & 3 & 7 & 1 \\ 6 & 4 & 8 & 2 \end{bmatrix}$$

be 4 × 4 matrix and

$$M_1 = M_5 = M_6 = \begin{bmatrix} 1 & 2 & 3 & 4 \\ 5 & 6 & 7 & 8 \\ 9 & 0 & 1 & 2 \\ 3 & 4 & 5 & 6 \end{bmatrix}$$

a 4 × 4 square matrix.

We see M is a 3-pseudo square false 7-matrix.

Now we give an example of a 3-pseudo mixed square false n-matrix.

***Example 2.3.16:*** Let $M = M_1 \cup M_2 \cup M_3 \cup \ldots \cup M_6$ be a 3-pseudo mixed rectangular matrix where

$$M_1 = M_2 = M_3 = \begin{bmatrix} 1 & 5 & 9 & 3 & 7 & 1 \\ 2 & 6 & 0 & 4 & 8 & 2 \\ 3 & 7 & 1 & 5 & 9 & 3 \\ 4 & 8 & 2 & 6 & 0 & 4 \end{bmatrix}$$

and

$$M_4 = M_5 = M_6 = \begin{bmatrix} 1 & 2 & 3 & 4 & 5 \\ 6 & 7 & 8 & 9 & 0 \\ 1 & 2 & 3 & 4 & 6 \\ 7 & 8 & 9 & 0 & 1 \end{bmatrix}.$$

Clearly one set of matrices is a 4 × 6 rectangular matrix where as another set of matrices is a 4 × 5 rectangular matrices. Thus M is a 3-pseudo mixed rectangular false 6-matrix.

***Example 2.3.17:*** Let $M = M_1 \cup M_2 \cup M_3 \cup \ldots \cup M_8$ where

$$M_1 = M_2 = M_5 = M_6 = M_8 = \begin{bmatrix} 1 & 5 & 6 & 7 \\ 2 & 8 & 9 & 0 \\ 3 & 1 & 2 & 3 \\ 4 & 5 & 6 & 7 \end{bmatrix}$$

and

$$M_3 = M_4 = M_7 = \begin{bmatrix} 4 & 3 & 2 & 1 & 9 \\ 6 & 7 & 8 & 9 & 0 \\ 1 & 3 & 5 & 7 & 4 \\ 2 & 4 & 6 & 8 & 5 \\ 1 & 0 & 8 & 7 & 3 \end{bmatrix}.$$

We see M is a 3-pseudo mixed square false 8-matrix. Now we present an example of a 3-pseudo mixed false n-matrix.

***Example 2.3.18:*** Let $M = M_1 \cup M_2 \cup M_3 \cup \ldots \cup M_9$ be a 3-pseudo 9 matrix where

$$\begin{bmatrix} 1 & 2 & 3 & 4 & 5 & 6 & 7 \\ 8 & 9 & 10 & 3 & 11 & 12 & 13 \\ 3 & 14 & 15 & 6 & 16 & 8 & 9 \end{bmatrix} = M_1 = M_2 =$$

$$M_4 = M_5 = M_7 = M_8$$

and

$$M_3 = M_6 = M_9 = \begin{bmatrix} 3 & 4 & 5 & 6 \\ 1 & 2 & 3 & 4 \\ 5 & 0 & 1 & 8 \\ 1 & 8 & 9 & 0 \end{bmatrix}$$

a 4 × 4 square matrix. M is a 3-pseudo mixed false 9-matrix.

***Example 2.3.19:*** Let $M = M_1 \cup M_2 \cup M_3 \cup \ldots \cup M_7$ where

$$M_1 = M_3 = M_5 = M_7 = \begin{bmatrix} 1 & 2 & 3 & 4 \\ 5 & 8 & 11 & 1 \\ 6 & 9 & 2 & 15 \\ 7 & 10 & 5 & 8 \end{bmatrix}$$

and

$$M_2 = M_4 = M_6 = \begin{bmatrix} 1 & 4 & 5 & 6 & 7 & 8 \\ 2 & 9 & 10 & 1 & 2 & 3 \\ 3 & 4 & 5 & 6 & 7 & 11 \end{bmatrix}$$

a 3 × 6 rectangular matrix. M is a 3-pseudo mixed false 7-matrix.

Now we proceed on to define the notion of m-pseudo false n-matrix m < n and (n ≥ 2m) m > 3.

**DEFINITION 2.3.10:** *Let $M = M_1 \cup M_2 \cup \ldots \cup M_n$ be a pseudo false n-matrix we call M to be a m-pseudo false n-matrix ($n \geq 2m$ and $m \geq 4$) if $M_1 = M_2 = \ldots = M_m = A$ and $M_{m+1} = M_{m+2} = \ldots = M_n = B$ with $A \neq B$. If in the m-pseudo false n-matrix ($n \geq 2m$ and $m \geq 4$) $M = M_1 \cup M_2 \cup \ldots \cup M_n$ if we have $M_1 = M_2 = \ldots = M_m = A$ be a $t \times t$ square matrix and $M_{m+1} = M_{m+2} = \ldots =$*

*$M_n$ = B be a t × t square matrix but A ≠ B then we define M to be a m-pseudo t × t square false n-matrix.*

We see if the conditions on m and n are not put we may not have the definition to be true.

Now we proceed on to give an illustration of the same.

***Example 2.3.20***: Let $M = M_1 \cup M_2 \cup \ldots \cup M_9$ be a false 9-matrix where $M_1 = M_3 = M_5 = M_7 = M_9 = A$ and $M_2 = M_4 = M_6 = M_8 = B$ (A ≠ B). We see M is a 4-pseudo false 9-matrix.

The following observations are important. (1) when the false n-matrix is a false 9-matrix we cannot have 5-pseudo 9-matrix to be defined we can have only 4-pseudo 9-matrix or 3-pseudo 9-matrix or 2-pseudo 9-matrix or 1-pseudo 9-matrix to be defined. Clearly 9 > 2.4 here n = 9 and m = 4.

***Example 2.3.21***: Let us consider $M = M_1 \cup M_2 \cup \ldots \cup M_{12}$ be a pseudo false 12-matrix. We see M can only be a maximum 6-pseudo false 12-matrix however it can be a m-pseudo false 12-matrix m = 1, 2, 3, 4, 5 and 6. Clearly M is not a 7-pseudo false 12-matrix. We see if $M_1 = M_3 = M_5 = M_7 = M_9 = M_{11} = A$ and $M_2 = M_4 = M_6 = M_8 = M_{10} = M_{12} = B$ (A ≠ B) then M is a 6-pseudo false 12-matrix.

Now we proceed on to define various classes of m-pseudo false n-matrices (n ≥ 2m, m < n).

**DEFINITION 2.3.11:** *Let $M = M_1 \cup M_2 \cup \ldots \cup M_n$ (n ≥ 2m and m < n) be a m-pseudo false n-matrix. If $M_1 = M_2 = M_5 = \ldots = M_m = A$ where A is a t × t square matrix and $M_{m+1} = M_{m+2} = \ldots = M_n = B$, B also a t × t square matrix with A ≠ B then we call M to be a m-pseudo square false n-matrix. If on the other hand in $M_1 = M_2 = \ldots = M_m = A$, A happens to be a p × q rectangular matrix (p ≠ q) and $M_{m+1} = M_{m+2} = \ldots = M_n = B$, (A ≠ B) but B is also a p × q rectangular matrix than we call M to be a m-pseudo rectangular false n-matrix. Suppose we have $M_1 = M_2 = \ldots = M_m = A$ to be a p × p square matrix and $M_{m+1} = M_{m+2} = \ldots = M_n = B$ to be a t × t square matrix (p ≠ t) then we call M to be a m-pseudo mixed square false n-matrix if we have in the false n matrix $M = M_1 \cup M_2 \cup \ldots \cup M_n$, $M_1 = M_2 = \ldots =$*

*$M_m$ = A, A a p × q rectangular matrix and $M_{m+1} = M_{m+2} = \ldots = M_n$ = B, a t × s rectangular matrix p ≠ t (or q ≠ s) (or not in the mutually exclusive sense) then we define M to be a m-pseudo mixed rectangular false n-matrix.*

*Let $M = M_1 \cup M_2 \cup \ldots \cup M_n$ be a m-pseudo false n-matrix we say M is a m-pseudo mixed false n-matrix if $M = M_1 \cup M_2 \cup \ldots \cup M_m$ = A is a p × p square matrix (or a t × s t ≠ s a rectangular matrix) and $M_{m+1} = M_{m+2} = \ldots = M_n$ = B is a t × s rectangular (or a p × p square matrix) then we call M to be a m-pseudo mixed false n-matrix.*

Now we illustrate these by the following examples:-

***Example 2.3.22:*** Let $M = M_1 \cup M_2 \cup \ldots \cup M_{10}$ be 4-pseudo false matrix where

$$M_1 = M_3 = M_4 = M_6 = M_7 = M_8 = \begin{bmatrix} 3 & 1 & 2 & 3 \\ 1 & 4 & 5 & 6 \\ 2 & 7 & 8 & 9 \\ 5 & 0 & 1 & 2 \end{bmatrix}$$

a 4 × 4 square matrix and

$$M_2 = M_5 = M_9 = M_{10} = \begin{bmatrix} 0 & 9 & 8 & 7 \\ 6 & 5 & 4 & 3 \\ 2 & 1 & 0 & 1 \\ 2 & 3 & 4 & 5 \end{bmatrix}$$

a 4 × 4 square matrix. Clearly M is a 4-pseudo false to matrix which is a 4-pseudo square false 10 matrix.

***Example 2.3.23:*** Let $M = M_1 \cup M_2 \cup M_3 \cup \ldots \cup M_{12}$ be a 6-pseudo false 12-matrix where

$$M_1 = M_3 = M_5 = M_7 = M_9 = M_{11} = \begin{bmatrix} 1 & 4 & 5 & 6 & 7 & 8 & 9 \\ 2 & 9 & 7 & 5 & 3 & 1 & 8 \\ 3 & 8 & 6 & 4 & 2 & 0 & 5 \end{bmatrix}$$

be a 3 × 7 rectangular matrix and

$$M_2 = M_4 = M_6 = M_8 = M_{10} = M_{12} = \begin{bmatrix} 9 & 6 & 5 & 4 & 3 & 2 & 1 \\ 8 & 1 & 3 & 5 & 7 & 9 & 8 \\ 7 & 2 & 4 & 6 & 8 & 0 & 5 \end{bmatrix}$$

be a 3 × 7 rectangular matrix. We see M is a 6-pseudo rectangular false 12-matrix.

We see a twelve matrix can maximum be only a 6-pseudo false 12-matrix and can never be a 7 or 8 or 9 or 10 or 11 pseudo false 12-matrix.

***Example 2.3.24:*** Consider the false n-matrix $M = M_1 \cup M_2 \cup \ldots \cup M_{15}$ where

$$M_1 = M_2 = M_3 = M_4 = M_5 = M_6 = M_{15} = \begin{bmatrix} 3 & 8 & 9 & 0 \\ 4 & 1 & 2 & 3 \\ 5 & 4 & 5 & 6 \\ 6 & 7 & 8 & 9 \end{bmatrix}$$

be a 4 × 4 square matrix and

$$M_7 = M_8 = M_9 = M_{10} = M_{11} = M_{12} = M_{13} = M_{14} = \begin{bmatrix} 3 & 2 & 1 \\ 0 & 5 & 6 \\ 9 & 8 & 7 \end{bmatrix}$$

a 3 × 3 square matrix. Clearly M is a 7-pseudo mixed square false 15-matrix.

***Example 2.3.25:*** Suppose $M = M_1 \cup M_2 \cup M_3 \cup \ldots \cup M_{17}$ be a false pseudo 17-matrix where

$$M_1 = M_2 = M_6 = M_7 = M_8 = M_9 = M_{10} = M_{17}$$

$$= \begin{bmatrix} 1 & 5 & 6 & 7 & 8 & 9 & 0 \\ 2 & 1 & 2 & 3 & 4 & 5 & 6 \\ 3 & 7 & 8 & 9 & 0 & 11 & 1 \\ 4 & 12 & 6 & 8 & 4 & 3 & 10 \end{bmatrix}$$

a 4 × 7 rectangular matrix and

$$M_3 = M_4 = M_5 = M_{11} = M_{12} = M_{13} = M_{14} = M_{15} = M_{16}$$

$$= \begin{bmatrix} 3 & 1 & 2 \\ 5 & 9 & 0 \\ 6 & 4 & 2 \\ 7 & 1 & 9 \\ 8 & 3 & 5 \end{bmatrix}$$

be a 5 × 3 rectangular matrix, then M is a 8-pseudo mixed rectangular false 17-matrix.

***Example 2.3.26:*** Let us consider the 18-matrix $M = M_1 \cup M_2 \cup M_3 \cup \ldots \cup M_{18}$ where

$$M_1 = M_2 = M_3 = \ldots = M_8 = M_9 = \begin{bmatrix} 1 & 6 & 7 & 8 & 9 \\ 2 & 10 & 9 & 8 & 7 \\ 3 & 6 & 5 & 4 & 3 \\ 4 & 2 & 1 & 11 & 9 \\ 5 & 8 & 7 & 6 & 5 \end{bmatrix}$$

a 5 × 5 square matrix and

$$M_{10} = M_{11} = \ldots = M_{18} = \begin{bmatrix} 3 & 0 & 1 & 2 & 3 & 4 & 5 & 6 & 7 \\ 1 & 8 & 9 & 1 & 9 & 2 & 8 & 7 & 3 \\ 4 & 4 & 5 & 9 & 8 & 7 & 6 & 5 & 4 \\ 5 & 3 & 2 & 1 & 4 & 5 & 6 & 7 & 8 \\ 6 & 0 & 1 & 3 & 0 & 7 & 0 & 9 & 0 \end{bmatrix}$$

a 5 × 9 rectangular matrix. We see M is a 9-pseudo mixed false 18-matrix.

***Example 2.3.27:*** Let $M = M_1 \cup M_2 \cup \ldots \cup M_9$ be a false 9-matrix, where

$$M_1 = M_3 = M_5 = M_7 = M_9 = \begin{bmatrix} 1 & 0 & 3 & 4 & 8 & 9 \\ 2 & 5 & 1 & 0 & 2 & 1 \\ 3 & 4 & 8 & 7 & 6 & 0 \\ 4 & 1 & 2 & 3 & 4 & 5 \end{bmatrix}$$

be a 4 × 6 rectangular matrix and

$$M_2 = M_4 = M_6 = M_8 = \begin{bmatrix} 3 & 6 & 1 \\ 4 & 7 & 0 \\ 5 & 10 & 9 \end{bmatrix}$$

be a a 3 × 3 square matrix. M is a 4-pseudo mixed false 9-matrix. Clearly M is not a 5-pseudo mixed false 9-matrix.

## 2.4 False n-codes ($n \geq 2$)

In this section we introduce yet another new class of codes called false n-codes which can find its applications in cryptography. If false n-codes are used it is very difficult to hack the secrets. They can also be used in the place of codes with ARQ protocols.

**DEFINITION 2.4.1:** *Let $C = C_1 \cup C_2$ be a bicode if $C_1 = C_2$ then we call C is a false bicode.*

***Example 2.4.1:*** Let $C = C_1 \cup C_2$ where $C_1$ is generated by

$$G_1 = \begin{bmatrix} 1 & 0 & 0 & 0 & 0 & 1 \\ 0 & 1 & 0 & 0 & 1 & 0 \\ 0 & 0 & 1 & 1 & 0 & 0 \end{bmatrix}$$

and $C_2$ generated by

$$G_2 = \begin{bmatrix} 1 & 0 & 0 & 0 & 0 & 1 \\ 0 & 1 & 0 & 0 & 1 & 0 \\ 0 & 0 & 1 & 1 & 0 & 0 \end{bmatrix}.$$

i.e. C = {(0 0 0 0 0 0), (1 0 0 0 0 1), (0 1 0 0 1 0), (0 0 1 1 0 0), (1 1 0 0 1 1), (1 0 1 1 0 1), (0 1 1 1 1 0), (1 1 1 1 1 1)} $\cup$ {(0 0 0 0 0 0), (1 0 0 0 0 1), (0 1 0 0 1 0), (0 0 1 1 0 0), (1 1 0 0 1 1); (1 0 1 1 0 1), (0 1 1 1 1 0), (1 1 1 1 1 1)}.

This false bicode has the following main advantages or purposes.

1. If they want ARQ messages they can use false bicode so that when the message is received if both vary one can take the better of the two or at times one may be a correct message and other a wrong message.
2. Also when they are not in a position to get or transform the channel as ARQ protocols then also these false bicodes may be helpful.
3. When false codes are used in crypto system it can easily mislead the hacker thereby maintaining security.

We will however illustrate this by one simple example.

***Example 2.4.2:*** Let $C = C_1 \cup C_1$ be a false bicode where $C_1$ = C (7.3) code $C_1$ is generated by the generator matrix

$$G_1 = \begin{bmatrix} 1 & 0 & 0 & 1 & 0 & 1 & 0 \\ 0 & 1 & 0 & 0 & 1 & 0 & 1 \\ 0 & 0 & 1 & 1 & 0 & 0 & 1 \end{bmatrix}.$$

The false bicode is given by

C = {(0 0 0 0 0 0 0), (1 0 0 1 0 1 0), (0 1 0 0 1 0 1), (0 0 1 1 0 0 1), (1 1 0 1 1 1 1), (1 0 1 0 0 1 1), (0 1 1 1 1 0 0), (1 1 1 0 1 1 0)} $\cup$ {(0 0 0 0 0 0 0), (1 0 0 1 0 1 0), (0 1 0 0 1 0 1), (0 0 1 1 0 0 1), (1 1 0 1 1 1 1), (1 0 1 0 0 1 1), (0 1 1 1 1 0 0), (1 1 1 0 1 1 0)}.

Suppose some message x = $(x_1\ x_2\ x_3 \ldots x_7) \cup (x_1\ x_2\ x_3 \ldots x_7) \in$ C the false bicode is sent and if y = (1 1 1 1 1 1 1) $\cup$ (1 1 0 1 1 1 1) is the received message, we see (1 1 1 1 1 1 1) $\notin C_1$ only (1 1 0 1 1 1 1) $\in C_1$ so the correct message is taken as (1 1 0 1 1 1 1) or one assumes (1 1 0 1 1 1 1) to be the sent message.

This false bicode save the user from unnecessarily spending money over sending the message or asking them to send the message several times that too when the ARQ protocols are impossible. We give another important use of these false bicodes. If the sender wants to maintain some confidentiality he would send from the false bicode C = $C_1 \cup C_1$ a bicode word x $\cup$ y $\in$ C where the sender knows only one message is relevant and the other is only to misguide the person who would try to know the confidentiality of the message. He would be certainly mislead. Thus the confidentiality can be preserved.

Now we proceed onto to define false tricode and false n-code.

**DEFINITION 2.4.2:** *Let $C = C_1 \cup C_2 \cup C_3$ be a tricode where $C_1 = C_2 = C_3$ then we call the tricode to be a false tricode.*

We illustrate this by the following example.

***Example 2.4.3:*** Let $C = C_1 \cup C_2 \cup C_3$ be a tricode where $C_1 = C_2 = C_3$ = C (7, 4) a linear code i.e. the generator trimatrix of C is

$$G = G_1 \cup G_2 \cup G_3$$

$$= \begin{bmatrix} 1 & 0 & 0 & 0 & 1 & 0 & 1 \\ 0 & 1 & 0 & 0 & 1 & 1 & 0 \\ 0 & 0 & 1 & 0 & 0 & 0 & 1 \\ 0 & 0 & 0 & 1 & 0 & 1 & 1 \end{bmatrix} \cup \begin{bmatrix} 1 & 0 & 0 & 0 & 1 & 0 & 1 \\ 0 & 1 & 0 & 0 & 1 & 1 & 0 \\ 0 & 0 & 1 & 0 & 0 & 0 & 1 \\ 0 & 0 & 0 & 1 & 0 & 1 & 1 \end{bmatrix} \cup$$

$$\begin{bmatrix} 1 & 0 & 0 & 0 & 1 & 0 & 1 \\ 0 & 1 & 0 & 0 & 1 & 1 & 0 \\ 0 & 0 & 1 & 0 & 0 & 0 & 1 \\ 0 & 0 & 0 & 1 & 0 & 1 & 1 \end{bmatrix}.$$

The false tricode

C = {(0 0 0 0 0 0 0), (1 0 0 0 1 0 1), (0 1 0 0 1 1 0), (0 0 1 0 0 0 1), (0 0 0 1 0 1 1), (1 1 0 0 0 1 1), (0 1 1 0 1 1 1), (0 0 1 1 0 1 0), (1 0 1 0 1 0 0), (0 1 0 1 1 1 1), (1 0 0 1 1 1 0), (1 1 1 0 0 1 0), (0 1 1 1 1 0 0), (1 1 0 1 0 0 0), (1 0 1 1 0 1 0), (1 1 1 1 0 0 1)} ∪ {(0 0 0 0 0 0 0), (1 0 0 0 1 0 1), (0 1 0 0 1 1 0), (0 0 1 0 0 0 1), (0 0 0 1 0 1 1), (1 1 0 0 0 1 1), (0 1 1 0 1 1 1), (0 0 1 1 0 1 0), (1 0 1 0 1 0 0), (0 1 0 1 1 1 1), (1 0 0 1 1 1 0), (1 1 1 0 0 1 0), (0 1 1 1 1 0 0), (1 1 0 1 0 0 0), (1 0 1 1 0 1 0), (1 1 1 1 0 0 1)} ∪ { (0 0 0 0 0 0 0), (1 0 0 0 1 0 1), (0 1 0 0 1 1 0), (0 0 1 0 0 0 1), (0 0 0 1 0 1 1), (1 1 0 0 0 1 1), (0 1 1 0 1 1 1), (0 0 1 1 0 1 0), (1 0 1 0 1 0 0), (0 1 0 1 1 1 1 1), (1 0 0 1 1 1 0), (1 1 1 0 0 1 0), (0 1 1 1 1 0 0), (1 1 0 1 0 0 0), (1 0 1 1 0 1 0), (1 1 1 1 0 0 1)}.

Any x, y, z ∈ C is of the form $x = x_1 \cup x_2 \cup x_3$, $y = y_1 \cup y_2 \cup y_3$, $z = z_1 \cup z_2 \cup z_3$, $z_i, x_i, y_i \in C_i$, $1 \leq i \leq 3$.

Now if x is sent and y is received where y = (1 1 1 1 0 0 1) ∪ (1 1 1 0 0 1 0) ∪ (0 1 1 1 1 0 0) = $y_1 \cup y_2 \cup y_3$ and if the messages in $x = x_1 \cup x_2 \cup x_3$ only $x_1$ is aimed to be the code to be sent or used where as $x_2$ and $x_3$ are used to mislead the hacker then while decoding this tricode the receiver will not bother about $y_2$ and $y_3$ he will only decode $y_1$. Thus when one uses such trick the data will be protected to some extent. If we increase the tricode from 3 to an arbitrarily large number say n then even pitching at the truly sent message would be impossible. To this end we define the false n code.

**DEFINITION 2.4.3:** *A code $C = C_1 \cup C_2 \cup \ldots \cup C_n$ where C is an n-code ($n \geq 4$) but in which $C_1 = C_2 = C_3 = \ldots = C_n$ then we call C to be a false n-code.*

We just give an example of a false n-code.

***Example 2.4.4:*** Let $C = C_1 \cup C_2 \cup C_3 \cup C_4 \cup C_5$ where each $C_i$ is a linear code generated by the generator matrix

$$G_i = \begin{bmatrix} 1 & 0 & 0 & 0 & 1 & 1 & 0 & 0 \\ 0 & 1 & 0 & 0 & 0 & 1 & 1 & 0 \\ 0 & 0 & 1 & 0 & 0 & 0 & 1 & 1 \\ 0 & 0 & 0 & 1 & 0 & 1 & 0 & 1 \end{bmatrix}$$

i = 1,2 ..., 5. The 5-code

$C = C_1 \cup C_2 \cup C_3 \cup C_4 \cup C_5$ ($C_1 = C_2 = 3 = C_4 = C_5$)
= {(0 0 0 0 0 0 0 0), (1 0 0 0 1 1 0 0), (0 1 0 0 0 1 1 0), (0 0 1 0 0 0 1 1), (0 0 0 1 0 1 0 1), (1 1 0 0 1 0 1 0), (0 1 1 0 0 1 0 1), (0 0 1 1 0 1 1 0), (1 0 1 0 1 1 1 1), (0 1 0 1 0 0 1 1), (1 0 0 1 1 0 0 1), (1 1 1 0 1 0 0 1), (0 1 1 1 0 0 0 0), (1 1 0 1 1 1 1 1), (1 0 1 1 1 0 1 0), (1 1 1 1 1 1 0 0)} ∪ {(0 0 0 0 0 0 0 0), (1 0 0 0 1 1 0 0), (0 1 0 0 0 1 1 0), (0 0 1 0 0 0 1 1), (0 0 0 1 0 1 0 1), (1 1 0 0 1 0 1 0), (0 1 1 0 0 1 0 1), (0 0 1 1 0 1 1 0), (1 0 1 0 1 0 1 1), (0 1 0 1 0 0 1 1), (1 0 0 1 1 0 0 1), (1 1 1 0 1 0 0 1), (0 1 1 1 0 0 0 0), (1 1 0 1 1 1 1 1), (1 0 1 1 1 0 0 0), (1 1 1 1 1 1 0 0)} ∪ {(0 0 0 0 0 0 0 0), (1 0 0 0 1 1 0 0), (0 1 0 0 0 1 1 0), (0 0 1 0 0 0 1 1), (0 0 0 1 0 1 0 1), (1 1 0 0 1 0 1 0), (0 1 1 0 0 1 0 1),

(0 0 1 1 0 1 1 0), (1 0 1 0 1 0 1 1), (0 1 0 1 0 0 1 1), (1 0 0 1 1 0 0 1), (1 1 1 0 1 0 0 1), (0 1 1 1 0 0 0 0), (1 1 0 1 1 1 1 1), (1 0 1 1 1 0 0 0), (1 1 1 1 1 1 0 0)} ∪ {(0 0 0 0 0 0 0 0), (1 0 0 0 1 1 0 0), (0 1 0 0 0 1 1 0), (0 0 1 0 0 0 1 1), (0 0 0 1 0 1 0 1), (1 1 0 0 1 0 1 0), (0 1 1 0 0 1 0 1), (0 0 1 1 0 1 1 0), (1 0 1 0 1 0 1 1), (0 1 0 1 0 0 1 1), (1 0 0 1 1 0 0 1), (1 1 1 0 1 0 0 1), (0 1 1 1 0 0 0 0), (1 1 0 1 1 1 1 1), (1 0 1 1 1 0 0 0), (1 1 1 1 1 1 0 0)} ∪ {(0 0 0 0 0 0 0 0), (1 0 0 0 1 1 0 0), (0 1 0 0 0 1 1 0), (0 0 1 0 0 0 1 1 ), (0 0 0 1 0 1 0 1), (1 1 0 0 1 0 1 0), (0 1 1 0 0 1 0 1), (0 0 1 1 0 1 1 0), (1 0 1 0 1 0 1 1 ), (0 1 0 1 0 0 1 1), (1 0 0 1 1 0 0 1), (1 1 1 0 1 0 0 1), (0 1 1 1 0 0 0 0), (1 1 0 1 1 1 1 1), (1 0 1 1 1 0 0 0), (1 1 1 1 1 1 0 0)}.

Note if x ∈ C then $x = x_1 \cup x_2 \cup x_3 \cup x_4 \cup x_5$ where $x_i \in C_1$, $1 \leq i \leq 5$; each of the $x_i$'s can be distinct or $x_i$'s can be the same. For instance let x = (0 1 1 1 0 0 0 0) ∪ (1 1 0 1 1 1 1 1) ∪ (1 0 1 1 1 0 0 0) ∪ (1 1 1 1 1 1 0 0) ∪ (1 0 0 0 1 1 0 0) we see each of the elements in x is distinct.

Let y = (0 1 1 1 0 0 0 0) ∪ (0 1 1 1 0 0 0 0) ∪ (0 1 1 1 0 0 0 0) ∪ (0 1 1 1 0 0 0 0), (0 1 1 1 0 0 0 0) ∈ $C_1$ we see each of the elements in y are one and the same. Both x, y ∈ $C = C_1 \cup C_1 \cup C_1 \cup C_1 \cup C_1$.

Now we proceed onto give the uses.

When one wants to use confidentiality in sending the message so that the intruder should not break the message in such cases we can use these false n-codes. That is in the false n-code $x = x_1 \cup x_2 \cup \ldots \cup x_n$ he can enclose the message say in one or more of the code and just sent the code x when he receives the message he need not decode or even verify whether the codes in other coordinates are correct code words but just decode only the codes which carries the message. Thus false codes play a major role in keeping the key intact i.e. the hacker cannot hack it. This can be ones private e-mail or any other thing in which confidentiality is needed.

Yet another use of this false n-code is when the message is sent and if the received message has error then they make an automatic repeat request this can be avoided instead the message after calculating the time period they can so place the length of the code so that we get the same message transmitted. This can be done any desired number of times. Use of this false n code would certainly be economically better when compared to using the method of ARQ protocols where it is at times impossible to capture the message after a time like photographs by moving satellites or missiles.

Now we proceed on to define the notion of 1-pseudo false n-code ($n \geq 3$). For when $n = 2$ we have the 1-pseudo false bicode coincides with the usual bicode which is not at all a false bicode.

**DEFINITION 2.4.4:** *Let $C = C_1 \cup C_2 \cup C_3$ be a tricode in which $C_1 = C_2$ and $C_3$ is different from $C_1$, then we define C to be a 1-pseudo false tricode.*

We illustrate this 1-pseudo false tricode by the following example.

***Example 2.4.5:*** Let us consider the 1-pseudo false tricode $C = C_1 \cup C_2 \cup C_2$ where $C_1 \neq C_2$ generated by the pseudo false trimatrix

$$G = G_1 \cup G_2 \cup G_2 =$$

$$\begin{bmatrix} 1 & 0 & 0 & 1 & 1 & 0 & 0 \\ 0 & 1 & 0 & 0 & 1 & 1 & 0 \\ 0 & 0 & 1 & 0 & 0 & 1 & 1 \end{bmatrix}$$

$$\cup \begin{bmatrix} 1 & 0 & 0 & 1 & 0 & 0 \\ 0 & 1 & 0 & 0 & 1 & 0 \\ 0 & 0 & 1 & 0 & 0 & 1 \end{bmatrix}$$

$$\cup \begin{bmatrix} 1 & 0 & 0 & 1 & 0 & 0 \\ 0 & 1 & 0 & 0 & 1 & 0 \\ 0 & 0 & 1 & 0 & 0 & 1 \end{bmatrix};$$

The 1-pseudo false tricode

C = {(0 0 0 0 0 0 0), (1 0 0 1 1 0 0), (0 1 0 0 1 1 0), (0 0 1 0 0 1 1), (1 1 0 1 0 1 0), (0 1 1 0 1 0 1), (1 0 1 1 1 1 1), (1 1 1 1 0 0 1)} $\cup$ {(0 0 0 0 0 0), (1 0 0 1 0 0), (0 1 0 0 1 0), (0 0 1 0 0 1), (1 1 0 1 1 0), (0 1 1 0 1 1), (1 0 1 1 0 1), (1 1 1 1 1 1)} $\cup$ {(0 0 0 0 0 0), (1 0 0 1 0 0), (0 1 0 0 1 0), (0 0 1 0 0 1), (1 1 0 1 1 0), (0 1 1 0 1 1), (1 0 1 1 0 1), (1 1 1 1 1 1)}.

Any $x = x_1 \cup x_2 \cup x_3$ is such that $x_1 \in C_1$ and $x_2 \in C_2$, $x_3 \in C_3$ i.e. x = (1 0 0 1 1 0 0) $\cup$ (1 0 0 1 0 0) $\cup$ (1 1 1 1 1 1) or x = (1 1 1 1 0 0 1) $\cup$ (0 1 1 0 1 1) $\cup$ (0 1 1 0 1 1).

Now we proceed onto define 1-pseudo false n-code n > 3.

**DEFINITION 2.4.5:** *Let $C = C_1 \cup C_2 \cup \ldots \cup C_i \cup \ldots \cup C_n$ be a union of n number of codes ($n \geq 4$). If $C_1 = C_2 = \ldots = C_{i-1} = C_{i+1} = \ldots = C_n$ and $C_i \neq C_j$ if $i \neq j$, j = 1, 2, 3, ..., i – 1, i + 1, ..., n then we call C to be a 1-pseudo false n-code i.e. in the set of n number of codes $C_1, C_2, \ldots, C_{i-1}, C_i, C_{i+1}, \ldots, C_n$ only $C_i$ is different and all the other n-1 codes are one and the same.*

Now we give an example of a 1-pseudo n-code, n > 3.

***Example 2.4.6:*** Let $C = C_1 \cup C_2 \cup C_1 \cup C_1 \cup C_1 \cup C_1$ be a 1-pseudo false six code generated by the 1-pseudo false 6-matrix

$$G = G_1 \cup G_2 \cup G_1 \cup G_1 \cup G_1 \cup G_1$$

where

$$G_1 = \begin{bmatrix} 1 & 0 & 1 & 1 & 0 & 1 & 1 \\ 0 & 1 & 0 & 1 & 1 & 0 & 1 \end{bmatrix}$$

and

$$G_2 = \begin{bmatrix} 1 & 0 & 0 & 1 & 0 & 0 & 1 \\ 0 & 1 & 0 & 1 & 1 & 1 & 0 \\ 0 & 0 & 1 & 1 & 0 & 0 & 1 \end{bmatrix}.$$

The 1-pseudo false 6-code

C = {(0 0 0 0 0 0 0), (1 0 1 1 0 1 1), (0 1 0 1 1 0 1), (1 1 1 0 1 1 0)} $\cup$ {(0 0 0 0 0 0 0), (1 0 0 1 0 0 1), (0 1 0 1 1 1 0), (0 0 1 1 0 0 1), (1 1 0 0 1 1 1), (0 1 1 0 1 1 1), (1 0 1 0 0 0 0), (1 1 1 1 1 1 0)} $\cup$ {(0 0 0 0 0 0 0), (1 0 1 1 0 1 1), (0 1 0 1 1 0 1), (1 1 1 0 1 1 0)} $\cup$ {(0 0 0 0 0 0 0), (1 0 1 1 0 1 1), (0 1 0 1 1 0 1), (1 1 1 0 1 1 0)} $\cup$ {( 0 0 0 0 0 0 0), (1 0 1 1 0 1 1), (0 1 0 1 1 0 1), (1 1 1 0 1 1 0)} $\cup$ {(0 0 0 0 0 0 0), (1 0 1 1 0 1 1), (0 1 0 1 1 0 1), (1 1 1 0 1 1 0)}.

Any element x in C would be of the form $x = x_1 \cup x_2 \cup x_3 \cup x_4 \cup x_5 \cup x_6$ where $x_2$, $x_1$, $x_3$, $x_4$, $x_5$ and $x_6$ are of length seven but $x_2 \in C_2$ and $x_i \in C_1$, i = 1, 3, 4, 5, 6, i.e. x = (1 0 1 1 0 1 1) $\cup$ (1 1 1 1 1 1 0) $\cup$ (0 1 0 1 1 0 1) $\cup$ (1 0 1 1 0 1 1) $\cup$ (1 1 1 0 1 1 0) $\cup$ (1 0 1 1 0 1 1) or if $y = y_1 \cup y_2 \cup y_3 \cup y_4 \cup y_5 \cup y_6$ = (1 0 1 1 0 1 1) $\cup$ (1 0 1 0 0 0 0) $\cup$ (1 0 1 1 0 1 1) $\cup$ (1 0 1 1 0 1 1) $\cup$ (1 0 1 1 0 1 1 ) $\cup$ (1 0 1 1 0 1 1) i.e., all the $y_i$ are equal to (1 0 1 1 0 1 1) for i = 1, 2, 3, 4, 5, 6 and $y_2$ = (1 0 1 0 0 0 0) $\in C_2$. Note it is interesting to mention when n = 3 we get the 1-pseudo false tricode.

*Note:* It is interesting to note that if $n \geq 2m$ we can have maximum only a m-pseudo false code we can never have m + r –pseudo false n-code where $r \geq 1$. Another interesting factor to make a note of, is if n = 2m + 1 then also we can have maximum only a m-pseudo false (2m + 1)-code. Here also we cannot have a m + r-pseudo false (2m + 1) code $r \geq 1$. This is the very property of t-pseudo false n codes. These class of t-pseudo false n-codes will be very useful when we want to send 2-repeated messages of different lengths we can use it to send 2-different messages of different length when we are not in a position to demand for ARQ protocols or when we want to send secret message which should be kept as very confidential ($n \geq 2t$; $t < n$).

Now we define yet another new class of pseudo false codes.

**DEFINITION 2.4.6:** *Let us consider a n code* $C = C_1 \cup C_2 \cup \ldots \cup C_n$ *where* $C_i = C(n_r, k_r)$; $C_j = C(n_j, k_j)$ $(i \neq j$ *and* $r \neq j)$ *and* $C_1 = C_2 = \ldots = C_{i-1} = C_{i+1} = \ldots = C_{j-1} = C_{j+1} = \ldots = C_n = (n_p, k_p)$ $p \neq j$, $p \neq r$ *and then we call C to be a (1, 1) pseudo false n-code,* $n \geq 4$.

We illustrate this by an example.

***Example 2.4.7:*** Let us consider a 6-code $C = C_1 \cup C_2 \cup C_3 \cup C_4 \cup C_5 \cup C_6$ where $C_2$ is a (6, 4) code generated by $G = G_1 \cup \ldots \cup G_6$ where

$$G_2 = \begin{bmatrix} 1 & 0 & 0 & 0 & 1 & 0 \\ 0 & 1 & 0 & 0 & 0 & 1 \\ 0 & 0 & 1 & 0 & 1 & 1 \\ 0 & 0 & 0 & 1 & 0 & 1 \end{bmatrix}.$$

$C_5$ is a (7, 3) code generated by

$$G_5 = \begin{bmatrix} 1 & 0 & 0 & 1 & 1 & 0 & 0 \\ 0 & 1 & 0 & 0 & 1 & 1 & 0 \\ 0 & 0 & 1 & 0 & 1 & 1 & 1 \end{bmatrix}.$$

$C_1$, $C_3$, $C_4$ and $C_6$ are (5, 2) codes generated by the generator matrix

$$G_1 = \begin{bmatrix} 1 & 0 & 1 & 1 & 0 \\ 0 & 1 & 1 & 0 & 1 \end{bmatrix} = G_3 = G_4 = G_6.$$

Thus the generator 6 matrix G is of the form

$$G = \begin{bmatrix} 1 & 0 & 1 & 1 & 0 \\ 0 & 1 & 1 & 0 & 1 \end{bmatrix} \cup$$

$$\begin{bmatrix} 1 & 0 & 0 & 0 & 1 & 0 \\ 0 & 1 & 0 & 0 & 0 & 1 \\ 0 & 0 & 1 & 0 & 1 & 1 \\ 0 & 0 & 0 & 1 & 0 & 1 \end{bmatrix} \cup \begin{bmatrix} 1 & 0 & 1 & 1 & 0 \\ 0 & 1 & 1 & 0 & 1 \end{bmatrix}$$

$$\cup \begin{bmatrix} 1 & 0 & 1 & 1 & 0 \\ 0 & 1 & 1 & 0 & 1 \end{bmatrix} \cup \begin{bmatrix} 1 & 0 & 0 & 1 & 1 & 0 & 0 \\ 0 & 1 & 0 & 0 & 1 & 1 & 0 \\ 0 & 0 & 1 & 0 & 1 & 1 & 1 \end{bmatrix}$$

$$\cup \begin{bmatrix} 1 & 0 & 1 & 1 & 0 \\ 0 & 1 & 1 & 0 & 1 \end{bmatrix}.$$

C = (5, 2) ∪ (6, 4) ∪ (5, 2) ∪ (5, 2) ∪ (7, 3) ∪ (5, 2) is a (1, 1) pseudo false 6-code given by

C = {(0 0 0 0 0), (1 0 1 1 0), (0 1 1 0 1), (1 1 0 1 1)} ∪ {(0 0 0 0 0 0), (1 0 0 0 1 0), (0 1 0 0 0 1), (0 0 1 0 1 1), (0 0 0 1 0 1), (1 1 0 0 1 1) (0 1 1 0 1 0), (0 0 1 1 1 0), (1 0 1 0 0 1), (0 1 0 1 0 0), (1 0 0 1 1 1), (1 1 1 0 0 0), (0 1 1 1 1 1), (1 1 0 1 1 0), (1 0 1 1 0 0), (1 1 1 1 0 1)} ∪ {(0 0 0 0 0), (1 0 1 1 0), (0 1 1 0 1), (1 1 0 1 1)} ∪ {(0 0 0 0 0), (1 0 1 1 0), (0 1 1 0 1), (1 1 0 1 1)} ∪ {(0 0 0 0 0 0 0), (1 0 0 1 1 0 0), (0 1 0 0 1 1 0), (0 0 1 0 1 1 1), (1 1 0 1 0 1 0), (0 1 1 0 0 0 1), (1 0 1 1 0 1 1), (1 1 1 1 1 0 1)} ∪ {(0 0 0 0 0) (1 0 1 1 0), (0 1 1 0 1), (1 1 0 1 1)}.

The main advantage of this code is it can be used in transmission when one is interested in sending 3 types of messages one a repeating one, where as the other two are just two different messages.

Another use of this code is when one wants to fool the hacker he will not be in any position to known which of the code in C that carries the true message so that it may be impossible for him to hack the information, thus this will be of immense use to a cryptologist.

Now we can like wise define any (t, t)-pseudo n-code when $t = 2$, $n \geq 5$ and (t, t)-pseudo n-code when $n \geq 2t + 1$ and $t < n$. This new class of codes will be very much helpful in several places.

**DEFINITION 2.4.7:** *Let $C = C_1 \cup C_2 \cup \ldots \cup C_n$ be any false n-code, we define C to be a (t, t)-pseudo false n-code ($n \geq 2t + 1$), $t \geq 2$ if $C_1 = C_2 = \ldots = C_t = A$; $C_{t+1} = C_{t+2} = \ldots = C_{2t} = B$ ($A \neq B$) $C_{2t+1} = C_{2t+2} = \ldots = C_n = C$; where $A \neq C$ and $B \neq C$.*

We now proceed on to give an example of it and proceed on to define new classes.

***Example 2.4.8:*** Let $C = C_1 \cup C_2 \cup \ldots \cup C_{16}$ be any false 16 code we define C to be a (5, 5) pseudo false 16 code ($16 \geq 10 + 1$) if

$$C_1 = C_2 = \ldots = C_5 = A = C(7, 3);$$

$$C_6 = C_7 = \ldots = C_{10} = B = C(9, 6); (A \neq B);$$

and

$$C_{11} = C_{12} = \ldots = C_{16} = D = C(7, 4);$$

where $A \neq B$, $D \neq B$ and $D \neq A$.

**DEFINITION 2.4.8:** *Let $C = C_1 \cup C_2 \cup \ldots \cup C_n$ be a pseudo false n-code we call C a (t, m)-pseudo false n-code ($t \neq m$) if $n > t + m + 1$ and $C_1 = C_2 = \ldots = C_t = A$, a $(n_t, k_t)$ code, $C_t = C_{t+1} = \ldots = C_{t+m} = C(n_m, k_m) = P$, $t \neq m$ and $C_{t+m+1} = \ldots = C_n = C(n_i, k_i) = B$, $i \neq t$, $i \neq m$ and $P \neq B$.*

This code will also be very useful to maintain secrecy of the message and it will be very difficult for the intruder to easily break open the key.

Now we illustrate it by an example.

***Example 2.4.9:*** Let $C = C_1 \cup C_2 \cup C_3 \cup \ldots \cup C_{19}$ be a (5, 8)-pseudo false 19-code where $C_1 = C_2 = \ldots = C_5$ is a C(8, 4) code $C_7 = C_8 = C_9 = C_{10} = C_{12} = C_{13} = C_{14} = C_{15}$ is a C(7, 3) code and $C_6 = C_{16} = C_{17} = C_{19}$ is a C(6, 3) code.

Now we show a n-code communication n-channel by figure 2.4.1.

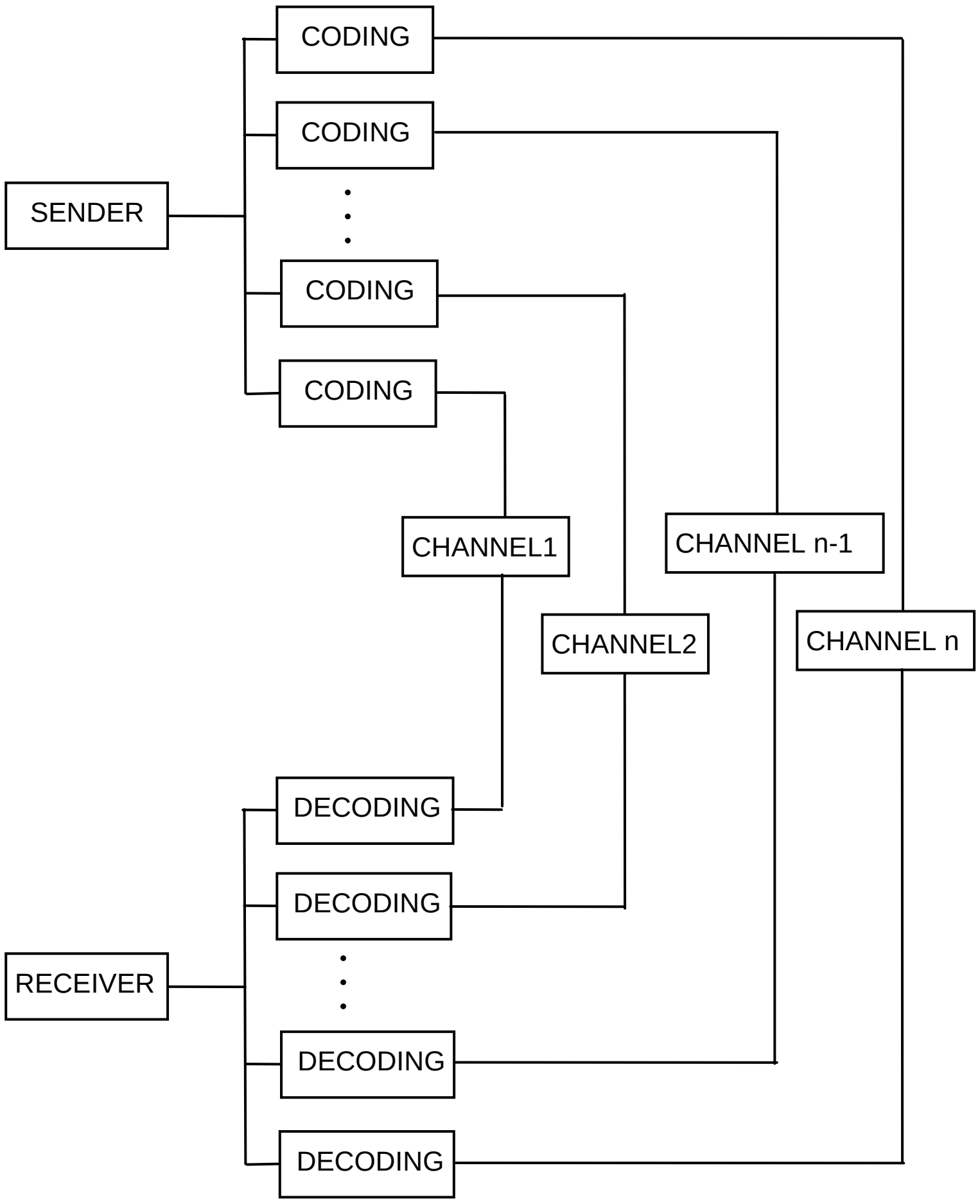

**FIGURE 2.4.1**

**Chapter Three**

# PERIYAR LINEAR CODES

This chapter has two sections. In section one we introduce a new class of codes called Periyar linear codes and in section two the applications of these new classes of codes are given.

## 3.1 Periyar Linear Codes and their Properties

In this section for the first time we introduce the new class of codes called Periyar Linear Codes and enumerate some of its properties. We proceed on to define a Periyar Linear Code. The authors to honour Periyar in his 125[th] birthday have named this code in his name.

**DEFINITION 3.1.1:** *A n-code $C = C_1 \cup C_2 \cup \ldots \cup C_n$ $(n \geq 2)$ is said to be a Periyar linear code (P-linear code) if it has a $(n_i, k_i)$ subcode $C_i \subseteq C$ such that no other subcode in C is a $(n_i, k_i)$ code, i fixed. C may have a subcode with $i \neq j$. Further $C \setminus C_i$ has no $(n_j, k_j)$ subcode (By subcode of C we mean any $C_i$ of C; $1 \leq i \leq n$).*

We first illustrate this by the following example.

***Example 3.1.1:*** Consider the bicode $C = C_1 \cup C_2$ where $C_1$ is a (5, 3) code and $C_2$ is a (7, 4) code. C is a Periyar linear code as C contains a subcode $C_1$ which is a (5, 3) and $C_1$ does not contain any (5, 3) code. Also $C_2$ is a (7, 4) subcode of C. C \ C2 does not contain any (7, 4) code. Thus this bicode is a Periyar linear code.

We cannot always say every bicode is a P-linear code to this we give the following example.

***Example 3.1.2:*** Consider the bicode $C = C_1 \cup C_2$ where $C_1$ and $C_2$ are both (7, 3) codes given by the parity check matrices $H_1$ and $H_2$ where

$$H_1 = \begin{bmatrix} 1 & 0 & 0 & 1 & 0 & 0 & 0 \\ 0 & 1 & 0 & 0 & 1 & 0 & 0 \\ 0 & 0 & 1 & 0 & 0 & 1 & 0 \\ 0 & 0 & 0 & 1 & 0 & 0 & 1 \end{bmatrix}$$

and

$$H_2 = \begin{bmatrix} 1 & 1 & 0 & 1 & 0 & 0 & 0 \\ 0 & 1 & 1 & 0 & 1 & 0 & 0 \\ 1 & 0 & 1 & 0 & 0 & 1 & 0 \\ 1 & 1 & 1 & 0 & 0 & 0 & 1 \end{bmatrix}.$$

$C_1$ = {(0 0 0 0 0 0 0), (1 0 0 1 1 1 1), (0 1 0 0 1 0 0), (0 0 1 0 0 1 0), (1 1 0 1 1 0 0), (1 0 1 1 0 1 0), (0 1 1 0 1 1 0), (1 1 1 1 1 1 0)} and $C_2$ = {(0 0 0 0 0 0 0), (1 0 0 1 0 1 1), (0 1 0 1 1 0 1), (0 0 1 0 1 1 1), (1 1 0 0 1 1 0), (0 1 1 1 0 1 0), (1 0 1 1 1 0 0), (1 1 1 0 0 0 1)}. We see both $C_1$ and $C_2$ are (7, 3) codes and C is a bicode but C is not a P-linear code for we cannot find a subcode other than (7, 3) subcodes.

We have a class of n-codes ($n \geq 2$) which are Periyar linear codes which is clear from the following.

**THEOREM 3.1.1:** *Every n-code $C = C_1 \cup C_2 \cup \ldots \cup C_n$ where $C_i$ is a $(n_i, k_i)$ code, i = 1, 2, …, n and $n_i \neq n_j$ if $i \neq j$ and (or) $k_i \neq k_j$ if $i \neq j$, $i \leq i, j \leq n$ is a P-linear code.*

*Proof:* Clearly we see if $n_i \neq n_j$ for $i \neq j$ then the n-code $C = C_1 \cup C_2 \cup \ldots \cup C_n = C(n_1, k_1) \cup C(n_2, k_2) \cup \ldots \cup C(n_n, k_n)$ is such that each $C_i = C(n_i, k_i)$ is distinct from $C_j = C(n_j, k_j)$. So C has $(n_i, k_i)$ subcodes hence C is a P-linear code. Here it may so happen $k_i = k_j$ even if $i \neq j$. Still $n_i \neq n_j$ for $i \neq j$ is sufficient to produce P-linear codes.

Now we have a n-code $C = C_1 \cup C_2 \cup \ldots \cup C_n = C_1(n_1, k_1) \cup C_2(n_2, k_2) \cup \ldots \cup C_n(n_n, k_n)$ where some $n_i = n_j$, $i \neq j$ but none of the $k_i = k_j$ if $i \neq j$ then also C is a P-linear code.

Clearly if both $n_i \neq n_j$ $(i \neq j)$ and $k_i \neq k_j$ $(i \neq j)$ then also C is a P-linear code.

*Note:* According to the definition of P-linear codes C we should not have two distinct subcodes B and B′ in C such that both B and B′ are $(n_t, k_t)$ linear codes with $B \not\subseteq B'$ and $B' \not\subseteq B$.

We will illustrate all the three cases by the following examples.

***Example 3.1.3:*** Let us consider a bicode $C = C_1 \cup C_2$ where $C_1$ is a (7, 3) linear code and $C_2$ is a (7, 4) linear code. We see both the codes $C_1$ and $C_2$ are of same length 7 but $k_1 \neq k_2$ hence C is a P-linear code.

***Example 3.1.4:*** Consider a bicode $C = C_1 \cup C_2$ where $C_1 = C(9, 5)$ and $C_2 = C(7, 5)$ both $C_1$ and $C_2$ are linear codes having same number of message symbols but have different lengths viz., 9 and 7. Hence C is a P-linear code.

***Example 3.1.5:*** Consider a bicode $C = C_1 \cup C_2$ where $C_1$ is a (7, 3) code and $C_2$ is a (8, 4) code then C is a P-linear code.

Now consider a tricode $C = C_1 \cup C_2 \cup C_3 = C_1(7, 5) \cup C_2(6, 3) \cup C_3(6, 3)$ with $C_2 \neq C_3$ still C is not a P-linear code; it is only a tricode.

Now we proceed on to define a new notion called weak Periyar linear code.

**DEFINITION 3.1.2:** *Let $C = C_1 \cup C_2 \cup \ldots \cup C_n$ be a linear n-code ($n \geq 3$) we say C is a weak Periyar linear code (weak P-linear code) if C has a subcode $C_j$ which is different from other subcodes $C_i$ i.e. C has atleast one subcode different from $C_j$.*

It is important to see that we may have several subcodes of same length and number of message symbols in $C = C_1 \cup C_2 \cup \ldots \cup C_n$ ($n \geq 3$). We still cannot say they are not P-linear codes but only they are weak P-linear codes.

***Example 3.1.6:*** Let $C = C_1 \cup C_2 \cup C_3 \cup C_4$ be a 4-code where $C_1$ = (7, 3) code $C_2$ = (7, 3) code, $C_3$ = (6, 3) code and $C_4$ = (6, 2) where $C_1$ is generated by

$$G_1 = \begin{bmatrix} 1 & 0 & 0 & 1 & 1 & 0 & 0 \\ 0 & 1 & 0 & 1 & 0 & 0 & 1 \\ 0 & 0 & 1 & 0 & 1 & 1 & 0 \end{bmatrix},$$

$C_2$ generated by

$$G_2 = \begin{bmatrix} 1 & 0 & 0 & 1 & 0 & 0 & 1 \\ 0 & 1 & 0 & 0 & 1 & 1 & 1 \\ 0 & 0 & 1 & 0 & 1 & 1 & 0 \end{bmatrix},$$

$C_3$ is generated by

$$G_3 = \begin{bmatrix} 1 & 0 & 0 & 1 & 1 & 0 \\ 0 & 1 & 0 & 1 & 0 & 1 \\ 0 & 0 & 1 & 0 & 1 & 1 \end{bmatrix}$$

and $C_4$ is generated by

$$G_4 = \begin{bmatrix} 1 & 0 & 1 & 1 & 1 & 0 \\ 0 & 1 & 0 & 1 & 1 & 1 \end{bmatrix}.$$

$C_1$ = {(0 0 0 0 0 0 0), (1 0 0 1 1 0 0), (0 1 0 1 0 0 1), (0 0 1 0 1 1 0), (1 1 0 0 1 0 1), (0 1 1 1 1 1 1), (1 0 1 1 0 1 0), (1 1 1 0 0 1 1)},

$C_2$ = {(0 0 0 0 0 0 0), (1 0 0 1 0 0 1), (0 1 0 0 1 1 1), (0 0 1 0 1 1 0), (1 1 0 1 1 1 0), (0 1 1 0 0 0 1), (1 0 1 1 1 1 1), (1 1 1 1 0 0 0)},

$C_3$ = {(0 0 0 0 0 0), (1 0 0 1 1 0), (0 1 0 1 0 1), (0 0 1 0 1 1), (1 1 0 0 1 1), (0 1 1 1 1 0), (1 0 1 1 0 1), (1 1 1 0 0 0)} and

$C_4$ = {(0 0 0 0 0 0), (1 0 1 1 1 0), (0 1 0 1 1 1), (1 1 1 0 0 1)}.

We see if we take $C_4$ or $C_3$ as a subcode of C then C is a P-linear code. But if we take $C_2$ or $C_1$ as a subcode then clearly C = $C_1 \cup C_2 \cup C_3 \cup C_4$ is not a P-linear but only a weak Periyar linear code. Thus from this example we see C is both a P-linear code as well as weak P-linear code. So we are forced to define a new concept because we see the notion of weak P-linear code and P-linear code does not divide C into two classes which are disjoint but we have a over lap.

So we define a new notion called duo Periyar linear code.

**DEFINITION 3.1.3:** *Let us consider a n-code* $C = C_1 \cup C_2 \cup \ldots \cup C_n$, *(*$n \geq 3$*) if C is both a P-linear code as well as a weak P-linear code then we call C to be a duo Periyar linear code (duo P-linear code).*

All P-linear nodes need not be duo P-linear codes or a all weak P-linear codes need not be duo P-linear code or a P-linear code need not be a weak P-linear code or vice versa.

**THEOREM 3.1.2:** *A P-linear code in general need not be a weak P-linear code. Likewise a weak P-linear code need not always be a P-linear code.*

*Proof:* To prove these above statements we give in support some examples.

Let C = C(6, 3) ∪ C(6 3) ∪ C(7, 4) ∪ C(7, 4) be a 4-code where C = $C_1 \cup C_2 \cup C_3 \cup C_4$ with C having the generator 4-matrix G = $G_1 \cup G_2 \cup G_3 \cup G_4$ where

$$G_1 = \begin{bmatrix} 1 & 0 & 0 & 0 & 0 & 1 \\ 0 & 1 & 0 & 0 & 1 & 0 \\ 0 & 0 & 1 & 1 & 0 & 0 \end{bmatrix},$$

$$G_2 = \begin{bmatrix} 1 & 0 & 0 & 1 & 1 & 0 \\ 0 & 1 & 0 & 1 & 0 & 1 \\ 0 & 0 & 1 & 0 & 1 & 1 \end{bmatrix},$$

$$G_3 = \begin{bmatrix} 1 & 0 & 0 & 0 & 1 & 0 & 1 \\ 0 & 1 & 0 & 0 & 0 & 1 & 0 \\ 0 & 0 & 1 & 0 & 0 & 1 & 1 \\ 0 & 0 & 0 & 1 & 1 & 1 & 0 \end{bmatrix}$$

and

$$G_4 = \begin{bmatrix} 1 & 0 & 0 & 0 & 0 & 1 & 0 \\ 0 & 1 & 0 & 0 & 1 & 1 & 1 \\ 0 & 0 & 1 & 0 & 0 & 1 & 1 \\ 0 & 0 & 0 & 1 & 1 & 1 & 0 \end{bmatrix}$$

where the code generated by $G_1$ is {(0 0 0 0 0 0), (1 0 0 0 0 1), (0 1 0 0 1 0), (0 0 1 1 0 0), (1 1 0 0 1 1), (0 1 1 1 1 0), (1 0 1 1 0 1), (1 1 1 1 1 1)} = C(6, 3) = $C_1$.

The code generated by $G_2$ is {(0 0 0 0 0 0), (1 0 0 1 1 0), (0 1 0 1 0 1), (0 0 1 0 1 1), (1 1 0 0 1 1), (0 1 1 1 1 0), (1 0 1 1 0 1), (1 1 1 0 0 0)} = C(6, 3) = $C_2$.

Clearly $C_1 \neq C_2$ The code generated by $G_3$ is {(0 0 0 0 0 0 0), (1 0 0 0 1 0 1), (0 1 0 0 0 1 0), (0 0 1 0 0 1 1), (0 0 0 1 1 1 0),

(1 1 0 0 1 1 1), (1 0 1 0 1 1 0), (0 1 1 0 0 0 1), (0 0 1 1 1 0 1), (0 1 0 1 1 0 0), (1 0 0 1 0 1 1), (1 1 1 0 1 0 0), (0 1 1 1 1 1 1), (1 1 0 1 0 0 0), (1 0 1 1 1 0 1), (1 1 1 1 0 1 1)} = $C_3$ a C(7, 3) code.

Now the code generated by $G_4$ is {(0 0 0 0 0 0 0), (1 0 0 0 0 1 0), (0 1 0 0 1 1 1), (0 0 1 0 0 1 1), (0 0 0 1 1 1 0), (1 1 0 0 1 0 1), (0 1 1 0 1 0 0), (0 0 1 1 1 0 1), (1 0 1 0 0 0 1), (0 1 0 1 0 0 1), (1 0 0 1 1 0 0), (1 1 1 0 1 1 0), (0 1 1 1 0 1 0), (1 1 0 1 0 0 1), (1 0 1 1 0 1 1), (1 1 1 1 0 0 0)} = C(7, 3) = $C_9$ Clearly $C_3 \neq C_4$. All the four codes are distinct, but C is not a P-linear code it is only a weak S linear code. Thus we see a code can be weak P-linear code but not a P-linear code.

***Example 3.1.7:*** Consider the 3 code C = $C_1 \cup C_2 \cup C_3$ = C(7, 3) $\cup$ C(6, 3) $\cup$ C(5, 3) where C(7, 3) is generated by the matrix

$$G_1 = \begin{bmatrix} 1 & 0 & 0 & 1 & 1 & 0 & 0 \\ 0 & 1 & 0 & 0 & 1 & 1 & 0 \\ 0 & 0 & 1 & 0 & 1 & 0 & 1 \end{bmatrix}.$$

The C(6, 3) code is generated by the matrix

$$G_2 = \begin{bmatrix} 1 & 0 & 0 & 0 & 0 & 1 \\ 0 & 1 & 0 & 0 & 1 & 0 \\ 0 & 0 & 1 & 1 & 0 & 0 \end{bmatrix}.$$

The code C(5, 3) is generated by the matrix

$$G_3 = \begin{bmatrix} 1 & 0 & 0 & 1 & 1 \\ 0 & 1 & 0 & 0 & 1 \\ 0 & 0 & 1 & 1 & 0 \end{bmatrix}.$$

The code generated by $G_1$ is {(0 0 0 0 0 0 0), (1 0 0 1 1 0 0), (0 1 0 0 1 1 0), (0 0 1 0 1 0 1), (1 1 0 1 0 1 0), (0 1 1 0 0 1 1), (1 0 1 1 0 0 1), (1 1 1 1 1 1 1)}. The code generated by $G_2$ is {(0 0 0 0 0 0), (1 0 0 0 0 1), (0 1 0 0 1 0), (0 0 1 1 0 0), (1 1 0 0 1 1), (0 1 1 1 1 0), (1 0 1 1 0 1), (1 1 1 1 1 1)}. The code generated by

$G_3$ is {(0 0 0 0 0), (1 0 0 1 1), (0 1 0 0 1), (0 0 1 1 0), (1 1 0 1 0), (0 1 1 1 1), (1 0 1 0 1), (1 1 1 0 0)}. We see this 3-code is a P-linear code and is not a weak P linear code. Thus we see we can have a n-code which is a P-linear code and not a weak P-linear code.

Now we give an example of a n-code (n > 2) which is both a P-linear code as well as a weak P-linear code.

***Example 3.1.8:*** Let $C = C_1 \cup C_2 \cup C_3 \cup C_4$ be a 4-linear code, where $C_1$ = C(7, 4) linear code, $C_2$ = C(5, 2) be another linear code, $C_3$ = C(5, 3) code and $C_4$ be a C(6, 3) linear code. It is easily verified C is a P-linear code for it has only one subcode C(7, 4) or C(5, 2). It is also a weak P-linear code for it has two subcodes of length 5. Thus $C = C_1 \cup C_2 \cup C_3 \cup C_4$ is a duo P-linear code.

Now having seen examples of all the types of codes defined. Now we proceed onto define P-linear repetition code.

**DEFINITION 3.1.4:** *Let $C = C_1 \cup C_2 \cup \ldots \cup C_n$ be a n-code where each $C_i$ is a repetition code of length $n_i$ with $n_i$-1 check symbols where the $n_i$-1 check symbols $x_2^i = x_3^i = \ldots = x_n^i$ are equal to $a_1$ ($a_1$ 'repeated' $n_i - 1$ times) then we may obtain a binary ($n_i$, 1) code with parity-check matrix*

$$H_i = \begin{bmatrix} 1 & 1 & 0 & \cdots & 0 \\ 1 & 0 & 1 & \cdots & 0 \\ \vdots & & & & \vdots \\ 1 & 0 & 0 & \cdots & 1 \end{bmatrix}.$$

*Clearly there are only two code words in the code C($n_i$, 1), namely (0 0 … 0) and (1 1 … 1) this is true for i = 1, 2, …, n. Thus we have n distinct repetition codes, since C is a n-code $C_i \neq C_j$ if $i \neq j$. We see every subcode is distinct so we call C a Periyar linear repetition code (P-linear repetition code).*

***Example 3.1.9:*** Let $C = C_1 \cup C_2 \cup C_3 \cup C_4$ be a 4-code, where $C_1$ = {(0 0 0 0), (1 1 1 1)}, $C_2$ = {(0 0 0 0 0 0), (1 1 1 1 1 1)}, $C_3$ = {(0 0 0 0 0 0 0), (1 1 1 1 1 1 1)} and $C_4$ = {(0 0 0 0 0 0 0 0 0),

(1 1 1 1 1 1 1 1 1)}. Each code $C_i$ is a repetition code which are distinct, this code is associated with the parity check 4-matrix H = $H_1 \cup H_2 \cup H_3 \cup H_4$.

$$H = \begin{bmatrix} 1 & 1 & 0 & 0 \\ 1 & 0 & 1 & 0 \\ 1 & 0 & 0 & 1 \end{bmatrix} \cup$$

$$\begin{bmatrix} 1 & 1 & 0 & 0 & 0 & 0 \\ 1 & 0 & 1 & 0 & 0 & 0 \\ 1 & 0 & 0 & 1 & 0 & 0 \\ 1 & 0 & 0 & 0 & 1 & 0 \\ 1 & 0 & 0 & 0 & 0 & 1 \end{bmatrix} \cup$$

$$\begin{bmatrix} 1 & 1 & 0 & 0 & 0 & 0 & 0 \\ 1 & 0 & 1 & 0 & 0 & 0 & 0 \\ 1 & 0 & 0 & 1 & 0 & 0 & 0 \\ 1 & 0 & 0 & 0 & 1 & 0 & 0 \\ 1 & 0 & 0 & 0 & 0 & 1 & 0 \\ 1 & 0 & 0 & 0 & 0 & 0 & 1 \end{bmatrix} \cup$$

$$\begin{bmatrix} 1 & 1 & 0 & 0 & 0 & 0 & 0 & 0 & 0 \\ 1 & 0 & 1 & 0 & 0 & 0 & 0 & 0 & 0 \\ 1 & 0 & 0 & 1 & 0 & 0 & 0 & 0 & 0 \\ 1 & 0 & 0 & 0 & 1 & 0 & 0 & 0 & 0 \\ 1 & 0 & 0 & 0 & 0 & 1 & 0 & 0 & 0 \\ 1 & 0 & 0 & 0 & 0 & 0 & 1 & 0 & 0 \\ 1 & 0 & 0 & 0 & 0 & 0 & 0 & 1 & 0 \\ 1 & 0 & 0 & 0 & 0 & 0 & 0 & 0 & 1 \end{bmatrix}.$$

Clearly C is a P-linear repetition code, where C = ⟨{(0 0 0 0), (1 1 1 1)} ∪ {(0 0 0 0 0 0), (1 1 1 1 1 1)} ∪ {(0 0 0 0 0 0 0), (1 1 1 1 1 1 1)} ∪ {(1 1 1 1 1 1 1 1 1), (0 0 0 0 0 0 0 0 0)}⟩.

Now we proceed on to define the Periyar linear parity check code.

**DEFINITION 3.1.5:** *Let $C = C_1 \cup C_2 \cup \ldots \cup C_m$ be a m-code. We define C to be a Periyar Parity check code if each $C_i$ is a $(m_i, m_i - 1)$ parity check binary code, i = 1, 2, …, m. Clearly if C is to be a m-code we need each $C_i$ must be distinct i.e. if $C_i = C(m_i, m_i - 1)$ then $m_i \neq m_j$ if $i \neq j$. C is associated with a parity check m-matrix $H = H \cup H_2 \cup \ldots \cup H_n$ where each $H_i$ is a $1 \times m_i$ unit vector i.e. $H_i = (1\ 1\ 1\ 1\ \ldots\ 1)$.*

We give an example of a P-linear parity check code.

***Example 3.1.10:*** Consider the 5 code $C = C_1 \cup C_2 \cup C_3 \cup C_4 \cup C_5$ where each $C_i$ is a parity check $(m_i, m_i - 1)$ binary linear code $H_1 = (1\ 1\ 1\ 1\ 1\ 1)$, $H_2 = (1\ 1\ 1\ 1\ 1)$, $H_3 = (1\ 1\ 1\ 1\ 1\ 1\ 1)$, $H_4 = (1\ 1\ 1\ 1\ 1\ 1\ 1\ 1\ 1)$ and $H_5 = (1\ 1\ 1\ 1)$, where $H = H_1 \cup H_2 \cup H_3 \cup H_4 \cup H_5$ is a parity check 5-matrix. Since each code $C_i$ is distinct for i = 1, 2, …, 5 we see C is a P-linear parity check code and not a weak P-linear parity check code.

**THEOREM 3.1.3:** *A linear binary repetition n-code ($n \geq 2$) is always a P-linear binary repetition code and never a weak P-linear binary repetition code and hence never a duo P-linear binary code.*

*Proof:* Given $C = C_1 \cup C_2 \cup \ldots \cup C_n$ is a linear binary repetition n-code ($n \geq 2$). Since C is a repetition n-code each code $C_i$ is distinct and is a repetition code for i = 1, 2, …, n. So C cannot have two subcodes of same length.

So C is not a weak P-linear code. Further C has n-distinct subcodes i.e. C has a subcode which has a length different from other subcodes, so C is a P-linear code. Since C is only a P-linear code and not a weak P-linear code so C cannot be a duo P-linear code. Hence the claim.

Next we prove the following interesting theorem about Parity Check n-code.

**THEOREM 3.1.4:** *Let $C = C_1 \cup C_2 \cup \ldots \cup C_m$ be linear binary m-parity check code. Then C is only a P-linear code and not a weak P-linear code; thus C is not a duo P-linear code.*

*Proof:* Given $C = C_1 \cup C_2 \cup \ldots \cup C_m$ is a linear binary m-parity check code; i.e. each $C_i$ is a parity check code of length $m_i$ and $m_i \neq m_j$ if $i \neq j$, $1 \leq i, j \leq m$. Thus each $C_i$ is distinct. Hence C has subcode which is unique i.e. C has subcode of length $m_j$ then C has no other subcode of length $m_j$. So C is a P-linear code. Now C has no two subcodes which are of same length, hence C is not a weak P-linear code.

Thus we have proved in case of every n-code if it is repetition n-code or a parity check n-code then it is only a P-linear code and never a weak P-linear code thus not a duo P-linear code.

Hence we have seen that all n-codes are not duo n-codes. We have a special class of n-codes which are also P-linear codes and not weak P-linear n-codes. We have already given examples of P-linear repetition code and P-linear parity check code.

We now proceed on to define the notion of binary Hamming P code.

**DEFINITION 3.1.6:** *Let $C = C_1 \cup C_2 \cup \ldots C_n$ be a n-code ($n \geq 2$). If each $C_i$ is a binary Hamming code of length $n_i = 2^{m_i} - 1$ ($m_i \geq 2$); $i = 1, 2, \ldots, n$. Thus any subcode of C is of length $n_i$ and clearly C has no subcode of length $n_i$; so C is a defined to be a Periyar linear Hamming code (P-linear Hamming code).*

***Example 3.1.11:*** Consider the bicode $C = C_1 \cup C_2$ where $C_1$ is a $n_1 = 2^3$-1 binary Hamming code and $C_2$ is a $n_2 = 2^4$-1 binary Hamming code. Then C is a P-linear Hamming code.

We expect the reader to prove the following theorem.

**THEOREM 3.1.5:** *Let $C = C_1 \cup C_2 \cup \ldots \cup C_n$ be a n-code ($n \geq 2$) such that each $C_i$ is a Hamming binary code of length $n_i = 2^i - 1$, $i = 1, 2, \ldots, n$. Then C is a P-linear Hamming code and C is never a weak P-linear code and hence is not a duo P-linear code.*

**DEFINITION 3.1.7:** *Let $C = C_1 \cup C_1^{\perp}$ be a bicode where $C_1$ is a (n, k) code and $C^{\perp}$ the dual code of C. We call C whole bicode. We see in case of whole bicode both of them are of same length but the number of message symbols in them is assumed to be different.*

We illustrate this by the following example.

***Example 3.1.12:*** Let $C = C_1 \cup C_1^{\perp}$ be a bicode where $C_1$ is generated by the generator matrix

$$G_1 = \begin{bmatrix} 1 & 0 & 0 & 0 & 1 & 1 & 0 \\ 0 & 1 & 0 & 1 & 0 & 1 & 0 \\ 0 & 0 & 1 & 1 & 1 & 0 & 1 \end{bmatrix}.$$

The generator matrix of $C_1^{\perp}$ is given by

$$G_2 = \begin{bmatrix} 0 & 1 & 1 & 1 & 0 & 0 & 0 \\ 1 & 0 & 1 & 0 & 1 & 0 & 0 \\ 1 & 1 & 0 & 0 & 0 & 1 & 0 \\ 0 & 0 & 1 & 0 & 0 & 0 & 1 \end{bmatrix}.$$

The code generated by $G_1$ is $C_1$ = {(0 0 0 0 0 0 0), (1 0 0 0 1 1 0), (0 1 0 1 0 1 0), (0 0 1 1 1 0 1), (1 1 0 1 1 0 0), (0 1 1 0 1 1 1), (1 0 1 1 0 1 1), (1 1 1 0 0 0 1)}. The code generated by $C_2$ i.e., $C_1^{\perp}$ = {(0 0 0 0 0 0 0), (0 1 1 1 0 0 0), (1 0 1 0 1 0 0), (1 1 0 0 0 1 0), (0 0 1 0 0 0 1), (1 1 0 1 1 0 0), (0 1 1 0 1 1 0), (1 1 1 0 0 1

1), (1 0 1 1 0 1 0), (1 0 0 0 1 0 1), (0 1 0 1 0 0 1), (0 0 0 1 1 1 0), (0 1 0 0 1 1 1), (1 0 0 1 0 1 1), (1 1 1 1 1 0 1), (0 0 1 1 1 1 1)}.

The whole bicode C = $C_1 \cup C_1^{\perp}$ where all the elements of $C_1$ and $C_1^{\perp}$ are given.

Now it is a natural question to ask can we define whole tricode the answer is no. But we can have a 4-whole code and in general a whole n-code if and only if n is an even number.

**DEFINITION 3.1.8:** *Let $C = C_1 \cup C_1^{\perp} \cup C_2 \cup C_2^{\perp} \ldots \cup C_n \cup C_n^{\perp}$ be a 2n-code where $C_i$ and its dual $C_i^{\perp}$ find their place in C then we call C a whole 2n-code ($n \geq 1$), i = 1, 2, …, n.*

**THEOREM 3.1.6:** *Let $C = C_1 \cup C_1^{\perp} \cup C_2 \cup C_2^{\perp} \cup \ldots \cup C_n \cup C_n^{\perp}$ be a whole 2n-code. C is not a P-linear code but C is always a weak P-linear code. Further C is not a duo P-linear code.*

*Proof:* We are given $C = C_1 \cup C_1^{\perp} \cup C_2 \cup C_2^{\perp} \cup \ldots \cup C_n \cup C_n^{\perp}$ to be a whole 2n-code. Thus C has atleast 2 subcodes of same length since if C has the subcode $C_i$ then $C_i^{\perp}$ is also a subcode as C is a whole 2n-code. Hence C is a weak P-linear code and never a P-linear code, thus the whole 2n code can never be a duo P-linear code. Hence the claim.

We have seen a class of n-codes which are P-linear codes and never weak P-linear codes. Now we have established that whole 2n-codes are weak P-linear codes and never a P-linear code.

Now we proceed onto give yet another class of duo P-linear codes to this end we define the notion of pseudo whole n-code ($n \geq 3$).

**DEFINITION 3.1.9:** *Let $C = C_1 \cup C_2 \cup \ldots \cup C_n$ (n an odd number) we call C to be a pseudo whole n-code if $C = C_1 \cup C_1^{\perp} \cup C_2 \cup C_2^{\perp} \cup \ldots \cup C_{\frac{n-1}{2}} \cup C_{\frac{n-1}{2}}^{\perp} \cup C_n$ as n is given to be an* *odd number where C = a union of the whole code plus another arbitrary code of length different from all other codes in C.*

We illustrate this by the following example.

***Example 3.1.13:*** Let $C = C_1 \cup C_1^{\perp} \cup C_2$ be a tricode where C is the union of code and its dual and another code $C_2$ where C is a pseudo whole code.

$C_1$ is generated by

$$G_1 = \begin{bmatrix} 1 & 0 & 0 & 1 & 1 & 1 & 1 \\ 0 & 1 & 0 & 0 & 0 & 1 & 0 \\ 0 & 0 & 1 & 0 & 1 & 1 & 0 \end{bmatrix}$$

and $C_1^{\perp}$ is generated by

$$G_2 = \begin{bmatrix} 1 & 0 & 0 & 1 & 0 & 0 & 0 \\ 1 & 0 & 1 & 0 & 1 & 0 & 0 \\ 1 & 1 & 1 & 0 & 0 & 1 & 0 \\ 1 & 0 & 0 & 0 & 0 & 0 & 1 \end{bmatrix}.$$

$C_2$ is a code generated by

$$G_3 = \begin{bmatrix} 1 & 0 & 0 & 1 & 0 & 0 \\ 0 & 1 & 0 & 0 & 1 & 0 \\ 0 & 0 & 1 & 0 & 0 & 1 \end{bmatrix}$$

C is a pseudo whole code.

**THEOREM 3.1.7:** *Let $C = C_1 \cup C_1^{\perp} \cup C_2 \cup C_2^{\perp} \cup \ldots \cup C_{\frac{n-1}{2}} \cup C_{\frac{n-1}{2}}^{\perp} \cup C_n$ (n-odd) be the pseudo whole n-code; then C is a duo P-linear whole code.*

*Proof:* Given $C = C_1 \cup C_1^{\perp} \cup C_2 \cup C_2^{\perp} \cup \ldots \cup C_{\frac{n-1}{2}} \cup C_{\frac{n-1}{2}}^{\perp} \cup C_n$ be the pseudo whole n-code (n-odd). C has a subcode $C_n$ which has no other subcode the same length. Hence C is a P-linear code, C has subcodes $C_i$, then C has same length subcode given by $C_i^{\perp}$ $1 \leq i \leq \frac{n-1}{2}$ . Thus C is a weak P-linear code. C is both a P-linear code as well as C is a weak- P-linear code, so C is a duo P-linear code which will be known as the duo P-linear pseudo whole code. Thus we have the class of pseudo whole codes which happen to be a duo P-linear code.

Now having seen the 3 distinct classes of codes we now proceed on to define the notion of Periyar cyclic code.

**DEFINITION 3.1.10:** *Let $C = C_1 \cup C_2 \cup \ldots \cup C_n$ be a n-code ($n \geq 2$) if each $C_i$ is a cyclic code then we define C to be a n-cyclic code. If each $C_i$ is of length $n_i$ and $n_i \neq n_j$ if $i \neq j$ for $1 \leq i, j \leq n$ then we define C to be a Periyar cyclic code. If C has more than one subcyclic code of length $n_i$, then we define C to be a weak Periyar cyclic code.*

First we will provide some examples of them so that it would make the reader understand the definition better.

***Example 3.1.14:*** Consider the cyclic 5-code $C = C_1 \cup C_2 \cup C_3 \cup C_4 \cup C_5$ where $C_1$ is a (6, 3) code $C_2$ a (7, 4) code, $C_3$ a (7, 3) code, $C_4$ a (8, 4) code and $C_5$ are (5, 4) cyclic codes. Let G be the associated generator 5-matrix of C.

$G = G_1 \cup G_2 \cup G_3 \cup G_4 \cup G_5$

$$= \begin{bmatrix} 1 & 0 & 0 & 1 & 0 & 0 \\ 0 & 1 & 0 & 0 & 1 & 0 \\ 0 & 0 & 1 & 0 & 0 & 1 \end{bmatrix} \cup$$

$$\begin{bmatrix} 1 & 1 & 0 & 1 & 0 & 0 & 0 \\ 0 & 1 & 1 & 0 & 1 & 0 & 0 \\ 0 & 0 & 1 & 1 & 0 & 1 & 0 \\ 0 & 0 & 0 & 1 & 1 & 0 & 1 \end{bmatrix} \cup$$

$$\begin{bmatrix} 1 & 1 & 1 & 0 & 1 & 0 & 0 \\ 0 & 1 & 1 & 1 & 0 & 1 & 0 \\ 0 & 0 & 1 & 1 & 1 & 0 & 1 \end{bmatrix} \cup$$

$$\begin{bmatrix} 1 & 0 & 0 & 0 & 1 & 0 & 0 & 0 \\ 0 & 1 & 0 & 0 & 0 & 1 & 0 & 0 \\ 0 & 0 & 1 & 0 & 0 & 0 & 1 & 0 \\ 0 & 0 & 0 & 1 & 0 & 0 & 0 & 1 \end{bmatrix} \cup$$

$$\begin{bmatrix} 1 & 1 & 0 & 0 & 0 \\ 0 & 1 & 1 & 0 & 0 \\ 0 & 0 & 1 & 1 & 0 \\ 0 & 0 & 0 & 1 & 1 \end{bmatrix}.$$

All the 5 codes $C_1$, $C_2$, $C_3$, $C_4$ and $C_5$ are cyclic. Thus C is a cyclic 5-code.

We see C is both a P-cyclic code as well as weak P-cyclic code. C is given by the following 5-codes

C = {(0 0 0 0 0 0), (1 0 0 1 0 0) (0 1 0 0 1 0), (0 0 1 0 0 1), (1 1 0 1 1 0), (0 1 1 0 1 1), (1 0 1 1 0 1), (1 1 1 1 1 1)} ∪ {(0 0 0 0 0 0 0), (1 1 0 1 0 0 0), (0 1 1 0 1 0 0), (0 0 1 1 0 1 0), (0 0 0 1 1 0 1), (1 0 1 1 1 0 0), (0 1 0 1 1 1 0), (0 0 1 0 1 1 1), (1 1 1 0 0 1 0), (0 1 1 1 0 0 1), (1 1 0 0 1 0 1), (1 0 0 0 1 1 0), (0 1 0 0 0 1 1), (1 0 1 0 0 0 1), (1 1 1 1 1 1 1), (1 0 0 1 0 1 1)} ∪ {(0 0 0 0 0 0 0), (1 1 1 0 1 0 0), (0 1 1 1 0 11 0), (0 0 1 1 1 0 1), (1 0 0 1 1 1 0),

(0 1 0 0 1 1 1), (1 1 0 1 0 0 1), (1 0 1 0 0 1 1)} $\cup$ {(0 0 0 0 0 0 0 0), (1 0 0 0 1 0 0 0), (0 1 0 0 0 1 0 0), (0 0 1 0 0 0 1 0), (0 0 0 1 0 0 0 1), (1 1 0 0 1 1 0 0), (0 1 1 0 0 1 1 0), (0 0 1 1 0 0 1 1), (1 0 1 0 1 0 1 0), (0 1 0 1 0 1 0 1), (1 0 0 1 1 0 0 1), (1 1 1 0 1 1 1 0), (0 1 1 1 0 1 1 1), (1 1 0 1 1 1 0 1), (1 0 1 1 1 0 1 1), (1 1 1 1 1 1 1 1)} $\cup$ {(0 0 0 0 0), (1 1 0 0 0), (0 1 1 0 0), (0 0 1 1 0), (0 0 0 1 1), (1 0 1 0 0), (0 1 0 1 0), (0 0 1 0 1), (1 1 1 1 0), (0 1 1 1 1), (1 1 0 1 1), (1 0 0 1 0), (0 1 0 0 1), (1 0 1 1 1), (1 1 1 0 1), (1 0 0 0 1)}

is a cyclic 5-code. For this has $C_5$ alone to be subcode likewise. $C_4$ is a subcode and has no other subcode of length 8; whereas C has $C_2$ and $C_3$ to be subcodes of length 7. Thus C is both a P-linear cyclic code as well as weak P-linear cyclic code.

***Example 3.1.15:*** Let $C = C_1 \cup C_2 \cup C_2 \cup C_4 \cup C_5$ where $C_1$ is a C(6, 3) code, $C_2$ is a C(6, 1) code, $C_3$ a (7, 3) code, $C_4$ a (7, 4) code and $C_5$ a (7, 3) code. Clearly C is not a P-linear cyclic code but only a weak P-linear cyclic code where C is generated by the generator 5-matrix given by

$$G = G_1 \cup G_2 \cup G_3 \cup G_4 \cup G_5 =$$

$$\begin{bmatrix} 1 & 0 & 0 & 1 & 0 & 0 \\ 0 & 1 & 0 & 0 & 1 & 0 \\ 0 & 0 & 1 & 0 & 0 & 1 \end{bmatrix} \cup$$

$$\begin{bmatrix} 1 & 1 & 0 & 0 & 0 & 0 \\ 0 & 1 & 1 & 0 & 0 & 0 \\ 0 & 0 & 1 & 1 & 0 & 0 \\ 0 & 0 & 0 & 1 & 1 & 0 \\ 0 & 0 & 0 & 0 & 1 & 1 \end{bmatrix} \cup$$

$$\begin{bmatrix} 1 & 0 & 1 & 1 & 1 & 0 & 0 \\ 0 & 1 & 0 & 1 & 1 & 1 & 0 \\ 0 & 0 & 1 & 0 & 1 & 1 & 1 \end{bmatrix} \cup$$

$$\begin{bmatrix} 1 & 0 & 1 & 1 & 0 & 0 & 0 \\ 0 & 1 & 0 & 1 & 1 & 0 & 0 \\ 0 & 0 & 1 & 0 & 1 & 1 & 0 \\ 0 & 0 & 0 & 1 & 0 & 1 & 1 \end{bmatrix} \cup$$

$$\begin{bmatrix} 1 & 1 & 1 & 0 & 1 & 0 & 0 \\ 0 & 1 & 1 & 1 & 0 & 1 & 0 \\ 0 & 0 & 1 & 1 & 1 & 0 & 1 \end{bmatrix}.$$

Clearly G generates a 5-cyclic code and G is not a P-linear cyclic code but only a weak P-cyclic code. For it has subcodes of length 6 and 7 respectively.

Now we have seen an example of a weak P-linear cyclic code which is not a P-cyclic code.

We proceed on to give an example of a P-cyclic code which is not a weak P-linear cyclic code.

***Example 3.1.16:*** Consider the cyclic 4-code $C = C_1 \cup C_2 \cup C_3 \cup C_4$ generated by the generator 4-matrix

$$G = G_1 \cup G_2 \cup G_3 \cup G_4 =$$

$$\begin{bmatrix} 1 & 1 & 1 & 0 & 1 & 0 & 0 \\ 0 & 1 & 1 & 1 & 0 & 1 & 0 \\ 0 & 0 & 1 & 1 & 1 & 0 & 1 \end{bmatrix} \cup$$

$$\begin{bmatrix} 1 & 1 & 0 & 0 & 0 & 0 \\ 0 & 1 & 1 & 0 & 0 & 0 \\ 0 & 0 & 1 & 1 & 0 & 0 \\ 0 & 0 & 0 & 1 & 1 & 0 \\ 0 & 0 & 0 & 0 & 1 & 1 \end{bmatrix} \cup$$

$$\begin{bmatrix} 1 & 1 & 0 & 0 & 0 \\ 0 & 1 & 1 & 0 & 0 \\ 0 & 0 & 1 & 1 & 0 \\ 0 & 0 & 0 & 1 & 1 \end{bmatrix} \cup$$

$$\begin{bmatrix} 1 & 0 & 0 & 0 & 1 & 0 & 0 & 0 \\ 0 & 1 & 0 & 0 & 0 & 1 & 0 & 0 \\ 0 & 0 & 1 & 0 & 0 & 0 & 1 & 0 \\ 0 & 0 & 0 & 1 & 0 & 0 & 0 & 1 \end{bmatrix}.$$

We see this cyclic 4-code has 4 subcodes each is of length different from the other i.e. C has a subcode of length 7, length 6, length 5 and length 8. So C is a P-linear code and not a weak P-linear code. Hence the claim.

Now we proceed on to give an example of a cyclic n-code which is a weak P-cyclic code and not a P-cyclic code.

***Example 3.1.17:*** Consider a cyclic 4 code C given by $C = C_1 \cup C_2 \cup C_3 \cup C_4$ generated by

$$G = G_1 \cup G_2 \cup G_3 \cup G_4$$

where

$$G_1 = \begin{bmatrix} 1 & 1 & 1 & 0 & 1 & 0 & 0 \\ 0 & 1 & 1 & 1 & 0 & 1 & 0 \\ 0 & 0 & 1 & 1 & 1 & 0 & 1 \end{bmatrix},$$

$$G_2 = \begin{bmatrix} 1 & 1 & 0 & 1 & 0 & 0 & 0 \\ 0 & 1 & 1 & 0 & 1 & 0 & 0 \\ 0 & 0 & 1 & 1 & 0 & 1 & 0 \\ 0 & 0 & 0 & 1 & 1 & 0 & 1 \end{bmatrix},$$

$$G_3 = \begin{bmatrix} 1 & 0 & 0 & 1 & 0 & 0 \\ 0 & 1 & 0 & 0 & 1 & 0 \\ 0 & 0 & 1 & 0 & 0 & 1 \end{bmatrix}$$

and

$$G_4 = \begin{bmatrix} 1 & 1 & 0 & 0 & 0 & 0 \\ 0 & 1 & 1 & 0 & 0 & 0 \\ 0 & 0 & 1 & 1 & 0 & 0 \\ 0 & 0 & 0 & 1 & 1 & 0 \\ 0 & 0 & 0 & 0 & 1 & 1 \end{bmatrix}.$$

We see the codes $C_1$ and $C_2$ are of same length 7 but have different number of message symbols viz. 3 and 4 respectively. C is a C(7, 3) linear code and $C_2$ is a (7, 4) linear code. Now $C_3$ and $C_4$ are linear codes of length 6, $C_3$ is a (6, 3) linear code and $C_4$ is a (6, 5) linear code. Thus C is a cyclic 4-code but not a P-linear cyclic code only a weak P-linear cyclic code as C has two subcodes of length 7 and two subcodes of length 6. It has no unique subcode.

Now we proceed onto define the notion of pseudo Periyar cyclic n-code.

**DEFINITION 3.1.11:** *Let $C = C_1 \cup C_2 \cup \ldots \cup C_n$ be a n-code. We call C to be a pseudo Periyar cyclic code if C has atleast one cyclic code i.e. atleast one of the $C_i$'s is a cyclic code ($n \geq 2$).*

We give example of a pseudo P-cyclic code.

***Example 3.1.18:*** Consider the 4-code $C = C_1 \cup C_2 \cup C_3 \cup C_4$ where C is generated by the generator matrix

$$G = G_1 \cup G_2 \cup G_3 \cup G_4$$

$$= \begin{bmatrix} 1 & 0 & 0 & 1 & 0 & 0 \\ 0 & 1 & 0 & 0 & 1 & 1 \\ 0 & 0 & 1 & 0 & 1 & 1 \end{bmatrix} \cup$$

$$\begin{bmatrix} 1 & 0 & 0 & 1 & 0 & 0 \\ 0 & 1 & 0 & 0 & 1 & 0 \\ 0 & 0 & 1 & 0 & 0 & 1 \end{bmatrix} \cup$$

$$\begin{bmatrix} 1 & 0 & 0 & 1 & 0 & 0 & 1 \\ 0 & 1 & 0 & 0 & 1 & 1 & 0 \\ 0 & 0 & 1 & 0 & 0 & 1 & 1 \end{bmatrix} \cup$$

$$\begin{bmatrix} 1 & 0 & 0 & 0 & 1 & 1 & 0 \\ 0 & 1 & 0 & 0 & 1 & 0 & 1 \\ 0 & 0 & 1 & 0 & 0 & 1 & 1 \\ 0 & 0 & 0 & 1 & 1 & 1 & 1 \end{bmatrix}.$$

The 4-code generated by G is given by C = {(0 0 0 0 0 0), (1 0 0 1 0 0), (0 1 0 0 1 1), (0 0 1 0 1 1), (1 1 0 1 1 1), (0 1 1 0 0 0), (1 0 1 1 1 1), (1 1 1 1 0 0)} ∪ {(0 0 0 0 0 0), (1 0 0 1 0 0), (0 1 0 0 1 0), (0 0 1 0 0 1), (1 1 0 1 1 0), (0 1 1 0 1 1), (1 0 1 1 0 1), (1 1 1 1 1 1)} ∪ {(0 0 0 0 0 0 0), (1 0 0 1 0 0 1), (0 1 0 0 1 1 0), (0 0 1 0 0 1 1), (1 1 0 1 1 1 1), (0 1 1 0 1 0 1), (1 0 1 1 0 1 0), (1 1 1 1 1 0 0)} ∪ {(0 0 0 0 0 0 0), (1 0 0 0 1 1 0), (0 1 0 0 1 0 1), (0 0 1 0 0 1 1), (0 0 0 1 1 1 1), (1 1 0 0 0 1 1), (0 1 1 0 1 1 0), (0 0 1 1 1 0 0), (1 0 1 0 1 0 1), (0 1 0 1 0 1 0), (1 0 0 1 0 0 1), (1 1 1 0 0 0 0), (0 1 1 1 0 0 1), (1 1 0 1 1 0 0), (1 0 1 1 0 1 0), (1 1 1 1 1 1 1)}. We see among the four codes only $C_2$ is a cyclic code so C is a pseudo P-cyclic code.

Now we proceed on to define the notion of pseudo strong Periyar cyclic code.

**DEFINITION 3.1.12:** *Let $C = C_1 \cup C_2 \cup \ldots \cup C_n$ be n-code if in C we have only n-1 of the codes among $C_1, C_2, \ldots, C_n$ to be cyclic then we define C to be a pseudo strong Periyar cyclic code.*

Now we proceed on to illustrate this by the following example.

***Example 3.1.19:*** Let $C = C_1 \cup C_2 \cup C_3 \cup C_4 \cup C_5 \cup C_6$ be a 6-code where C is generated by the 6-matrix

$$G = G_1 \cup G_2 \cup G_3 \cup G_4 \cup G_5 \cup G_6 =$$

$$= \begin{bmatrix} 1 & 1 & 0 & 1 & 0 & 0 & 0 \\ 0 & 1 & 1 & 0 & 1 & 0 & 0 \\ 0 & 0 & 1 & 1 & 0 & 1 & 0 \\ 0 & 0 & 0 & 1 & 1 & 0 & 1 \end{bmatrix} \cup$$

$$\begin{bmatrix} 1 & 1 & 1 & 0 & 1 & 0 & 0 \\ 0 & 1 & 1 & 1 & 0 & 1 & 0 \\ 0 & 0 & 1 & 1 & 1 & 0 & 1 \end{bmatrix} \cup$$

$$\begin{bmatrix} 1 & 0 & 0 & 1 & 0 & 0 \\ 0 & 1 & 0 & 0 & 1 & 0 \\ 0 & 0 & 1 & 0 & 0 & 1 \end{bmatrix} \cup$$

$$\begin{bmatrix} 1 & 1 & 0 & 0 & 0 & 0 \\ 0 & 1 & 1 & 0 & 0 & 0 \\ 0 & 0 & 1 & 1 & 0 & 0 \\ 0 & 0 & 0 & 1 & 1 & 0 \\ 0 & 0 & 0 & 0 & 1 & 1 \end{bmatrix} \cup$$

$$\begin{bmatrix} 1 & 0 & 0 & 0 & 1 & 0 & 0 & 0 \\ 0 & 1 & 0 & 0 & 0 & 1 & 0 & 0 \\ 0 & 0 & 1 & 0 & 0 & 0 & 1 & 0 \\ 0 & 0 & 0 & 1 & 0 & 0 & 0 & 1 \end{bmatrix} \cup$$

$$\begin{bmatrix} 1 & 0 & 0 & 0 & 1 & 1 & 1 & 0 \\ 0 & 1 & 0 & 0 & 0 & 0 & 0 & 1 \\ 0 & 0 & 1 & 0 & 1 & 1 & 0 & 0 \\ 0 & 0 & 0 & 1 & 0 & 0 & 1 & 1 \end{bmatrix}.$$

Clearly the codes $C_1$, $C_2$, …, $C_5$ are cyclic only $C_6$ is a non cyclic code so C is a pseudo strong P-cyclic code. We see in this example $C_1$ is a (7, 4) cyclic code generated by the polynomial $x^3 + x + 1$ and $x^3 + x + 1 / x^7 + 1$. $C_2$ is a (7, 3) cyclic code generated by the polynomial $x^4 + x^2 + x + 1$ and $x^4 + x^2 + x + 1 / x^7 + 1$. $C_3$ is a (6, 3) cyclic code generated by the polynomial $x^3 + 1$ and $x^3 + 1 / x^6 + 1$. $C_4$ is a (6, 5) cyclic code generated by the polynomial x+1, $C_5$ is a (8, 4) cyclic code generated by a generator polynomial $x^4+1$ and $x^4+1 / x^8+1$ and $C_6$ is not a cyclic code. So C is a pseudo strong P-cyclic 6-code.

We see all pseudo strong P-cyclic codes are trivially pseudo P-cyclic codes. But in general a pseudo P-cyclic code need not be a strong pseudo P-cyclic code.

Now we would define the notion of weak Periyar false codes.

**DEFINITION 3.1.13:** *$C = C_1 \cup C_2 \cup \ldots \cup C_n$ is a false n-code if atleast one of the subcodes is different from the other codes then we define C to be a weak Periyar false code.*

Clearly a false code is never a weak P-false code.

**THEOREM 3.1.8:** *The class of m-pseudo false n-codes ($n \geq 2m$, $2 \leq m < n$) are weak P-false codes.*

*Proof:* From the very definition of m-pseudo false n-codes $C = C_1 \cup C_2 \cup \ldots \cup C_n$ we see $C_1 = C_2 = \ldots = C_m = A$ and $C_{m+1} = C_{m+2} = \ldots = C_n = B$; $A \neq B$. Thus C has atleast one subcode is different from other codes.

Hence C is a weak P-false code.

**DEFINITION 3.1.14:** *A false n-code $C = C_1 \cup C_2 \cup \ldots \cup C_n$ ($n \geq 2$) is said to be Periyar false code if and only if C is a 1-pseudo false code i.e. C has one and only one subcode $C_i$ different from other subcodes of C.*

**THEOREM 3.1.9:** *Every 1-pseudo false n-code is a P-false code.*

*Proof:* Follows from the very definition of P-false codes.

We illustrate by an example.

***Example 3.1.20:*** Consider the 1-pseudo false n-code $C = C_1 \cup C_2 \cup C_3 \cup C_4 \cup C_5$ given by $C_1 = C(5, 3)$, $C_2 = C(7, 4)$, $C_3 = C(7, 4) = C_4 = C_5$. Clearly C has one and only subcode $C_1 = C(5, 3)$ and other subcodes of C are identical with C(7, 4) thus C is a pseudo P-false code.

We proceed onto define yet another new type of P-linear codes built using the set of different types of codes like repetition code, parity check code, cyclic code and any other general code. This type of new codes will be known as P-linear mixed codes.

**DEFINITION 3.1.15:** *Let $C = C_1 \cup C_2 \cup \ldots \cup C_n$ be a n-code if some of the codes $C_i$ are parity check codes, some of the codes $C_j$ are repetition codes and some others $C_k$ are cyclic codes then we call C to be a Periyar linear mixed code where $1 \leq i, j, k \leq n$.*

We illustrate this by the following example.

***Example 3.1.21:*** Let $C = C_1 \cup C_2 \cup C_3 \cup C_4 \cup C_5$ be a P-linear mixed code where $C_1$ is a repetition code, $C_2$ a parity check code, $C_3$ a cyclic code and $C_4$ and $C_5$ a general code. Suppose H = $H_1 \cup H_2 \cup H_3 \cup H_4 \cup H_5$.

$$= \begin{bmatrix} 1 & 1 & 0 & 0 & 0 & 0 \\ 1 & 0 & 1 & 0 & 0 & 0 \\ 1 & 0 & 0 & 1 & 0 & 0 \\ 1 & 0 & 0 & 0 & 1 & 0 \\ 1 & 0 & 0 & 0 & 0 & 1 \end{bmatrix} \cup [1\ 1\ 1\ 1\ 1\ 1\ 1] \cup$$

$$\begin{bmatrix} 1 & 0 & 0 & 1 & 0 & 0 \\ 0 & 1 & 0 & 0 & 1 & 0 \\ 0 & 0 & 1 & 0 & 0 & 1 \end{bmatrix} \cup \begin{bmatrix} 1 & 1 & 0 & 0 & 0 & 1 \\ 0 & 1 & 0 & 0 & 1 & 0 \\ 0 & 1 & 1 & 1 & 0 & 0 \end{bmatrix} \cup$$

$$\begin{bmatrix} 1 & 0 & 0 & 0 & 1 & 1 & 0 \\ 0 & 1 & 0 & 0 & 0 & 1 & 1 \\ 0 & 0 & 1 & 0 & 1 & 1 & 0 \\ 0 & 0 & 0 & 1 & 1 & 0 & 1 \end{bmatrix}$$

H is a mixed parity check matrix of the P-linear mixed codes.

## 3.2 Application of these New Classes of Codes

In this section we give a few applications of these new classes of codes.

1. The m-pseudo false n-codes (t, t) pseudo false n-codes and (t, m) pseudo false n-codes can be used in cryptography. The true message can be sent in the m codes and just the remaining (n – m) codes will serve to mislead the intruder. Infact this code can also be used for m persons simultaneously and each one will read his message from the n-code i.e., if $C = C_1 \cup C_2 \cup \ldots \cup C_n$ here $C_{i_1} = C_{i_2} = \ldots = C_{i_m} = A$, $1 \leq i_1, i_2, \ldots, i_m \leq n$, $i_j \neq i_t$ if $j \neq t$. Then $i_1$ can be assigned the first customer, $i_2$ to the second and $i_m$ to the $m^{th}$ customer so they know the place in which their message is present. For instance $i_1 = 5$ the first customer will only decode the code word $x_5$ from $x = x_1 \cup \ldots \cup x_n$ and he will not bother to pay any attention to other codes. Thus it is very difficult for the intruder to hack this way of message transmission. In case they have two sets of customers say t and m in number then (t, m) pseudo n code can be used.

2. These codes can be used in networking of computers when in the work place some m of them work on the same data and t on another set of data. Infact the (t, m) pseudo n-code can be extended to any $(t_1, \ldots, t_r)$ pseudo n code without any difficultly. Further the transmission or working with them is not difficult as we have super computer to cater to our needs.

3. These n-codes or n-false codes or m pseudo false n codes or (m, t) pseudo false n codes can be easily used when we cannot use ARQ protocols, so n- false codes and m-pseudo false n codes will serve best in those cases.

4. These new classes of codes can also be used in data storage. The authors strongly feel these codes are best suited in any networking. Finally they can be used with highest degree of security in defence department of any nation.

5. Periyar code finds its applications in channels where all the n information which are passed simultaneously are distinct. In such channels Periyar codes are best suited.

# FURTHER READING

# INDEX

### A

## B

## C

## D

## E

## F

## G

## H

## I

## L

## M

## N

**O**

**P**

## R

## S

**T**

**W**

# ABOUT THE AUTHORS

**Dr.W.B.Vasantha Kandasamy** is an Associate Professor in the Department of Mathematics, Indian Institute of Technology Madras, Chennai. In the past decade she has guided 11 Ph.D. scholars in the different fields of non-associative algebras, algebraic coding theory, transportation theory, fuzzy groups, and applications of fuzzy theory of the problems faced in chemical industries and cement industries.

She has to her credit 640 research papers. She has guided over 57 M.Sc. and M.Tech. projects. She has worked in collaboration projects with the Indian Space Research Organization and with the Tamil Nadu State AIDS Control Society. This is her 33rd book.

On India's 60th Independence Day, Dr.Vasantha was conferred the Kalpana Chawla Award for Courage and Daring Enterprise by the State Government of Tamil Nadu in recognition of her sustained fight for social justice in the Indian Institute of Technology (IIT) Madras and for her contribution to mathematics. (The award, instituted in the memory of Indian-American astronaut Kalpana Chawla who died aboard Space Shuttle Columbia). The award carried a cash prize of five lakh rupees (the highest prize-money for any Indian award) and a gold medal.

She can be contacted at vasanthakandasamy@gmail.com
You can visit her on the web at: http://mat.iitm.ac.in/~wbv or: http://www.vasantha.net

---

**Dr. Florentin Smarandache** is an Associate Professor of Mathematics at the University of New Mexico in USA. He published over 75 books and 100 articles and notes in mathematics, physics, philosophy, psychology, literature, rebus. In mathematics his research is in number theory, non-Euclidean geometry, synthetic geometry, algebraic structures, statistics, neutrosophic logic and set (generalizations of fuzzy logic and set respectively), neutrosophic probability (generalization of classical and imprecise probability). Also, small contributions to nuclear and particle physics, information fusion, neutrosophy (a generalization of dialectics), law of sensations and stimuli, etc.

He can be contacted at smarand@unm.edu